\documentclass[10pt]{article}

\usepackage{graphicx} % Add graphics capabilities
\usepackage{amsgen,amsmath,amstext,amsbsy,amsopn, amssymb,amsfonts}%\usepackage[colorlinks=true,linkcolor=red]{hyperref} % Hyperlink capabilities
\usepackage{mathrsfs}

\usepackage{memhfixc} % This package is required to resolve incompatibilities
\usepackage{pdfsync}
\usepackage{amsgen,amsmath,amstext,amsbsy,amsopn, amssymb,amsfonts}
\usepackage{float}
\usepackage{epsfig}

\usepackage{fullpage}
\DeclareMathAlphabet{\mathpzc}{OT1}{pzc}{m}{it}

\def\N{\mathbb N}
\def\R{\mathbb R}
\def\Z{\mathbb Z}
\def\C{\mathbb C}
\def\B{\mathbb B}
\def\SS{\mathbb S}

\def\BM{\mathcal{M}}
\def\BN{\mathcal{N}}

\def\BA{\mathcal{A}}

\def\BJ{\mathcal{J}}

\def\BR{\mathcal{R}}
\def\BL{\mathcal{L}}

\def\BY{\mathcal{Y}}
\def\BX{\mathcal{X}}
\def\BD{\mathcal{D}}
\def\BF{\mathcal{F}}
\def\BH{\mathcal{H}}

\def\BC{\mathcal{C}}
\def\BO{\mathcal{O}}
\def\BW{\mathcal{W}}
\def\BU{\mathcal{U}}
\def\BV{\mathcal{V}}
\def\BZ{\mathcal{Z}}

\def\vp{\varphi}

\def\e{\varepsilon}

\def\Chi{\raise .3ex \hbox{\large $\chi$}} 
\def\vp{\varphi}
\def\s{\sigma}
\def\d{\delta}

\def\a{\alpha}

\def\p{\partial}

\def\OO{\Omega}

\def\oo{\omega}

\def\abs{\Big |}

\def\carre{\ \hfill $\Box$}
\def\ii{\infty}

\def\n{|\!|}
\def\no{\Big|\!\Big|}

\def\ds{\displaystyle}
\def\bs{\boldsymbol}

\def\t{\tilde}

\def\cs{{\mathfrak c}_s}
\def\ex{{\sf e}}

\def\ini{{\rm in}}

\def \be{\begin{equation}}
\def \ee{\end{equation}}

\newtheorem{corollary}{Corollary}

\newtheorem{definition}{Definition}
\newtheorem{remark}{Remark}
\newtheorem{remarks}{Remark}[section]
\newtheorem{lemma}{Lemma}
\newtheorem{lemmas}{Lemma}[section]
\newtheorem{proposition}{Proposition}
\newtheorem{propositions}{Proposition}[section]
\newtheorem{theorem}{Theorem}

\title{Stability and instability for subsonic\\ travelling waves 
of the Nonlinear Schr\"odinger Equation\\ in dimension one}
\author{D. Chiron\footnote{Laboratoire J.A. Dieudonn\'e, Universit{\'e} 
de Nice-Sophia Antipolis, Parc Valrose, 06108 Nice Cedex 02, France. 
\quad \quad \quad \quad \quad \quad
{\sf e-mail}: chiron@unice.fr.} }

\date{}

\begin{document}

\maketitle

\begin{abstract}
We study the stability/instability of the subsonic travelling waves of 
the Nonlinear Schr\"odinger Equation in dimension one. 
Our aim is to propose several methods for showing instability 
(use of the Grillakis-Shatah-Strauss theory, proof of existence 
of an unstable eigenvalue via an Evans function) or stability. 
For the later, we show how to construct in a systematic way 
a Liapounov functional for which the travelling wave is a local minimizer.
%in the stable case of the Grillakis-Shatah-Strauss theory. 
These approaches allow 
to give a complete stability/instability analysis in the energy space 
including the critical case of the kink solution. We also treat 
the case of a cusp in the energy-momentum diagram.
\end{abstract}


%%%%%%%%%%%%%%%%%%%%%%%%%%%%%%%%%%%%%%%%%%%%%
%%%%%%%%%%%%%%%%%%%%%%%%%%%%%%%%%%%%%%%%%%%%%
\section{Introduction}

\ \quad This paper is a continuation of our previous work \cite{C1d}, where 
we consider the one dimensional Nonlinear Schr\"odinger Equation
\be
\tag{NLS}
i \frac{\p \Psi}{\p t} + \p_x^2 \Psi + \Psi f(|\Psi|^2) = 0 . 
\ee
This equation appears as a relevant model in condensed matter physics: 
Bose-Einstein condensation and superfluidity (see \cite{RB}, \cite{GP}, 
\cite{G}, \cite{AHMNPTB}); Nonlinear Optics (see, for instance, the 
survey \cite{KL}). Several nonlinearities may be encountered in physical 
situations: $ f(\varrho) = \pm \varrho $ gives rise to the focusing/defocusing 
cubic (NLS); $ f(\varrho) = 1- \varrho$ to the so called Gross-Pitaevskii 
equation; $ f(\varrho)= - \varrho^2 $ (see \cite{KNSQ} for Bose-Einstein 
condensates); more generally a pure power; the ``cubic-quintic'' (NLS) 
(see \cite{BP}), where
$$ f(\varrho) = - \a_1 + \a_3 \varrho - \a_5 \varrho^2 $$
and $\a_1 $, $\a_3 $ and $\a_5$ are positive constants such 
that $f $ has two positive roots; and in Nonlinear Optics, we may 
take (see \cite{KL}):
\be
\label{nonlin}
f(\varrho) = - \alpha \varrho^\nu - \beta \varrho^{2\nu}, 
\quad 
f(\varrho) = - \frac{\varrho_0 }{2} 
\Big( \frac{1}{ (1+ \frac{1}{\varrho_0} )^\nu } 
- \frac{1}{(1+ \frac{\varrho}{\varrho_0} )^\nu} \Big), 
\quad 
f(\varrho) = - \alpha \varrho \Big( 1 + \gamma \, {\rm tanh} 
( \frac{\varrho^2 - \varrho_0^2}{\sigma^2} ) \Big) ,
\ee
where $\alpha$, $\beta$, $\gamma $, $\nu$,  $\sigma > 0$ are 
given constants (the second one, for instance, takes into account 
saturation effects), etc. As a consequence, as in our work \cite{C1d}, 
we shall consider a rather general nonlinearity $f$, with $f$ of 
class $ \BC^2 $. In the context of Bose-Einstein condensation or 
Nonlinear Optics, the natural condition at infinity appears to be
$$ | \Psi |^2 \to r_0^2 \quad \quad \quad {\rm as} \quad |x| \to +\ii ,$$
where $ r_0 > 0 $ is such that $f(r_0^2)=0$.\\

For solutions $\Psi$ of (NLS) which do not vanish, we may use the 
Madelung transform
$$ \Psi = A \exp ( i \phi ) $$
and rewrite (NLS) as an hydrodynamical system with an additional 
quantum pressure
\be
\label{MadTWk}
\left\{\begin{array}{ll}
\ds{ \p_t A + 2 \p_{x} \phi \p_{x} A 
+ A \p_{x}^2 \phi } = 0 \\ \ \\ 
\ds{ \p_{t} \phi + (\p_{x} \phi)^2 - f ( A^2 ) 
- \frac{ \p_{x}^2 A }{A} } = 0 
\end{array}\right. 
\quad \quad {\rm or} \quad \quad
\left\{\begin{array}{ll}
\ds{ \p_t \rho + 2 \p_{x} (\rho u ) } = 0 \\ \ \\ 
\ds{ \p_{t} u + 2 u \p_{x} u - \p_x (f ( \rho )) 
- \p_x \Big( \frac{ \p_{x}^2 (\sqrt{\rho}) }{\sqrt{\rho} } \Big) } = 0 ,
\end{array}\right. 
\ee
with $(\rho,u) \equiv (A^2,\p_x \phi)$. When neglecting the quantum 
pressure and linearizing this Euler system around the 
particular trivial solution $\Psi = r_0$ (or $(A,u)= (r_0,0)$), 
we obtain the free wave equation
$$ \left\{\begin{array}{ll}
\ds{ \p_t \bar{A} + r_0 \p_{x} \bar{U} } = 0 \\ \ \\ 
\ds{ \p_{t} \bar{U} - 2 r_0 f'(r_0^2) \p_x \bar{A} } = 0 
\end{array}\right. $$
with associated speed of sound
$$ \cs \equiv \sqrt{- 2 r_0^2 f'(r_0^2)} > 0 $$
provided $f$ satisfies the defocusing assumption $f'(r_0^2) < 0$ (that 
is the Euler system is hyperbolic in the region $\rho \simeq r_0^2$), 
which we will assume throughout the paper. Concerning the rigorous 
justification of the free wave regime for the Gross-Pitaevkii equation 
(in arbitrary dimension), see \cite{BDS}. The speed of sound $ \cs $ 
enters in a crucial way in the question of existence of travelling waves 
for (NLS) with modulus tending to $ r_0 $ at infinity (see, {\it e.g.}, 
\cite{C1d}).\\

The Nonlinear Schr\"odinger equation formally preserves the energy
$$ E(\psi) \equiv \int_\R |\p_x \psi |^2 + F (|\psi|^2) \ dx , $$
where $ F(\varrho) \equiv \displaystyle \int_{\varrho}^{r_0^2} f$. 
Since $ F(\varrho) \sim \ds{ \frac{ \cs^2}{8r_0^2} ( \varrho - r_0^2 )^2 } 
\sim \ds{ \frac{ \cs^2}{2} ( \sqrt{\varrho} - r_0 )^2 }$ when 
$\varrho \to r_0^2 $, it follows that 
the natural energy space turns out to be the space
$$ \mathcal{Z} \equiv \Big\{ \psi \in L^\ii(\R) , \ 
\p_x \psi \in L^2(\R), \ |\psi| - r_0 \in L^2(\R) \Big\} \subset \BC_b(\R,\C) , $$
endowed with the distance
$$ 
d_{\mathcal{Z}} ( \psi , \tilde{\psi} ) \equiv \n \p_x \psi - \p_x \tilde{\psi} \n_{L^2(\R)} 
+ \n \, | \psi | - | \tilde{\psi} | \, \n_{L^2(\R)} 
+ \abs  \psi(0) - \tilde{\psi}(0) \abs .
$$
The Cauchy problem has been shown to be locally well posed in the 
Zhidkov space $ \{ \psi \in L^\ii(\R) , \ \p_x \psi \in L^2(\R) \} $ by 
P. Zhidkov \cite{Z} (see also the work by C. Gallo \cite{GaZhi}). 
For global well-posedness results, see \cite{Ga} and \cite{GG}. 
More precisely, the local well-posedness we shall use is the following.

%%%%%%%%%%%%%%%%%%%%%
\begin{theorem} [\cite{Z}, \cite{GaZhi}]
\label{Cauchy} 
Let $\Psi^\ini \in \mathcal{Z} $. Then, there exists $ T_* > 0 $ and a unique 
solution $ \Psi $ to {\rm (NLS)} such that $\Psi_{|t=0} = \Psi^\ini $ and 
$ \Psi - \Psi^\ini \in \BC([0,T_*), H^1(\R))$. Moreover, $E(\Psi(t))$ does 
not depend on $t$.
\end{theorem}
%%%%%%%%%%%%%%%%%%%%%

The other quantity formally conserved by the Schr\"odinger flow, 
due to the invariance by translation, is the momentum. 
The momentum is not easy to define in dimension one 
for maps that vanish somewhere (see \cite{BGSsurvey}, \cite{BGSS}). 
However, if $\psi $ does not vanish, we may lift $ \psi = A \ex^{i\phi} $, 
and then the correct definition of the momentum is given by \cite{KY}:
$$ P(\psi) \equiv \int_\R \langle i \psi | \p_x \psi \rangle 
\Big( 1 - \frac{r_0^2}{|\psi|^2} \Big) \ dx = 
\int_\R (A^2 - r_0^2) \p_x \phi \ dx , $$
where $\langle \cdot | \cdot \rangle$ denotes the real scalar 
product in $\C$. We define
$$ 
\mathcal{Z}_{\rm hy} \equiv \{ v \in \mathcal{Z}, \inf_\R | v | > 0 \} ,
$$
which is the open subset of $ \mathcal{Z} $ in which we have lifting 
and where the hydrodynamical formulation \eqref{MadTWk} of (NLS) is 
possible through the Madelung transform. It turns out that if the initial 
datum belongs to $ \mathcal{Z}_{\rm hy} $, the solution of (NLS) provided 
by Theorem \ref{Cauchy} remains in $ \mathcal{Z}_{\rm hy} $ for small times, 
and that the momentum is indeed conserved on this time interval 
(see \cite{GaZhi}).

%%%%%%%%%%%%%%%%%%%%%%%%%%%%%%%%%%%%%%%%%%%%%%%%%%%%%%%%%%%%%%
\subsection{The travelling waves and energy-momentum diagrams}
\label{sextravelo}

The travelling waves with speed of propagation $c$ are 
special solutions of (NLS) of the form
$$ \Psi (t,x) = U(x-c t) . $$
The profile $U$ has then to solve the ODE
\be
\tag{TW$_c$}
\p_x^2 U + U f(|U|^2) = i c \p_{x} U 
\ee
together with the condition  $|U(x)| \to r_0$ as $x\to \pm \ii$. 
These particular solutions play an important role in the long 
time dynamics of (NLS) with nonzero condition at infinity. 
Possibly conjugating (TW$_c$), we see that we may assume 
that $c \geq 0$ without loss of generality. Moreover, we 
shall restrict ourselves to travelling waves which belong 
to the energy space $\mathcal{Z}$ (so that $ |U| \to r_0 $ at $\pm \ii$ 
by Sobolev embedding $H^1(\R) \hookrightarrow 
\BC_0(\R,\C) \equiv \{ h \in \BC(\R,\C), \ \lim_{\pm \ii } h = 0 \}$). 
For travelling waves $ U_c \in \mathcal{Z} $ that do not vanish in $\R$, hence that may 
be lifted $ U_c = A_c \ex^{i\phi_c}$, the ODE (TW$_c$) can be transformed 
(see, {\it e.g.}, \cite{C1d}) into the system
$$
\p_x \phi_c = \frac{c}{2} \cdot \frac{\eta_c}{\eta_c + r_0^2} , 
\quad \quad \quad 
2 \p^2_x \eta_c + \BV_c'(\eta_c) = 0 ,
\quad \quad {\rm with} \quad \quad
\eta_c \equiv A_c^2 - r_0^2 ,
$$
and where the function $\BV_c$ is related to $f$ by the formula
$$
 \BV_c(\xi ) \equiv c^2 \xi^2 - 4 ( r_0^2 + \xi) F ( r_0^2 + \xi) .
$$
To a nontrivial travelling wave $U_c$ is associated (see \cite{C1d}) 
some $ \xi_c \geq -r_0^2 $ such that $ \BV_c(\xi_c ) = 0 \not = \BV_c' (\xi_c ) $ 
and $ \BV_c $ is negative between $ \xi_c $ and $ -r_0^2 $, and 
$ \eta_c $ varies between $0$ and $\xi_c$, that is there holds 
$ \{ \inf_\R |U_c| , \sup_\R |U_c| \} = \{ r_0 , \sqrt{ r_0^2 + \xi_c} \} $. 
Moreover, the only travelling wave solution (if it exists) that vanishes 
somewhere is for $c=0$ and is called the kink: it is an odd solution 
(up to a space translation) and then $\xi_0 = 0 $.

We have also seen in \cite{C1d} that any travelling wave in 
$ \mathcal{Z} $ with speed $ c > \cs $ is constant; and also 
that any nonconstant travelling wave in $ \mathcal{Z} $ of speed 
$ c_* \in ( 0, \cs )$ belongs to a unique (up to the natural invariances: 
phase factor and translation) local branch $ c \mapsto U_c $ defined for $c$ 
close to $ c_* $.

\bigskip

In \cite{C1d}, we have investigated the qualitative behaviours of 
the travelling waves for (NLS) with nonzero condition at infinity 
for a general nonlinearity $f$. A particular attention has been payed 
in \cite{C1d} to the transonic limit, where we have an asymptotic 
behaviour governed by the Korteweg-de Vries or the generalized Korteweg-de 
Vries equation. In order to illustrate the 
very different situations we may encounter when we allow a general 
nonlinearity $f$, we give now some energy-momentum diagrams we have 
obtained (one is taken from the appendix in \cite{CS}, where we 
have performed numerical simulations in dimension two for the 
model cases we have studied in \cite{C1d}).\\

\begin{figure}[H]
\begin{center}
\includegraphics[viewport=10 10 500 450,width=6cm,height=7cm]{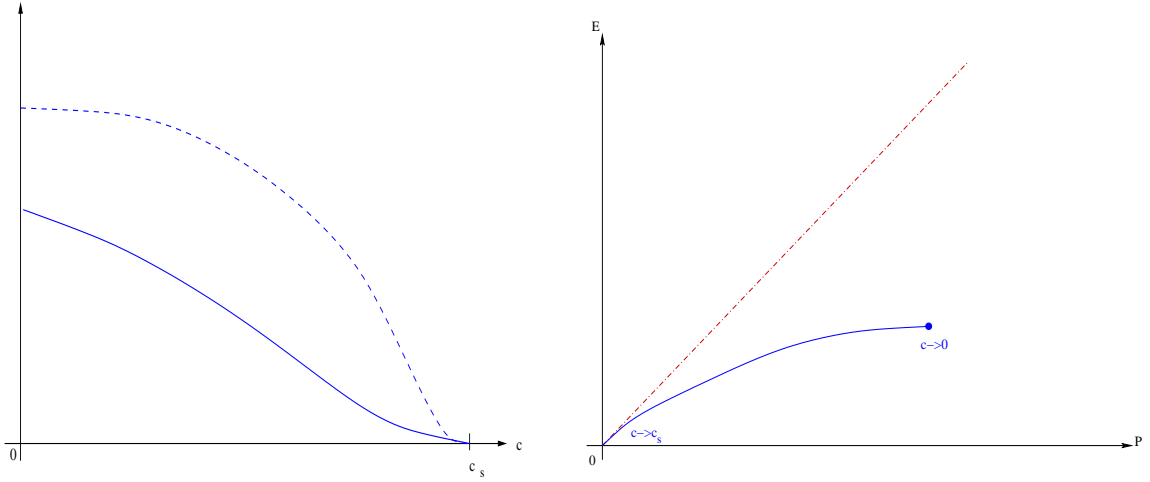} 
\hspace{1.5cm} 
\includegraphics[viewport=10 0 500 450,width=6cm,height=6cm]{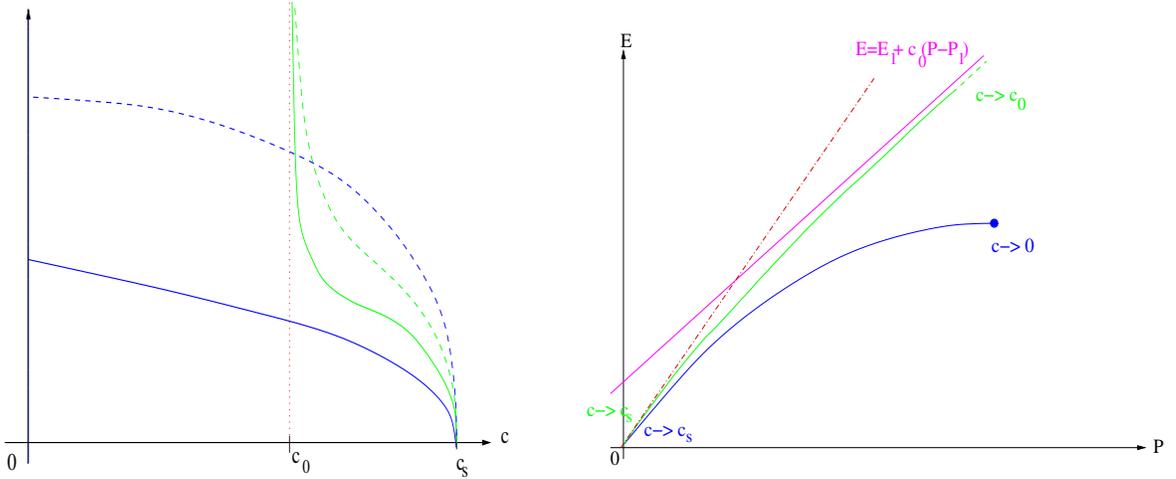}
\end{center}
\caption{(a) Energy (dashed curve) and momentum (full curve) vs. speed;  
(b) $(E,P)$ diagram}
\label{zGP}
\end{figure}

\noindent $\bullet$ The Gross-Pitaevskii nonlinearity: 
$ f( \varrho ) = 1 - \varrho  $ (see figure \ref{zGP}).

\begin{figure}[H]
\begin{center}
\includegraphics[viewport=10 10 500 450,width=6cm,height=7cm]{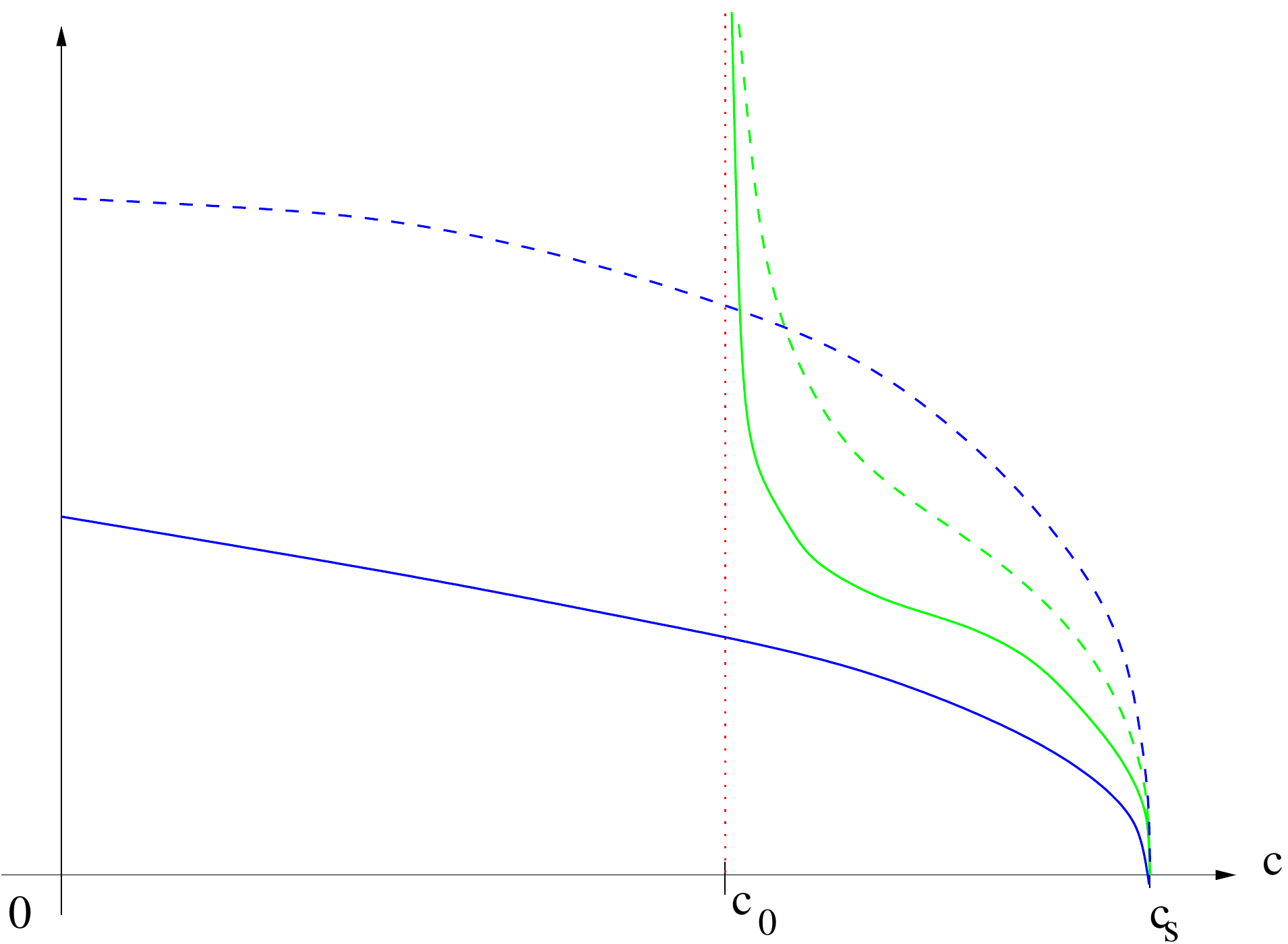} 
\hspace{1.5cm}  
\includegraphics[viewport=10 0 500 450,width=6cm,height=6cm]{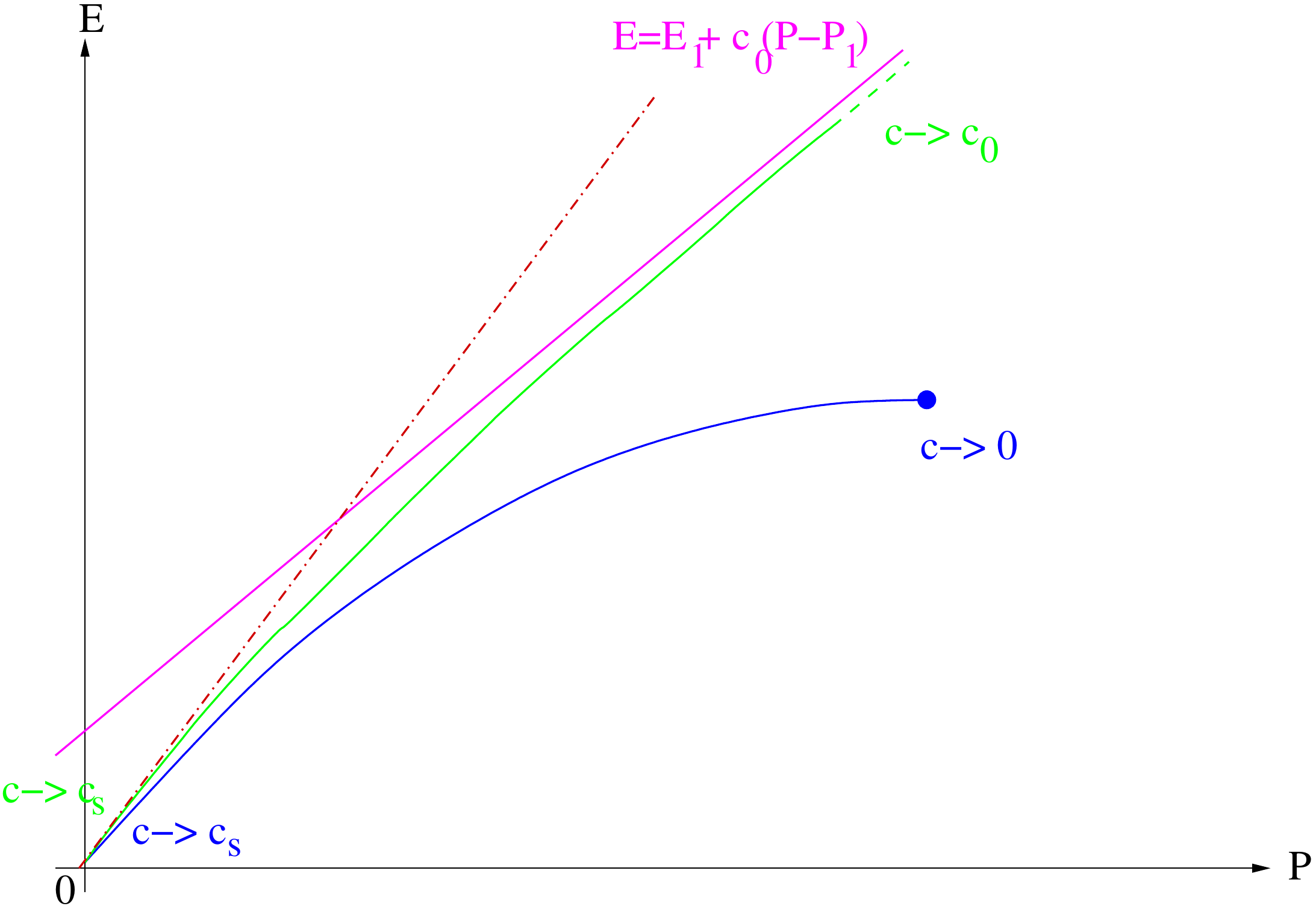}
\end{center}
\caption{(a) Energy (dashed curve) and momentum (full curve) vs. speed;  
(b) $(E,P)$ diagram}
\label{zcqsI}
\end{figure}

\noindent $\bullet$ A cubic-quintic-septic nonlinearity (I): $ f( \varrho ) = - (\varrho-1) 
+ \frac{3}{2} (\varrho-1)^2 - \frac{3}{2} (\varrho-1)^3 $ (see figure \ref{zcqsI}).

\begin{figure}[H]
\begin{center}
\includegraphics[viewport=10 10 500 450,width=6cm,height=7cm]{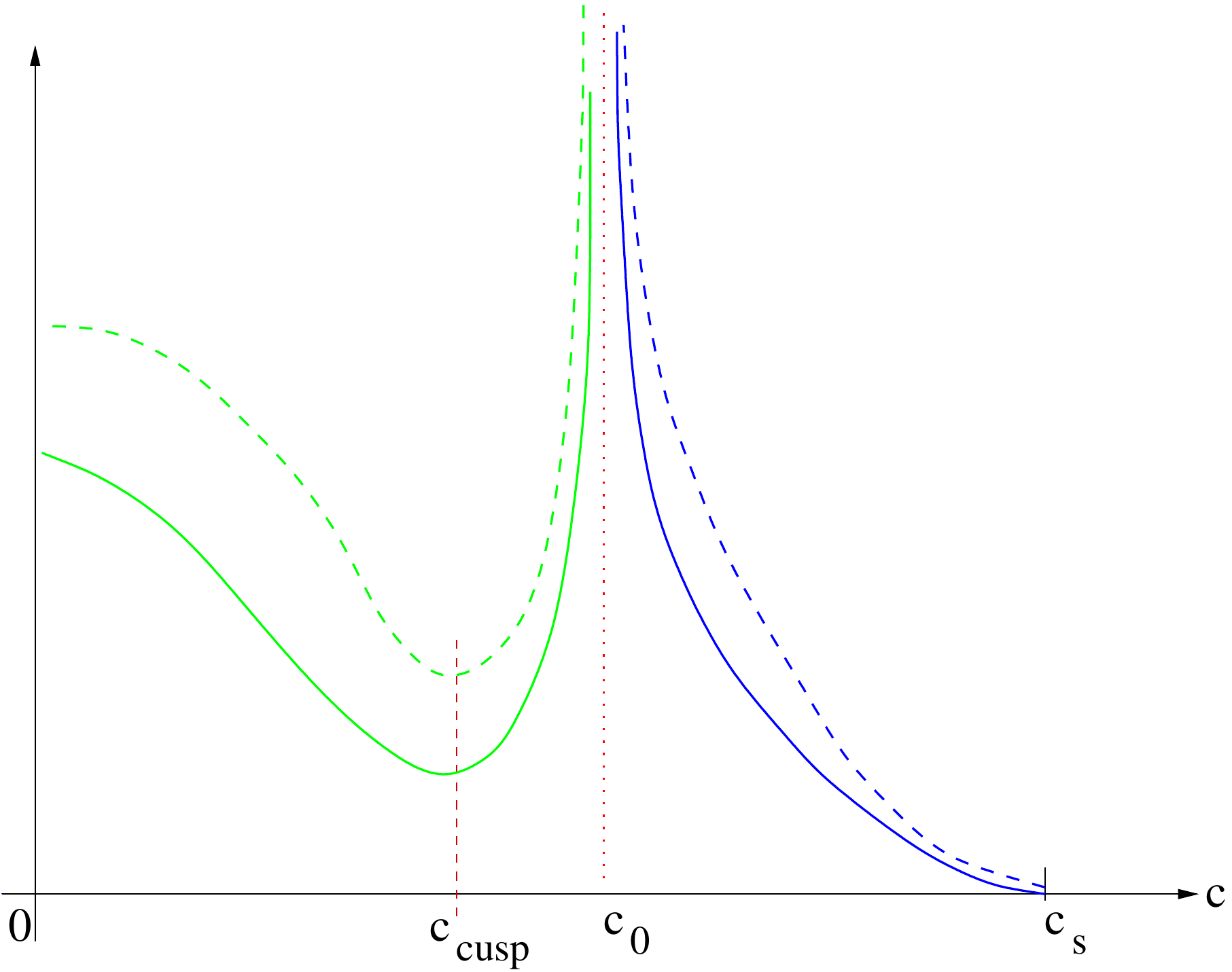} 
\end{center}
\caption{Energy (dashed curve) and momentum (full curve) vs. speed}
\label{zcqsII}
\end{figure}

\begin{figure}[H]
\begin{center}
\includegraphics[viewport=10 10 500 450,width=4.5cm,height=6cm]{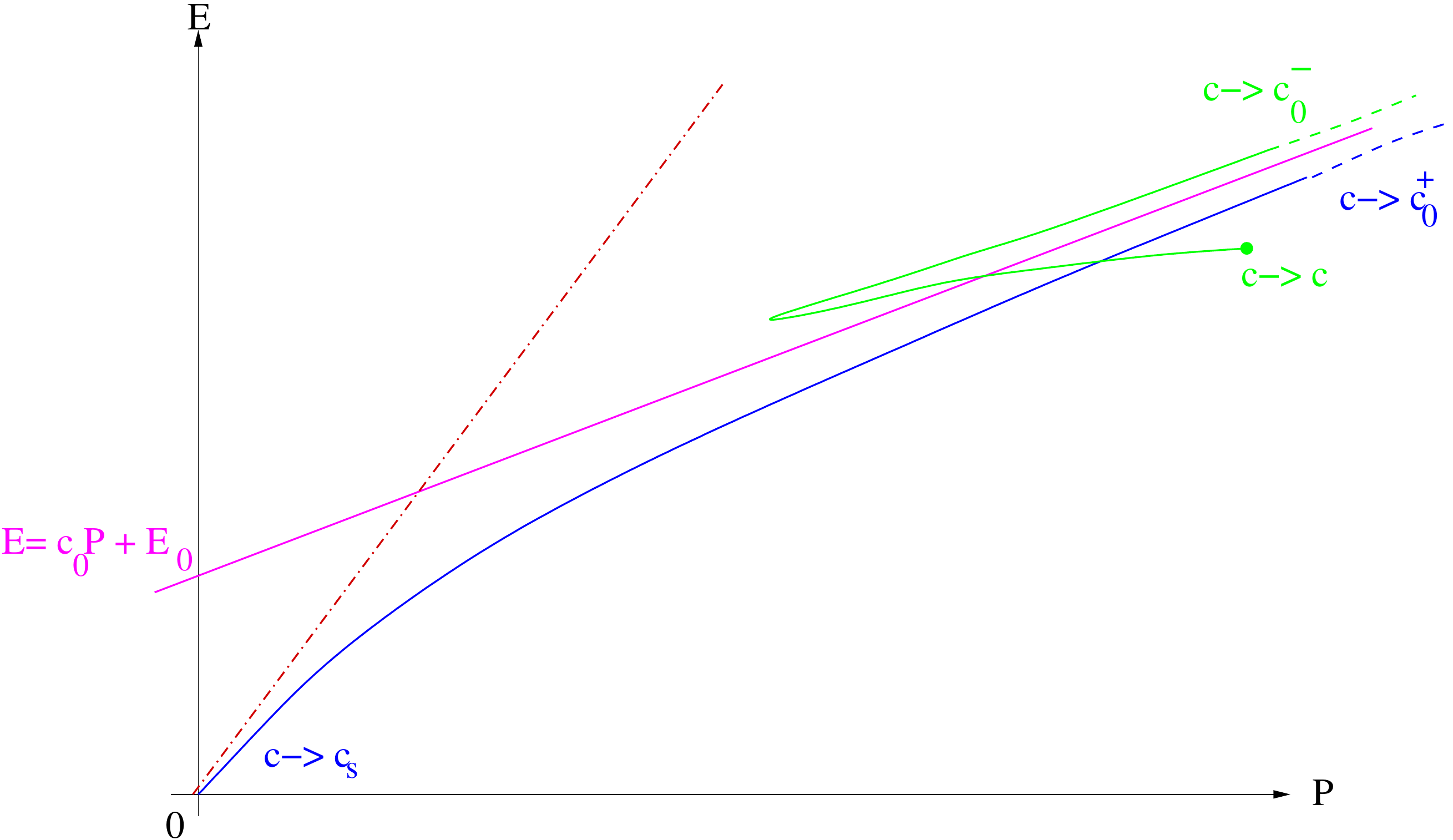} 
\hspace{3cm} 
\includegraphics[viewport=10 10 500 450,width=4.5cm,height=6cm]{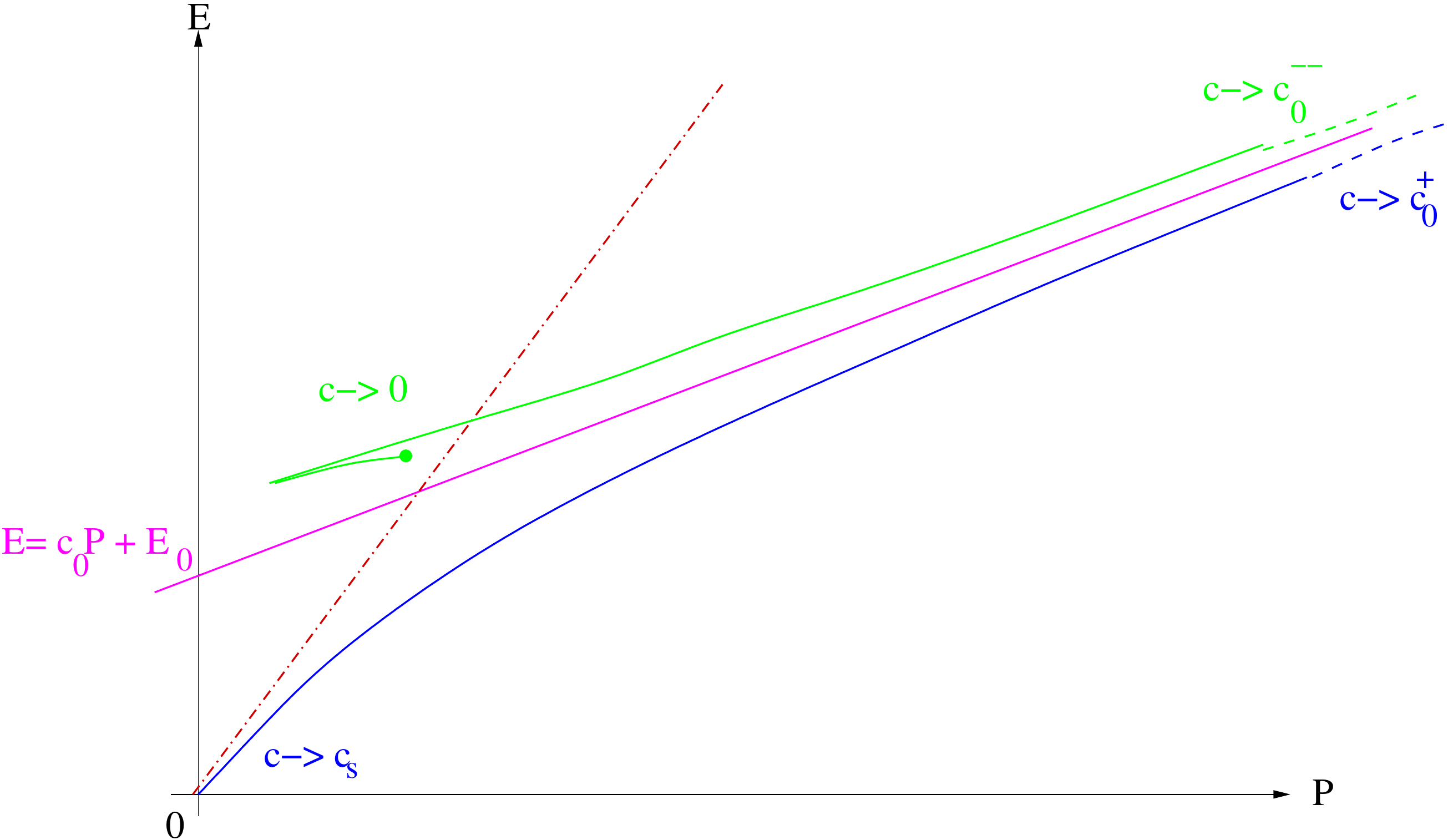}
\end{center}
\caption{The two $(E,P)$ diagrams}
\label{zzcqsII}
\end{figure}

\noindent $\bullet$ A cubic-quintic-septic nonlinearity (II): 
$ f(\varrho) \equiv - 4 (\varrho-1) - 36 (\varrho-1)^3 $ or 
$ f(\varrho) \equiv - 4 (\varrho-1) - 60 (\varrho-1)^3 $. For 
these two nonlinearities, the graph of $E$ and $P$ vs. speed $c$ 
is given in figure \ref{zcqsII}, but the $(E,P)$ diagrams are 
respectively those in figure \ref{zzcqsII}.

\begin{figure}[H]
\begin{center}
\includegraphics[viewport=10 10 500 450,width=6cm,height=7cm]{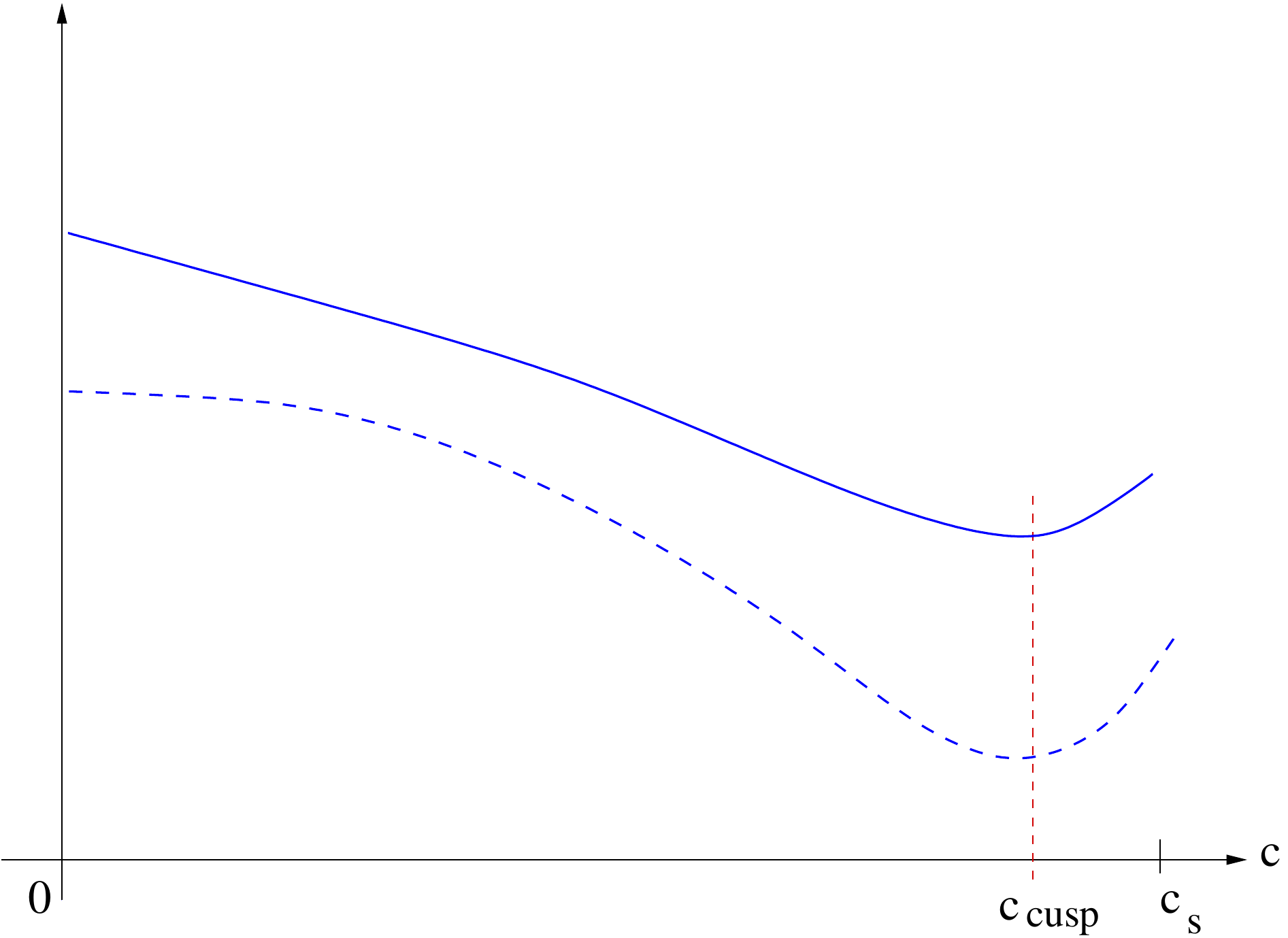} 
\hspace{1.5cm}  
\includegraphics[viewport=10 0 500 450,width=5cm,height=6cm]{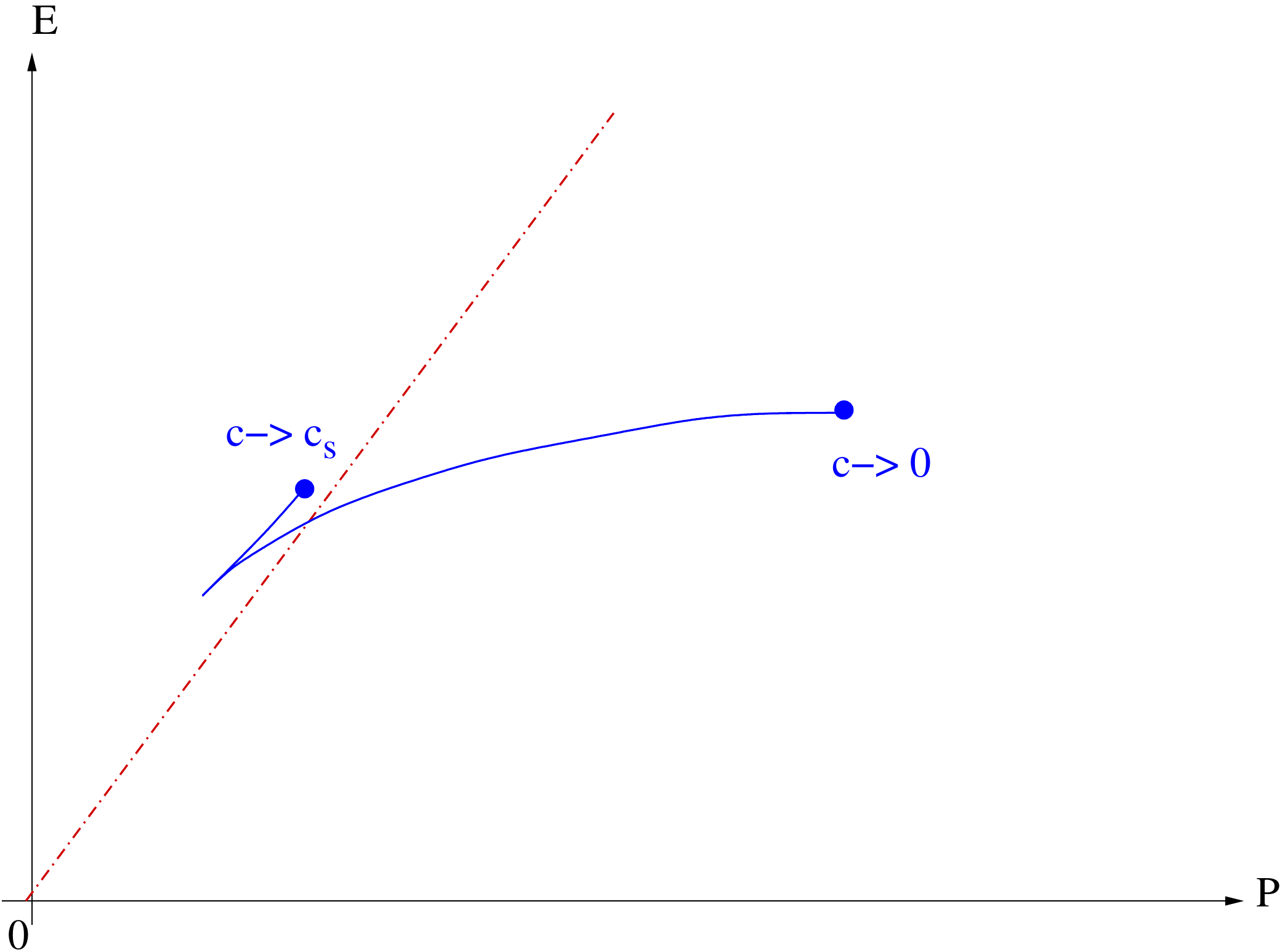}
\end{center}
\caption{(a) Energy (dashed curve) and momentum (full curve) vs. speed; 
(b) $(E,P)$ diagram}
\label{zcqsIII}
\end{figure}

\noindent $\bullet$ A cubic-quintic-septic nonlinearity (III): 
$ f(\varrho) \equiv - \frac12 (\varrho-1) + \frac34 (\varrho-1)^2 - 2 (\varrho-1)^3 $ 
(see figure \ref{zcqsIII}).

\begin{figure}[H]
\begin{center}
\includegraphics[viewport=10 20 500 450,width=6cm,height=7cm]{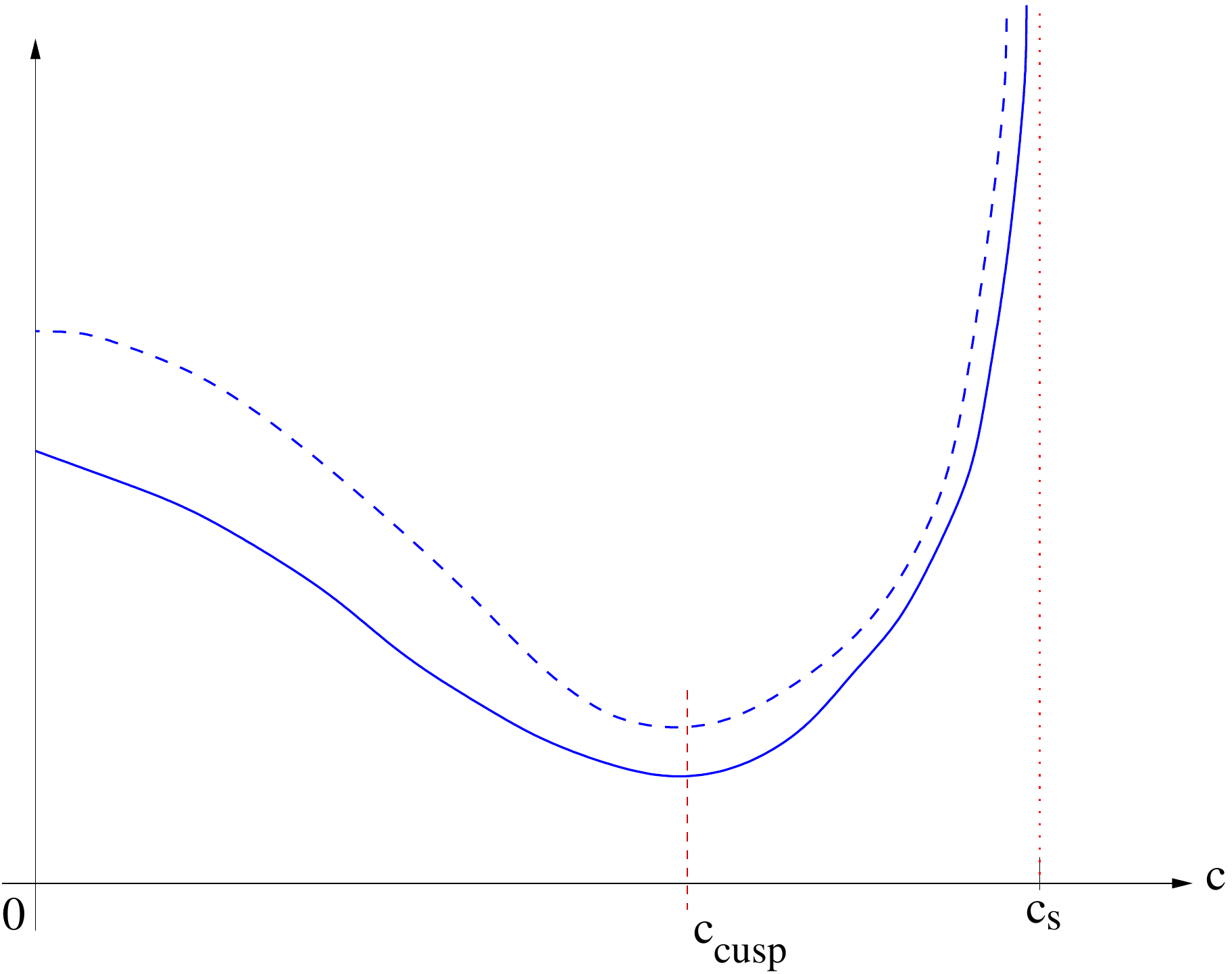} 
\hspace{1.5cm}  
\includegraphics[viewport=10 0 500 450,width=5cm,height=6cm]{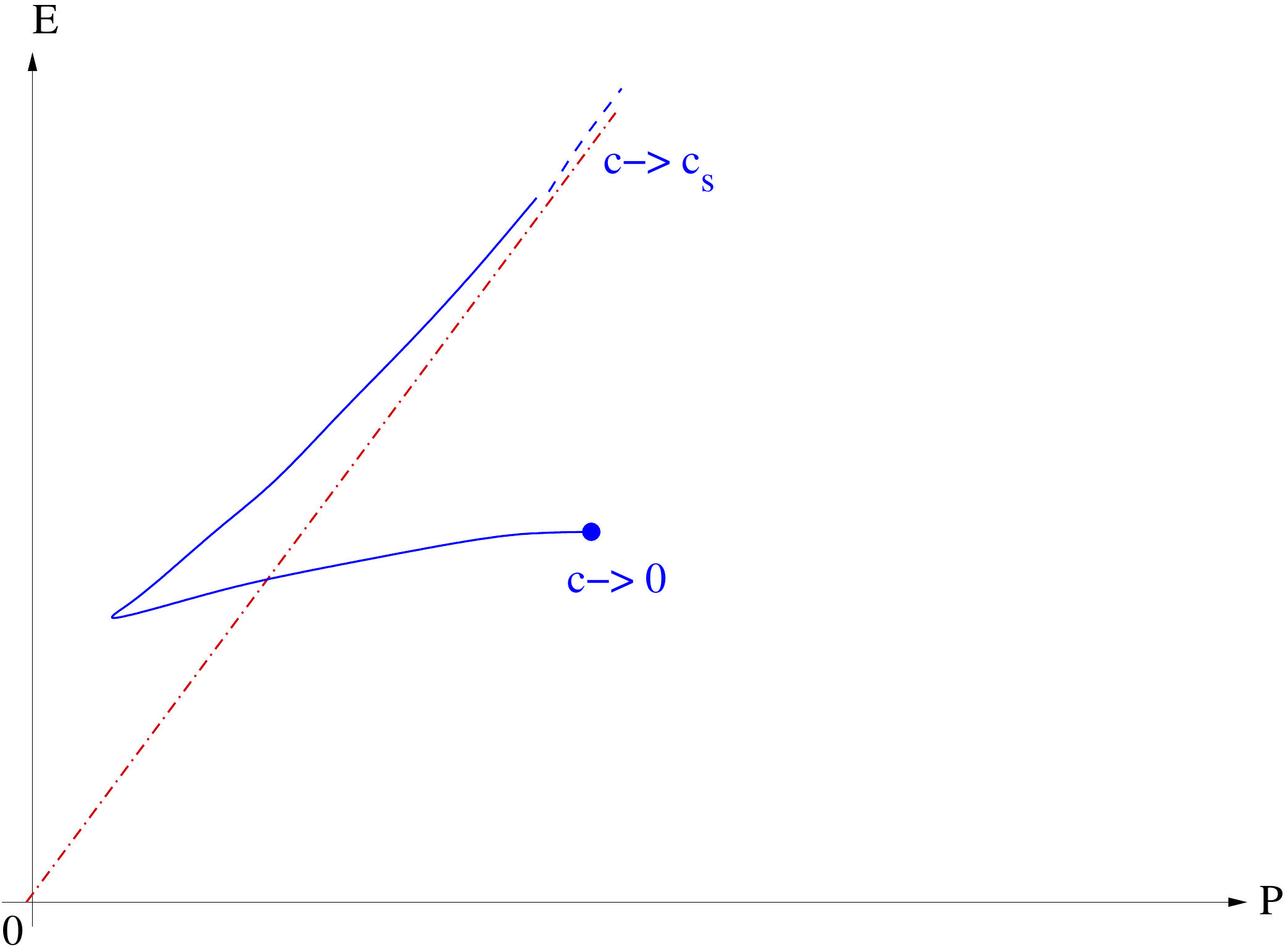}
\end{center}
\caption{(a) Energy (dashed curve) and momentum (full curve) vs. speed; 
(b) $(E,P)$ diagram}
\label{zdege}
\end{figure}

\noindent $\bullet$ A degenerate case: 
$ f(\varrho) \equiv - 2(\varrho - 1) + 3 (\varrho - 1)^2 
- 4 (\varrho - 1)^3 + 5(\varrho - 1)^4 - 6 (\varrho - 1)^5 $ (see figure \ref{zdege}).

\begin{figure}[H]
\begin{center}
\includegraphics[viewport=10 15 500 450,width=6cm,height=7cm]{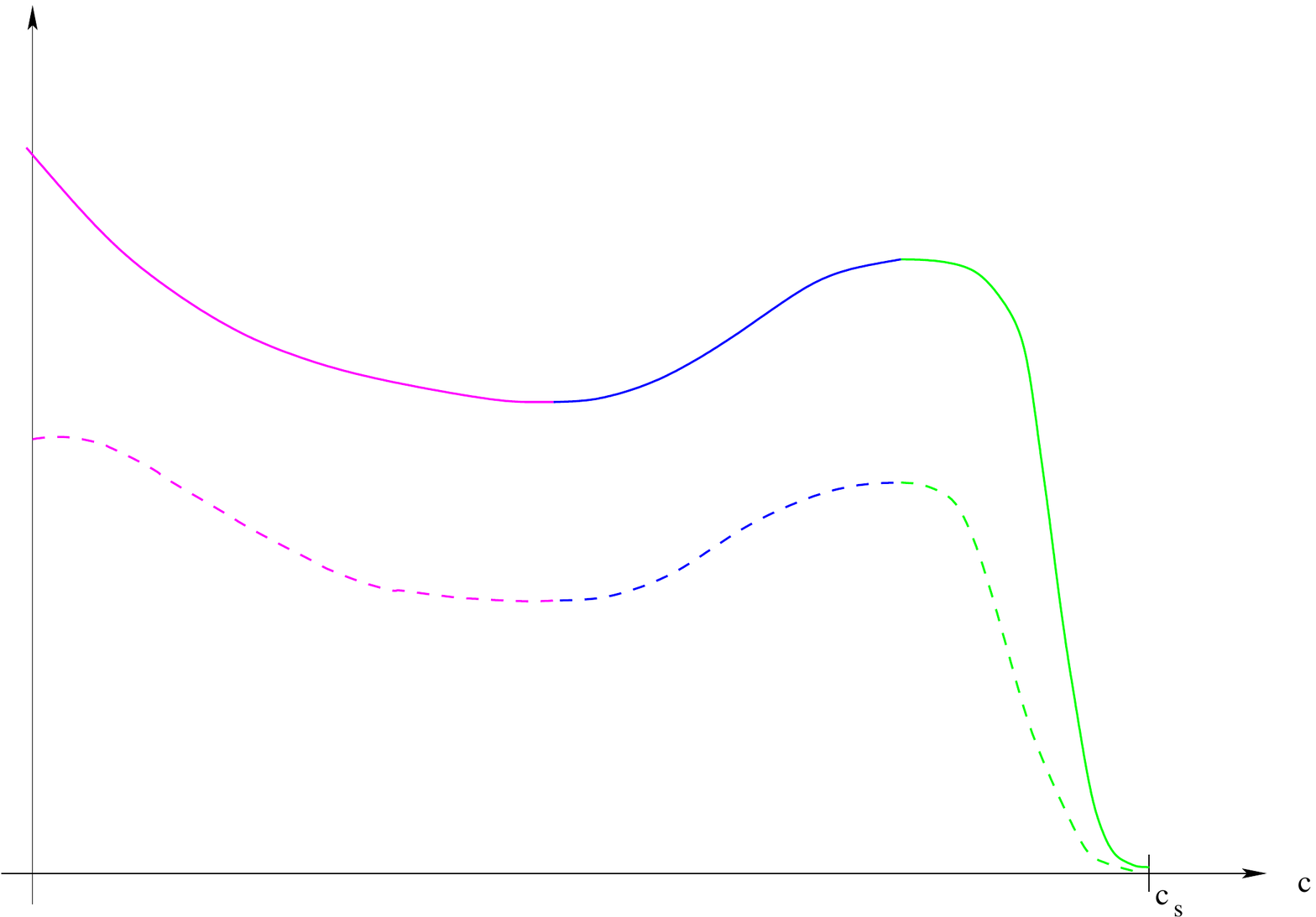} 
\hspace{2cm}  
\includegraphics[viewport=10 10 500 450,width=5cm,height=6cm]{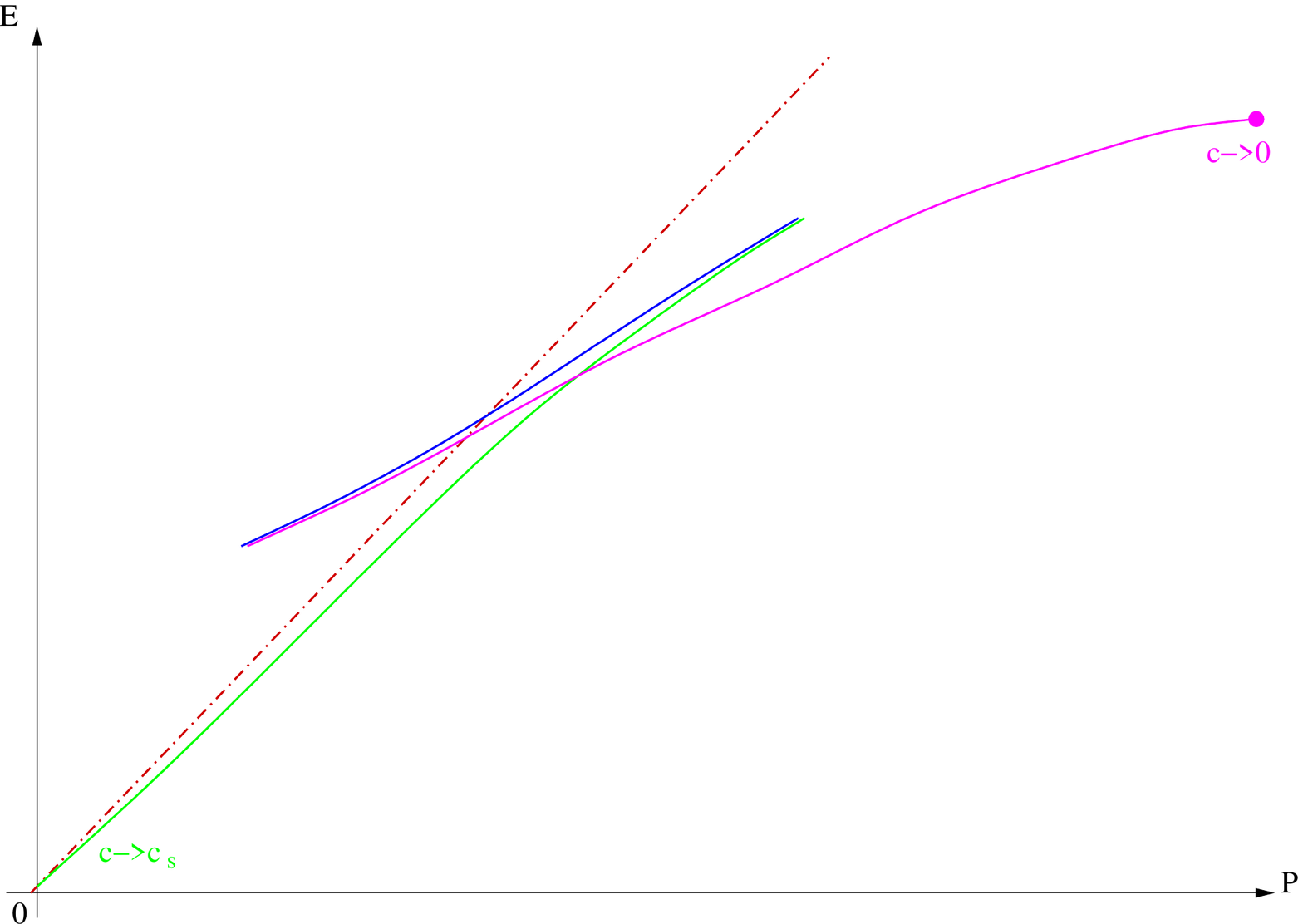}
\end{center}
\caption{(a) Energy (*) and momentum (+) vs. speed; 
(b) qualitative $(E,P)$ diagram}
\label{zdegeper}
\end{figure}

\noindent $\bullet$ A perturbation of the previous degenerate case: 
$ f(\varrho) \equiv - 2(\varrho - 1) + (3- 10^{-3}) (\varrho - 1)^2 
- 4 (\varrho - 1)^3 + 5(\varrho - 1)^4 - 6 (\varrho - 1)^5 $ (see figure \ref{zdegeper}).

\begin{figure}[H]
\begin{center}
\includegraphics[viewport=10 0 500 450,width=6cm,height=5.5cm]{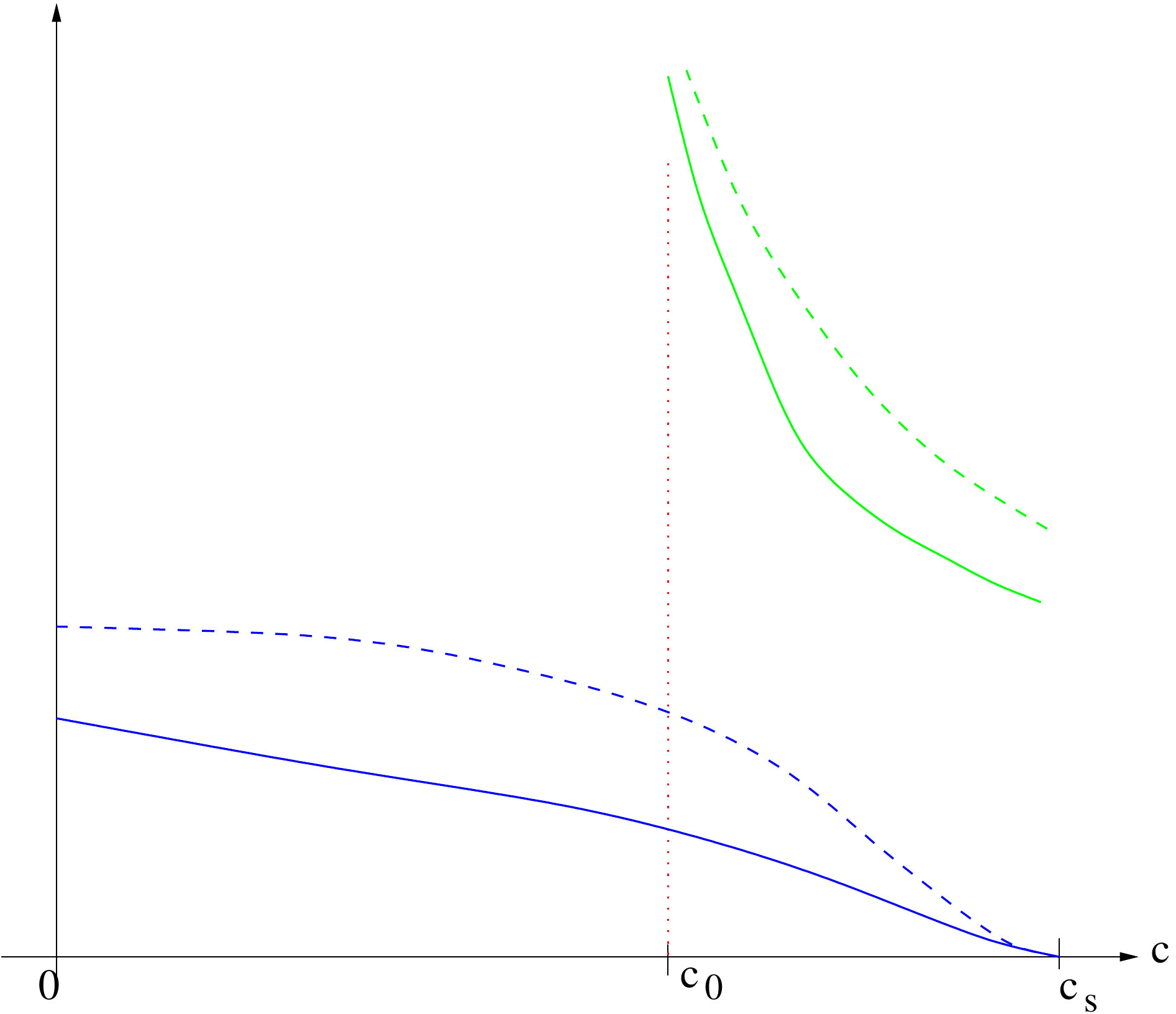} 
\hspace{1.5cm}  
\includegraphics[viewport=10 0 500 450,width=5cm,height=6cm]{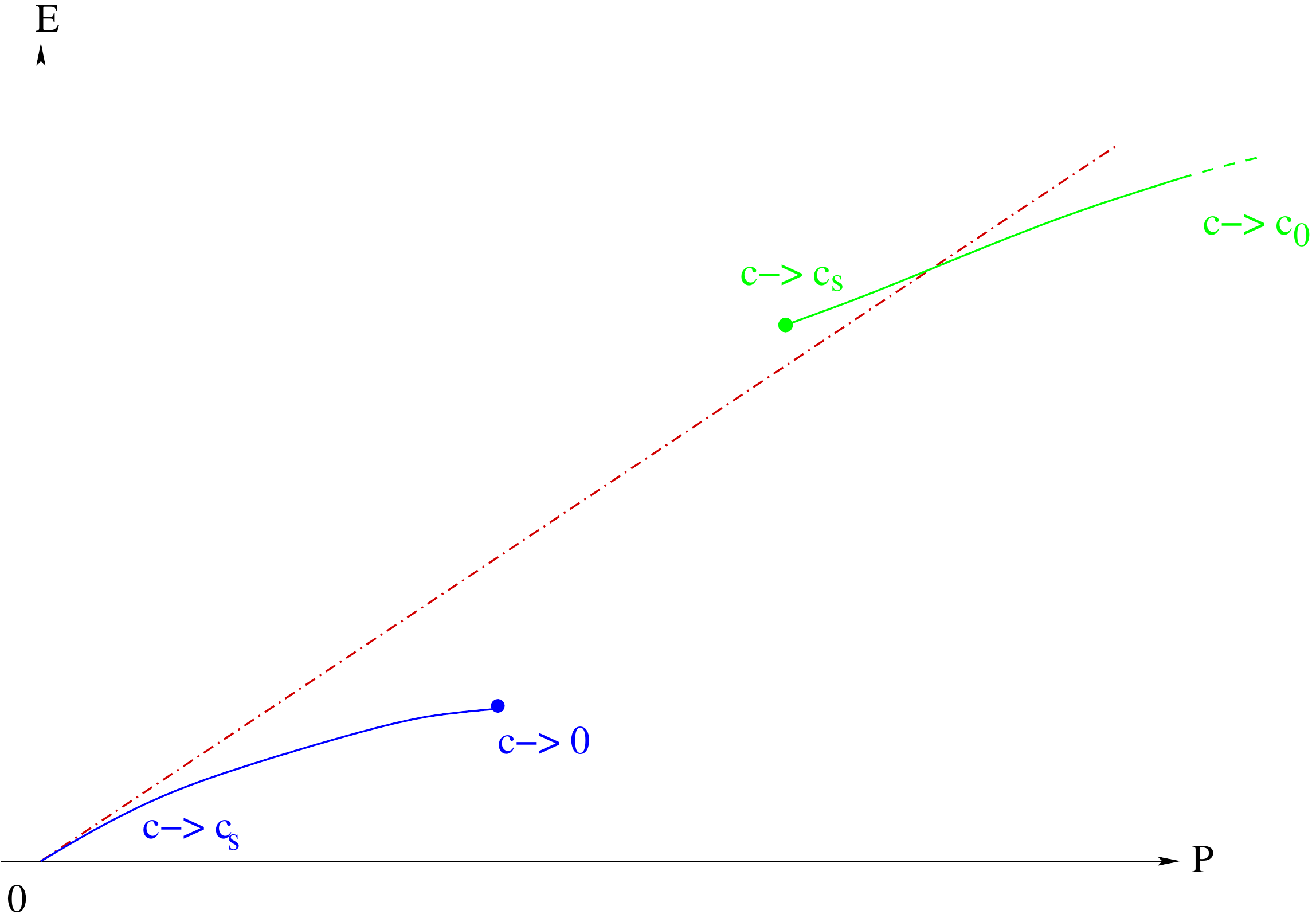}
\end{center}
\caption{(a) Energy (dashed curve) and momentum (full curve) vs. speed; 
(b) $(E,P)$ diagram}
\label{zsatur}
\end{figure}

\noindent $\bullet$ A saturated (NLS):
$ f(\varrho) \equiv \ds{\exp \Big( \frac{1 - \varrho}{\varrho_0} \Big) - 1 } $ 
with $\varrho_0 = 0.4$ (see figure \ref{zsatur}).

\begin{figure}[H]
\begin{center}
\includegraphics[viewport=10 10 500 450,width=6cm,height=7cm]{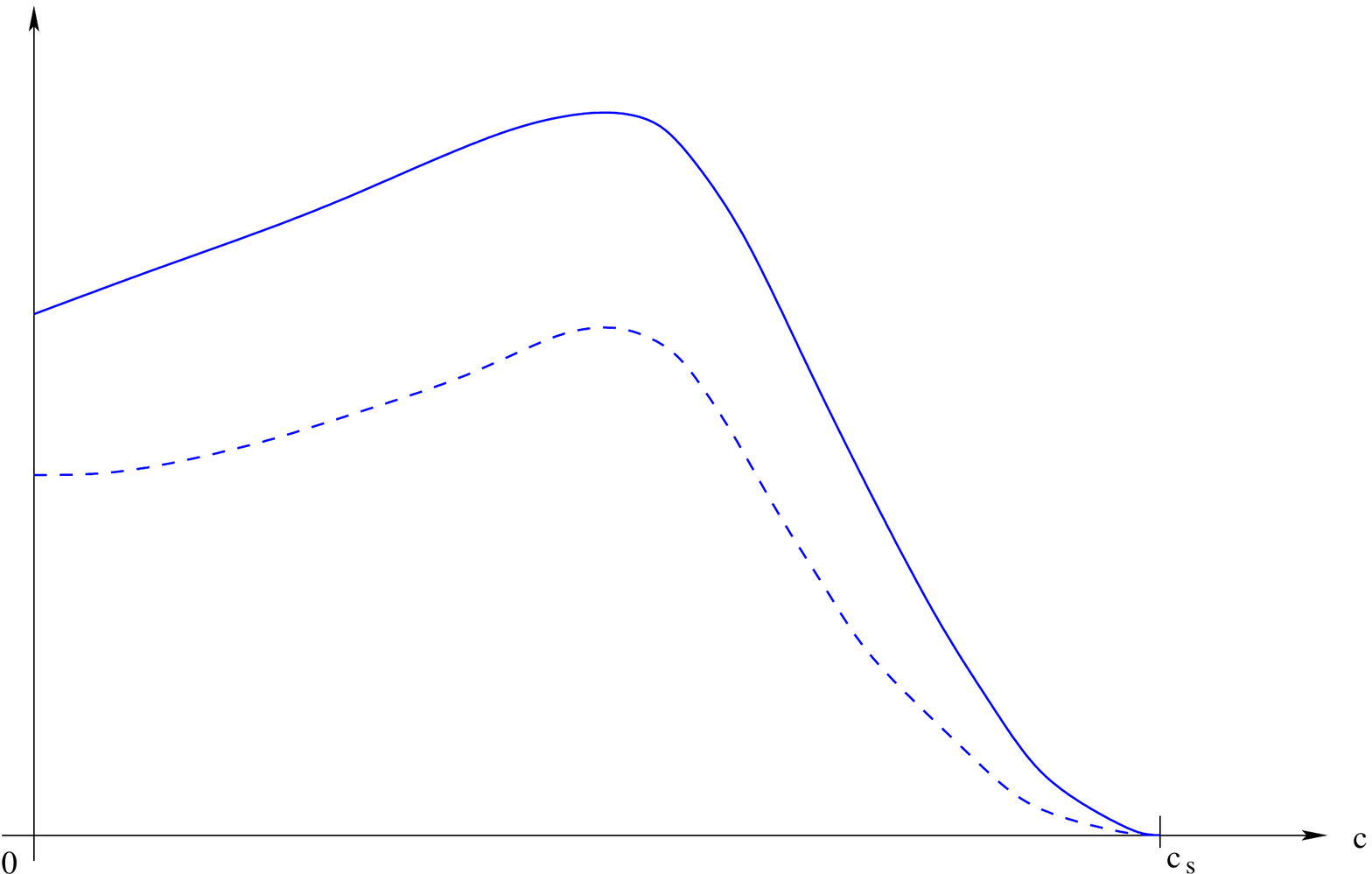} 
\hspace{2.5cm}  
\includegraphics[viewport=10 10 500 450,width=5cm,height=6cm]{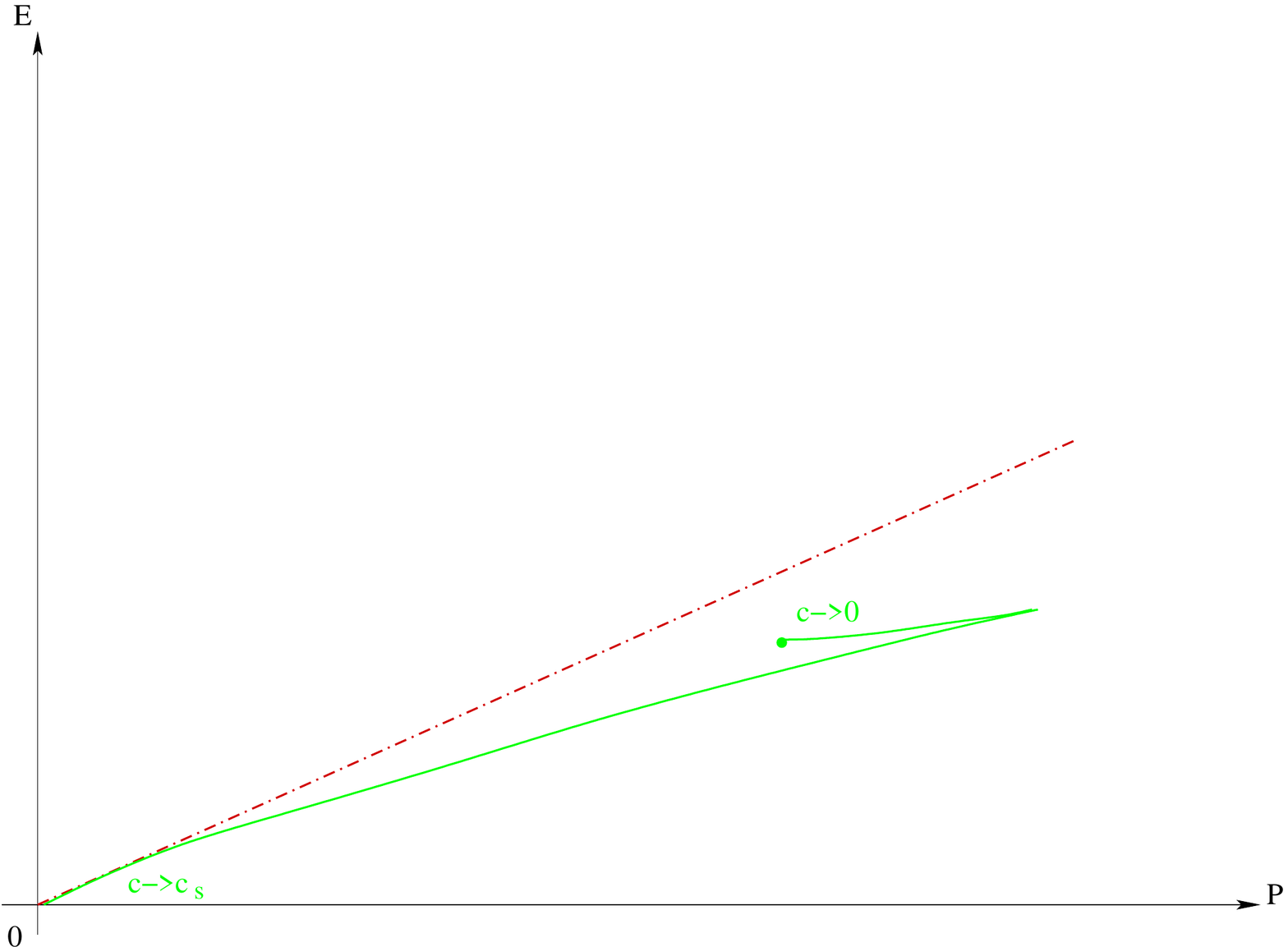}
\end{center}
\caption{(a) Energy (dashed curve) and momentum (full curve), 
(b) $(E,P)$ diagram}
\label{zsaturtrop}
\end{figure}

\noindent $\bullet$ Another saturated (NLS): 
$ f(\varrho) \equiv \ds{\frac{\varrho_0}{2} \Big( \frac{1}{(1+\varrho/\varrho_0)^2} 
- \frac{1}{(1+1/\varrho_0)^2} \Big) } $, with $ \varrho_0 = 0.08$ (see figure \ref{zsaturtrop}).

\begin{figure}[H]
\begin{center}
\includegraphics[viewport=10 10 500 450,width=6cm,height=7cm]{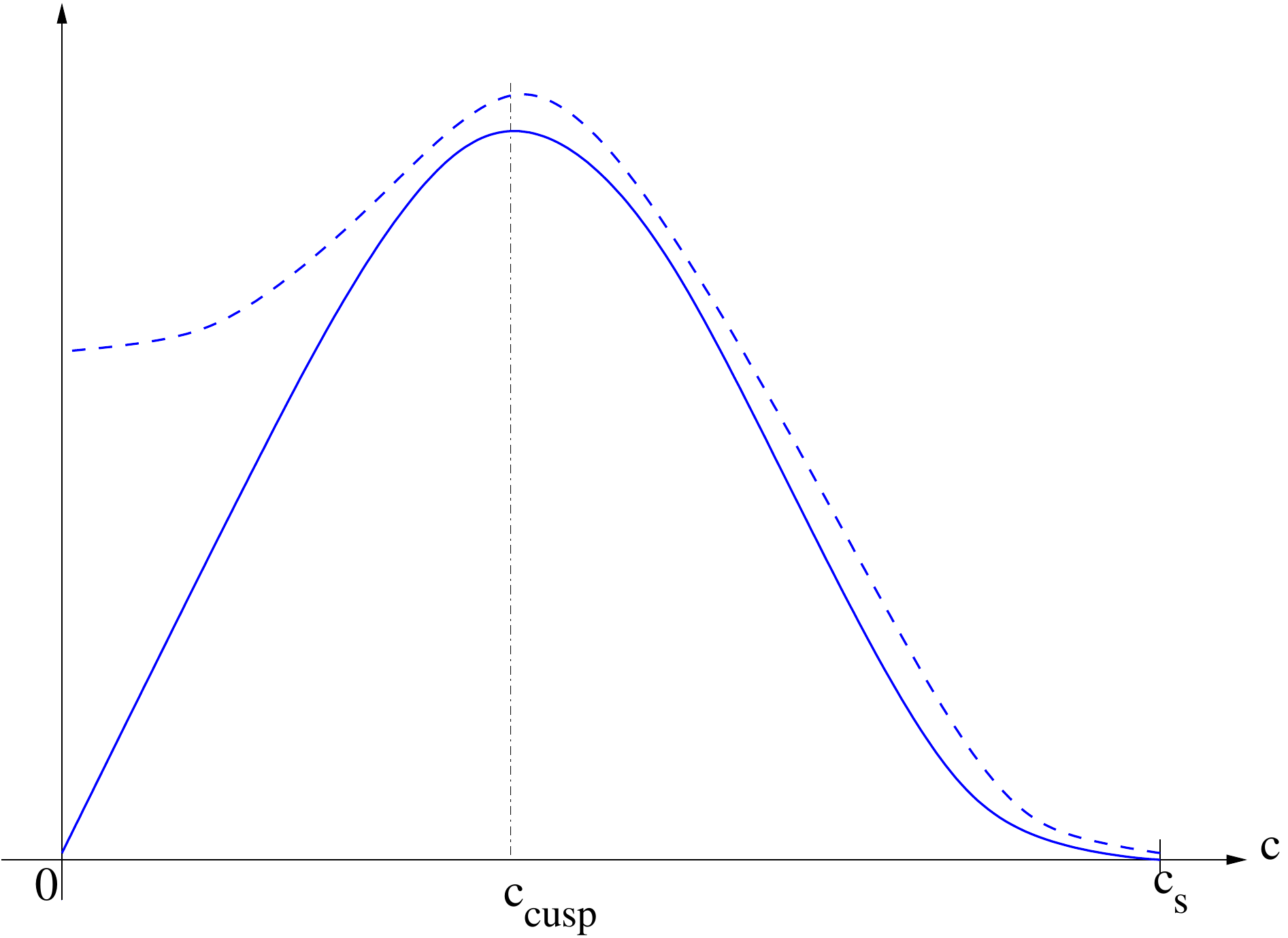} 
\hspace{1.5cm}  
\includegraphics[viewport=10 0 500 450,width=5cm,height=6cm]{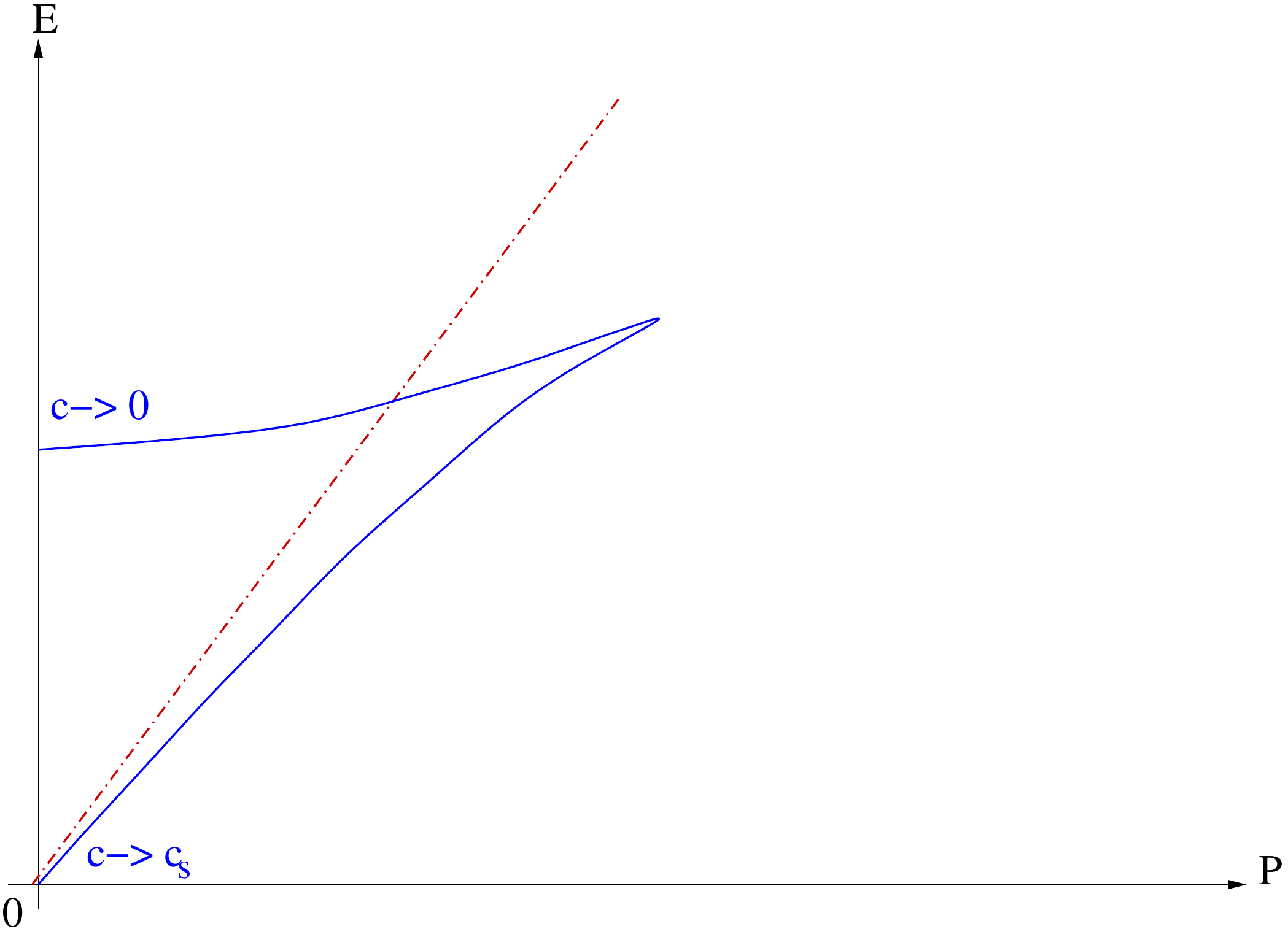}
\end{center}
\caption{(a) Energy (dashed curve) and momentum (full curve), 
(b) $(E,P)$ diagram}
\label{zcq}
\end{figure}

\noindent $\bullet$ The cubic-quintic nonlinearity: 
$ f(\varrho) \equiv - (\varrho-1) - 3 (\varrho-1)^2 $ (see figure \ref{zcq}).

Through the study (in \cite{C1d}) of these model cases, we have shown 
that if the energy-momentum diagram is well-known for the Gross-Pitaevskii 
equation, the qualitative properties of the travelling waves solutions 
can {\it not} be easily deduced from the global shape of the nonlinearity $f$. 
In particular, even if we restrict ourselves to smooth and decreasing 
nonlinearities (as is the Gross-Pitaevkii one), we see that we may have a 
great variety of behaviours: multiplicity of solutions; branches with diverging 
energy and momentum; nonexistence of travelling wave for some 
$ c_0 \in (0,\cs)$; branches of solutions that cross; existence of sonic 
travelling wave; transonic limit governed by the (mKdV) or more generally 
by the (gKdV) solitary wave equation instead of the usual (KdV) one; 
existence of cusps... \\

We investigate now the behaviour at infinity of the nontrivial travelling 
waves, which depend whether $c = \cs$ or not. We consider 
for $m \in \N $ the following assumption:\\

\noindent $ ( \BA_m ) $ $f$ is of class $\BC^{m+3}$ near $r_0^2$. Moreover, 
for $ 1 \leq j < m+ 2 $ , we have 
$ \ds{\frac{ f^{(j)} (r_0^2) }{(j+1)!} r_0^{2j}} = (-1)^{j+1} \ds{\frac{ \cs^2 }{4} } $ 
but $ \ds{\frac{ f^{(m+2)} (r_0^2) }{(m+3)!} r_0^{2(m+2)}} \not 
= (-1)^{m+3} \ds{\frac{ \cs^2 }{4} } $ (note that for $j=1$, equality always 
holds by definition of the speed of sound $ \cs = \sqrt{- 2 r_0^2 f'(r_0^2)}$). \\

%%%%%%%%%%%%%%%%%%%
\begin{proposition}
\label{decroiss} 
Let $ U_c \in \mathcal{Z} $ be a non constant travelling wave of speed 
$ 0 \leq c \leq \cs $. \\
% be fixed, and assume that $ U_c \in \mathcal{Z} $ 
% is a non constant travelling wave of speed $ c $.\\
%, with $ |U_c| \to r_0 $ as $ | x | \to + \ii $. \\
$(i)$ If $ c = 0 $, then there exist $ \phi_0 \in \R $ such that 
$ \ex^{i \phi_0} U_0 $ is a real-valued function and there exist 
two real constants $ M_0 \not = 0 $ (depending only on $f$ and $ \xi_0 $) 
and $x_0 $ such that, as $ x \to \pm \ii $,
$$
 \ex^{i \phi_0} U_0 (x) \mp r_0 \sim M_0 \exp \Big( - \cs | x - x_0 | \Big) 
\quad {\it if} \ \xi_0 = - r_0^2 , 
\quad \quad %\quad 
 \ex^{i \phi_0} U_0 (x) - r_0 \sim M_0 \exp \Big( - \cs | x - x_0 | \Big) 
\quad {\it if} \ \xi_0 \not = - r_0^2 .
$$
$ (ii) $ If $ 0 < c < \cs $, then $ U_c $ does not vanish, hence can 
be lifted: $ U_c = A_c \ex^{i \phi_c} $. Furthermore, there exist four 
real constants $ M_c $, $ \Theta_c $ (depending only on $f$, $c$ and 
$\xi_c$), $x_0 $ and $ \phi_0 $ such that, as $ x \to \pm \ii $,
$$
|U_c(x) |^2 - r_0^2 = \eta_{c} (x) \sim \frac{ 2 r_0^2 }{c} \p_x \phi (x) 
\sim M_c \exp \Big( - \sqrt{\cs^2 - c^2} | x - x_0 | \Big) ,
$$
and
$$
\phi (x) - \phi_0 \mp \Theta_c \sim 
- {\rm sgn} (x) \frac{c M_c }{2 r_0^2 \sqrt{\cs^2 - c^2} } \exp 
\Big( - \sqrt{\cs^2 - c^2} | x - x_0 | \Big) . 
$$
$ (iii) $ If $ c = \cs $ then $ U_{\cs} $ does not vanish, hence can 
be lifted: $ U_{\cs} = A_{\cs} \ex^{i \phi_{\cs}} $. We assume that there 
exists $m \in \N$ such that $ ( \BA_m ) $ is verified and denote
$$
\Lambda_m \equiv 
\frac{4}{r_0^{2(m+1)}} \Big[ \frac{ r_0^{2(m+2)} }{(m+3)!} f^{(m+2)} ( r_0^2) 
+ (-1)^{m+2} \frac{ \cs^2 }{4} \Big] \not = 0 .
$$
Then, we have, as $ x \to \pm \ii $,
$$
|U_{\cs}(x) |^2 - r_0^2 = \eta_{\cs}(x) \sim \frac{2r_0^2}{\cs} \p_x \phi (x) 
\sim {\rm sgn}(\xi_{\cs}) 
\Big( \frac{4}{(m+1)^2 | \Lambda_m | x^2 } \Big)^{\frac{1}{m+1}} $$
and
$$
\phi (x) \sim \frac{\cs {\rm sgn}(\xi_{\cs}) }{2r_0^2} 
\Big( \frac{4}{(m+1)^2 | \Lambda_m | } \Big)^{\frac{1}{m+1}} 
\left\{\begin{array}{ll}
\ds{ {\rm sgn}(x) \ln |x| } & \quad {\it if} \ m = 1 \\ \ \\ 
\ds{ \frac{m+1}{m-1} {\rm sgn}(x) |x|^\frac{m-1}{m+1} } & \quad {\it if} \ m \geq 2 ,
\end{array}\right.
$$
and if $m=0$, there exists $ \Theta_{\cs} \in \R $ and $ \phi_0 \in \R $ 
such that
$$
\phi (x) - \phi_0 \mp \Theta_{\cs} \sim {\rm sgn}(\xi_{\cs}) 
\frac{2 \cs}{ r_0^2 | \Lambda_0 | x} .
$$
In particular, since we impose $ U_{\cs} \in \mathcal{Z} $, we 
must have $ m \in \{ 0 , 1 , 2 \} $.
\end{proposition}
%%%%%%%%%%%%%%%%%

For the Gross-Pitaevkii nonlinearity ($f ( \varrho ) = 1 - \varrho $), 
we may compute explicitly the travelling waves for $0 < c < \cs = \sqrt{2}$ 
(see \cite{Tsu}, \cite{BGSsurvey})
$$
U_c (x) = \sqrt{\frac{2-c^2}{2}}\, {\rm tanh} 
\Big( x \frac{\sqrt{2-c^2}}{2} \Big) - i \frac{c}{\sqrt{2}} ,
$$
up to the invariances of the problem: translations and multiplications 
by a phase factor. On this explicit formula, the decay of the phase 
and modulus can be checked. In particular, as $ x \to \pm \ii $, we 
have
$$ 
U_c (x) \to \pm \ds{ \sqrt{1 - \frac{c^2}{\cs}} - i \frac{c}{\cs} } .
$$

%%%%%%%%%%%%%%
\begin{remark} \rm In the above statements, the constants $ \phi_0 $ 
and $ x_0 $ reflect the gauge and translation invariance. In the spirit 
of the model cases proposed in \cite{C1d}, for
$$
f ( \varrho ) \equiv - 2 ( \varrho - 1 ) + 3 ( \varrho - 1 )^2 
- 4 ( \varrho - 1 )^3 + 5 ( \varrho - 1 )^4 - 12 ( \varrho - 1 )^5 ,
$$
we obtain a smooth decreasing nonlinearity tending to $ - \ii $ 
at $+ \ii$ (thus qualitatively similar to the Gross-Pitaevskii nonlinearity) 
for which we have $ r_0 = 1$, $ \cs = 2 $, and 
$ \BV_{\cs} (\xi) = - 4\xi^4 - 8 \xi^5 $. For this nonlinearity $f$, 
there exists a nontrivial sonic travelling wave of infinite energy 
(corresponding to $\xi_{\cs} = - 1/2 $), since $m=3$.
\end{remark}
%%%%%%%%%%%%

The aim of this paper is to investigate the stability of the 
travelling waves for the one dimensional (NLS). We recall the definition 
of orbital stability in a metric space $(\BX,d_\BX)$ for which we have 
a local in time existence result.

%%%%%%%%%%%%%%%%%%
\begin{definition}
\label{defstabi} 
Let $ 0 \leq c \leq \cs $ and $ U_c \in \mathcal{Z} $ be 
a nontrivial travelling wave of speed $ c $. We say that $ U_c $ 
is orbitally stable in $(\BX,d_\BX)$, where $\BX \subset \BZ $, if 
for any $ \epsilon > 0 $, there exists $ \d > 0 $ such that for any 
initial datum $ \Psi^\ini \in \BX $ such that $ d_\BX( \Psi^\ini , U_c ) \leq \d $, 
any solution $ \Psi $ to {\rm (NLS)} with initial datum $ \Psi^\ini $ is global 
in $ \BX $ and
$$
\sup_{t\geq 0} \inf_{ \scriptsize{\begin{array}{c} y \in \R \\ \theta \in \R \end{array}}} 
d_\BX ( \Psi(t) , \ex^{i\theta} U_c (\cdot - y ) ) \leq \epsilon .
$$ 
\end{definition}
%%%%%%%%%%%%%%%%%%

In the sequel, $U_c$ will always stands for a nontrivial travelling 
wave, and we freeze the translation invariance by imposing that 
$|U_c|$ is even. Moreover, the solutions of (NLS) we consider will always 
be those given by Theorem \ref{Cauchy}.

%%%%%%%%%%%%%%%%%%%%%%%%%%%%%%%%%%%%%%%%%%%%%%%%%%%%%%%%%%%%%%%%%%%%%%%%
\subsection{Stability and instability in the case $ \bs{ 0 < c < \cs }$}
\label{sex0cs}

%%%%%%%%%%%%%%%%%%%%%%%%%%%%%%%%%%%%%%%%%%%%%%%%%%%%%%%%%%%%%%%%%%%%%%%%%
\subsubsection{Stability for the hydrodynamical and the energy distances}
\label{sex0cshydroener}

The first stability result for the travelling waves for (NLS) 
with nonzero condition at infinity is due to Z. Lin \cite{Lin}. 
The analysis relies on the hydrodynamical form of (NLS), which is 
valid for solutions that never vanish. The advantage is to work 
with a fixed functional space since 
$ (\eta , u ) = ( A^2 - r_0^2 , \p_x \phi ) \in H^1(\R) \times L^2 (\R) $, 
whereas the travelling waves have a limit $ r_0 \ex^{\pm i \Theta_c} $ 
(up to a phase factor) at $\pm \ii $ depending on the speed $c$. The 
result of Z. Lin \cite{Lin} establishes rigorously the stability criterion 
found in \cite{BKK}, \cite{Ba}.

%%%%%%%%%%%%%%%%%%%%%%%%%%%%
\begin{theorem} [\cite{Lin}]
\label{stab}
Assume that $ 0 < c_* < \cs $ is such that there exists a nontrivial 
travelling wave $ U_{c_*} $. Then, there exists some small $\s>0$ 
such that $U_{c_*}$ belongs to a locally unique continuous branch of 
nontrivial travelling waves $U_c$ defined for $ c_* - \s \leq c \leq c_* + \s$.\\
$(i)$ Assume
$$  \frac{dP(U_c)}{dc}_{|c=c_*} < 0 . $$
Then, $U_{c_*} = A_* \ex^{i\phi_*}$ is orbitally stable in the sense 
that for any $ \epsilon > 0$, there exists $\d>0$ such that if 
$ \Psi^\ini = A^\ini \ex^{i\phi^\ini} \in \mathcal{Z} $ 
verifies
$$ \n A^\ini - A_* \n_{H^1(\R)} + \n \p_x \phi^\ini - \p_x \phi_* \n_{L^2(\R)} 
\leq \d , $$
then the solution $\Psi$ to {\rm (NLS)} such that $\Psi_{|t=0} = \Psi^\ini$ 
never vanishes, can be lifted $\Psi = A \ex^{i \phi}$, and we have
$$ \sup_{t \geq 0} \inf_{y \in \R} 
\Big\{ 
\n A(t) - A_* ( \cdot - y ) \n_{H^1(\R)} 
+ \n \p_x \phi(t) - \p_x \phi_* ( \cdot - y ) \n_{L^2(\R)} \Big\} \leq \epsilon . $$
$(ii)$ Assume
$$  \frac{dP(U_c)}{dc}_{|c=c_*} > 0 . $$
Then, $U_{c_*} = A_* \ex^{i\phi_*}$ is orbitally unstable in the sense 
that there exists $ \epsilon > 0 $ such that, for any $\d>0$, there 
exists $ \Psi^\ini = A^\ini \ex^{i\phi^\ini} \in \mathcal{Z} $ verifying 
$$ \n A^\ini - A_* \n_{H^1(\R)} + \n \p_x \phi^\ini - \p_x \phi_* \n_{L^2(\R)} 
\leq \d , $$
but such that if $\Psi$ denotes the solution to {\rm (NLS)} with 
$\Psi_{|t=0} = \Psi^\ini $, then there exists $t > 0$ such that 
$\Psi$ does not vanish on the time interval $[0,t]$ but
$$ \inf_{y \in \R} \Big\{ \n A(t) - A_* ( \cdot - y ) \n_{H^1(\R)} 
+ \n \p_x \phi(t) - \p_x \phi_* ( \cdot - y ) \n_{L^2(\R)} \Big\} 
\geq \epsilon . $$
\end{theorem}
%%%%%%%%%%%%%

By one dimensional Sobolev embedding $ H^1(\R) \hookrightarrow \BC_0(\R) $, 
it is clear that since $ U_{c_*} $ does not vanish in $ \R $, by imposing 
$ |\!|\, | \Psi^\ini | - | U_{c_*} | \, |\!|_{H^1(\R)} 
= |\!| A^\ini - A_* |\!|_{H^1(\R)}$ small, $ \Psi^\ini $ does not 
vanish in $\R $ and thus can be lifted.

%%%%%%%%%%%%%%
\begin{remark} \rm We point out that C. Gallo in \cite{GaZhi} fills 
two gaps in the proof of \cite{Lin}: the first one concerns the local 
in time existence for the hydrodynamical system (see \eqref{hydrohamilto} 
in section \ref{sexremarkpf}) and the second one is about the conservation 
of the energy and the momentum. Furthermore, we make two additional remarks 
on the proof of \cite{Lin} in section \ref{sexremarkpf}.
\end{remark}
%%%%%%%%%%%%

Theorem \ref{stab} is stability or instability in the open set 
$ \mathcal{Z}_{\rm hy} \subset \mathcal{Z} $ for the hydrodynamical distance
$$ 
d_{\rm hy} ( \psi , \tilde{\psi} ) \equiv \n A - \tilde{A} \n_{H^1(\R)} 
+  \n \p_x \phi - \p_x \tilde{ \phi} \n_{L^2(\R)} 
+ \abs {\rm arg} \Big( \frac{\psi(0)}{\tilde{\psi}(0)} \Big) \abs , 
\quad \quad \quad \psi = A \ex^{i \phi } , \quad \tilde{\psi} 
= \tilde{A} \ex^{i \tilde{\phi} } ,
$$
which is not the energy distance. Here, $ {\rm arg} : \C^* \to ( - \pi , + \pi ] $ 
is the principal argument. For the stability, it suffices to consider the 
phase $ \theta \in \R $ such that 
$ \ds{ \arg \Big( \frac{\Psi(t)}{ {\sf e}^{i \theta} U_{c_*} ( \cdot - y ) } \Big) } $ 
is zero at $ x=0 $, where $y$ is the translation parameter. For the 
instability, no matter what is the phase $ \theta \in \R $. The result 
of \cite{Lin} is based on the application of the Grillakis-Shatah-Strauss 
theory \cite{GSS}  (see also \cite{BSS}, \cite{SS}) to the hydrodynamical 
formulation of (NLS) (see section \ref{sexremarkpf}). One difficulty is 
to overcome the fact that the Hamiltonian operator $\p_x$ is not onto.

On the energy-momentum diagrams, the stability can be checked either 
on the graphs of $E$ and $P$ with respect to $c$, either on the 
concavity of the curve $P \mapsto E$. Indeed, we have seen 
in \cite{C1d} that the so called Hamilton group relation
$$ c = \frac{dE}{dP} , \quad \quad \quad {\rm or} 
\quad \quad \quad \frac{dE(U_c)}{dc} = c \frac{dP(U_c)}{dc} , $$
holds, where the derivative is computed on the local branch. 
Therefore,
$$ \frac{d^2E}{dP^2} = \frac{d}{dP} \frac{dE}{dP} = \frac{dc }{dP} . $$
This means that we have stability when $ P \mapsto E $ is 
concave, that is $ \ds{\frac{d^2E}{dP^2}} < 0 $, and instability if $ P \mapsto E $ 
is convex, {\it i.e.} $ \ds{ \frac{d^2E}{dP^2} > 0 }$.

Actually, the proof of \cite{GSS}, \cite{Lin} provides an 
explicit control, as shown in the following lemma.

%%%%%%%%%%%%%%
\begin{lemma} \label{souskontrol} Under the assumptions of Theorem \ref{stab} 
and in the case $(i)$ of stability, we have, provided $ d_{\rm hy} ( \Psi^\ini , U_{c_*} ) $ 
is small enough,
\begin{align}
\label{Kontr}
\sup_{t\geq 0} \inf_{y \in \R} 
\Big\{ \n A(t) - A_* ( \cdot - y ) \n_{H^1(\R)} 
+ & \, \n \p_x \phi(t) - \p_x \phi_* ( \cdot - y ) \n_{L^2(\R)} \Big\} 
\nonumber \\
& \, 
\leq K \sqrt{ | E(\Psi^\ini) - E ( U_{c_*} ) | + | P(\Psi^\ini) - P ( U_{c_*} ) | } ,
\end{align}
as well as the control
\be
\label{Kontrlin}
\sup_{t\geq 0} \inf_{ \scriptsize{\begin{array}{c} y \in \R \\ \theta \in \R \end{array}}} 
d_{\rm hy} ( \Psi(t) , \ex^{i\theta} U_{c_*} (\cdot - y ) ) \leq K d_{\rm hy} ( \Psi^\ini , U_{c_*} ) .
\ee
\end{lemma}
%%%%%%%%%%%%

%%%%%%%%%%%%%%
\begin{remark} \rm The second estimate \eqref{Kontrlin} is not a simple 
consequence of the control \eqref{Kontr}, but relies on a comparison to 
$ U_c $ for some $c$ close to $c_*$ instead on a comparison to $ U_{c_*} $ 
(this idea has also been used in \cite{Wein}). It follows that, in the 
definition of stability for $ U_ {c_*} $, one can take $ \d = \BO( \e )$.
\end{remark}
%%%%%%%%%%%%

Let us lay the emphasis on the fact that Theorem \ref{stab} of Z. Lin is 
given in the hydrodynamical distance $ d_{\rm hy} $, which is not the energy distance 
$ d_\mathcal{Z} $. As a matter of fact, the Madelung transform
$$
\mathscr{M} : ( \mathcal{Z}_{\rm hy} , d_{ \mathcal{Z}} ) \ni U \mapsto 
\Big( \eta , u , \frac{ U(0)}{|U(0)| } \Big) \in H^1(\R,\R) \times L^2(\R,\R) \times \SS^1 ,
$$
where $ U = A \ex^{i\phi} $, $ \eta = A^2 - r_0^2 $ and $u = \p_x \phi $ is 
not so well behaved.

%%%%%%%%%%%%%
\begin{lemma} 
\label{cestfaux} 
$(i)$ The mapping $ \mathscr{M} : ( \mathcal{Z}_{\rm hy} , d_{ \mathcal{Z}} ) 
\to H^1(\R,\R) \times L^2(\R,\R) \times \SS^1 $ is an homeomorphism.\\
$(ii)$ There exists $ \phi_* \in \BC^2 (\R, \R )$ such that $ \p_x \phi_* \in L^2(\R) $ 
and a sequence $ (\phi_n )_{n \geq 1 } $ of functions in $ H^1 (\R, \R) $ such that, 
when $ n \to + \ii $,
$$
0 < d_{\rm hy} ( \ex^{i \phi_*} , \ex^{i \phi_*} \ex^{i\phi_n} ) \to 0 
\quad \quad \quad but \quad \quad \quad
\frac{ d_\mathcal{Z} ( \ex^{i \phi_*} , \ex^{i \phi_*} \ex^{i\phi_n} ) }
{ d_{\rm hy} ( \ex^{i \phi_*} , \ex^{i \phi_*} \ex^{i\phi_n} ) } \to +\ii .
$$
\end{lemma}
%%%%%%%%%%

Therefore, $ \mathscr{M}^{-1} $ is not locally Lipschitz continuous in general. 
However, for the stability issues, we compare the $ d_{\mathcal{Z}} $ and the 
$ d_{\rm hy} $ distances to some fixed travelling wave $ U_* $, which enjoys some 
nice decay properties at infinity. Let us now stress the link between the two 
distances $d_{\rm hy}$ and $ d_{\mathcal{Z}} $ in this case.

%%%%%%%%%%%%%
\begin{lemma} 
\label{equivdistance}
Let $ 0 \leq c_* \leq \cs $ and assume that $ U_* \in \mathcal{Z} $ is a 
non constant travelling wave with speed $c_*$ that does not vanish. 
If $ c_* = \cs $, we further assume that assumption $(\BA_0)$ is verified. 
%for some $m \in \{ 0, 1 , 2 \}$, assumption $(\BA_m)$ is verified. 
Then, there exists some constants $K$ and $\d > 0 $, depending only on $U_*$, 
such that for any $ \psi \in \BZ $ verifying $d_{\mathcal{Z}} ( \psi , U_* ) \leq \d$, 
we have
$$
\frac{1}{K} \, d_{\rm hy} ( \psi , U_* ) \leq d_{\mathcal{Z}} ( \psi , U_* ) \leq 
K d_{\rm hy} ( \psi , U_* ) .
$$
\end{lemma}
%%%%%%%%%%

An immediate corollary of Lemma \ref{equivdistance} is that Theorem 
\ref{stab} is also a stability/instability result in the energy distance. 
If one wishes only a stability/instability result, it is sufficient to 
invoke the fact that the mapping $ \mathscr{M} $ is an homeomorphism. 
However, the use of Lemma \ref{equivdistance} provides a stronger explicit 
control similar to the one obtained in Remark \ref{souskontrol} (see \eqref{Kontr}). 
In particular, in the definition of stability for $ U_ {c_*} $ in 
$ ( \mathcal{Z}, d_\mathcal{Z}) $, one can take $ \d = \BO( \e )$.

%%%%%%%%%%%%%%%%%
\begin{corollary} 
\label{stabenergy}
Assume that $ 0 < c_* < \cs $ is such that there exists a nontrivial 
travelling wave $ U_{c_*} $. Then, there exists some small $\s>0$ 
such that $U_{c_*}$ belongs to a locally unique continuous branch of 
nontrivial travelling waves $U_c$ defined for $ c_* - \s \leq c \leq c_* + \s$.\\
$(i)$ If $ \ds{ \frac{dP(U_c)}{dc}_{|c=c_*} < 0 } $, 
%$$  \frac{dP(U_c)}{dc}_{|c=c_*} < 0 , $$
then $U_{c_*} = A_* \ex^{i\phi_*}$ is orbitally stable in 
$( \mathcal{Z} , d_{\mathcal{Z} } ) $. Furthermore, if $ \Psi (t) $ is 
the (global) solution to {\rm (NLS)} with initial datum $ \Psi^\ini $, 
then we have, for some constant $K$ depending only on $U_{c_*} $ and provided 
$ d_{\mathcal{Z} } ( \Psi^\ini , U_{ c_*} ) $ is sufficiently small,
$$
\sup_{t \geq 0} \inf_{ \scriptsize{\begin{array}{c} y \in \R \\ \theta \in \R \end{array}}} 
d_{\mathcal{Z} } ( \Psi(t) , \ex^{i\theta} U_{ c_*} (\cdot - y) ) \leq K 
\sqrt{ | E(\Psi^\ini) - E ( U_{c_*} ) | + | P(\Psi^\ini) - P ( U_{c_*} ) | } ,
$$
as well as the control
$$
 \sup_{t \geq 0} \inf_{ \scriptsize{\begin{array}{c} y \in \R \\ \theta \in \R \end{array}}} 
d_{\mathcal{Z} } ( \Psi(t) , \ex^{i\theta} U_{ c_*} (\cdot - y) ) \leq 
K d_{\mathcal{Z} } ( \Psi^\ini , U_{ c_*} ).
$$
$(ii)$ If $ \ds{ \frac{dP(U_c)}{dc}_{|c=c_*} > 0 } $, 
%$$ \frac{dP(U_c)}{dc}_{|c=c_*} > 0 , $$
then $U_{c_*} = A_* \ex^{i\phi_*}$ is orbitally unstable in 
$( \mathcal{Z} , d_{\mathcal{Z} } ) $.
\end{corollary}
%%%%%%%%%%%%%%%

% 
% %%%%%%%%%%%%%%
% \begin{remark} \rm 
% If one wishes only a stability result, it is sufficient to invoke 
% the fact that the Madelung mapping
% $$
% \mathfrak{M} : ( \mathcal{Z}_{\rm hy} , d_{ \mathcal{Z}} ) \ni U \mapsto 
% \Big( \eta , u , \frac{U(0)}{|U(0)| } \Big) \in H^1(\R,\R) \times L^2(\R,\R) \times \SS^1 
% $$
% is an homeomorphism, where $ U = A \ex^{i\phi} $, $ \eta = A^2 - r_0^2 $ 
% and $u = \p_x \phi $. However, the use of Lemma \ref{equivdistance} provides 
% a stronger explicit control similar to the one obtained (see \eqref{Kontr}) 
% in Remark \ref{souskontrol}. In particular, in the definition of stability 
% for $ U_ {c_*} $ in $ ( \mathcal{Z}, d_\mathcal{Z}) $, one can take $ \d = \BO( \e )$.
% \end{remark}
% %%%%%%%%%%%%

For the Gross-Pitaevskii nonlinearity ($f( \varrho ) = 1 - \varrho$), 
the stability (for the energy distance $d_\mathcal{Z}$) of the travelling 
waves with speed $ 0 < c < \cs $ was proved by F. B\'ethuel, P. Gravejat and 
J.-C. Saut in \cite{BGSsurvey} through the 
variational characterization that these solutions are minimizers of the energy 
under the constraint of fixed momentum. However, in view of the energy momentum 
diagrams in section \ref{sextravelo}, this constraint minimization 
approach can not be used in the general setting we consider here. Indeed, 
this method provides only stability, but there may exist unstable travelling 
waves. Moreover, it follows from the proof of Theorem \ref{stab} that stable 
waves are local minimizers of the energy at fixed momentum but not necessarily 
global minimizers. Finally, we emphasize that the spectral methods allow to 
derive an explicit (Lipschitz) control in case of stability.

%%%%%%%%%%%%%%%%%%%%%%%%%%%%%%%%%%%%%%%%%%%%%%%%%%%%
\subsubsection{Stability via a Liapounov functional}
\label{sex0csLiap}

Another way to prove the orbital stability is to find a Liapounov 
functional. By Liapounov functional, we mean a functional which is 
conserved by the (NLS) flow and for which the travelling wave $ U_c $ 
is a local minimum (for instance, a critical point with second 
derivative $ \geq \d $ Id for some $ \d > 0 $). Such a Liapounov 
functional always exists in the Grillakis-Shatah-Strauss theory 
when $ \ds{ \frac{dP(U_c)}{dc}_{|c=c_*} < 0 } $, as shown by Theorem A 
in Appendix A. Its direct application to our problem leads us to 
define the functional in $ \mathcal{Z}_{\rm hy} $
$$
 \mathscr{L} ( \psi ) \equiv E ( \psi ) - c_* P( \psi ) 
+ \frac{M}{2} \Big( P( \psi ) - P( U_{c_*} ) \Big)^2 ,
$$
where $M$ is some positive parameter. It turns out that $ \mathscr{L} $ 
is such a Liapounov functional when $M$ is sufficiently large. Since 
the proof relies on the Grillakis-Shatah-Strauss framework, we have to 
work in the hydrodynamical variables. However, by Lemma \ref{equivdistance}, 
we recover the case of the energy distance. 

%%%%%%%%%%%%%%%
\begin{theorem} 
\label{Liapoupou} 
Assume that for some $ c_* \in ( 0 ,\cs ) $ and $ \s > 0 $ small, 
$ (0, \cs )\supset [ c_* - \s , c_* + \s ] \ni c \mapsto U_c \in \mathcal{Z} $ 
is a continuous branch of nontrivial travelling waves with 
$ \ds{ \frac{dP(U_c)}{dc}_{|c=c_*} < 0 } $. 
%$$ \frac{dP(U_c)}{dc}_{|c=c_*} < 0 . $$
If
$$ M > \frac{1}{ - \ds{ \frac{dP(U_c)}{dc}_{|c=c_*} } } > 0 , $$
there exists $\epsilon > 0 $ and $K$, depending only on $ U_{c_*} $, such that 
for any $  \psi \in \mathcal{Z} $ with 
$ \ds{ \inf_{ \scriptsize{\begin{array}{c} y \in \R \\ \theta \in \R \end{array}}} }  
d_{\rm hy} ( \psi , \ex^{i \theta} U_{c_*} (\cdot -y ) ) \leq \epsilon $, we have
$$ 
\inf_{ \scriptsize{\begin{array}{c} y \in \R \\ \theta \in \R \end{array}}} 
d_{\rm hy}^2 ( \psi , \ex^{i \theta} U_{c_*} (\cdot -y ) ) 
\leq K \Big( \mathscr{L} ( \psi ) - \mathscr{L} ( U_{c_*} ) \Big)  
$$
and analogously with $ d_{\rm hy} $ replaced by $ d_\mathcal{Z} $. 
Consequently, $ U_{c_*} = A_* \ex^{i\phi_*} $ is orbitally stable in  
$( \mathcal{Z}_{\rm hy} , d_{\rm hy} ) $ and in  $( \mathcal{Z} , d_{\mathcal{Z} } ) $. 
Furthermore, if $ \Psi (t) $ is the (global) solution to {\rm (NLS)} with 
initial datum $ \Psi^\ini $, then we have
$$
\sup_{t \geq 0} \inf_{ \scriptsize{\begin{array}{c} y \in \R \\ \theta \in \R \end{array}}} 
d_{\rm hy} ( \Psi(t) , \ex^{i\theta} U_{ c_*} (\cdot - y) ) \leq K 
\sqrt{ \mathscr{L} ( \Psi^\ini ) - \mathscr{L} ( U_{c_*} ) } 
\leq K d_{\rm hy} ( \Psi^\ini , U_{ c_*} ) .
$$
provided $ d_{\rm hy} ( \Psi^\ini , U_{ c_*} ) $ is sufficiently small, 
and analogously with $ d_{\rm hy} $ replaced by $ d_\mathcal{Z} $.
\end{theorem}
%%%%%%%%%%%%%

For the travelling waves for (NLS) in dimension one, this type of 
Liapounov functional appears for the first time in the paper 
\cite{Ba} by I. Barashenkov. However, in \cite{Ba}, the problem 
is treated directly on the wave function $ \Psi $, whereas the 
correct proof holds on the hydrodynamical variables, in particular 
because of the gauge invariance $ ( \theta , \Psi ) \mapsto \ex^{i\theta} \Psi $. 
For instance, the work \cite{Ba} suggests that we have stability for 
$H^1$ perturbations, whereas it holds only for perturbations in the 
energy space. Finally, we fill some gaps in the proof of \cite{Ba}.

%%%%%%%%%%%%%%%%%%%%%%%%%%%%%%%%%%%%%%%%%%%%%%%%%%%%%%%%%%%%%%%%%%%%%%%
\subsubsection{Instability via the existence of an unstable eigenvalue}
\label{sex0cseigen}

In the Grillakis-Shatah-Strauss theory \cite{GSS}, the instability is not 
shown by proving the existence of a unstable eigenvalue for the linearized (NLS) 
and then a nonlinear instability result (see however \cite{GSS2} when the 
Hamiltonian skew-adjoint operator is onto). There exists, however, some 
general results that prove the existence of unstable eigenvalues. 
For the instability of bound states for (NLS) (and also for nonlinear Klein-Gordon 
equation), that is solutions of the form $ \ex^{i \oo t } U_\oo (x) $, the work 
\cite{Gri} by M. Grillakis shows that the condition $\ds{ \frac{d}{d\oo} \Big( 
\int_{ \R^d} |U_\oo|^2 \Big)_{|\oo=\oo_*} > 0 }$ is sufficient for the existence 
of such an unstable eigenvalue. However, the proof relies on the fact that the 
bound states are real valued functions (up to a phase factor) and it is 
not clear whether it extends to the case of travelling waves we 
are studying. Indeed, since we have to work in hydrodynamical variables in 
order to have a fixed functional space, the linearized operator does not have 
(for $ c \not = 0 $) the structure required for the application of \cite{Gri}. 
Another general result is due to O. Lopes \cite{Lopes} but it assumes 
that the linearized equation can be solved using a semigroup. This is not 
the case for our problem once it is written in hydrodynamical variables (see 
below). Finally, the paper \cite{Lin2} by Z. Lin proposes an alternative 
approach for the existence of unstable eigenvalues. The method has the advantage 
of allowing pseudo-differential equations (like the Benjamin-Ono equation). 
However, the results are given for three model equations involving a scalar 
unknown, and it is not clear whether the proof can be extended to the case 
of systems.

The linearization of (NLS) near the travelling wave $ U_{c_*} $ 
in the frame moving with speed $ c_* $ is
\be
\label{NLSlin}
i \frac{\p \psi}{\p t} - i c_* \p_x \psi + \p^2_x \psi 
+ \psi f(|U_{c_*}|^2 ) + 2 \langle \psi | U_{c_*} \rangle f'(|U_{c_*}|^2 ) U_{c_*} = 0 ,
\ee
and thus, searching for exponentially growing modes 
$ \psi (t,x) = \ex^{\lambda t} w(x) $ leads to the eigenvalue problem
\be
\label{pbvp}
i \lambda w - i c_* \p_x w + \p^2_x w 
+ w f(|U_{c_*}|^2 ) + 2 \langle w | U_{c_*} \rangle f'(|U_{c_*}|^2 ) U_{c_*} = 0 ,
\ee
with $ {\rm Re}(\lambda) > 0 $ and $ w \not = 0 $. For one dimensional 
problems, the linear instability is commonly shown through the use of 
Evans functions (see the classical paper \cite{PeWe} and also the review article 
\cite{San}). For our problem, we look for an unstable eigenvalue for the equation 
written in hydrodynamical variables, namely we look for exponentially 
growing solutions $ ( \underline{\eta}, \underline{u} ) $ of the linear problem 
(written in the moving frame)
\be
\label{EKlin}
\left\{\begin{array}{ll}
\ds{ \p_t \underline{\eta} - c_* \p_x \underline{\eta} 
+ 2 \p_{x} ( ( r_0^2 + \eta_*) \underline{u} + \underline{\eta} u_* ) } = 0 \\ \ \\ 
\ds{ \p_t \underline{u} - c_* \p_x \underline{u} + 2 \p_{x} ( u_* \underline{u} ) 
- \p_x (f' ( r_0^2 + \eta_* ) \underline{\eta} ) 
- \p_x \Big\{ \frac{1}{2 \sqrt{r_0^2 + \eta_*} } \p_{x}^2 \Big( 
\frac{ \underline{\eta} }{ \sqrt{r_0^2 + \eta_*} } \Big) 
- \frac{ \underline{\eta} \p_{x}^2 (\sqrt{r_0^2 + \eta_*}) }{2 (r_0^2 + \eta_*)^{3/2} } \Big\} } = 0 ,
\end{array}\right. 
\ee
where $( \eta_*, u_*)$ is the reference solution. The advantage is here again 
to work with a fixed functional space in variables $( \eta, u ) $.  Due to the term 
$ \ds{ \p_x \Big\{ \frac{1}{2 \sqrt{r_0^2 + \eta_*} } \p_{x}^2 \Big( 
\frac{ \underline{\eta} }{ \sqrt{r_0^2 + \eta_*} } \Big) \Big\} } $, this equation 
can not be solved using a semigroup, except in the trivial case where $ \eta_* $ 
is constant, hence the result of \cite{Lopes} does not apply. However, system 
\eqref{EKlin} is a particular case of the Euler-Korteweg system for capillary fluids 
(see \cite{Benzo} for a survey on this model). We may then use a linear instability 
result already shown for the Euler-Korteweg system with the Evans function 
method, as in \cite{Zu} by K. Zumbrun for a simplified system, 
and more recently in \cite{BeGa} by S. Benzoni-Gavage for the complete 
Euler-Korteweg system.

%%%%%%%%%%%%%%%
\begin{theorem} 
\label{vpinstable} 
Assume that for some $ c_* \in ( 0 ,\cs ) $ and $ \s > 0 $ small, 
$ (0, \cs )\supset [ c_* - \s , c_* + \s ] \ni c \mapsto U_c \in \mathcal{Z} $ 
is a continuous branch of nontrivial travelling waves with
$$ \frac{dP(U_c)}{dc}_{|c=c_*} > 0 . $$
Then, there exists exactly one unstable eigenvalue $ \gamma_0 \in \{ {\rm Re} > 0 \} $ 
for \eqref{pbvp} and $ \gamma_0 \in ( 0 , + \ii ) $, that is {\rm (NLS)} is (spectrally) 
linearly unstable.
\end{theorem}
%%%%%%%%%%%%%

Once we have shown the existence of an unstable eigenvalue for the 
linearized (NLS) equation \eqref{NLSlin}, we can prove a nonlinear 
instability result as in \cite{HPW}, \cite{dBbubble}. Note that here, we 
no longer work in the hydrodynamical variables, where the high order 
derivatives involve nonlinear terms, but on the semilinear (NLS) equation. 
%Note first that (NLS) is locally well-posed in 
%$ U_{c_*} + H^1 (\R,\C) $ since we have a semilinear problem in 
%dimension 1. Therefore, we can prove a nonlinear instability result 
%as in \cite{HPW}, \cite{dBbubble}. 

%%%%%%%%%%%%%%%%%
\begin{corollary} 
\label{instabH1} 
Under the assumptions of Theorem \ref{vpinstable}, $ U_{c_*} $ is unstable 
in $ U_{c_*} + H^1(\R , \C ) $ (endowed with the natural $H^1$ distance): 
there exists $ \epsilon $ such that for any $ \d > 0 $, there exists 
$ \Psi^\ini \in U_{c_*} + H^1 (\R) $ such that 
$ \n \Psi^\ini - U_{c_*} \n_{H^1(\R)} \leq \d $ but if 
$ \Psi \in U_{c_*} + \BC ( [0, T^*) , H^1(\R) ) $ denotes the maximal 
solution of {\rm (NLS)}, then there exists $ 0 < t < T^* $ such that 
$ \n \Psi (t) - U_{c_*} \n_{H^1(\R)} \geq \epsilon $.
\end{corollary}
%%%%%%%%%%%%%%%

Since the proof is very similar to the one in \cite{HPW}, \cite{dBbubble}, we 
omit it. We may actually prove a stronger instability result, since the above 
one is not proved by tracking the exponentially growing mode. In 
\cite{dMG}, a spectral mapping theorem is shown and used to show the 
nonlinear instability by tracking this exponentially growing mode, which 
is a natural mechanism of instability. In Appendix B, we show that this 
spectral mapping theorem holds for a wide class of Hamiltonian equation. 
The direct application of Corollary B.2 in Appendix B gives the following 
nonlinear instability result.

%%%%%%%%%%%%%%%%%
\begin{corollary} 
\label{instabamelioree} 
We make the assumptions of Theorem \ref{vpinstable}, so that there exists 
an unstable eigenmode $( \gamma_0 , w ) \in ( 0 , +\ii) \times H^1(\R) $, 
$ |\!| w |\!|_{H^1} = 1 $. There exists $ M > 0 $ such that 
for any solution $\psi \in \BC(\R_+,H^1(\R,\C))$ of the linearized equation 
\eqref{NLSlin}, we have the growth estimate of the semigroup
$$
\forall t \geq 0 \quad \quad \quad 
\n \psi(t) \n_{H^1 (\R)} \leq M \ex^{ \gamma_0 t } \n \psi(0) \n_{H^1 (\R)} .
$$
Moreover, $ U_{c_*} $ has also the following instability property: there 
exists $K > 0$, $\d> 0 $ and $\e_0 >0 $, such that for any $ 0 < \d < \d_0 $, 
the solution $\Psi(t) $ to {\rm (NLS)} with initial datum $ \Psi^\ini = U_{c_*} + \d w 
\in U_{c_*} + H^1 (\R) $ exists at least on $ [ 0 , \gamma_0^{-1} \ln( 2 \e_0/ \d ) ] $ 
and verifies
$$
 \n \Psi(t) - U_{c_*} - \d \ex^{\gamma_0 t} w \n_{H^1(\R)} \leq K \d^2 \ex^{ 2 \gamma_0 t} .
$$
In particular, for $ t =  \gamma_0^{-1} \ln( 2 \e_0/ \d ) $ and $ \epsilon \equiv \frac{\e_0}{K} $, 
we have
$$
 \inf_{ y \in \R } \n |\Psi (t)| - |U_{c_*} | ( \cdot - y ) \n_{L^2(\R)} \geq \epsilon 
\quad \quad \quad {\it and } \quad \quad \quad 
 \inf_{ y \in \R } \n |\Psi (t)| - |U_{c_*} | ( \cdot - y ) \n_{L^\ii(\R)} \geq \epsilon ,
$$
which implies
$$
 \inf_{ \scriptsize{\begin{array}{c} y \in \R \\ \theta \in \R \end{array}} } 
\n \Psi(t) - \ex^{i \theta} U_{c_*} ( \cdot - y ) \n_{H^1(\R)} \geq \epsilon $$
as well as
$$
 \inf_{ \scriptsize{\begin{array}{c} y \in \R \\ \theta \in \R \end{array}} } 
d_{\rm hy} ( \Psi(t) , \ex^{i \theta} U_{c_*} ( \cdot - y ) ) \geq \epsilon 
\quad \quad \quad  {\it and } \quad \quad \quad 
 \inf_{ \scriptsize{\begin{array}{c} y \in \R \\ \theta \in \R \end{array}} } 
d_\mathcal{Z} ( \Psi(t) , \ex^{i \theta} U_{c_*} ( \cdot - y ) ) \geq \epsilon .
$$
%$$
% \inf_{ y \in \R } \n |\Psi (t)| - |U_{c_*} | ( \cdot - y ) \n_{L^\ii(\R)} \geq \epsilon 
% \quad \quad \quad {\it and} \quad \quad \quad 
%  \inf_{ \scriptsize{\begin{array}{c} y \in \R \\ \theta \in \R \end{array}}} 
%  \n \p_x \Psi (t) -  \ex^{i \theta} \p_x U_{c_*} ( \cdot - y ) \n_{L^2(\R)} \geq \epsilon .
%$$
%
%$$
% \inf_{ \scriptsize{\begin{array}{c} y \in \R \\ \theta \in \R \end{array}}} 
% \n \Psi(t) - \ex^{i \theta} U_{c_*} ( \cdot - y ) \n_{H^1(\R)} \geq \frac{1}{K} \d \ex^{ \gamma_0 t} 
% - K \d^2 \ex^{ 2 \gamma_0 t} ,
%$$
% 
% $$
% \quad  \quad \quad {\it and} \quad  \quad \quad 
% \n \Psi(t) - U_{c_*} \n_{H^1(\R)}
%$$
%
%$ \epsilon $ such that 
%for any $ \d > 0 $, there exists $ \Psi^\ini \in U_{c_*} + H^1 (\R) $ 
%such that $ \n \Psi^\ini - U_{c_*} \n_{H^1(\R)} \leq \d $ but if 
%$ \Psi \in U_{c_*} + \BC ( [0, T^*) , H^1(\R) ) $ denotes the maximal 
%solution of {\rm (NLS)}, then there exists $ 0 < t < T^* $ such that 
%$ \n |\Psi (t)| - |U_{c_*} | \n_{L^2(\R)} \geq \epsilon $, 
%$ \n |\Psi (t)| - |U_{c_*} | \n_{L^\ii(\R)} \geq \epsilon $ and 
%$ \n \p_x \Psi (t) - \p_x U_{c_*} \n_{L^2(\R)} \geq \epsilon $.
\end{corollary}
%%%%%%%%%%%%%%%

With the above result, we then show the nonlinear instability also in the 
energy space, and thus recover the instability result of Z. Lin but 
this time by tracking the unstable growing mode.

%%%%%%%%%%%%%%%%%%%%%%%%%%%%%%%%%%%%%
\subsubsection{Instability at a cusp}
\label{sex0cscusp}

In this section, we investigate the question of stability in the degenerate 
case $ \ds{ \frac{dP}{dc} = 0 } $. In \cite{GSS} (see also \cite{GSS2}), a 
stability result for the wave of speed $c_*$ is shown when the action 
$ c \mapsto S(c) = E( U_c ) - c P ( U_c ) $ (on the local branch) 
is such that, for instance, $ \ds{ \frac{d^2 S}{dc^2} = - \frac{dP}{dc} } $ 
is positive for $ c \not = c_* $ but vanishes for $ c = c_* $. In the 
energy-momentum diagrams of section \ref{sextravelo}, the situation 
is different since $ \ds{ \frac{dP}{dc} } $ changes sign {\it at} 
the cusps, or, equivalently, the action 
$ c \mapsto S(c) = E( U_c ) - c P ( U_c ) $ (on the local branch) 
changes its concavity {\it at} the cusp. The work \cite{CoPe} by 
A. Comech and D. Pelinovsky shows that for the Nonlinear Schr\"odinger 
equation, a bound state associated with a cusp in the energy-charge diagram 
is unstable. The proof relies on a careful analysis of the linearized equation, 
which is spectrally stable, but linearly unstable (with polynomial growth 
for the linear problem). A similar technique was used by A. Comech, 
S. Cuccagna and D. Pelinovsky in \cite{CoCuPe} for the generalized 
Korteweg-de Vries equation. 
Then, M. Ohta in \cite{Ohta} also proves the nonlinear 
instability of these ``bound states'' using a Liapounov functional as 
in \cite{GSS}. However, in \cite{Ohta}, it is assumed that $ J = T'(0) $ and 
is onto, which are both not true here (and there are further restrictions due 
to the introduction of an intermediate Hilbert space). Recently, M. Maeda 
has extended in \cite{Maeda} the above instability result, removing some 
assumptions in \cite{Ohta}. We show the instability of travelling waves 
associated with a cusp in the energy-momentum diagram in the generic case 
where $ \ds{ \frac{d^2P}{dc^2} \not = 0 } $. Our approach follows the 
lines of \cite{Maeda}, but with some modifications since our problem 
does not fit exactly the general framework of this paper. In particular, 
we can not find naturally a space ``$Y$'', and some functions appearing 
in the proof do not lie in the range of the skew-adjoint operator $ \p_x $ 
involved in the Hamiltonian formalism. We overcome this difficulty using 
an approximation argument (similar to the one used in \cite{Lin}).

%%%%%%%%%%%%%%%
\begin{theorem} 
\label{stabcusp} 
Assume that for some $ c_* \in ( 0 ,\cs ) $ and $ \s > 0 $ small, 
$ (0, \cs )\supset [ c_* - \s , c_* + \s ] \ni c \mapsto U_c \in \mathcal{Z} $ 
is a continuous branch of nontrivial travelling waves with
$$ \frac{dP(U_c)}{dc}_{|c=c_*} = 0 \not = \frac{d^2P(U_c)}{dc^2}_{|c=c_*} , $$
and assume in addition that $f$ is of class $\BC^2$. Then, $ U_{c_*} $ is orbitally 
unstable in $ ( \mathcal{Z} , d_\mathcal{Z} ) $.
\end{theorem}
%%%%%%%%%%%%%

%%%%%%%%%%%%%%%%%%%%%%%%%%%%%%%%%%%%%%%%%%%%%%%%%%%
\subsection{Stability in the case $ \bs{ c = 0 } $}
\label{sex0}

%%%%%%%%%%%%%%%%%%%%%%%%%%%%%%%%%%%%%%%%
\subsubsection{Instability for the bubbles}
\label{sex0bubble}

When $ c=0 $, we have two types of stationary waves: 
the bubbles ,when $ \xi_0 > - r_0^2 $, are even functions 
(up to a translation) that do not vanish, and the kinks, when $ \xi_0 = - r_0^2 $, 
are odd functions (up to a translation). The instability 
of stationary bubbles has been shown by A. de Bouard 
\cite{dBbubble} (and is true even in higher dimension). The proof 
in \cite{dBbubble} relies on the proof of the existence of an 
unstable eigenvalue for the linearized (NLS), and then the proof 
of a nonlinear instability result. An alternative proof of the 
linear instability of the bubbles is given in \cite{PeKe} (Theorem 3.11 $(ii)$).

%%%%%%%%%%%%%%%%%%%%%%%%%%%%%%%%%
\begin{theorem} [\cite{dBbubble}]
\label{unstabdeBouard}
Assume that there exists a bubble, that is a nontrivial stationary 
($c=0$) wave $U_0$ which does not vanish. Then, $U_0$ 
is (linearly and nonlinearly) unstable in $ U_0 + H^1(\R) $ (endowed 
with the natural $H^1$ metric), that is there exists $ \epsilon $ such 
that for any $ \d > 0 $, there exists $ \Psi^\ini \in U_0 + H^1 (\R) $ 
such that $ \n \Psi^\ini - U_0 \n_{H^1(\R)} \leq \d $ but if 
$ \Psi \in U_0 + \BC ( [0, T^*) , H^1(\R) ) $ denotes the maximal 
solution of {\rm (NLS)}, then there exists $ 0 < t < T^* $ such that 
$ \n \Psi (t) - U_0 \n_{H^1(\R)} \geq \epsilon $.
\end{theorem}
%%%%%%%%%%%%%%%%%%%%%

Actually, in the same way that Corollary \ref{instabamelioree} is 
a better instability result than Corollary \ref{instabH1}, we have 
the following stronger instability result, which is a direct consequence 
of Corollary B.2 in Appendix B.

%%%%%%%%%%%%%%%%%%%
\begin{proposition} 
\label{unstabbubble}
Assume that there exists a bubble, that is a nontrivial stationary 
($c=0$) wave $U_0$ which does not vanish. Then, $U_0$ is (nonlinearly) 
unstable in $ U_0 + H^1(\R) $, $( \mathcal{Z} , d_{\mathcal{Z}} ) $ and 
$( \mathcal{Z}_{\rm hy} , d_{\rm hy} ) $ in the same sense as in Corollary 
\ref{instabamelioree}.
\end{proposition}
%%%%%%%%%%%%%%%%%

Finally, we would like to emphasize that we may recover the instability 
result for bubbles from the proof of Theorem \ref{stab}, relying on the 
hydrodynamical form of (NLS), which holds true here since bubbles do not 
vanish. Our result holds in the energy space and for the hydrodynamical 
distance.

%%%%%%%%%%%%%%%
\begin{theorem} 
\label{unstab}
Assume that there exists a bubble, that is a nontrivial stationary 
($c=0$) wave $U_0$ which does not vanish. Then, there exists 
some small $\s>0$ such that $U_{0}$ belongs to a locally unique 
continuous branch of nontrivial travelling waves $U_c$ defined 
for $ 0 \leq c \leq \s$. Then, $ c \mapsto P(U_c) $ has a derivative 
at $ c = 0$,
$$ \frac{dP(U_c)}{dc}_{|c=0} > 0 $$
and $U_0 = A_* \ex^{i\phi_*}$ is orbitally unstable for the distances 
$d_{\mathcal{Z}} $ and $ d_{\rm hy} $.
\end{theorem}
%%%%%%%%%%%%%

\noindent {\it Proof.} We give a proof based on the argument of Z. Lin 
\cite{Lin}, which is possible since $U_0$ is a bubble hence does not 
vanish and the spectral decomposition used in \cite{Lin} still holds when $c=0$. 
Moreover, it is clear that the mapping $ c \mapsto ( \eta_c, u_c ) \in H^1 \times L^2 $ 
is smooth up to $c=0$, using the uniform exponential decay at infinity 
near $c=0$ and arguing as in \cite{C1d}. Therefore, it suffices to show that 
$ \frac{dP(U_c)}{dc}_{|c=0} > 0 $. From the expression of the momentum 
given in subsection 1.2 in \cite{C1d}, we have, for $ 0 \leq c \leq \s $,
$$ P(U_c) = c \, {\rm sgn}(\xi_c) \int_0^{\xi_c} \frac{\xi^2}{r_0^2 + \xi} 
\frac{d\xi}{\sqrt{ - \BV_c (\xi) } } 
= c \Big| \int_0^{\xi_0} \frac{\xi^2}{r_0^2 + \xi} 
\frac{d\xi}{\sqrt{- \BV_{0} (\xi)} } \Big| + o(c) $$
since $ \xi_0 > - r_0^2 $. Indeed, we are allowed to pass to the limit 
in the integral once it is written with the change of variables 
$ \xi = t \xi_c $:
$$
\int_0^{\xi_c} \frac{\xi^2}{r_0^2 + \xi} 
\frac{d\xi}{ - \BV_c(\xi)}  
= 
\int_0^{1} \frac{\xi_c^3 t^2 }{r_0^2 + t \xi_c } 
\frac{d t}{\sqrt{ - \BV_c( t \xi_c )} } ,
$$
since $ \xi_0 > - r_0^2 $. Therefore,
$$ \frac{dP(U_c)}{dc}_{|c=0} = \Big| \int_0^{\xi_0} \frac{\xi^2}{r_0^2 + \xi} 
\frac{d\xi}{\sqrt{- \BV_0(\xi)} } \Big| > 0 $$
since $ \xi_0 \not = 0 $ ($U_0$ is not trivial). The conclusion 
follows then from the proof of Theorem \ref{stab}. \carre \\

When we know that $ \ds{ \frac{dP(U_c)}{dc}_{|c=0} > 0 } $, we may also 
use the Evans function as in Theorem \ref{vpinstable} to show the 
existence of an unstable eigenmode. However, due to the fact that 
the kink $ U_0 $ is real-valued, we can use the arguments 
in \cite{dBbubble}, \cite{PeKe}.

%%%%%%%%%%%%%%%%%%%%%%%%%%%%%%%%%%%%%%%%%%%%%
\subsubsection{Stability analysis for the kinks}
\label{sex0kink}

We now turn to the case of the kinks ($ \xi_0 = - r_0^2 $ and 
$ U_0 $ is odd up to a translation). Since $U_0$ vanishes at 
the origin, the hydrodynamical form of (NLS) can not be used. 
The stability of the kink as attracted several recent works. 
L. Di Menza and C. Gallo \cite{dMG} have investigated the linear 
stability through the Vakhitov-Kolokolov function {\rm VK}, defined by
$$ {\rm VK}(\lambda) \equiv \int_\R \Big( [- \p^2_x - f(U_0^2) - \lambda ]^{-1} 
(\p_x U_0) \Big) (\p_x U_0) \ dx , $$
where $U_0$ is the kink, for $ \lambda \in ( \lambda_* , 0 ) $ 
for some $ \lambda_* < 0 $. They show that the Vakhitov-Kolokolov 
function {\rm VK} has a limit VK$_0$ when $ \lambda \to 0^- $. If 
VK$_0 > 0 $, then the linearization of (NLS) around the kink 
has an unstable real positive eigenvalue. When VK$_0 < 0 $, 
the linearization of (NLS) around the kink has a spectrum 
included in $i \R$ (spectral stability). Note that the approach 
of \cite{Lin} (extending \cite{GSS}) does not give directly unstable 
eigenvalues in the case $ \ds{ \frac{dP}{dc} >0} $. Recently, the link 
between the quantity $ \ds{ \frac{dP}{dc} } $ and the 
sign of VK$_0$ has been given by D. Pelinovsky and 
P. Kevrekidis \cite{PeKe} (proof of Lemma 3.10 there, the factor $ \sqrt{2} $ 
coming from the coefficients of the (NLS) equation in \cite{PeKe}):
\be
\label{jolie}
 2 \sqrt{2} {\rm VK}_0 = \lim_{c \to 0} \frac{dP(U_c)}{dc} ;
\ee
and they also prove, in a different way from \cite{dMG} 
that we have spectral stability when $ \ds{ \lim_{c \to 0} \frac{dP}{dc} < 0 } $ 
and existence of an unstable eigenvalue (in $ \R_+^*$) 
if $ \ds{\lim_{c \to 0} \frac{dP}{dc} > 0 } $. It is shown in \cite{PeKe} 
that the limit $ \lim_{c \to 0} \frac{dP(U_c)}{dc} $ does exist. Actually, 
they prove that the function $ [ 0, c_0 ) \ni c \mapsto P(U_c) $ 
is of class $ \BC^1 $ and that the derivative at $c=0$ is also given by 
(see \eqref{jolie})
\begin{align}
\label{Angelina}
\lim_{c \to 0} \frac{dP(U_c)}{dc} = & \, 2 \sqrt{2} {\rm VK}_0 = 2 \sqrt{2} \lim_{\lambda \to 0^- }
\int_\R \Big( [- \p^2_x - f(U_0^2) - \lambda ]^{-1} (\p_x U_0) \Big) (\p_x U_0) \ dx 
\nonumber \\
= & \, 2 \sqrt{2} \int_\R  {\rm Im} \Big( \frac{\p U_c}{\p c}_{|c=0} \Big) \p_x U_0 \ dx .
\end{align}

Our next lemma gives an explicit formula of the expression \eqref{Angelina}, 
involving only the nonlinearity $f$. 

%%%%%%%%%%%%%
\begin{lemma}
\label{prolongement} Assume that $ U_0 $ is a kink. Then, there exists 
$ c_0 \in ( 0 , \cs ) $ such that $ U_0 $ belongs to the (locally) unique 
branch $ [ 0, c_0 ) \ni c \mapsto U_c \in \mathcal{Z} $. Moreover, 
$ P(U_c) \to r_0^2 \pi $ as $ c \to 0 $ and the continuous extension 
$ [ 0, c_0 ) \ni c \mapsto P(U_c) $ has a derivative at $ c= 0$ 
given by
$$
 \frac{d P(U_c)}{dc}_{|c=0} = 
- \frac{8 r_0^3}{ 3 \sqrt{F(0)} } 
+ \frac12 \int_0^{r_0^2} \frac{ (\varrho-r_0^2)^2 }{ \varrho^{3/2} } 
\Big( \frac{1}{ \sqrt{ F( \varrho) }} - \frac{1}{\sqrt{ F(0) } } 
\Big) \ d \varrho . 
$$
\end{lemma}
%%%%%%%%%%%

The advantage of the formula given in Lemma \ref{prolongement} 
compared to \eqref{Angelina} is that it allows a direct computation of 
$ \ds{ \frac{d P(U_c)}{dc}_{|c=0} } $ when $f$ is known, which does not 
require to compute numerically 
$U_0$ and $ \ds{\frac{\p U_c}{\p c}_{|c=0} } $. For instance, it is 
quite well adapted to the stability analysis as in \cite{FaCa}. 
Let us observe that it may happen that a kink is unstable (see \cite{KivKro}, 
\cite{dMG}).\\

In the case of linear instability, \cite{dMG} shows that then, 
nonlinear instability holds. Actually, C. Gallo and L. Di Menza 
prove in \cite{dMG} a stronger result, where they show that the 
$L^\ii$ norm (and not only the $ H^1 $ norm) does not remain small.

%%%%%%%%%%%%%%%%%%%%%%%%%%%%%%%%%
\begin{theorem} [\cite{dMG}]
\label{unstabdMG}
Assume that there exists a kink, that is a nontrivial stationary 
($c=0$) wave $U_0$ vanishing somewhere, and satisfying 
$ \ds{ïŸ\frac{dP(U_c)}{dc}_{|c = 0} > 0}ïŸ$. Then, $U_0$ 
is (linearly and nonlinearly) unstable in the sense that there 
exists $ \epsilon $ such that for any $ \d > 0 $, there exists 
$ \Psi^\ini \in U_0 + H^1 (\R) $ such that 
$ \n \Psi^\ini - U_0 \n_{H^1(\R)} \leq \d $ but if 
$ \Psi \in U_0 + \BC ( [0, T^*) , H^1(\R) ) $ denotes the maximal 
solution of {\rm (NLS)}, then there exists $ 0 < t < T^* $ such that 
$ \n \Psi (t) - U_0 \n_{L^\ii(\R)} \geq \epsilon $.
\end{theorem}
%%%%%%%%%%%%%%%%%%%%%

The proof in \cite{dMG} relies on the tracking of the exponentially 
growing eigenmode. One may actually improve slightly the result 
as this was done in Corollary \ref{instabamelioree}. As a matter 
of fact, this was the result in Theorem \ref{unstabdMG} that has 
motivated us for Corollary \ref{instabamelioree}. 

We focus now on the nonlinear stability issue when there is linear 
(spectral) stability, that is when $ \ds{ \frac{d P(U_c)}{dc}_{|c=0} < 0 } $. 
Concerning the Gross-Pitaevskii nonlinearity ($f(\varrho) = 1 - \varrho$), 
for which we have $ \ds{ \frac{d P(U_c)}{dc}_{|c=0} < 0 } $, 
we quote two papers on this question. The first one is the work of P. G\'erard 
and Z. Zhang \cite{GZ} where the stability is shown by inverse scattering, 
hence in a space of functions sufficiently decaying at infinity. The analysis 
then relies on the integrability of the one-dimensional (GP) equation. 
The other work is by F. B\'ethuel, P. Gravejat, J.-C. Saut and D. Smets 
\cite{BGSS}. They prove the orbital stability of the kink of the 
Gross-Pitaevskii equation by showing that the kink is a global minimizer of the 
energy under the constraint that a variant of the momentum is fixed 
(recall that the definition of the momentum has to be clarified for 
an arbitrary function in the energy space), and that the corresponding 
minimizing sequences are compact (up to space translations and phase 
factors). In this approach, it is crucial (see \cite{BGSsurvey}, \cite{BGSS}) 
that $ E_{\rm kink} < \cs P_{\rm kink} = \cs r_0^2 \pi $ in order to prevent 
the dichotomy case for the minimizing sequences. However, since the
energy of the kink is equal to
$$
E_{\rm kink} = 4 \int_{-r_0^2}^0 \frac{ F( r_0^2 + \xi) }{\sqrt{- \BV_0(\xi) }} \ d \xi 
= 2 \int_{-r_0^2}^0 \sqrt{\frac{ F( r_0^2 + \xi) }{r_0^2+\xi }} \ d \xi 
= 2 \int_0^{r_0^2} \sqrt{\frac{ F( \varrho) }{\varrho }} \ d \varrho ,
$$
whereas its momentum is always equal to $ r_0^2 \pi $, it is clear that 
the condition $ E_{\rm kink} < \cs P_{\rm kink} =\cs r_0^2 \pi $ does not 
hold in general, as shown in the following example.\\

\noindent {\bf Example.} For $ \kappa \geq 0 $, consider
$$ 
f( \varrho) \equiv 1 - \varrho + \kappa ( 1 - \varrho )^3 ,
$$
which is smooth and decrease to $-\ii$ as the Gross-Pitaevskii 
nonlinearity. We have $r_0 =1 $, $ \cs = \sqrt{2} $, 
$ F( \varrho) = ( 1 - \varrho )^2 / 2 + \kappa ( 1 - \varrho )^4 / 4 $ and
$$
 E_{\rm kink} = 2 \int_0^{r_0^2} \sqrt{\frac{ F( \varrho) }{\varrho }} \ d \varrho 
 = 2 \int_0^{r_0^2} \sqrt{\frac{ 2 ( 1 - \varrho )^2 + \kappa ( 1 - \varrho )^4 }{4 \varrho }} \ d \varrho 
 > \cs r_0^2 \pi =  \pi \sqrt{2}  
$$
for $\kappa$ large (the left-hand side tends to $+\ii$), and numerical 
computations show that it is the case for $ \kappa \geq 14 $. Furthermore, 
Lemma \ref{prolongement} gives
\be
\label{cestpaslemoment}
 \sqrt{F(0)} \frac{d P(U_c)}{dc}_{|c=0} = 
- \frac{8}{ 3 } 
+ \frac12 \int_0^{1} \frac{ (\varrho - 1)^2 }{ \varrho^{3/2} } 
\Big( \sqrt{ \frac{F(0)}{F( \varrho) }} - 1 \Big) \ d \varrho . 
\ee
Since $ \frac{F(0)}{F( \varrho) } = \ds{ \frac{2+\kappa}{ 2 (\varrho - 1)^2 + \kappa (\varrho - 1)^4} } $, 
it can be easily checked that the right-hand side of \eqref{cestpaslemoment} 
is a decreasing function of $ \kappa $ tending to
$$
- \frac{8}{ 3 } 
+ \frac12 \int_0^{1} \frac{ (\varrho - 1)^2 }{ \varrho^{3/2} } 
\Big( \frac{1}{( \varrho -1)^2 } - 1 \Big) \ d \varrho = - 1 
$$
when $ \kappa \to +\ii $ (by monotone convergence). In particular, for any 
$ \kappa \geq 0 $, we have $ \ds{\frac{dP(U_c)}{dc}_{|c=0} < 0 } $, that 
is the kink is always (linearly) stable. The energy-momentum diagram for 
this type of nonlinearity with $\kappa $ large is as in on the right of 
figure \ref{zzcqsII} (the left part correspond to $\kappa$ smaller).\\

In comparison with the constraint minimization approach as in 
\cite{BGSsurvey, BGSS}, which allows to establish a global minimization result, the 
spectral methods as in \cite{GSS, Lin} allow to put forward locally minimizing 
properties, which turn out to be useful for the stability analysis in dimension 1. 

In the stability analysis of the kink, one issue is the definition of the momentum $P$, 
which was up to now given only for maps in $\mathcal{Z}_{\rm hy} $, that is for 
maps that never vanish, but the kink vanishes at the origin. In \cite{BGSS}, the 
notion of momentum was extended to the whole energy space $ \mathcal{Z} $, 
hence including maps vanishing somewhere, as a quantity defined mod $2\pi$, 
and was called ``untwisted momentum''. This notion will be useful for our stability result.

%%%%%%%%%%%%%
\begin{lemma} [\cite{BGSS}] If $ \psi \in \mathcal{Z}$, the limit 
$$
\mathfrak{P} (\psi) \equiv \lim_{ R \to +\ii } \Big[ 
\int_{-R}^{+R} \langle i \psi | \p_x \psi \rangle \ dx 
- r_0^2 \Big( \arg(\psi(+R)) -  \arg(\psi(-R)) \Big) \Big]
$$
exists in $ \R / ( 2 \pi r_0^2 \Z ) $. The mapping 
$ \mathfrak{P} : \mathcal{Z} \to \R / ( 2 \pi r_0^2 \Z ) $ 
is continuous and if $ \psi \in \mathcal{Z}$ verifies $ \inf_\R |\psi| > 0 $ 
(i.e. $ \psi \in \mathcal{Z}_{\rm hy} $), then $ \mathfrak{P}(\psi) = P(\psi) \mod 2\pi r_0^2 $. 
Finally, if $ \Psi \in \BC ( [0,T) , \mathcal{Z} ) $ is a solution 
to {\rm (NLS)}, then $ \mathfrak{P} (\Psi(t) ) $ does not depend on $t$.
\end{lemma}
%%%%%%%%%%%

\noindent {\it Proof.} For sake of completeness, we recall the 
proof of \cite{BGSS}. Let $ \psi \in \mathcal{Z} $ and let us verify 
the Cauchy criterion. Since $ |\psi| \to r_0 > 0 $ at $ \pm \ii $, 
we may lift $ \psi = A_\pm \ex^{ i \phi_\pm } $ in $ ( -\ii , - R_0 )$ 
and in $( + R_0 , + \ii )$ for some $ R_0 $ sufficiently large. For  
$ R' > R > R_0 $, we thus have in $ \R / ( 2 \pi r_0^2 \Z ) $
\begin{align*}
\Big[ 
\int_{-R'}^{+R'} \langle i \psi | \p_x \psi \rangle \ dx 
- & \, r_0^2 \Big( \arg(\psi(+R')) -  \arg(\psi(-R')) \Big) \Big] \\
& \quad \quad - \Big[ 
\int_{-R}^{+R} \langle i \psi | \p_x \psi \rangle \ dx 
- r_0^2 \Big( \arg(\psi(+R)) -  \arg(\psi(-R)) \Big) \Big] \\ 
& = \int_{R}^{R'} \langle i \psi | \p_x \psi \rangle \ dx 
+ \int_{- R'}^{- R} \langle i \psi | \p_x \psi \rangle \ dx 
\\ & \quad \quad \quad 
- r_0^2 \Big( \arg(\psi(R')) -  \arg(\psi(R)) \Big) 
+ r_0^2 \Big( \arg(\psi( - R')) -  \arg(\psi(-R)) \Big) \\
& = \int_{R}^{R'} A_+^2 \p_x \phi_+ \ dx + \int_{- R'}^{- R} A_-^2 \p_x \phi_- \ dx 
\\ & \quad \quad \quad 
- r_0^2 \Big( \phi_+(R') - \phi_+(R) \Big) 
+ r_0^2 \Big( \phi_-(- R') - \phi_-(-R) \Big) \\ 
& = \int_{R}^{R'} ( A_+^2 - r_0^2 ) \p_x \phi_+ \ dx 
+ \int_{- R'}^{- R} ( A_-^2 - r_0^2 ) \p_x \phi_- \ dx .
\end{align*}
The absolute value of each term is $ \leq K 
\int_{ \pm x \geq \pm R} |\p_x \psi|^2 + (|\psi| - r_0)^2 \ dx $ thus 
tends to zero if $ R \to + \ii $. Thus, $ \mathfrak{P} (\psi) $ is 
well-defined. 
The proof of the continuity follows the same lines, and allows to show 
that $ \mathfrak{P} $ is actually locally Lipschitz continuous. Let 
$ \psi \in \mathcal{Z} $. If $ \tilde{\psi} \in \mathcal{Z} $ and 
$ d_\mathcal{Z} ( \tilde{\psi} , \psi ) $ is small enough, we have 
$ |\!| |\tilde{\psi}| - | \psi | |\!|_{L^\ii} $ as small as we want. In 
particular, if $ R_0 > 0 $ is large enough so that $ |\psi| \geq 3r_0 / 4$ 
for $|x| \geq R $, we have $ |\tilde{\psi}| \geq r_0 / 2 $ for $|x| \geq R_0 $. 
As a consequence, writing $ \psi = A_\pm \ex^{ i \phi_\pm } $ 
and $ \tilde{\psi} = \tilde{A}_\pm \ex^{ i \tilde{\phi}_\pm } $ in 
$ ( -\ii , -R_0 )$ and in $ ( + R_0 , + \ii )$, we have, in $\R / (2 \pi r_0^2 \Z ) $ 
and for $ R > R_0 $,
\begin{align*}
\Big[ 
\int_{-R}^{+R} & \, \langle i \psi | \p_x \psi \rangle \ dx 
- r_0^2 \Big( \arg(\psi(+R)) - \arg(\psi(-R)) \Big) \Big] \\
& \quad \quad - \Big[ 
\int_{-R}^{+R} \langle i \tilde{\psi} | \p_x \tilde{\psi} \rangle \ dx 
- r_0^2 \Big( \arg(\tilde{\psi}(+R)) -  \arg(\tilde{\psi}(-R)) \Big) \Big] \\ 
& = \int_{-R_0}^{+R_0} \langle i (\psi - \tilde{\psi} ) | \p_x \psi \rangle 
+ \langle i \tilde{\psi} | \p_x ( \psi - \tilde{\psi} ) \rangle \ dx 
\\ & \quad \quad \quad 
+ \int_{ R_0 }^R A_+^2 \p_x \phi_+ - \tilde{A}_+^2 \p_x \tilde{\phi}_+ \ dx 
- r_0^2 \Big( \phi_+(+R) - \tilde{\phi}_+(+R) \Big) 
\\ & \quad \quad \quad 
+ \int_{-R}^{ - R_0 } A_-^2 \p_x \phi_- - \tilde{A}_-^2 \p_x \tilde{\phi}_- \ dx 
+ r_0^2 \Big( \phi_-(- R ) - \tilde{\phi}_-(-R) \Big) \\
& = \int_{-R_0}^{+R_0} \langle i (\psi - \tilde{\psi} ) | \p_x \psi \rangle 
+ \langle i \tilde{\psi} | \p_x ( \psi - \tilde{\psi} ) \rangle \ dx 
\\ & \quad \quad \quad 
+ r_0^2 \Big( \phi_+(+R_0) - \tilde{\phi}_+(+R_0) \Big) 
+ r_0^2 \Big( \phi_-(-R_0) - \tilde{\phi}_-(-R_0) \Big) 
\\ & \quad \quad \quad 
+ \int_{ R_0 }^R ( A_+^2 -r_0^2 ) \p_x \phi_+ - ( \tilde{A}_+^2 -r_0^2 ) \p_x \tilde{\phi}_+ \ dx 
+ \int_{ - R}^{- R_0 } ( A_-^2 -r_0^2 ) \p_x \phi_- - ( \tilde{A}_-^2 -r_0^2 ) \p_x \tilde{\phi}_- \ dx .
\end{align*}
We now estimate all the terms. For the last line, we use Cauchy-Schwarz to get
$ | \int_{ R_0 }^R ( A_+^2 - r_0^2 ) \p_x \phi_+ \ dx | \leq K (\psi) 
|\!| A_+ - r_0 |\!|_{L^2(\R)} |\!| A_+ \p_x \phi_+ |\!|_{L^2(\R)} 
\leq K(\psi) d_\mathcal{Z} (\psi , \tilde{\psi} )$, and similarly for the other terms. 
Moreover, using that $ ( \psi - \tilde{\psi} ) (x) = (\psi - \tilde{\psi})(0) 
+ \int_0^x \p_x ( \psi - \tilde{\psi} ) $, we get by Cauchy-Schwarz 
$ |\!| \psi - \tilde{\psi} |\!|_{\BC^0([-R_0,+R_0]) } 
\leq | (\psi - \tilde{\psi})(0) | + \sqrt{ R_0} |\!| \p_x \psi - \p_x \tilde{\psi} |\!|_{L^2(\R)} 
\leq K(R_0) d_\mathcal{Z} (\psi , \tilde{\psi} )$. Thus, the terms 
of the second line can be estimated by $ K(\psi,R_0) d_\mathcal{Z} (\psi , \tilde{\psi} )$, 
and for those of the first line, they can also be bounded by 
$ \leq K(\psi,R_0) d_\mathcal{Z} (\psi , \tilde{\psi} ) $. Passing to the limit 
as $R \to + \ii $ then gives
$$
\Big| \mathfrak{P}(\psi) - \mathfrak{P}(\tilde{\psi}) \quad {\rm mod} \, 2 \pi r_0^2 \Big| 
\leq K(\psi,R_0) d_\mathcal{Z} (\psi , \tilde{\psi} ) .
$$
This completes the proof for the definition of $ \mathfrak{P} $. To show that 
$ \mathfrak{P} $ is constant under the (NLS) flow, we use that $\Psi \in \Psi(0) + \BC([0,T),H^1) $ 
and the approximation by smoother solutions (see Proposition 1 in \cite{BGSS}). \carre \\

For the stability of the kink, we can no longer use the 
Grillakis-Shatah-Strauss theory applied to the hydrodynamical 
formulation of (NLS), since the kink vanishes at the origin. 
Therefore, it is natural to consider the Liapounov functional 
$ \mathscr{L} $ introduced in section \ref{sex0csLiap}, which 
becomes in the stationary case $ c = 0 $:
$$
\mathscr{L} (\psi) = E (\psi) + \frac{M}{2} \Big( P(\psi) - P(U_0) \Big)^2 .
$$
Since the momentum $P$ is not well-defined in $\mathcal{Z}$, we 
have to replace it by the untwisted momentum $ \mathfrak{P} $, which 
is defined modulo $ 2 \pi r_0^2 $. Consequently, it is natural to define 
the functional in $ \mathcal{Z} $
$$
\mathscr{K} (\psi) \equiv E (\psi) + 2M r_0^4 \sin^2 \Big( \frac{\mathfrak{P}(\psi) - r_0^2 \pi}{2r_0^2} \Big) ,
$$
which is well-defined and continuous in $ \mathcal{Z} $ since 
$ \sin^2 $ is $ \pi $-periodic. In addition, $ \mathscr{K} $ is conserved 
by the (NLS) flow as $E$ and $ \mathfrak{P} $.

%%%%%%%%%%%%%%%
\begin{theorem} 
\label{statiomini}
Assume that there exists a kink, that is a nontrivial stationary 
($c=0$) wave $U_0$ which is odd. Assume also that
$$ \frac{dP(U_c)}{dc}_{|c=0} < 0 . $$
Then, there exists some small $ \mu_* > 0 $ such that $ U_0 $ 
is a {\rm local} minimizer of $ \mathscr{K} $. More precisely, 
denoting 
$$
\mathscr{V}_{\mu_*} \equiv \{ \psi \in \mathcal{Z} , \ \inf_\R |\psi| < \mu_* \} ,
$$
we have, for any $ \psi \in \mathscr{V}_{\mu_*} \setminus \{ \ex^{i\theta} 
U_0( \cdot - y ) , \theta \in \R , y \in \R \} $,
$$
\mathscr{K} ( \psi ) > \mathscr{K} ( U_0 ) = E ( U_0 ) . 
$$
\end{theorem}
%%%%%%%%%%%%%

The crucial point in this result is to prove that the functional 
$ \mathscr{K} (\psi) $ controls the infimum $ \inf_\R |\psi| $. From 
this locally minimizing property of the kink when 
$ \ds{ \frac{dP(U_c)}{dc}_{|c=0} < 0 }$, we infer its orbital stability, 
provided we can prove some compactness on the minimizing sequences. 
Our method allows to infer a control on the distance of the solution 
to (NLS) to the orbit of the kink, but it is much weaker than those 
obtained by spectral methods in Lemma \ref{souskontrol} or 
Corollary \ref{stabenergy} for instance.

%%%%%%%%%%%%%%%
\begin{theorem}
 \label{stabmini} 
Assume that there exists a kink, that is an odd nontrivial stationary 
($c=0$) wave $U_0$ and that
$$ \frac{dP(U_c)}{dc}_{|c=0} < 0 . $$
Then, $U_0$ is orbitally stable in $ ( \mathcal{Z} , d_\mathcal{Z} ) $. 
Moreover, if $\Psi(t) $ is the (global) solution to {\rm (NLS)} with 
initial datum $ \Psi^\ini $, we have the control
$$
\sup_{t \geq 0} \inf_{ \scriptsize{\begin{array}{c} y \in \R \\ \theta \in \R \end{array}}} 
d_\mathcal{Z} ( \Psi(t) , \ex^{i\theta} U_{0} (\cdot - y) ) \leq 
K \sqrt[8]{\mathscr{K}(\Psi^\ini) - E(U_0)} \leq K \sqrt[4]{d_{\mathcal{Z}} ( \Psi^\ini , U_0 )} 
$$
provided that the right-hand side is sufficiently small.
\end{theorem}
%%%%%%%%%%%%%

This result settles the nonlinear stability under the condition 
$ \ds{ \frac{dP(U_c)}{dc}_{|c=0} < 0ïŸ} $ for a general nonlinearity $f$. 
In particular, it may be applied to the nonlinearity $f$ given in the 
example above. It shows that the stability of the kink holds with 
$ \d = \BO (\e^4 )$. We do not claim that the exponent $1/8$ is optimal.

\bigskip

For a complete study of the stability of the travelling waves, it would remain 
to investigate the case of the sonic ($ c = \cs $) travelling waves 
(when they exist). The methods we have developed do not apply directly, 
and we give in section \ref{sexsonic} some of the difficulties associated 
with this critical situation.

%%%%%%%%%%%%%%%%%%%%%%%%%%%%%%%%%%%%%%%%%%%%%%
\section{Decay at infinity (proof of Proposition \ref{decroiss})}

For simplicity, we shall denote
$$
\BV ( \xi ) \equiv \BV_{\cs}(\xi) = \cs^2 \xi^2 - 4 ( r_0^2+ \xi ) F ( r_0^2+ \xi ) .
$$
We freeze the invariance by translation by imposing $ |U_c| $ 
(hence also $\p_x \phi$) even, so that we can use the formulas in \cite{C1d}. 
In particular, it suffices to show the asymptotics for $ x \to + \ii $: 
the case $x \to - \ii$ follows by symmetry. 
We start with the proof of case $(iii)$. Under assumption 
$ ( \BA_m ) $ and since $ F' = - f $, we infer the Taylor expansion
\begin{align*}
\BV ( \xi ) = & \, \cs^2 \xi^2 + 4 ( r_0^2+ \xi ) \Big( \frac{1}{2!} f' ( r_0^2) \xi^2 
+ ... + \frac{1}{(m+2)!} f^{(m+1)} ( r_0^2) \xi^{m+2} 
+ \frac{1}{(m+3)!} f^{(m+2)} ( r_0^2) \xi^{m+3} + \BO( \xi^{m+4} ) \Big)
\\ = & \, 
\frac{4 r_0^2}{(m+3)!} f^{(m+2)} ( r_0^2) \xi^{m+3} 
+ \frac{4}{(m+2)!} f^{(m+1)} ( r_0^2) \xi^{m+3} 
+ \BO( \xi^{m+4} )
\\ = & \, 
\frac{4}{r_0^{2(m+1)}} \Big[ \frac{ r_0^{2(m+2)} }{(m+3)!} f^{(m+2)} ( r_0^2) 
+ (-1)^{m+2} \frac{ \cs^2 }{4} \Big] \xi^{m+3} 
+ \BO( \xi^{m+4} ) 
= \Lambda_m \xi^{m+3} 
+ \BO( \xi^{m+4} ) 
\end{align*}
since when $ ( \BA_m ) $ holds, all the terms $\BO(\xi^{m+2} )$ cancel out. 
The coefficient $ \Lambda_m$ is not zero by assumption. Note that the existence 
of a nontrivial sonic wave, which depends on the global behaviour of $\BV$, imposes that 
%either $m$ is even, either ($m$ is odd and 
%$ \ds{\frac{ f^{(m+2)} (r_0^2) }{(m+3)!} r_0^{2(m+2)}} < (-1)^{m+3} \ds{\frac{ \cs^2 }{4} } $), 
%and in all cases, we have 
$ \Lambda_m \xi^{m+3} < 0 $ when $ \xi $ is small and has the sign 
of $ \xi_{\cs} $. Therefore, from the formula (following from the Hamiltonian 
equation $ 2 \p^2_x \eta_c + \BV'(\eta_c) = 0 $, see \cite{C1d} for example)
\begin{align*}
x = - {\rm sgn}(\xi_{\cs}) \int_{\xi_{\cs}}^{\eta_{\cs}(x)} \frac{d\xi}{\sqrt{-\BV(\xi)}}
\end{align*}
and since there holds, as $\eta \to 0 $ (with the sign of $\xi_{\cs}$),
\begin{align*}
\int_{\xi_{\cs}}^{\eta} \frac{d\xi}{\sqrt{-\BV(\xi)}}
= & \, \int_{\xi_{\cs}}^{\eta} \frac{d\xi}{ \sqrt{- \Lambda_m \xi^{m+3 } }}  
+ \int_{\xi_{\cs}}^{\eta} \frac{\BV(\xi) - \Lambda_m \xi^{m+3 }}{ 
\sqrt{-\BV(\xi)} \sqrt{- \Lambda_m \xi^{m+3 } } \Big[\sqrt{-\BV(\xi)} + 
\sqrt{- \Lambda_m \xi^{m+3 } } \Big] } \ d\xi 
\\ = & \, - \frac{2 \ {\rm sgn}(\xi_{\cs}) }{m+1} \Big( \frac{1}{ \sqrt{- \Lambda_m \eta^{m+1 }}} 
- \frac{1}{ \sqrt{- \Lambda_m \xi_{\cs}^{m+1 }}} \Big) 
+ \left\{\begin{array}{ll}
\ds{ \BO(1) } & \quad {\rm if} \ m = 0 \\ \ \\ 
\ds{ \BO( | \ln | \eta | | ) } & \quad {\rm if} \ m = 1 \\ \ \\ 
\ds{ \BO ( \eta^{- \frac{m-1}{2}} ) } & \quad {\rm if} \ m \geq 2
\end{array}\right. 
\end{align*}
(here, we use that the last integrand is $ \BO(\xi^{-(m+1)/2}) $ as $\xi \to 0 $), 
it follows that, as $ x \to + \ii $,
$$
\eta_{\cs}(x) = {\rm sgn}(\xi_{\cs}) \Big( \frac{4}{(m+1)^2 |\Lambda_m| } \Big)^{\frac{1}{m+1}} 
\frac{1}{ x^{\frac{2}{m+1}} } 
+ \left\{\begin{array}{ll}
\ds{ \BO \Big( \frac{1}{x^3} \Big) } & \quad {\rm if} \ m = 0 \\ \ \\ 
\ds{ \BO \Big( \frac{ \ln (x) }{x^2} \Big) } & \quad {\rm if} \ m = 1 \\ \ \\ 
\ds{ \BO \Big( \frac{ 1}{x^\frac{4}{m+1} } \Big) } & \quad {\rm if} \ m \geq 2 .
\end{array}\right.
$$
This shows the asymptotics for the modulus, or $ \eta_{\cs} $. The 
asymptotic expansion for $ \p_x \phi_{\cs} $ is easily deduced from the 
equation on the phase $ 2 \p_x \phi_{\cs} = \cs \eta_{\cs} / ( r_0^2 + \eta_{\cs} ) $, 
and the phase $ \phi_{\cs} $ is then computed by integration, which completes 
the proof of case $(iii)$.

The proof of $(ii)$ is easier. Indeed, in this case, the function 
$\BV$ has the expansion
$$
\BV ( \xi ) = \cs^2 \xi^2 - 4 ( r_0^2+ \xi ) F ( r_0^2+ \xi ) = \BO ( \xi^3 ) ,
$$
hence
$$
\BV_c ( \xi ) = \BV ( \xi ) - ( \cs^2 - c^2) \xi^2 = - ( \cs^2 - c^2) \xi^2 + \BO ( \xi^3 ) . 
$$
As a consequence, the result follows from the expansion, for $ \eta \to 0 $
\begin{align*}
\int_{\xi_{c}}^{\eta} \frac{d\xi}{\sqrt{-\BV(\xi)}} = & \, 
\int_{\xi_{c}}^{\eta} \frac{d\xi}{ \sqrt{ ( \cs^2 - c^2) \xi^2 }}  
+ \int_{\xi_c}^{\eta} \frac{\BV(\xi) }{ 
\sqrt{-\BV_c(\xi)} \sqrt{ ( \cs^2 - c^2) \xi^2 } \Big[\sqrt{-\BV_c(\xi)} 
+ \sqrt{ ( \cs^2 - c^2) \xi^2 } \Big] } \ d\xi 
\\ = & \, 
{\rm sgn (\xi_c)} \frac{ \ln ( \eta / \xi_c ) }{ \sqrt{ \cs^2 - c^2 } } 
+ \int_{\xi_c}^{0} \frac{\BV(\xi) }{ 
\sqrt{-\BV_c(\xi)} \sqrt{ ( \cs^2 - c^2) \xi^2 } \Big[\sqrt{-\BV_c(\xi)} 
+ \sqrt{ ( \cs^2 - c^2) \xi^2 } \Big] } \ d\xi 
+ \BO (\eta) 
\end{align*}
since the integrand for the last integral is continuous at $ \xi = 0 $. 
This yields the desired expansion for the modulus:
\begin{align*}
\eta_c ( x) = & \, \xi_c \exp \Big( - x \sqrt{ \cs^2 - c^2 } 
- \int_{\xi_c}^{0} \frac{\BV(\xi) }{ 
\sqrt{- \xi^2 \BV_c(\xi)} \Big[\sqrt{-\BV_c(\xi)} 
+ \sqrt{ ( \cs^2 - c^2) \xi^2 } \Big] } \ d\xi \Big) 
+ \BO \Big[ \exp \Big( - 2 x \sqrt{ \cs^2 - c^2 } \Big) \Big] 
\\ = & \, 
M_c \exp \Big( - x \sqrt{ \cs^2 - c^2 } \Big) 
+ \BO \Big[ \exp \Big( - 2 x \sqrt{ \cs^2 - c^2 } \Big) \Big] ,
\end{align*}
with
$$
M_c \equiv \xi_c \exp \Big( - \int_{\xi_c}^{0} \frac{\BV(\xi) }{ 
\sqrt{- \xi^2 \BV_c(\xi)} \Big[\sqrt{-\BV_c(\xi)} 
+ \sqrt{ ( \cs^2 - c^2) \xi^2 } \Big] } \ d\xi \Big) \not = 0,
$$
and hence for the phase by similar computations to those above. 
The proof of case $(i)$ is similar, separating the case $\xi_0 = - r_0^2 $ of 
the kink (even solution) from the case $ \xi_0 \not = - r_0^2 $ of the bubble (odd 
solution) and is omitted. \carre

%%%%%%%%%%%%%%%%%%%%%%%%%%%%%%%%%%%%%%%%%%%%%%%%%%%%%%%%%%%%%%%%%%%%%%%%%%%%%%%%
\section{Stability results deduced from the hydrodynamical formulation of (NLS)}
\label{proofhydro}

%%%%%%%%%%%%%%%%%%%%%%%%%%%%%%%%%%%%%%%%%%%
\subsection{Proof of Lemma \ref{cestfaux}}
\label{sexcestfaux}

\noindent {\bf (i) The mapping $\bs{\mathscr{M}}$ is an homeomorphism.} 
Let $ \psi = A \ex^{i \phi} $, $ (\psi_n = A_n \ex^{i \phi_n} )_n \in \mathcal{Z} $ 
such that $ \psi_n \to \psi $ for $ d_{\rm hy} $. Then, $ A_n - A \to 0 $ in $H^1$, 
$ \p_x \phi_n \to \p_x \phi $ in $L^2$ and we may assume (possibly adding some 
multiple of $2\pi$ to $\phi_n$, that $ \phi_n (0) \to \phi (0) $. We write, using 
the embedding $H^1(\R) \hookrightarrow L^\ii(\R)$ for the before last line,
\begin{align}
\label{Lipsch}
 d_\mathcal{Z} (\psi_n , \psi ) = & \, 
\n \p_x \psi_n - \p_x \psi \n_{L^2} + \n | \psi_n | - | \psi | \n_{L^2} + | \psi_n(0) - \psi (0) | 
\nonumber \\ 
= & \, 
\n \ex^{i \phi_n} \p_x A_n + i A_n \ex^{i \phi_n} \p_x \phi_n 
- \ex^{i \phi} \p_x A - i A \ex^{i \phi} \p_x \phi \n_{L^2} + \n A_n - A \n_{L^2} 
+ | A_n(0) \ex^{i\phi_n(0)} - A(0) \ex^{i\phi(0)} | \nonumber  \\ 
\leq & \, 
\n ( \ex^{i \phi_n} - \ex^{i \phi}) \p_x A \n_{L^2} 
+ \n \ex^{i \phi} ( \p_x A_n - \p_x A ) \n_{L^2} 
+ \n ( A_n - A) \ex^{i \phi_n} \p_x \phi_n \n_{L^2} 
+ \n A ( \ex^{i \phi_n} - \ex^{i \phi} ) \p_x \phi \n_{L^2} \nonumber \\
& \quad + \n A \ex^{i \phi_n} ( \p_x \phi - \p_x \phi_n ) \n_{L^2} 
+ \n A_n - A \n_{L^2} 
+ | (A_n(0) - A(0) ) \ex^{i\phi_n(0)} | + | A(0) ( \ex^{i\phi_n(0)} - \ex^{i\phi(0)} ) | \nonumber \\ 
\leq & \, 
\n ( \ex^{i \phi_n} - \ex^{i \phi}) \p_x A \n_{L^2} 
+ K \n A_n - A \n_{H^1} 
+ \n ( A_n - A) \n_{H^1} \n \p_x \phi_n \n_{L^2} 
+ \n A \n_{L^\ii} \n ( \ex^{i \phi_n} - \ex^{i \phi} ) \p_x \phi \n_{L^2} \nonumber \\
& \quad + \n A \n_{L^\ii} \n \p_x \phi - \p_x \phi_n \n_{L^2} 
+ \n A \n_{L^\ii} | \ex^{i\phi_n(0)} - \ex^{i\phi(0)} | \\ 
= & \, \n ( \ex^{i \phi_n} - \ex^{i \phi}) \p_x A \n_{L^2} 
+ \n A \n_{L^\ii} \n ( \ex^{i \phi_n} - \ex^{i \phi} ) \p_x \phi \n_{L^2} + o_{n \to + \ii }(1), \nonumber 
\end{align}
from the convergences we have. Now observe that $ \phi_n(x) = 
\phi_n(0) + \int_0^x \p_x \phi_n(t) \ dt \to \phi(0) + \int_0^x \p_x \phi(t) \ dt = \phi(0) $ 
pointwise, hence, by the Dominated Convergence Theorem, 
$ |\!| ( \ex^{i \phi_n} - \ex^{i \phi}) \p_x A |\!|_{L^2} \to 0 $ and similarly 
for the other term. Therefore, $ d_\mathcal{Z} (\psi_n , \psi ) \to 0 $ as wished.

Let now $ \psi = A \ex^{i \phi} $, $ (\psi_n = A_n \ex^{i \phi_n} )_n \in \mathcal{Z} $ 
such that $ \psi_n \to \psi $ for $ d_\mathcal{Z} $. Then, $ A_n - A = 
| \psi_n | - | \psi | \to 0 $ in $L^2$, $ \p_x \psi_n \to \p_x \psi $ in $L^2$ and 
$ \psi_n (0) \to \psi (0) $. Since $ |\cdot | $ is $1$-Lipschitz continuous, we infer 
for the modulus
$$
 \n \p_x A_n - \p_x A \n_{L^2} = \n \p_x | \psi_n | - \p_x | \psi | \n_{L^2} 
\leq \n \p_x \psi_n - \p_x \psi \n_{L^2}  .
$$
Moreover, $ \psi_n(0) \to \psi(0) $ and this implies $ {\rm arg}( \psi_n(0) / \psi(0) ) \to 0 $. 
Therefore, it suffices to show that $ \p_x \phi_n \to \p_x \phi $ in $ L^2 $. We use 
the formula $ A^2 \p_x \phi = \langle i \psi |\p_x \psi \rangle $, which yields
$$
 \p_x \phi_n - \p_x \phi = 
 \frac{\langle i \psi_n |\p_x \psi_n \rangle}{A_n^2} - \frac{\langle i \psi |\p_x \psi \rangle}{A^2} 
 = \langle i \psi_n |\p_x \psi_n \rangle \Big( \frac{1}{A_n^2} - \frac{1}{A^2} \Big)
 + \frac{\langle i \psi_n |\p_x (\psi - \psi_n) \rangle}{A^2} 
  - \frac{\langle i ( \psi - \psi_n ) |\p_x \psi \rangle}{A^2} ,
$$
%\begin{align*} 
% \p_x \phi_n - \p_x \phi = &\, 
% \frac{\langle i \psi_n |ïŸæµ¬p_x \psi_n \rangle}{A_n^2} - \frac{\langle i \psi |ïŸæµ¬p_x \psi \rangle}{A^2} 
% %\\
% %= & \, 
% = \langle i \psi_n |ïŸæµ¬p_x \psi_n \rangle \Big( \frac{1}{A_n^2} - \frac{1}{A^2} \Big)
% %\frac{(A_n-A) (A_n+A)}{A^2A_n^2} 
% + \frac{\langle i \psi_n |ïŸæµ¬p_x (\psi - \psi_n) \rangle}{A^2} 
%  - \frac{\langle i ( \psi - \psi_n ) |ïŸæµ¬p_x \psi \rangle}{A^2} ,\\
% = & \, 
%\end{align*}
hence
\be
\label{Lipp}
 \n \p_x \phi_n - \p_x \phi \n_{L^2} \leq 
 \frac{|\!| \psi_n |\!|_{L^\ii} |\!| A |\!|_{L^\ii}}{ ( \inf_\R A^2 ) ( \inf_\R A_n^2 )} \n A_n - A \n_{L^2} 
 + \frac{|\!| \psi_n |\!|_{L^\ii}}{\inf_\R A^2 } \n \p_x \psi - \p_x \psi_n \n_{L^2} 
 + \frac{1}{\inf_\R A^2 } \n | \psi_n - \psi | \p_x \psi \n_{L^2} .
 \ee
The first two terms tend to zero as $ n \to + \ii $. For the last term, we use here 
again the Dominated Convergence Theorem since $ \psi_n(x) = 
\psi_n(0) + \int_0^x \p_x \psi_n(t) \ dt \to \psi(0) + \int_0^x \p_x \psi(t) \ dt = \psi(0) $ 
pointwise. This concludes for $(i)$.\\

\noindent {\bf Proof of (ii).} Let us define 
$ \phi_* : \R \to \R $ by $ \phi_*(x) \equiv 
\ds{ \frac12 (\ln x)^2} {\bf 1}_{x\geq 1} $. Then, straightforward 
computations give $ \p_x \phi_* (x) = \ds{ \frac{ \ln x}{x} } {\bf 1}_{x\geq 1} 
\in L^2(\R) $ and, for $ X \geq \ex $, by monotonicity of $ \p_x \phi_* $,
\be
\label{lo}
 \int_X^{2X} ( \p_x \phi_* )^2 \ d x \geq X \frac{ \ln^2(2X)}{(2X)^2} 
\geq \frac{ (\ln X)^2}{4 X}. 
\ee
We now consider $ \phi_n  : \R \to \R $ defined by 
$ \phi_n(x) = 0 $ if $ x \leq 0 $ or $ x \geq 3 n \pi $, 
$ \phi_n(x) = x / n $ if $ 0 \leq x  \leq n \pi $, 
$ \phi_n(x) = \pi $ if $ n \pi \leq x \leq 2 n \pi $ and 
$ \phi_n(x) = 3 \pi - x / n $ if $ 2 n \pi \leq x \leq 3 n \pi $. 
Then, we easily obtain
$$
 d_{\rm hy} ( \ex^{i \phi_*} , \ex^{i \phi_* + i \phi_n} ) = 
 \n \p_x \phi_n \n_{L^2} = \sqrt{ 2 \times \frac{\pi n}{n^2} } = \sqrt{ \frac{2\pi}{n} } 
\to 0 .
$$
Moreover,
$$
 d_\mathcal{Z} ( \ex^{i \phi_*} , \ex^{i \phi_* + i \phi_n} ) 
= \n \p_x \phi_* \ex^{i \phi_*} - ( \p_x \phi_* + \p_x \phi_n ) \ex^{i \phi_* + i \phi_n} \n_{L^2} 
\geq \n \p_x \phi_* ( \ex^{i \phi_n} - 1 ) \n_{L^2} - \n \p_x \phi_n \n_{L^2} ,
$$
and, by our choice of $ \phi_n $ and using \eqref{lo},
$$
\n \p_x \phi_* ( \ex^{i \phi_n} - 1 ) \n_{L^2}^2 
\geq \int_{n \pi}^{2 n \pi} 4 ( \p_x \phi_*)^2 \ dx 
\geq \frac{ (\ln X)^2}{X}_{|X=n\pi} \sim \frac{ (\ln n)^2}{n \pi} .
$$
Since $ \ds{ \frac{ \ln n}{ \sqrt{n \pi} } \gg \sqrt{\frac{2\pi}{n} } 
= d_{\rm hy} ( \ex^{i \phi_*} , \ex^{i \phi_* + i \phi_n} ) } $, it follows that, 
as wished
$$
d_\mathcal{Z} ( \ex^{i \phi_*} , \ex^{i \phi_* + i \phi_n} ) \geq 
\frac{ \ln n }{ \sqrt{n \pi} } ( 1 + o(1) ) \gg \sqrt{\frac{2\pi}{n} } 
= d_{\rm hy} ( \ex^{i \phi_*} , \ex^{i \phi_* + i \phi_n} ) .
$$

We do not know whether the mapping $ \mathscr{M} $ is locally Lipschitz, 
but is probably not.

%%%%%%%%%%%%%%%%%%%%%%%%%%%%%%%%%%%%%%%%%%%%%%%
\subsection{Proof of Lemma \ref{equivdistance}}
\label{pfequiv}

Note first that since $U_*$ does not vanish, if 
$\d$ is sufficiently small and $d_{\mathcal{Z}} ( \psi , U_* ) \leq \d$, then 
$ \n | \psi | - | U_* | \n_{L^\ii} \leq (1/2) \inf_{\R} | U_* | $, hence 
$ | \psi | \geq (1/2) \inf_{\R} | U_* | > 0 $ in $\R$, thus  $\psi$ does 
not vanish, may be lifted $ \psi = A \exp( i \phi) $, and we may further 
assume $ \phi (0) - \phi_*(0) \in ( - \pi , + \pi ] $. In \eqref{Lipsch}, we can 
easily check that the terms leading to the "$o(1)$" are indeed controlled by 
$ K(\psi) d_{\rm hy} ( \psi_n , \psi ) $. In other words, we have
$$
  d_\mathcal{Z} ( \psi , U_* ) \leq 
  \n ( \ex^{i \phi} - \ex^{i \phi_* }) \p_x A_* \n_{L^2} 
+ \n A_* \n_{L^\ii} \n ( \ex^{i \phi} - \ex^{i \phi_*} ) \p_x \phi_* \n_{L^2} 
+ K(U_*) d_{\rm hy} ( \psi , U_* ) ,
$$
provided $d_{\rm hy} ( \psi , U_* )$ is small enough. In order to bound the
two remaining terms, we write, for $ x \in \R $,
$$
\phi (x) - \phi_*(x) = \phi (0) - \phi_*(0) + \int_0^x \p_x \phi (y) - \p_x \phi_*(y) \ d y ,
$$
which implies, using that $ \R \ni \theta \mapsto \ex^{i \theta} $ is $1$-Lipschitz 
and the Cauchy-Schwarz inequality, 
\be
\label{integ}
\Big | 1 - \ex^{ i( \phi_* (x) - \phi (x)) } \Big| \leq \Big | \phi (0) - \phi_*(0) \Big| 
+ \sqrt{|x|} \n u - u_* \n_{L^2} .
\ee
Consequently, 
$$
  \n ( \ex^{i \phi} - \ex^{i \phi_* }) \p_x A_* \n_{L^2} \leq 
\Big | \phi (0) - \phi_*(0) \Big| \n \p_x A_* \n_{L^2} 
+ \n u - u_* \n_{L^2} \n  \sqrt{|x|} \p_x A_* \n_{L^2} 
$$
and
$$
 \n ( \ex^{i \phi} - \ex^{i \phi_*} ) \p_x \phi_* \n_{L^2} \leq 
 \Big | \phi (0) - \phi_*(0) \Big| \n \p_x \phi_* \n_{L^2} 
+ \n u - u_* \n_{L^2} \n \sqrt{|x|} \p_x \phi_* \n_{L^2} .
$$
Both terms are $ \leq K(U_*) d_{\rm hy} ( \psi , U_* ) $. 
Indeed, $U_* \in \mathcal{Z} $ is a travelling wave, hence 
$ A_* $, $ \p_x A_* $, $ \p_x \phi_* $ are bounded functions 
which decay at infinity exponentially if $ 0 \leq c < \cs $ ({\it cf.} 
Proposition \ref{decroiss} $(i)$ or $(ii)$). If $ c = \cs $, 
since assumption $(\BA_0)$ is satisfied, we invoke 
Proposition \ref{decroiss} $(iii)$, which ensures that  $ \p_x \phi_* $ and 
$ \p_x A_* $ decay at the rate $ \BO ( |x|^{- 2} ) $ ($ \p_x A_* $ decays 
faster actually). Therefore, $  \sqrt{|x|} \p_x \phi_* \in L^2 $. Gathering these 
estimates provides
$$ 
d_{\mathcal{Z}} (  \psi , U_* ) \leq K( U_* ) d_{\rm hy} (  \psi , U_* ) .
$$

On the other hand, from \eqref{Lipp} and the estimate 
$ |\!| A - A_* |\!|_{H^1} \leq d_{\mathcal{Z}} (  \psi , U_* ) $ (see 
the proof of $(i)$), we infer
$$
d_{\rm hy} (  \psi , U_* ) \leq K(U_*) d_{\mathcal{Z}} (  \psi , U_* ) 
+ \frac{1}{\inf_\R A^2 } \n | \psi - U_* | \p_x U_* \n_{L^2} .
$$
Using here again the estimate $ | \psi (x) - U_* (x) | \leq | \phi (0) - \phi_*(0) | 
+ \sqrt{|x|} \n \p_x \psi - U_* \n_{L^2} $, we deduce
$$
d_{\rm hy} (  \psi , U_* ) \leq K(U_*) d_{\mathcal{Z}} (  \psi , U_* ) .
$$
The proof is complete. \carre 

%%%%%%%%%%%%%%%%%%%%%%%%%%%%%%%%%%%%%%%%%%%%%%%%%%%%%%%%%%%
\subsection{Two remarks on the proof of Theorem \ref{stab}}
\label{sexremarkpf}

We would like to point out two minor points concerning the 
proof of Theorem \ref{stab} by Z. Lin. We recall that the proof 
of Z. Lin \cite{Lin} relies on the Grillakis-Shatah-Strauss 
theory \cite{GSS} once we have written (NLS) under the 
hydrodynamical form \eqref{MadTWk}, denoting $ \psi = A \ex^{i \phi}$, 
$ (\rho, u) \equiv (|\psi|^2 = A^2 , \p_x \phi )$:
$$
\left\{\begin{array}{ll}
\ds{ \p_t \rho + 2 \p_{x} (\rho u ) } = 0 \\ \ \\ 
\ds{ \p_{t} u + 2 u \p_{x} u - \p_x (f ( \rho )) 
- \p_x \Big( \frac{ \p_{x}^2 (\sqrt{\rho}) }{\sqrt{\rho} } \Big) } = 0 ,
\end{array}\right. 
$$
or more precisely, with $ \eta \equiv \rho - r_0^2 = |\psi|^2  - r_0^2 $ and 
denoting $\frac{\d E}{\d \eta}$, $\frac{\d E}{\d u} $ the variational derivative,
\be
\label{hydrohamilto}
\frac{\p}{\p t} \left( \begin{array}{c}
\eta \\ u 
\end{array}\right) = J \left( \begin{array}{c}
\ds{\frac{\d E}{\d \eta}} \\ \ds{ \frac{\d E}{\d u} }
\end{array}\right) ,
\quad \quad \quad 
J \equiv \left( \begin{array}{cc}
0 & \p_x \\ \p_x & 0 
\end{array}\right) .
\ee
The first remark is that 
the scalar product in the Hilbert space $ X = H^1\times L^2 $ can 
not be $ ( (\eta, u ) , ( \tilde{\eta} , \tilde{u} ) )_{H^1\times L^2} 
= \int_\R \eta \tilde{\eta} + u \tilde{u} \ dx $ as used in \cite{Lin}, 
but the natural one is $ ( (\eta, u ) , ( \tilde{\eta} , \tilde{u} ) )_{H^1\times L^2} 
= \int_\R \eta \tilde{\eta} + \p_x \eta \p_x \tilde{\eta} + u \tilde{u} \ dx $. 
This requires to make some minor changes in the proof, 
especially not to identify $(H^1)^*$ with $H^1$. For instance, 
a linear mapping $B$ is associated with the momentum through 
the formula
$$
P_{\rm hy} ( \eta , u ) \equiv 
\int_\R \eta u \ dx = \frac12 ( B (\eta, u ) , ( \eta , u ) )_{H^1\times L^2} 
\quad \quad \quad {\rm with} \quad B \equiv \left( \begin{array}{cc} 
0 & 1 \\ 1 & 0 
\end{array}\right)
$$
for the (non hilbertian) scalar product 
$ ( (\eta, u ) , ( \tilde{\eta} , \tilde{u} ) )_{H^1\times L^2} = 
\int_\R \eta \tilde{\eta} + u \tilde{u} \ dx $. The correct definition 
is actually
$$
P_{\rm hy} ( \eta , u ) 
=  \int_\R \eta u \ dx = \frac12 \langle B (\eta, u ) , ( \eta , u ) \rangle_{X^*, X} 
\quad \quad \quad {\rm with} \quad B \equiv \left( \begin{array}{cc} 
0 & \iota^* \\ \iota & 0 
\end{array}\right) ,
$$
where $ \iota : H^1 \hookrightarrow L^2 $ is the canonical injection. 
As already mentioned in \ref{sex0cs}, the two points in the proof 
of \cite{Lin} that have been completed by C. Gallo in \cite{GaZhi} are 
that: \cite{Lin} uses a local in time existence for the hydrodynamical system 
\eqref{hydrohamilto} in $ H^1 \times L^2 $, and not only in 
$ \{ \rho \in L^\ii,\ \p_x \rho \in L^2\} \times L^2 $; and that the 
energy and the momentum are indeed conserved for the local solution 
if the initial datum does not vanish. \\

The second point is that in the proof of stability (theorem 3.5 in \cite{GSS}), 
it is made use of the fact that if $ \BU \in X $ and $ ( \BU_n )_{n \in \N} \in X $ is 
a sequence such that $ E(\BU_n) \to E(\BU) $ and $ P_{\rm hy} (\BU_n) \to P_{\rm hy}(\BU) $, 
then there exists a sequence $ ( \tilde{\BU}_n )_{n \in \N} \in X $ such that 
$ \BU_n - \tilde{\BU}_n \to 0 $ in $X$, $ E(\tilde{\BU}_n) \to E(\BU) $ and 
$ P_{\rm hy}(\tilde{\BU}_n) = P_{\rm hy}(\BU) $. In the context of bound states, 
the existence of such a sequence $ ( \tilde{\BU}_n )_{n \in \N} \in X $ follows 
by simple scaling in space, since then the momentum or charge is simply 
$ \int_{\R^d} \BU_n^2 \ d x $. However, for the one dimensional travelling waves 
for (NLS), the momentum $P$ is scaling invariant. We do not know if the existence 
of such a sequence holds in a general framework, but for the problem 
we are studying, we can rely on the following lemma, which is 
an adaptation of Lemma 6 in \cite{BGSsurvey} (see also Lemma in \cite{BGSS}). 
% The computations being very similar, the proof is omitted.

%%%%%%%%%%%%%%
\begin{lemmas}
\label{raccord} 
There exists $ \mathfrak{p}_0 > 0 $ and $ K > 0 $, depending only on 
$f$ such that for any $ \mathfrak{p} \in ( -  \mathfrak{p}_0 , +  \mathfrak{p}_0 ) $ 
and $ \mu \in \R $ with $ | \mu | \leq | \mathfrak{p} | $, there exists 
$ w = a \ex^{i \vp } \in H^1 ( [0 , 1 / ( 2 | \mathfrak{p}| ) ) , \C ) $ 
verifying
$$
w(0) = w \Big( \frac{1}{2 | \mathfrak{p}| } \Big) , 
\quad \quad \quad 
|w(0)| = r_0 + \mu ,
\quad \quad \quad 
\int_0^{1 / ( 2 | \mathfrak{p}| )} ( a^2 - r_0^2 ) \p_x \vp \ dx = \mathfrak{p} 
$$
and
$$
 \int_0^{1 / ( 2 | \mathfrak{p}| )} | \p_x w |^2 + F(|w|^2) \ d x \leq K | \mathfrak{p}| 
$$
\end{lemmas}
%%%%%%%%%%%%

\noindent {\it Proof.} If $ \mathfrak{p} = 0 $, we simply take $ w = r_0 $. 
We then assume $ 0 <  \mathfrak{p} <  \mathfrak{p}_0 $, since the case 
$ - \mathfrak{p}_0 <  \mathfrak{p} < 0 $ will follow by complex conjugation. 
We then define, for some small $ \d $ to be determined later,
$$
w(x) \equiv \sqrt{ r_0^2 - \d + 2\mathfrak{p} \Big( 1 - | 8 \mathfrak{p} x - 1| \Big)_+ } 
\exp \Big[ i \Big( 1 - | 4 \mathfrak{p} x - 1| \Big)_+ \Big] = a \ex^{i \vp }.
$$
It is clear that $w \in H^1( [0 , 1 / ( 2 \mathfrak{p} ) ] , \C ) $ and that 
$ w(0) = w ( 1/ (2 \mathfrak{p} ) ) = \sqrt{ r_0^2 - \d } $, thus 
$ |w(0)| = r_0 + \mu $ provided we choose $ \d =  - \mu^2 - 2 r_0 \mu = \BO(|\mu|) $. 
Moreover, using that the phase $ \vp $ has compact support 
$ [ 0 ,  1 / ( 2 \mathfrak{p} ) ] $, 
\begin{align*}
\int_0^{1 / ( 2 \mathfrak{p} )} ( a^2 - r_0^2 ) \p_x \vp \ dx 
= & \, 
\int_0^{1 / ( 2 \mathfrak{p} )} \Big\{ 
- \d + 2\mathfrak{p} \Big( 1 - | 8 \mathfrak{p} x - 1| \Big)_+ \Big\} \p_x 
\Big( 1 - | 4 \mathfrak{p} x - 1| \Big)_+ \ dx 
\\ = & \, 
2\mathfrak{p} 
\int_0^{1 / ( 2 \mathfrak{p} )} \Big( 1 - | 8 \mathfrak{p} x - 1| \Big)_+ 
\p_x \Big( 1 - | 4 \mathfrak{p} x - 1| \Big)_+ \ dx 
\\ = & \, 
2\mathfrak{p} 
\int_0^{1 / ( 2 \mathfrak{p} )} \Big( 1 - | 8 \mathfrak{p} x - 1| \Big)_+ 
\p_x \Big( 1 - | 4 \mathfrak{p} x - 1| \Big)_+ \ dx .
\end{align*}
For the last integral, the first factor is equal to $0$ if 
$ x \geq 1 / ( 4 \mathfrak{p} ) $ and the second factor is equal to 
$ 4 \mathfrak{p} $ when $ x \leq 1 / ( 4 \mathfrak{p} ) $. Hence, 
direct computation gives 
$$
\int_0^{1 / ( 2 \mathfrak{p} )} ( a^2 - r_0^2 ) \p_x \vp \ dx 
= 2 \mathfrak{p} \int_0^{1 / ( 4 \mathfrak{p} )} 
\Big( 1 - | 8 \mathfrak{p} x - 1| \Big)_+ \times 4 \mathfrak{p}  \ dx 
= \mathfrak{p} .
$$
For the energy part, notice first that
$$
| a^2 - r_0^2 | = \Big| - \d + 2 \mathfrak{p} \Big( 1 - | 8 \mathfrak{p} x - 1| \Big)_+  \Big| 
\leq | \d | + 2 \mathfrak{p}_0 
$$
is as small as we want if $ | \d | $ and $ \mathfrak{p}_0 $ are chosen 
sufficiently small. Therefore,
$$
F( | w |^2 ) \leq K ( a^2 - r_0^2 )^2 .
$$
By simple computations, we have
\begin{align*}
\int_0^{1 / ( 2 \mathfrak{p} )} & \, | \p_x w |^2 + F(|w|^2) \ d x \\ 
& \leq K \int_0^{1 / ( 2 \mathfrak{p} )} 
\mathfrak{p}^2 \Big| \p_x  \Big( 1 - | 8 \mathfrak{p} x - 1| \Big)_+ \Big|^2 
+ \Big| \p_x \Big( 1 - | 4 \mathfrak{p} x - 1| \Big)_+ \Big|^2
+ \Big( - \d + 2 \mathfrak{p} \Big( 1 - | 8 \mathfrak{p} x - 1| \Big)_+ \Big)^2 \ d x \\ 
&  \leq K \mathfrak{p}^3 + K \mathfrak{p} + K \frac{\d^2 + \mathfrak{p}^2}{ \mathfrak{p} }
\leq K \mathfrak{p} 
\end{align*}
since $ \d = \BO(|\mu|) = \BO ( \mathfrak{p}) $, which concludes the proof. \carre \\

We then consider a sequence $ \BU_n = (\eta_n , u_n ) \in X = H^1 \times L^2 $ and 
show the existence of the desired sequence 
$ \tilde{\BU}_n = ( \tilde{\eta}_n , \tilde{u}_n ) \in X $. We recall that 
$ \BU_n $ (resp. $\BU$) is associated with a mapping $ \psi_n \in \mathcal{Z} $ 
(resp. $U_*$) that does not vanish. We have $ P_{\rm hy} ( \BU_n ) = P ( \psi_n ) \to P(U_*) $, 
thus for $ n $ large enough, $ | P ( \psi_n ) - P(U_*) | \leq \mathfrak{p}_0 $. 
For $ n $ fixed, we now pick $ R_n > 0 $ large enough so that
$$
\int_{R_n}^{+ \ii } | \p_x  \psi_n |^2 + ( | \psi_n | - r_0 )^2 \ dx 
\leq [ P ( \psi_n ) - P(U_*) ]^2 . 
$$
In particular, by Sobolev embedding,
$$ 
\big| \, | \psi_n | (R_n) - r_0 \, \big| \leq \n  | \psi_n | - r_0 \n_{L^\ii( [ R_n,+\ii ))}
\leq \sqrt{ \int_{R_n}^{+ \ii } | \p_x | \psi_n | |^2 + ( | \psi_n | - r_0 )^2 \ dx } 
\leq | P ( \psi_n ) - P(U_*) | .
$$
We are now in position to apply (for $n$ large) Lemma \ref{raccord} with 
$ ( \mathfrak{p} , \mu ) = ( P(U_*) - P ( \psi_n ) , | \psi_n | (R_n) - r_0 ) $. 
This provides the mapping $ w_n \in H^1 ( [0 , 1 / ( 2 | \mathfrak{p}| ) ) , \C ) $. 
Since $ | \psi_n | (R_n) - r_0 \to 0 $, for $n$ large enough, there exists 
$\theta_n \in \R$ such that $ \psi_n (R_n) = \ex^{ i \theta_n } | \psi_n | (R_n) 
= \ex^{ i \theta_n } ( r_0 + \mu ) = \ex^{ i \theta_n } w_n(0) $. We then consider 
the mapping $ \tilde{\psi}_n \in \mathcal{Z} $ defined by
$$
\tilde{\psi}_n (x) \equiv \left\{\begin{array}{ll}
\ds{ \psi_n(x) } & \quad {\rm if} \ x \leq R_n \\ \ \\ 
\ds{ \ex^{ i \theta_n } w_n (x -R_n) } & \quad {\rm if} \ R_n \leq x 
\leq R_n + \frac{1}{ 2 | P ( \psi_n ) - P(U_*) | } \\ \ \\ 
\ds{ \psi_n \Big( x - \frac{1}{ 2 | P ( \psi_n ) - P(U_*) | } \Big) } & \quad
 {\rm if} \ x \geq R_n + \frac{1}{ 2 | P ( \psi_n ) - P(U_*) | } .
\end{array}\right. 
$$
From the construction of $w_n$ and the phase factor $ \theta_n $, 
$\tilde{\psi}_n $ is well-defined and continuous. It is clear that
$$
P ( \tilde{\psi}_n ) = P ( \psi_n ) + \int_0^{1 / (2\mathfrak{p} )} 
( a_n^2 - r_0^2 ) \p_x \vp_n \ dx 
= P ( \psi_n ) + \mathfrak{p} = P(U_*) 
$$
for every (large) $n$, and that
$$
E ( \tilde{\psi}_n ) = E ( \psi_n ) + \int_0^{1 / ( 2 \mathfrak{p} )} 
| \p_x w |^2 + F(|w|^2) \ dx 
= E ( U_* ) + o(1) + \BO( |\mathfrak{p}| ) 
= E ( U_* ) + o(1) + \BO ( | P(U_*) - P ( \psi_n ) | ) 
$$
converges to $ P(U_*) $ as $ n \to + \ii $. Denoting $ \tilde{\BU}_n \in X $ 
the hydrodynamical expression of $ \tilde{\psi}_n $, it remains to show 
that $ \BU_n - \tilde{\BU}_n \to 0 $ in $X = H^1 \times L^2 $. We thus 
compute, with the definition of $ \tilde{\psi}_n $,
\begin{align*}
\n \BU_n - \tilde{\BU}_n \n_{X}^2 = & \, 
\int_{R_n}^{+\ii} \Big| \p_x | \psi_n | - \p_x | \tilde{\psi}_n | \Big|^2 
+ ( | \psi_n | - | \tilde{\psi}_n | )^2 + ( u_n - \tilde{u}_n )^2 \ dx 
\\ \leq & \, 2 \int_{R_n}^{+\ii} | \p_x | \psi_n | |^2 + | \p_x | \tilde{\psi}_n | |^2 
+ ( | \psi_n | - r_0) ^2 + ( | \tilde{\psi}_n | - r_0 )^2 + u_n^2 + \tilde{u}_n^2 \ dx 
\\ \leq & \, 4 K \int_{R_n}^{+\ii} | \p_x  \psi_n |^2 + ( | \psi_n | - r_0 )^2 \ dx 
+ 2  \int_0^{1 / ( 2 | P(U_*) - P ( \psi_n )| )} | \p_x w_n |^2 + (|w_n|^2 - r_0^2 )^2 \ d x 
\\ \leq & \, 4 K [ P ( \psi_n ) - P(U_*) ]^2 + K | P ( \psi_n ) - P(U_*) | \to 0 .
\end{align*}
For the before last inequality, we have used that for $|x| \geq R_n $, 
$\psi_n$ has modulus uniformly close to $r_0$, hence 
$ | \p_x | \psi_n | |^2 + u_n^2 \leq K | \p_x  \psi_n |^2 $. Note that 
the construction still holds for the energy distance, the computations 
being similar. 

%%%%%%%%%%%%%%%%%%%%%%%%%%%%%%%%%%%%%%%%%%%%%%%%%%%%%%%%%%%
\subsection{Proof of Lemma \ref{souskontrol}}
\label{pfsouskontrol}

\noindent {\bf Proof of estimate $\bs{\eqref{Kontr}}$.} Instead of 
concluding the stability proof as in \cite{GSS}, we can notice that 
we have actually the bound
\be
\label{minicontr}
E_{\rm hy} (\BU) - E_{\rm hy}( \BU_{c_*}) \geq 
\frac{1}{K} \inf_{y\in \R } |\!| \BU - \BU_{c_*} ( \cdot - y) |\!|^2
\ee
as soon as $ P(U_{c_*}) = P_{\rm hy}(\BU_{c_*} ) = P_{\rm hy}( \BU )$ and 
$ \BU \in \mathscr{O}_\e \equiv \Big\{ \BV \in X , \ 
\inf_{ y \in \R } \n \BV - \BU_{c_*}( \cdot - y ) \n_X < \e \Big\} $ 
for some small $ \e $. If $ \Psi^\ini $ does not have momentum equal 
to $ P_{\rm hy}(\BU_{c_*} )$, we use Lemma \ref{raccord} to infer that there 
exists $ \tilde{\Psi}(t) $, with momentum equal to $ P_{\rm hy} (\BU_{c_*} ) = P(U_{c_*})$, 
and such that $ E (\tilde{\Psi}(t)) - E (\Psi(t)) = \BO( | P( \Psi(t) ) - P(U_{c_*} ) | ) $ 
and
$ d_{\rm hy}( \Psi(t), \tilde{\Psi}(t) ) \leq \BO( \sqrt{| P( \Psi(t) ) - P(U_{c_*} ) |} ) $. 
Therefore, for $ t \geq 0 $, denoting $ \Psi_{\rm hy}(t) \in X $ and 
$ \tilde{\Psi}_{\rm hy} (t) \in X $ the hydrodynamical variables for 
$ \Psi $ and $\tilde{\Psi} (t)$,
\begin{align*}
\inf_{ y \in \R } |\!| \Psi_{\rm hy} (t) - \BU_{c_*} (\cdot - y ) |\!| 
\leq & \, 
\inf_{ y \in \R } \Big[ |\!| \tilde{\Psi}_{\rm hy} (t) - \BU_{c_*} (\cdot - y ) |\!| 
+ |\!| \Psi_{\rm hy} (t) - \tilde{\Psi}_{\rm hy} (t) |\!| \Big] \\
\leq & \, \sqrt{K} \sqrt{E (\tilde{\Psi}(t)) - E (U_{c_*})} 
+ \BO( \sqrt{| P( \Psi(t) ) - P(U_{c_*} ) |} ) 
\\ \leq & \, 
K \Big[ \sqrt{ | E( \Psi(t) ) - E ( U_{c_*} ) |
+ | P( \Psi^\ini ) - P(U_{c_*} ) |}
+ \sqrt{| P( \Psi^\ini ) - P(U_{c_*} ) |} \Big] ,
\end{align*}
which yields \eqref{Kontr}.\\
 
The above estimate is optimal when $ P(\Psi^\ini) = P ( U_{c_*} ) $ since 
$ U_{c_*} $ is a critical point of the action $ E - c_* P $. This 
bound shows that, in the definition of stability, one has to take 
$ \d = \BO( \e^2 )$ in general. The estimate $\bs{\eqref{Kontr}}$ shows 
that one can actually take $ \d = \BO( \e )$.\\

\noindent {\bf Proof of estimate $\bs{\eqref{Kontrlin}}$.} The 
point is to compare $ \Psi(t) $ to $ U_{c} $ with $ c \simeq c_* $ 
such that $ P(U_{c}) = P(\Psi^\ini) $ instead of comparing 
to $ U_{c_*} $. In other words, we replace $ \tilde{\Psi}(t)$ 
by $ U_{c} $. Note first that since $ \ds{\frac{dP}{dc}_{ | c = c_*} < 0 }$, 
there exists, by the implicit function theorem, such a 
$ c \simeq c_* $. We then proceed as follows. Let $ \Psi^\ini \in \mathcal{Z} $ 
be close to $ U_{c_*} $. Then, there exists $ c = c(\Psi^\ini) \simeq c_* $ 
such that $ P( U_{c}) = P(\Psi^\ini) $. Moreover, since 
$ \ds{\frac{dP}{dc}_{ | c = c_*} \not = 0 }$, it follows
\be
\label{crapule}
|\!| \BU_c - \BU_{c_*}  |\!| \leq K | c - c_* | 
\leq K | P( U_c ) - P(U_{c_*} ) | = K | P(\Psi^\ini) - P(U_{c_*} ) | 
\leq K |\!| \Psi_{\rm hy}^\ini - \BU_{c_*} |\!| \leq K d_{\rm hy} (\Psi^\ini , U_{c_*} ) .
\ee
From \eqref{minicontr}, it comes 
$$
E_{\rm hy} (\BU) - E_{\rm hy} ( \BU_{c}) \geq 
\frac{1}{K} \inf_{y\in \R } |\!| \BU - \BU_{c} ( \cdot - y) |\!|^2
$$
as soon as $ P_{\rm hy} (\BU) = P_{\rm hy} ( \BU_{c} ) $. The fact that 
the constant $K$ can be taken uniform with respect to $c$ for $c$ close 
to $c_*$ comes directly from the proof in \cite{GSS}. Therefore, 
for $ t \geq 0 $,
\begin{align*}
\inf_{ y \in \R } |\!| \Psi_{\rm hy} (t) - \BU_{c_*} (\cdot - y ) |\!| 
\leq & \, 
\inf_{ y \in \R } \Big[ |\!| \Psi_{\rm hy} (t) - \BU_{c} (\cdot - y ) |\!| 
+ |\!| \BU_{c} (\cdot - y ) - \BU_{c_*} (\cdot - y ) |\!| \Big] \\
\leq & \, \sqrt{K} \sqrt{E ( \Psi (t)) - E (U_{c})} 
+ \BO( | P( \Psi^\ini ) - P(U_{c_*} ) | ) .
\end{align*}
Using that $ P (\Psi(t)) = P_{\rm hy}(\Psi_{\rm hy}(t)) = P_{\rm hy} ( \BU_{c} ) $ and that 
$ \BU_{c} $ is a critical point of the action $ E_{\rm hy} - c P_{\rm hy} $, we infer 
$ E ( \Psi (t)) - E (\BU_{c}) =  [E_{\rm hy}-cP_{\rm hy}] 
(\Psi_{\rm hy}^\ini) - [E_{\rm hy} - c P_{\rm hy}] (\BU_{c}) 
= \BO ( |\!| \Psi_{\rm hy}^\ini - \BU_{c} |\!|^2 ) $. Consequently,
$$
\inf_{ y \in \R } |\!| \Psi_{\rm hy}(t) - \BU_{c_*} (\cdot - y ) |\!| 
\leq K \Big( |\!| \Psi_{\rm hy}^\ini - \BU_{c} |\!| + |\!| \Psi_{\rm hy}^\ini - \BU_{c_*} |\!| \Big) 
\leq K d_{\rm hy} ( \Psi^\ini , U_{c_*} ) + K |\!| \BU_c - \BU_{c_*} |\!| 
\leq K d_{\rm hy} ( \Psi^\ini , U_{c_*} ) ,
$$
by \eqref{crapule}. This gives \eqref{Kontrlin}.

%%%%%%%%%%%%%%%%%%%%%%%%%%%%%%%%%%%%%%%%%%%%%%%%%%%%%%%%%%%%%%%%%%%%%%%
\section{Instability result for cusps: proof of Theorem \ref{stabcusp}}

In this section, we set $ \BF_c \equiv E_{\rm hy} - c P_{\rm hy} $ 
and we assume
$$ 
- \frac{d^2 \BF_c(U_c)}{dc^2}_{|c=c_*} = \frac{d P(U_c)}{dc}_{|c=c_*} = 0 
\quad \quad \quad {\rm and } \quad \quad \quad 
0 \not = \ddot{P}_* \equiv \frac{d^2 P(U_c)}{dc^2}_{|c=c_*} 
= - \frac{d^3 \BF_c(U_c)}{dc^3}_{|c=c_*} .
$$
The approach is reminiscent to the proof of M. Maeda \cite{Maeda}. 
Several modifications are necessary since for the skew-adjoint 
operator $ J = \p_x $, we can not find the required Hilbert space $Y$. 
More degenerate cases can probably be considered as in \cite{Maeda}. 

We shall denote $ \mathbb{I} : X \to X^* $ and $ \mathbb{I}_{H^1} : H^1 \to (H^1)^* $ 
the Riesz isomorphisms and $ \BU = ( \eta , u )^t \in X = H^1(\R,\R) \times L^2 (\R,\R) $ 
and $ H \equiv L^2(\R,\R) \times L^2 (\R,\R) $, endowed with its canonical 
scalar product. They are the corresponding Hilbert spaces needed in \cite{Maeda}. 
We consider the symmetric matrix
$$
\mathbb{B} \equiv \left( \begin{array}{cc} 
0 & 1 \\ 1 & 0 
\end{array}\right) ,
$$
which is such that $ \mathbb{B}^2 = {\rm I}_2 $ and 
$ 2 P_{\rm hy} ( \BU ) = ( \mathbb{B} \BU , \BU)_{H} $.

Our assumption $ \ds{ \frac{d P(U_c)}{dc}_{|c=c_*} = 0 \not = \frac{d^2 P(U_c)}{dc^2}_{|c=c_*} } $ 
will simplify a little the computations in \cite{Maeda}. The functions $ \eta_1 $ 
and $ \eta_2 $ used in \cite{Maeda} become now
$$
 \eta_1 ( \gamma ) = \BF_{c_* + \gamma } ( \BU_{c_* + \gamma } ) 
- \BF_{c_* } ( \BU_{c_* } ) - \gamma \frac{d \BF_c(U_c)}{dc}_{|c=c_*} 
\sim - \frac{\gamma^3}{6} \ddot{P}_* $$
and
$$
\eta_2 (\gamma ) = \frac{d \eta_1 }{d \gamma } 
= - P(U_{c_* + \gamma}) + P(U_{c_*} ) 
\sim - \frac{\gamma^2}{2} \ddot{P}_* .
$$
In order to clarify the dualities used in \cite{Maeda}, we provide some 
elements of the proof adapted to our context.

%%%%%%%%%%%%%%
\begin{lemmas} 
\label{Maeda1} 
There exists $ \gamma_0 > 0 $ small and $ \s : ( - \gamma_0 , + \gamma_0 ) \to \R $ 
with $ \s ( \gamma) \sim - \ds{ \frac{\gamma^2 \ddot{P}_* }{2\n \BU_{* } \n^2_H } } $ 
and such that, for any $ \gamma \in ( - \gamma_0 , + \gamma_0 ) $,
$$ P_{\rm hy} ( \BU_{c_* + \gamma } + \s ( \gamma) \mathbb{B} \BU_{c_* + \gamma } ) 
= P_{\rm hy} ( \BU_* ) . $$
\end{lemmas}
%%%%%%%%%%%%

\noindent {\it Proof.} We have
$$ P_{\rm hy} ( \BU_{c_* + \gamma } + \s \mathbb{B} \BU_{c_* + \gamma } ) 
= \frac12 ( \mathbb{B} \BU_{c_* + \gamma } + \s \BU_{c_* + \gamma } , 
\BU_{c_* + \gamma } + \s \mathbb{B} \BU_{c_* + \gamma } )_H 
= P_{\rm hy} ( \BU_{c_* + \gamma } ) 
+ \s \n \BU_{c_* + \gamma } \n^2_H 
+ \s^2 P_{\rm hy} ( \BU_{c_* + \gamma } ). $$
Since $ \n \BU_{* } \n^2_H \not = 0 $, the conclusion follows 
from an easy implicit function argument near $\s = \gamma = 0 $. In \cite{Maeda}, 
the linear mapping $B$ is seen from $X$ to $X^*$, but here, there is no confusion 
to define 
$ \BU_{c_* + \gamma } + \s \mathbb{B} \BU_{c_* + \gamma } \in H = L^2 \times L^2 $. \carre \\

We define, for $ \gamma \in ( - \gamma_0 , + \gamma_0 ) $,
$$
 \BW ( \gamma ) \equiv \BU_{c_* + \gamma } + \s ( \gamma) \mathbb{B} \BU_{c_* + \gamma } ,
$$
which then verifies by construction $ P_{\rm hy} (\BW ( \gamma ) ) =  P_{\rm hy} ( \BU_* ) $.

%%%%%%%%%%%%%%
\begin{lemmas} 
\label{Maeda2} 
As $ \gamma \to 0 $, there holds $ \BF_{c_* } ( \BW( \gamma ) ) - \BF_{c_* } ( \BU_{c_* } ) 
\sim - \ds{ \frac{\gamma^3}{6} \ddot{P}_* } $.
\end{lemmas}
%%%%%%%%%%%%

\noindent {\it Proof.} Using that $ \BF'_{c_* + \gamma } ( \BU_{c_* + \gamma } ) = 0 $, 
$ P_{\rm hy} (\BW ( \gamma ) ) =  P_{\rm hy} ( \BU_{c_*} ) = - \ds{\frac{d \BF_c(U_c)}{dc}_{|c=c_*}} $ 
and $ \s (\gamma ) = \BO (\gamma^2)$, we have by Taylor expansion
\begin{align*}
\BF_{c_* } ( \BW( \gamma ) ) - \BF_{c_* } ( \BU_{c_* } ) 
& \, = 
\BF_{c_* + \gamma} (  \BU_{c_* + \gamma } + \s ( \gamma) \mathbb{B} \BU_{c_* + \gamma } ) 
- \BF_{c_* } ( \BU_{c_* } ) + \gamma P_{\rm hy} (\BW ( \gamma ) ) \\
& \, = \BF_{c_* + \gamma} (  \BU_{c_* + \gamma } ) - \BF_{c_* } ( \BU_{c_* } ) 
- \gamma \frac{d \BF_c(U_c)}{dc}_{|c=c_*} + \BO(\gamma^4) 
\sim - \frac{\gamma^3}{6} \ddot{P}_* ,
\end{align*}
as wished. \carre \\

We recall that we have defined the tubular neighbourhood 
$ \mathscr{O}_\e = \{ \BV \in X , \ \inf_{ y \in \R } \n \BV - \BU_{*}( \cdot - y ) \n_X < \e \} $. 

%%%%%%%%%%%%%%
\begin{lemmas} 
\label{Maeda3} 
For $ \e > 0 $ small enough, there exist four $\BC^1$ mappings 
$ \bar{\gamma} $, $ \alpha $, $ \bar{y} : \mathscr{O}_\e \to \R $ 
and $ \vartheta : \mathscr{O}_\e \to X $, satisfying, for 
$ \BU \in \mathscr{O}_\e $,
$$
 \BU ( \cdot - \bar{y}( \BU) ) = \BW( \bar{\gamma}( \BU ) ) 
+ \vartheta ( \BU ) + \alpha ( \BU ) \mathbb{B} \BU_{c_* + \bar{\gamma}( \BU ) }
$$
and the orthogonality relations
$$
 (\vartheta ( \BU ) , \p_x \BU_{c_* + \bar{\gamma}( \BU )} )_H 
= (\vartheta ( \BU ) , [\p_c \BU_c]_{|c = c_* + \bar{\gamma}( \BU )} )_H 
= (\vartheta ( \BU ) , \mathbb{B} \BU_{c_* + \bar{\gamma}( \BU ) } )_H = 0 .
$$
Finally, $ \mathbb{I}^{-1} \bar{\gamma}' \in H^2 \times H^1 $ and 
$ \mathbb{I}_{H^1}^{-1} \frac{ \p \bar{\gamma} }{ \p \eta } \in H^4 $.
\end{lemmas}
%%%%%%%%%%%%

\noindent {\it Proof.} We consider the mapping 
$ G : X \times \R \times ( - \gamma_0 , + \gamma_0 ) \times \R \to \R^3 $ defined by
$$
 G ( \BU , y , \gamma , \alpha ) \equiv \left( \begin{array}{c} 
( \BU( \cdot -y) - \BW(\gamma) - \alpha \mathbb{B} \BU_{c_* + \gamma } , \p_x \BU_{c_*+ \gamma} )_H \\ 
( \BU( \cdot -y) - \BW(\gamma) - \alpha \mathbb{B} \BU_{c_* + \gamma }, [\p_c \BU_c]_{|c = c_* + \gamma } )_H \\ 
( \BU( \cdot -y) - \BW(\gamma) - \alpha \mathbb{B} \BU_{c_* + \gamma }, \mathbb{B} \BU_{c_* + \gamma } )_H 
\end{array}\right) .
$$
Then $ G ( \BU_* , 0 , 0, 0 ) = 0 $ since $ \BW( 0 ) = \BU_* $. In order to show 
that $G$ is of class $ \BC^1 $, we have to pay attention to the translation term 
$ \BU ( \cdot - y ) $, since differentiation in $y$ requires $ \BU \in H^1 \times H^1 $ 
whereas we only have $ \BU \in X = H^1 \times L^2 $. It thus suffices to write
$$
 G ( \BU , y , \gamma , \alpha ) = \left( \begin{array}{c} 
( \BU , \p_x \BU_{c_*+ \gamma} ( \cdot +y))_H 
- ( \BW(\gamma) + \alpha \mathbb{B} \BU_{c_* + \gamma } , \p_x \BU_{c_*+ \gamma} )_H \\ 
( \BU , [\p_c \BU_c]_{|c = c_* + \gamma }( \cdot +y) )_H 
- ( \BW(\gamma) + \alpha \mathbb{B} \BU_{c_* + \gamma }, [\p_c \BU_c]_{|c = c_* + \gamma } )_H \\ 
( \BU , \mathbb{B} \BU_{c_* + \gamma }( \cdot +y) )_H 
- (\BW(\gamma) + \alpha \mathbb{B} \BU_{c_* + \gamma }, \mathbb{B} \BU_{c_* + \gamma } )_H 
\end{array}\right) 
$$
to see that $G$ is indeed of class $\BC^1$ on $ X \times \R \times ( - \gamma_0 , + \gamma_0 ) \times \R $ 
since $ c \mapsto \BU_{c} $ is smooth. Moreover, using that 
$ \p_\gamma \BW_{|\gamma = 0} = [\p_c \BU_c]_{|c=c_* } $, we infer
$$
 \frac{\p G }{\p (y , \gamma , \alpha )} ( \BU_* , 0 , 0 , 0) = 
\left( \begin{array}{ccc} 
 ( \BU_* , \p^2_x \BU_* )_H & 
- ( [\p_c \BU_c]_{|c=c_* } , \p_x \BU_* )_H & 
- ( \mathbb{B} \BU_* , \p_x \BU_* )_H \\ 
- ( \BU_* , \p_x [\p_c \BU_c]_{|c=c_*} )_H & - |\!| [\p_c \BU_c]_{|c=c_* } |\!|^2_H & 
- ( \mathbb{B} \BU_* , [\p_c \BU_c]_{|c=c_*} )_H \\ 
- ( \BU_* , \mathbb{B} \p_x \BU_* )_H & - ( [\p_c \BU_c]_{|c=c_* } ,\mathbb{B} \BU_* )_H & 
- |\!| \mathbb{B} \BU_* |\!|_H^2 
\end{array}\right) .
$$
Here, \cite{Maeda} uses assumption 2 $(iii)$, that reads for us 
$ ( [\p_c \BU_c]_{|c=c_* } , \p_x \BU_* )_H = - ( \p_x [\p_c \BU_c]_{|c=c_* } , \BU_* )_H = 0 $. 
It is indeed the case since $ \BU_c $ is chosen even for any $c$ (close to $c_*$). 
Furthermore, $ ( \BU_* , \p^2_x \BU_* )_H = - |\!| \p_x \BU_* |\!|_H^2 $ 
by integration by parts, $ ( \mathbb{B} \BU_* , \p_x \BU_* )_H 
= ( \BU_* , \mathbb{B} \p_x \BU_* )_H = 0 $ since $ \mathbb{B} \p_x = J $ 
is skew-adjoint, and $ ( \mathbb{B} \BU_* , [\p_c \BU_c]_{|c=c_*} )_H = 
\p_c [ P_{\rm hy} (\BU_c) ]_{|c=c_*} = 0 $ by hypothesis. Therefore,
$$
 \frac{\p G }{\p (y , \gamma , \alpha )} ( \BU_* , 0 , 0 , 0 ) 
= \left( \begin{array}{ccc} 
 - |\!| \p_x \BU_{c_* } |\!|^2_H & 0 & 0 \\ 
0 & - |\!| [\p_c \BU_c]_{|c=c_* } |\!|^2_H & 0 \\
0 & 0 & - |\!| \mathbb{B} \BU_* |\!|_H^2 
\end{array}\right) 
$$
is invertible, thus the implicit function theorem provides three real-valued 
functions $ \underline{y} $, $ \underline{\gamma} $ and $ \underline{\alpha} $, 
defined near $ \BU_*$ (in $X$) and with 
$ \underline{y} (\BU_*) = \underline{\gamma} (\BU_*) = \underline{\alpha}(\BU_*)= 0 $, 
such that $ G ( \BU , \underline{y} (\BU) , \underline{\gamma} (\BU) , \underline{\alpha}(\BU)) = 0 $. 
These functions are extended to $ \mathscr{O}_\e $ (for $\e $ small enough) 
by the formulas $ \bar{y} ( \BU ) \equiv \underline{y} (\BU( \cdot - y)) + y $, 
$ \bar{\gamma} ( \BU ) \equiv \underline{\gamma} (\BU( \cdot - y)) $ and 
$ \alpha ( \BU ) \equiv \underline{\alpha} (\BU( \cdot - y)) $ for any 
$ y \in \R $ such that $ \BU( \cdot - y) $ lies in the neighbourhood of 
$ \BU_*$ where $ \underline{y} $, $ \underline{\gamma} $ and $ \underline{\alpha} $ 
are defined. Consequently, the mapping
$$
 \vartheta ( \BU ) \equiv \BU ( \cdot - \bar{y}( \BU) ) 
- \BW( \bar{\gamma}( \BU ) ) - \bar{\alpha}( \BU ) \mathbb{B} \BU_{c_* + \bar{\gamma }( \BU )} 
$$
is orthogonal in $H$ to $ \p_x \BU_{c_* + \bar{\gamma}( \BU )} $, 
$ [\p_c \BU_c]_{|c = c_* + \bar{\gamma}( \BU )} $ and $ \mathbb{B} \BU_{c_* + \bar{\gamma}( \BU ) } $, 
as desired. Since $ f $ is assumed of class $\BC^2$, we have $ \BU_c \in H^4 $ 
and the regularities $ \mathbb{I}^{-1} \bar{\gamma}' \in H^2 \times H^1 $ and 
$ \mathbb{I}_{H^1}^{-1} \frac{ \p \bar{\gamma} }{ \p \eta } \in H^4 $ follow easily.\carre \\

%%%%%%%%%%%%%%%
\begin{remarks} \rm 
We would like to point out that in \cite{Maeda} (Lemma 3), it is claimed that ``$w(u)$'' 
is orthogonal to ``$ \p_\oo \phi_{\oo + \Lambda(u)} $'' (we refer to the notations 
in \cite{Maeda}). However, since ``$ T(\theta(u)) - \Psi(\Lambda(u) ) $'' is already 
orthogonal to ``$ \p_\oo \phi_{\oo + \Lambda(u)} $'' by construction, this is equivalent to 
``$ \langle B \phi_{\oo + \Lambda(u)} , \p_\oo \phi_{\oo + \Lambda(u)} \rangle = 0 $'', 
or ``$ \p_{\oo'} [ Q( \phi_{\oo'}) ] = 0 $'' at ``$ \oo' = \oo + \Lambda(u) $''. We 
have not understood why this should happen since in general, for the function 
$ \oo' \mapsto Q( \phi_{\oo'} ) $, the point $ \oo $ is the only local critical point. 
For this reason, we have added a component to the original mapping $G$ in \cite{Maeda}. 
Let us observe that then, Lemma 3 in \cite{Maeda} uses the assumption ``$ d''(\oo) = 0 $``. 
On the other hand, the derivative of $G$ in \cite{Maeda} assumes ``$ u \in D(T'(0)) $'', 
for otherwise the expression ``$ G_{1,1} ( u, \theta , \Lambda ) = \langle T'(0) T(\theta) u , 
T'(0) \phi_{\oo + \Lambda} \rangle $'', for instance, is meaningless. We have 
therefore given some details showing clearly the smoothness of $G$.
\end{remarks}
%%%%%%%%%%%%%%%

We now prove a lemma which shows that the quadratic functional associated with 
$ \BF_*'' $ gives a good control on $ \vartheta ( \BU ) $ thanks to the 
orthogonality conditions on this function. This result is in the spirit of 
Lemma 7 in \cite{Ohta}.

%%%%%%%%%%%%%%
\begin{lemmas} 
\label{Maeda4} 
There exist $ 0 < \gamma_1 \leq \gamma_0 $ and $ K_0 > 0 $ such that if 
$ \gamma \in ( - \gamma_1 , + \gamma_1 ) $ and if $ \vartheta \in X $ verifies
$$
 ( \vartheta , \p_x \BU_{c_* + \gamma } )_H 
= ( \vartheta , [\p_c \BU_c]_{|c = c_* + \gamma } )_H 
= ( \vartheta , \mathbb{B} \BU_{c_* + \gamma } )_H = 0 ,
$$
then $ \langle \BF_{c_* + \gamma }'' ( \BU_{c_* + \gamma } ) \vartheta , \vartheta \rangle_{X^*,X} 
\geq K_0 |\!| \vartheta |\!|_X^2 $.
\end{lemmas}
%%%%%%%%%%%%

\noindent {\it Proof.} As a first step, we prove that if $ \vartheta \in X $ 
verifies $ \vartheta \not = 0 $,
$$
 ( \vartheta , \p_x \BU_* )_H = ( \vartheta , [\p_c \BU_c]_{|c = c_* } )_H 
= ( \vartheta , \mathbb{B} \BU_{*} )_H = 0 ,
$$
then $ \langle \BF_*'' ( \BU_*) \vartheta , \vartheta \rangle_{X^*,X} > 0 $. 
Indeed, assume that $ \langle \BF_*'' ( \BU_*) \vartheta , \vartheta \rangle_{X^*,X} \leq 0 $. 
Let $ \chi \in X $ be a negative eigenvector of $ \BF_*''$. We claim that can not 
have $ ( ( \vartheta , \chi )_H , ( [\p_c \BU_c]_{|c=c_*} , \chi )_H ) = ( 0,0 ) $. 
For otherwise, $ ( \vartheta , \chi )_H = 0 $ implies that $ \vartheta $ is 
$L^2$-orthogonal to $ \chi $, which is the eigenvector associated with the only 
negative eigenvalue $ - \mu_0 $ of $ \BF_*'' $ seen as an unbounded operator on $L^2$, 
thus $ \langle \BF_*'' ( \BU_*) \vartheta , \vartheta \rangle_{X^*,X} \geq 0 $, 
and since we assume equality, this means that $ \vartheta $ belongs to 
the kernel of $ \BF_*'' ( \BU_*) $, which is spanned by $ \BU_{*} = \p_x \BU_{*} $, 
but the condition $ ( \vartheta , \p_x \BU_* )_H = 0 $ then implies $ \vartheta = 0 $, 
a contradiction. Therefore, there exists $ (a, b) \in \R^2 $ such that 
$ (a,b) \not = (0,0)$ and $ ( a [\p_c \BU_c]_{|c=c_*} + b \vartheta , \chi )_H = 0 $. 
The nonzero vector $ p \equiv a [\p_c \BU_c]_{|c=c_*} + b \vartheta $ then 
verifies $ ( p , \chi )_H = 0 $ and $ ( p , J \BU_* )_H = a ( [\p_c \BU_c]_{|c=c_*},   J \BU_* )_H 
+ b ( \vartheta , J \BU_* )_H = 0 $, so that 
$ \langle \BF_*'' ( \BU_*) p , p \rangle_{X^*,X} > 0 $. Here, we have used 
once again that $ ( [\p_c \BU_c]_{|c=c_*}, J \BU_* )_H = 0 $ since the left vector 
is an even function and the right vector an odd function. However, in view 
of the equality $ \langle \BF_*'' ( \BU_*) [\p_c \BU_c]_{|c=c_*} , \phi \rangle_{X^*,X} 
= ( \mathbb{B} \BU_* , \phi )_H $, valid for any $ \phi \in X $ (which follows 
from differentiation of $ E'_{\rm hy} (\BU_c) = c P_{\rm hy}'(\BU_c) 
= c ( \mathbb{B} \BU_c , \cdot )_H $ at $c=c_*$)), we have
$$
 \langle \BF_*'' (\BU_*) [\p_c \BU_c]_{|c=c_*} , \vartheta \rangle_{X^*,X} 
= ( \mathbb{B} [\p_c \BU_c]_{|c=c_*} , \vartheta )_H = 0 .
$$
As a consequence,
\begin{align*}
 0 < \langle \BF_*'' ( \BU_*) p , p \rangle_{X^*,X} 
= &\,
a^2 \langle \BF_*'' (\BU_*) [\p_c \BU_c]_{|c=c_*} , [\p_c \BU_c]_{|c=c_*} \rangle_{X^*,X} 
+ b^2 \langle \BF_*'' (\BU_*) \vartheta , \vartheta \rangle_{X^*,X} \\
= &\, 
a^2 ( \mathbb{B} \BU_* , [\p_c \BU_c]_{|c=c_*} )_H 
+ b^2 \langle \BF_*'' (\BU_*) \vartheta , \vartheta \rangle_{X^*,X} 
= b^2 \langle \BF_*'' (\BU_*) \vartheta , \vartheta \rangle_{X^*,X} , 
\end{align*}
since 
$ ( \mathbb{B} \BU_* , [\p_c \BU_c]_{|c=c_*} )_H = \p_c [ P_{\rm hy} (\BU_c) ]_{|c=c_*} = 0 $ 
in our situation. We reach a contradiction since the right-hand side is 
supposed $ \leq 0 $.

We now prove the lemma by contradiction, and then assume that there 
exists sequences $ ( \vartheta_n)_{ n\geq 1 } \in X $ and $ ( \gamma_n)_{n\geq 1} \in 
(0 , \gamma_0 ) $ such that $ \gamma_n \to 0 $, $ |\!| \vartheta_n |\!|_X^2 = 1 $ and
\be
\label{orthogo}
 ( \vartheta_n , \p_x \BU_{c_* + \gamma_n } )_H 
= ( \vartheta_n , [\p_c \BU_c]_{|c = c_* + \gamma_n } )_H 
= ( \vartheta_n , \mathbb{B} \BU_{c_* + \gamma_n } )_H = 0 ,
\ee
but $ \langle \BF_{c_* + \gamma_n }'' ( \BU_{c_* + \gamma_n } ) \vartheta_n , \vartheta_n \rangle_{X^*,X} 
\to 0 $. Possibly passing to a subsequence, we may assume the existence of 
some $ \vartheta = ( \zeta , \upsilon ) \in X $ such that 
$ \vartheta_n \equiv ( \zeta_n , \upsilon_n ) \rightharpoonup \vartheta $ in $ X = H^1 \times L^2 $. 
We then show the lower semicontinuity of $ \langle \BF_*'' (\BU_*) \vartheta , \vartheta \rangle_{X^*,X} $. 
This is roughly a verification of part of assumption (A3) in \cite{Ohta}, used in Lemma 7 
there. By compact Sobolev embedding, we may assume 
$ \zeta_n \to \zeta $ in $ L^\ii _{\rm loc} (\R) $. A straightforward computation gives
\begin{align*}
\langle \BF_{c_* + \gamma }'' ( \BU_{c_* + \gamma } ) \vartheta , \vartheta \rangle_{X^*,X} 
= \int_\R & \, \frac{ (\p_x \zeta )^2 }{ 2(r_0^2 + \eta_{c_* + \gamma }) } 
- \frac{\p_x \zeta \p_x \eta_{c_* +\gamma }  }{ (r_0^2 + \eta_{c_* +\gamma } )^2 } 
+ \frac{ \zeta^2 (\p_x \eta_{c_* +\gamma } )^2 }{4 (r_0^2 + \eta_{c_* +\gamma } )^3 } \\
& + 2( r_0^2 + \eta_{c_* +\gamma }  ) \upsilon^2 
+ 2 ( 2 u_{c_* +\gamma }  - (c_* +\gamma ) ) \upsilon \zeta  
- f' ( r_0^2 + \eta_{c_* +\gamma } ) \zeta^2 \ dx .
\end{align*}
Since $ r_0^2 + \eta_{c_* +\gamma_n } $ remains bounded away from zero uniformly 
and $ \eta_{c_* +\gamma_n } \to \eta_{c_* } $ in $W^{1,\ii}(\R) \cap H^1(\R) $ as 
$ n \to + \ii $, the weak convergence $ \zeta_n \rightharpoonup \zeta $ in $ H^1 $ implies
\begin{align}
\int_\R \frac{ (\p_x \zeta )^2 }{ 2(r_0^2 + \eta_{c_* }) } 
& \, - \frac{\p_x \zeta \p_x \eta_{c_* }  }{ (r_0^2 + \eta_{c_* } )^2 } 
+ \frac{ \zeta^2 (\p_x \eta_{c_*} )^2 }{4 (r_0^2 + \eta_{c_* } )^3 } \ dx \nonumber \\
\label{diri}
& \leq 
\varliminf_{n \to +\ii} \int_\R \frac{ (\p_x \zeta_n )^2 }{ 2(r_0^2 + \eta_{c_* + \gamma_n }) } 
- \frac{\p_x \zeta_n \p_x \eta_{c_* +\gamma_n }  }{ (r_0^2 + \eta_{c_* +\gamma_n } )^2 } 
+ \frac{ \zeta_n^2 (\p_x \eta_{c_* +\gamma_n } )^2 }{4 (r_0^2 + \eta_{c_* +\gamma_n } )^3 } 
\ dx .
\end{align}
For the remaining terms, we write, for some $R> 0$ to be determined later,
\begin{align*}
  \int_\R & \, 2( r_0^2 + \eta_{c_* +\gamma_n }  ) \upsilon_n^2 
+ 2 ( 2 u_{c_* +\gamma_n }  - (c_* +\gamma_n ) ) \upsilon_n \zeta_n   
- f' ( r_0^2 + \eta_{c_* +\gamma_n  } ) \zeta_n ^2 \ dx \\ 
=  & \, \int_\R 2 \Big[ ( r_0^2 + \eta_{c_* +\gamma_n } )^{1/2} \upsilon_n 
+ \frac{ ( 2 u_{c_* +\gamma_n }  - (c_* +\gamma_n ) ) \zeta_n }{2 ( r_0^2 + \eta_{c_* +\gamma_n } )^{1/2} }\Big]^2 \ dx \\
& \quad + \int_{|x| \leq R} +  \int_{|x| \geq R } 
\frac12 \Big( - \frac{ ( 2 u_{c_* +\gamma_n }  - (c_* +\gamma_n ) )^2}{ r_0^2 + \eta_{c_* +\gamma_n } } 
- 2 f' ( r_0^2 + \eta_{c_* +\gamma_n  } )  \Big) \zeta_n ^2 \ dx .
\end{align*}
For the first integral, we may use that $( \zeta_n , \upsilon_n ) \rightharpoonup 
( \zeta , \upsilon ) $ in $L^2 \times L^2$ and the fact that 
$ ( \eta_{c_* +\gamma_n } , u_{c_* +\gamma_n } ) $ converges to $ ( \eta_* , u_*) $ 
uniformly to deduce
\be
\label{cachalot}
( r_0^2 + \eta_{c_* +\gamma_n } )^{1/2} \upsilon_n 
+ \frac{ ( 2 u_{c_* +\gamma_n }  - (c_* +\gamma_n ) ) \zeta_n }{2 ( r_0^2 + \eta_{c_* +\gamma_n } )^{1/2}}
\rightharpoonup ( r_0^2 + \eta_* )^{1/2} \upsilon + \frac{(2u_* - c_* ) \zeta }{2 ( r_0^2 + \eta_* )^{1/2}} 
\quad \quad \quad {\rm in} \quad L^2 ,
\ee
hence,
\begin{align}
\label{baleine}
\int_\R 2 & \, 
\Big[( r_0^2+\eta_* )^{1/2} \upsilon + \frac{ (2u_* -c_* )\zeta}{2 ( r_0^2 + \eta_* )^{1/2}} \Big]^2 \ dx 
\nonumber \\ 
& \leq \varliminf_{n \to +\ii} \int_\R 2 \Big[ ( r_0^2 + \eta_{c_* +\gamma_n } )^{1/2} \upsilon_n 
+ \frac{ ( 2 u_{c_* +\gamma_n}-(c_* +\gamma_n ) ) \zeta_n }{2 ( r_0^2 + \eta_{c_* +\gamma_n } )^{1/2}} \Big]^2 \ dx .
\end{align}
Since $ \zeta_n \to \zeta $ in $ L^\ii([-R,+R]) $ and 
$ ( u_{c_* +\gamma_n } , \eta_{c_* +\gamma_n } ) \to ( u_* , \eta_*) $ uniformly, 
it follows that
\begin{align*}
\int_{|x| \leq R} & \, \frac12 \Big( - \frac{( 2 u_*  - c_*)^2 }{ r_0^2 + \eta_*}
- 2 f' ( r_0^2 + \eta_* )  \Big) \zeta^2 \ dx \\
& = \lim_{n \to +\ii} 
\int_{|x| \leq R} \frac12 \Big( - \frac{ ( 2 u_{c_* +\gamma_n } 
- (c_* +\gamma_n ) )^2 }{ r_0^2 + \eta_{c_* +\gamma_n }}
- 2 f' ( r_0^2 + \eta_{c_* +\gamma_n  } )  \Big) \zeta_n^2 \ dx .
\end{align*}
For the last integral, we have to use the decay at infinity of 
$ \eta_{c_* + \gamma} $ and $ u_{c_* + \gamma} $ uniformly for 
$ | \gamma | $ small. This gives
$$
- \frac{ ( 2 u_{c_* +\gamma_n }  - (c_* +\gamma_n ) )^2 }{ r_0^2 + \eta_{c_* +\gamma_n }  }
- 2 f' ( r_0^2 + \eta_{c_* +\gamma_n  } ) \to \frac{\cs^2 - c_*^2 }{r_0^2} 
$$
as $ |x| \to + \ii $, uniformly in $n$. Since $ 0 < c_* < \cs $, there exist some 
small $\d > 0 $ and some $ R> 0$ large such that, for any $n$ and 
any $x$ with $|x| \geq R$,
$$
- \frac{ ( 2 u_{c_* +\gamma_n }  - (c_* +\gamma_n ) )^2 }{ r_0^2 + \eta_{c_* +\gamma_n }  }
- 2 f' ( r_0^2 + \eta_{c_* +\gamma_n  } ) \geq \d .
$$
In particular, since $ \zeta_n \rightharpoonup \zeta $ in $ L^2 $,
$$
{\bf 1}_{|x| \geq R} \Big( - \frac{ ( 2 u_{c_* +\gamma_n }  - (c_* +\gamma_n ) )^2 }
{ r_0^2 + \eta_{c_* +\gamma_n }} 
- 2 f' ( r_0^2 + \eta_{c_* +\gamma_n  } ) \Big)^{1/2} \zeta_n \rightharpoonup 
{\bf 1}_{|x| \geq R} 
\Big(- \frac{ ( 2 u_* - c_*)^2 }{ r_0^2 + \eta_* }- 2 f' ( r_0^2 + \eta_* ) \Big)^{1/2} \zeta ,
$$
in $ L^2 $, thus
\begin{align}
\label{faiblounet}
\int_{|x| \geq R} & \, \frac12 \Big( - \frac{( 2 u_*  - c_*)^2 }{ r_0^2 + \eta_*}
- 2 f' ( r_0^2 + \eta_* )  \Big) \zeta^2 \ dx \nonumber \\
& \leq \varliminf_{n \to +\ii} 
\int_{|x| \geq R} \frac12 \Big( - \frac{ ( 2 u_{c_* +\gamma_n } 
- (c_* +\gamma_n ) )^2 }{ r_0^2 + \eta_{c_* +\gamma_n }}
- 2 f' ( r_0^2 + \eta_{c_* +\gamma_n  } )  \Big) \zeta_n^2 \ dx .
\end{align}
Combining these three $ \varliminf $ inequalities, we deduce
\be
\label{semi}
\langle \BF_*'' ( \BU_*) \vartheta , \vartheta \rangle_{X^*,X} 
\leq \varliminf_{n \to +\ii}\langle \BF_{c_* + \gamma_n }'' ( \BU_{c_* + \gamma_n } ) \vartheta_n , 
\vartheta_n \rangle_{X^*,X} = 0 .
\ee

Turning back to our sequence $( \vartheta_n , \gamma_n )$, we may pass 
to the limit in \eqref{orthogo}:
$$
 ( \vartheta , \p_x \BU_* )_H 
= ( \vartheta , [\p_c \BU_c]_{|c = c_* } )_H 
= ( \vartheta , \mathbb{B} \BU_* )_H = 0 .
$$
Comparing with \eqref{semi}, we deduce from our first claim that $ \vartheta = 0 $. 
This means that we must have equality in all the above $\varliminf$ inequalities. 
In particular, the weak convergence \eqref{faiblounet} is actually strong, 
thus $ \zeta_n \to \zeta = 0 $ in $ L^2 (\R) $ (the strong convergence 
in $ \{ |x| \leq R \} $ being already known since $ \zeta_n \to \zeta $ in 
$ L^\ii_{\rm loc} (\R) $). Going back to the equality in \eqref{diri} thus provides
$ \p_x \zeta_n \to \p_x \zeta = 0 $ in $ L^2(\R) $, since $ r_0^2 + \eta_{ c_* + \gamma_n} $ 
remains uniformly bounded away from zero and by weak convergence, 
$ \ds{ 0 = \int_\R \frac{ \zeta^2 (\p_x \eta_{c_*} )^2 }{4 (r_0^2 + \eta_{c_* } )^3 } \ dx 
= \lim_{n \to +\ii} \int_\R 
\frac{\p_x \zeta_n \p_x \eta_{c_* +\gamma_n}}{ (r_0^2 + \eta_{c_* +\gamma_n } )^2 } 
\ dx } $. Finally, the equality in \eqref{baleine} means that \eqref{cachalot} 
is actually a strong convergence, that is $ \upsilon_n \to \upsilon = 0 $ in $ L^2 $ 
since $ \zeta_n \to \zeta $ in $ L^2 $. The contradiction then follows: 
$ 1 = |\!| \vartheta_n |\!|_X^2 = |\!| \zeta_n |\!|_{L^2}^2 + |\!| \p_x \zeta_n |\!|_{L^2}^2 
+ |\!| \upsilon_n |\!|_{L^2}^2 \to 0 $. \carre

%%%%%%%%%%%%%%%
\begin{remarks} \rm 
This result is also Lemma 7 in \cite{Maeda}, and is said to be Lemma 7 
in \cite{Ohta}. However, the hypothesis of Lemma 7 in \cite{Ohta} are 
not verified, and in particular assumption (B3) there. It is natural to believe that 
this assumption is verified is most physical situations, but it is not clear 
whether it always holds true in the general framework of \cite{Maeda} 
without further hypothesis.
\end{remarks}
%%%%%%%%%%%%%

The next lemma provides a control for $ \alpha( \BU ) $.

%%%%%%%%%%%%%%
\begin{lemmas} 
\label{Maeda5} 
Assume $ \e > 0 $ small enough. Then, there exists $ K > 0 $ such that for any 
$ \BU \in \mathscr{O}_\e $ satisfying $ P_{\rm hy}( \BU ) = P_{\rm hy}( \BU_* ) $, 
there holds
$$
 | \alpha( \BU ) | \leq K \Big( \bar{\gamma}^2(\BU) \n \vartheta ( \BU ) \n_X 
+ \n \vartheta ( \BU ) \n_X^2 \Big) .
$$
\end{lemmas}
%%%%%%%%%%%%

\noindent {\it Proof.} It is the same as in \cite{Maeda}, Lemma 8, but give it for 
completeness. We expand and use that $ \mathbb{B}^2 = {\rm Id}_2 $ and 
the definition 
$ \BW ( \gamma ) \equiv \BU_{c_* + \gamma } + \s ( \gamma) \mathbb{B} \BU_{c_* + \gamma } $ 
for the second line:
\begin{align*}
 P_{\rm hy}( \BU_* ) = & \, P_{\rm hy}( \BU ) = P_{\rm hy}( \BU ( \cdot - \bar{y}( \BU) ) )  
= P_{\rm hy} \Big( \BW( \bar{\gamma}( \BU ) ) 
+ \vartheta ( \BU ) + \alpha ( \BU ) \mathbb{B} \BU_{c_* + \bar{\gamma}( \BU ) } \Big) \\ 
= & \, P_{\rm hy}( \BW( \bar{\gamma}( \BU ) ) ) + P_{\rm hy}( \vartheta ( \BU ) ) 
+ \alpha^2(\BU) P_{\rm hy}(\mathbb{B} \BU_{c_* + \bar{\gamma}( \BU ) } ) \\
& \, + \alpha(\BU) ( \mathbb{B} \vartheta ( \BU ) , \mathbb{B} \BU_{c_* + \bar{\gamma}( \BU ) } )_H 
+ ( \mathbb{B} \BU_{c_* + \bar{\gamma}( \BU ) } , \vartheta ( \BU ) )_H 
+ \alpha(\BU) ( \mathbb{B} \BU_{c_*+\bar{\gamma}(\BU)}, 
    \mathbb{B} \BU_{c_*+\bar{\gamma}( \BU ) } )_H \\
& \, + \sigma(\bar{\gamma}(\BU)) ( \BU_{c_*+\bar{\gamma}( \BU ) } , \vartheta (\BU) )_H 
+ \sigma(\bar{\gamma}(\BU)) \alpha(\BU) (\BU_{c_*+\bar{\gamma}(\BU)} , 
    \mathbb{B} \BU_{c_*+\bar{\gamma}(\BU)} )_H .
\end{align*}
Since $ P_{\rm hy}( \BW( \bar{\gamma}( \BU ) ) ) = P_{\rm hy}( \BU_* ) $, we infer 
$$
 - \alpha(\BU) \Big[ \n \BU_{*} \n_H^2 + o(1) \Big] = 
\sigma(\bar{\gamma}(\BU)) ( \BU_{c_*+\bar{\gamma}( \BU ) } , \vartheta (\BU) )_H 
+ P_{\rm hy}( \vartheta ( \BU ) ) 
$$
and the conclusion follows since $ \s(\gamma) = \BO(\gamma^2)$ by Lemma \ref{Maeda1}.\carre \\

Now, we give a lemma useful to estimate $ \vartheta( \BU ) $.

%%%%%%%%%%%%%%
\begin{lemmas} 
\label{Maeda6} 
Assume $ \e > 0 $ small enough. Then, there exists $ K > 0 $ such that for any 
$ \BU \in \mathscr{O}_\e $ satisfying $ P_{\rm hy}( \BU ) = P_{\rm hy}( \BU_* ) $ 
and $ \BF_* ( \BU ) - \BF_* ( \BU_* ) < 0 $, there holds
$$
 |\!| \vartheta(\BU) |\!|_X^2 \leq K |\bar{\gamma}(\BU) |^3 .
$$
In particular, $ | \alpha( \BU ) | \leq K |\bar{\gamma}(\BU) |^3 $. 
\end{lemmas}
%%%%%%%%%%%%

\noindent {\it Proof.} It is the same as in \cite{Maeda}, Lemma 9. Note first that 
the last assertion is a direct consequence of the first one and Lemma \ref{Maeda5}. 
Next, we argue by contradiction and assume that there exists a sequence 
$ \BU_n \to \BU_* $ in $X$ such that $ \BF_* ( \BU_n ) - \BF_* (\BU_*) < 0 $ and 
$ |\!| \vartheta(\BU_n) |\!|_X^2 \gg |\bar{\gamma}(\BU_n) |^3 $. For simplicity, 
we denote $ \bar{\gamma}_n = \bar{\gamma}(\BU_n)$, $ \vartheta_n = \vartheta ( \BU_n ) $, 
$ \alpha_n = \alpha ( \BU_n ) $. Then, by Lemma \ref{Maeda5}, we have 
$ | \alpha_n | 
\leq K ( \bar{\gamma}^2_n |\!| \vartheta_n \n_X + |\!| \vartheta_n |\!|_X^2 ) 
\leq K ( |\!| \vartheta_n |\!|_X^{7/3} + |\!| \vartheta_n |\!|_X^2 )= \BO( |\!| \vartheta_n |\!|_X^2) $. 
Therefore, by Taylor expansion and Lemma \ref{Maeda3}, it holds
\begin{align}
 \BF_* (\BU_n) - \BF(\BU_*) 
 = & \, 
 \BF_* (\BU_n(\cdot - \bar{y}_n )) - \BF(\BU_*) 
 = \BF_* ( \BW( \bar{\gamma}_n )) 
+ \vartheta_n + \alpha_n \mathbb{B} \BU_{c_* + \bar{\gamma}_n } ) - \BF(\BU_*) 
\nonumber \\ \label{Ford}
 = & \, 
 \BF_* ( \BW( \bar{\gamma}_n ) ) - \BF(\BU_*) 
+ \langle \BF'_* ( \BW( \bar{\gamma}_n) ) , 
		\vartheta_n + \alpha_n \mathbb{B} \BU_{c_* + \bar{\gamma}_n } \rangle_{X^*,X} \\
& \, + \frac12 \langle \BF''_* ( \BW( \bar{\gamma}_n )) \vartheta_n ,  \vartheta_n \rangle_{X^*,X} 
+ o( \n \vartheta (\BU_n) \n_X^2 ) .
\nonumber
\end{align}
However, by Lemma \ref{Maeda2}, $ \BF_{c_* } ( \BW( \gamma ) ) - \BF_{c_* } ( \BU_{c_* } ) 
= \BO ( | \gamma|^3) $, and since $ \BF'_* ( \BW( \gamma ) ) = 
\BF'_* ( \BW(0) ) + o(1) = \BF'_* ( \BU_* ) + o(1) = o(1) $, there holds 
$ \langle \BF'_* ( \BW( \bar{\gamma}_n ) ) , \alpha_n \mathbb{B} \BU_{c_* + \bar{\gamma}_n } \rangle_{X^*,X} 
= o ( |\!| \vartheta_n |\!|_X^2 ) $. Furthermore, using 
$ \BF'_* = \BF'_{c_* +  \bar{\gamma}_n} +  \bar{\gamma}_n \mathbb{B} $, 
the third orthogonality condition in Lemma \ref{Maeda3} and that 
$ \s(\gamma ) = \BO(\gamma^2 ) $, we deduce 
\begin{align*}
\langle \BF'_*& \, ( \BW( \bar{\gamma}_n ) ) , \vartheta_n \rangle_{X^*,X} 
= 
\langle \BF'_{c_* +  \bar{\gamma}_n} ( \BW( \bar{\gamma}_n)) , \vartheta_n \rangle_{X^*,X} 
+ \bar{\gamma}_n ( \mathbb{B} \BW( \bar{\gamma}_n) , \vartheta_n )_{H} 
\\ 
= & \, 
\langle \BF'_{c_* +  \bar{\gamma}_n} ( \BU_{c_* + \bar{\gamma}_n} )
+ \s( \bar{\gamma}_n) \mathbb{B} \BU_{c_* + \bar{\gamma}_n} , \vartheta_n \rangle_{X^*,X} 
+ \bar{\gamma}_n \s ( \bar{\gamma}_n ) ( \BU_{c_*+ \bar{\gamma}_n} , 
\vartheta_n )_{H} \\ 
= & \, \BO( \bar{\gamma}^2_n \n \vartheta_n \n_X ) 
+ \BO( | \bar{\gamma}_n |^3 \n \vartheta_n \n_X ) 
= o( \n \vartheta_n \n_X^{7/3 } ) = o ( \n \vartheta_n \n_X^2 ) .
\end{align*}
For the last line, we have used another Taylor expansion with 
$ \BF'_{c_* +  \bar{\gamma}_n} ( \BU_{c_* + \bar{\gamma}_n}  ) = 0 $. Finally, 
Lemma \ref{Maeda4} yields $ \langle \BF_{c_* + \bar{\gamma}_n }'' 
(\BU_{c_*+\bar{\gamma}_n}) \vartheta_n , \vartheta_n \rangle_{X^*,X} 
\geq K_0 |\!| \vartheta_n |\!|_X^2 $. Reporting these expansions in \eqref{Ford} 
yields
$$
 \BF_* (\BU_n) - \BF(\BU_*) \geq \frac{K_0}{4} \n \vartheta_n \n_X^2 
 + o ( \n \vartheta_n \n_X^2 ) \geq \frac{K_0}{8} \n \vartheta_n \n_X^2 
$$
for $n$ sufficiently large, which contradicts our assumption. \carre \\

We now need to find an extension of the functionals "$A$" and "$P$" used 
in \cite{Maeda} (and also in \cite{Ohta}). In these works, these functionals are 
built on what should be here "$ J^{-1} \p_c \BU_c = \mathbb{B} \p_x^{-1} \p_c \BU_c $",  
but unfortunately, $ \p_c \phi_c $ does not have vanishing integral over $\R$ 
(for instance, $\p_c \eta_c$ has constant sign). 
We rely instead on a construction of a suitable approximation of 
"$ J^{-1} \p_c \BU_c $". A similar construction is used in \cite{Lin}.

%%%%%%%%%%%%%%
\begin{lemmas} 
\label{Maeda7} 
For any $ 0 < \kappa < 1 $, there exists a $ \BC^2 $ mapping 
$ \Upsilon_\kappa : ( - \gamma_1, + \gamma_1 ) \to X $ such that, 
for any $ \gamma \in ( - \gamma_1, + \gamma_1 ) $, 
$ \Upsilon_\kappa (\gamma ) \in H^2 \times H^1 $ is an odd function verifying
$$
 \n J \Upsilon_\kappa (\gamma) - [ \p_c \BU_c]_{|c = c_*+ \gamma} \n_X 
 \leq \kappa .
$$
\end{lemmas}
%%%%%%%%%%%%

\noindent {\it Proof.} We fix an even function $ \Theta_0 \in \BC_c^\ii (\R ) $ 
such that $ \int_\R \Theta_0 \ dx = 1 $. For $ T > 0 $ to be fixed later, but 
independent of $\gamma$ and $\kappa$, we set 
$ t_\kappa \equiv T / \kappa^2 > 0 $ and
$$
 \Upsilon_\kappa ( \gamma ) (x) \equiv \mathbb{B} 
 \int_0^x \Big[ [ \p_c \BU_c]_{|c = c_*+ \gamma} (y) 
 - \frac{1}{t_\kappa} \Theta_0 \Big( \frac{y}{t_\kappa} \Big) 
 \int_\R [\p_c \BU_c]_{|c = c_*+ \gamma} (z) \ dz \Big] \ dy .
$$
It is clear that $ \Upsilon_\kappa ( \gamma ) \in \BC^1 (\R) $ and that, 
since $ J = \p_x \mathbb{B} $ and $  \mathbb{B}^2 = {\rm Id}_2 $,
$$
 J \Upsilon_\kappa ( \gamma ) - [ \p_c \BU_c]_{|c = c_*+ \gamma}  = 
 \frac{1}{t_\kappa} \Theta_0 \Big( \frac{\cdot}{t_\kappa} \Big) 
 \int_\R [\p_c \BU_c]_{|c = c_*+ \gamma} (z) \ dz
$$
In particular,
$$ \n J \Upsilon_\kappa (\gamma ) - \p_c [ \BU_c]_{|c = c_*+ \gamma} \n_X^2 
= \Big[ \frac{1}{t_\kappa} \n \Theta_0 \n_{L^2}^2 + \frac{1}{t_\kappa^3} \n \p_x \Theta_0 \n_{L^2}^2 \Big] 
\Big( \int_\R [\p_c \BU_c]_{|c = c_*+ \gamma} (z) \ dz \Big)^2 \leq \kappa^2 
$$
if we choose $ T = T( c_*, \BU_* , \Theta_0 ) > 0 $ sufficiently large 
%$ T \equiv (1/2) |\!| \Theta_0 |\!|_{L^2}^2 (  \int_\R [\p_c \BU_c]_{|c = c_*} (z) \ dz )^2 $ 
and $ \gamma_1 $ smaller if necessary. Moreover, $  \Upsilon_\kappa ( \gamma ) $ 
is odd since $ \BU_c $ and $ \Theta_0 $ are even. In addition, the even function 
$ y \mapsto [ \p_c \BU_c]_{|c = c_*+ \gamma} (y) 
 - \frac{1}{t_\kappa} \Theta_0 \Big( \frac{y}{t_\kappa} \Big) 
 \int_\R [\p_c \BU_c]_{|c = c_*+ \gamma} (z) \ dz $ decays exponentially at 
infinity (since $ \Theta_0 $ has compact support and $ \p_c \BU_c $ decays 
exponentially), and has zero integral (since $ \Theta_0 $ has integral equal to one), hence 
$$
\Upsilon_\kappa ( \gamma ) (x) = - \mathbb{B} 
 \int_x^{+\ii} \Big[ [ \p_c \BU_c]_{|c = c_*+ \gamma} (y) 
 - \frac{1}{t_\kappa} \Theta_0 \Big( \frac{y}{t_\kappa} \Big) 
 \int_\R [\p_c \BU_c]_{|c = c_*+ \gamma} (z) \ dz \Big] \ dy 
$$
and decays exponentially at infinity. It follows easily from these two equalities 
that $ \gamma \mapsto \Upsilon_\kappa (\gamma) \in L^2 \times L^2 $ is well 
defined and continuous, hence also $ \gamma \mapsto \Upsilon_\kappa (\gamma) \in H^2 \times H^1 $. 
%, thus belongs to $ H = L^2 \times L^2 $. 
By the same type of arguments,
$$
\frac{ \p \Upsilon_\kappa}{\p \gamma } ( \gamma ) (x) = \mathbb{B} 
 \int_0^x \Big[ [ \p^2_c \BU_c]_{|c = c_*+ \gamma} (y) 
 - \frac{1}{t_\kappa} \Theta_0 \Big( \frac{y}{t_\kappa} \Big) 
 \int_\R [\p^2_c \BU_c]_{|c = c_*+ \gamma} (z) \ dz \Big] \ dy 
$$
is well-defined and is a continuous function of $\gamma$ with values into $H^2 \times H^1$, 
and similarly for the second derivative. \carre \\

We now define, in the tubular neighbourhood $ \mathscr{O}_\e $ of $ \BU_* $, 
the functional (corresponding to "$A$" in \cite{Maeda})
$$
 \OO_\kappa ( \BU) \equiv ( \BU ( \cdot - \bar{y}(\BU) ) , 
 \Upsilon_\kappa ( \bar{\gamma}(\BU) ) )_H 
 = ( \BU , \Upsilon_\kappa ( \bar{\gamma}(\BU) ) ( \cdot + \bar{y}(\BU) ) )_H
 $$
depending on $ \kappa \in (0,1 ) $, which will be determined later. 
The first properties of $ \OO_\kappa $ are given below.

%%%%%%%%%%%%%%
\begin{lemmas} 
\label{Maeda8}  
For any $ 0 < \kappa < 1 $, $ \OO_\kappa : \mathscr{O}_\e \to \R $ is 
of class $ \BC^1 $. In addition, there exists some bounded mapping 
$ \BN_{\bar{\gamma}} : \mathscr{O}_\e \to X $ such that if 
$ \Psi_{\rm hy} \in \BC^1 ( [ 0, T) , X ) $ is a solution to 
\eqref{hydrohamilto} that remains in $ \mathscr{O}_\e $, then 
$$
 \frac{d}{dt} \OO_\kappa ( \Psi_{\rm hy} (t) ) = \Xi_\kappa ( \Psi_{\rm hy} (t) ) ,
$$
where $ \Xi_\kappa : \mathscr{O}_\e \to \R $ is defined by
$$
 \Xi_\kappa (\BU) \equiv 
- \Big\langle \BF_{c_*+ \bar{\gamma}(\BU)} ' (\BU) , 
\Big\{ J \Upsilon_\kappa ( \bar{\gamma}(\BU)) (\cdot + \bar{y}(\BU)) 
+ ( \BU, \p_\gamma \Upsilon_\kappa (\bar{\gamma}(\BU)) (\cdot + \bar{y}(\BU) ) )_H 
\BN_{\bar{\gamma}} (\BU) \Big\} \Big\rangle_{X^*,X} .
$$
\end{lemmas}
%%%%%%%%%%%%

\noindent {\it Proof.} The fact that $ \OO_\kappa $ is of class $ \BC^1 $ follows 
directly from the second expression and the fact that $ \bar{y} $ and $ \bar{\gamma} $ 
are $ \BC^1 $ (in \cite{Maeda} formula (3.11), the same 
remark as for the smoothness of $G$ after Lemma \ref{Maeda3} holds, since 
it requires "$ u \in D(T'(0)) $"). If $ \Psi_{\rm hy} = ( \eta , u) \in \BC^1 ([0,T),X)$ 
is a solution to \eqref{hydrohamilto} that remains in $ \mathscr{O}_\e $, we therefore 
have, denoting $ \bar{\gamma}(t) = \bar{\gamma}(\Psi_{\rm hy} (t)) $ and 
$ \bar{y}( t) = \bar{y}(\Psi_{\rm hy} (t)) $,
\begin{align}
\label{OOdot}
 \frac{d}{dt} \OO_\kappa ( \Psi_{\rm hy} (t) ) 
 = & \, ( \p_t \Psi_{\rm hy} (t) , \Upsilon_\kappa ( \bar{\gamma}(t) ) 
 ( \cdot + \bar{y}(t) ) )_H 
\nonumber \\ 
& \, + 
( \Psi_{\rm hy} (t) , \p_x \Upsilon_\kappa (\bar{\gamma}(t)) (\cdot + \bar{y}(t) ) )_H 
\langle \bar{y}' (\Psi_{\rm hy}(t)) , \p_t \Psi_{\rm hy}(t) \rangle_{X^*,X} \\ 
& \, + 
( \Psi_{\rm hy} (t) , \p_\gamma \Upsilon_\kappa ( \bar{\gamma}(t) ) 
(\cdot + \bar{y}(t) ) )_H 
\langle \bar{\gamma}' (\Psi_{\rm hy} (t)) , \p_t \Psi_{\rm hy} (t) \rangle_{X^*,X} .
\nonumber 
\end{align}
We now observe that the invariance of $ \OO_\kappa $ by translation 
provides by differentiation the equality, for $ \BU \in \mathscr{O}_\e $,
\begin{align}
\label{trans} 
 0 = & \, \frac{d}{dy} \OO_\kappa (\BU ( \cdot - y) )_{|y=0} 
= ( \BU , \p_x \Upsilon_\kappa ( \bar{\gamma}(\BU) ) ( \cdot + \bar{y}(\BU) ) )_H 
= ( \mathbb{B} \BU , J \Upsilon_\kappa ( \bar{\gamma}(\BU) ) ( \cdot + \bar{y}(\BU) ) )_H 
\nonumber \\ 
= & \, \langle P_{\rm hy}' (\BU) , J \Upsilon_\kappa ( \bar{\gamma}(\BU) ) ( \cdot + \bar{y}(\BU) ) 
\rangle_{X^*,X} .
\end{align}
In particular, the second term in \eqref{OOdot} vanishes. In addition, since 
$ \Psi_{\rm hy} = ( \eta , u) \in \BC^1 ( [ 0, T) , X ) $ a solution to 
\eqref{hydrohamilto} that remains in $ \mathscr{O}_\e $, there holds, 
denoting by $ \ds{ \frac{\d E_{\rm hy}}{\d \Psi} } $ the variational derivative,
\begin{align*}
( \p_t & \Psi_{\rm hy} (t) , \Upsilon_\kappa ( \bar{\gamma}(t)) ( \cdot + \bar{y}(t)) )_H 
= \Big( J \frac{\d E_{\rm hy}}{\d \Psi_{\rm hy}}(\Psi_{\rm hy}(t)) , 
\Upsilon_\kappa ( \bar{\gamma}(t)) ( \cdot + \bar{y}(t)) \Big)_H 
\\ & = - \Big( \frac{\d E_{\rm hy}}{\d \Psi_{\rm hy}}(\Psi_{\rm hy}(t)) , 
J \Upsilon_\kappa ( \bar{\gamma}(t)) ( \cdot + \bar{y}(t)) \Big)_H 
= 
- \Big\langle E'_{\rm hy}(\Psi_{\rm hy}(t)) , 
J \Upsilon_\kappa ( \bar{\gamma}(t)) ( \cdot + \bar{y}(t)) \Big\rangle_{X^*,X} 
\\ & = 
- \langle \BF_{c_* + \bar{\gamma}(t)}' (\Psi_{\rm hy}(t)) , J \Upsilon_\kappa ( \bar{\gamma}(t) ) 
 ( \cdot + \bar{y}(t) ) \rangle_{X^*,X} 
- ( c_*+ \bar{\gamma}(t) ) \Big\langle P'_{\rm hy}(\Psi_{\rm hy}(t)) , 
J \Upsilon_\kappa ( \bar{\gamma}(t)) ( \cdot + \bar{y}(t)) \Big\rangle_{X^*,X} 
\\ & = 
- \langle \BF_{c_* + \bar{\gamma}(t)}' (\Psi_{\rm hy}(t)) , J \Upsilon_\kappa ( \bar{\gamma}(t) ) 
 ( \cdot + \bar{y}(t) ) \rangle_{X^*,X} ,
\end{align*}
by \eqref{trans}. In addition, since 
$ \ds{ \frac{\d P_{\rm hy}}{\d \Psi_{\rm hy}} = \mathbb{B} \Psi_{\rm hy} } $ and 
$ J \mathbb{B} = \p_x $,
\begin{align*}
\langle \bar{\gamma}' (\Psi_{\rm hy}(t)) ,& \, \p_t \Psi_{\rm hy}(t) \rangle_{X^*,X} 
= \Big\langle \bar{\gamma}' (\Psi_{\rm hy}(t)) , 
J \frac{\d E_{\rm hy}}{\d \Psi_{\rm hy}}(\Psi_{\rm hy}(t)) \Big\rangle_{X^*,X} \\
= & \, \Big\langle \bar{\gamma}' (\Psi_{\rm hy}(t)) , 
J \frac{\d \BF_{c_*+ \bar{\gamma}(\BU)}}{\d \Psi_{\rm hy}}(\Psi_{\rm hy}(t)) \Big\rangle_{X^*,X} 
+ ( c_*+ \bar{\gamma}(\BU) ) \Big\langle \bar{\gamma}' (\Psi_{\rm hy}(t)) , 
\p_x \Psi_{\rm hy}(t) \Big\rangle_{X^*,X} .
\end{align*}
The second term vanishes since $ \bar{\gamma} $ is invariant by 
translation (by definition, see the proof of Lemma \ref{Maeda3}). As a consequence,
\begin{align*}
\langle \bar{\gamma}' (\Psi_{\rm hy}(t)) , \p_t \Psi_{\rm hy}(t) \rangle_{X^*,X} 
= & \, \Big( J \frac{\d \BF_{c_*+ \bar{\gamma}(t)}}{\d \Psi_{\rm hy}}(\Psi_{\rm hy}(t)) , 
    \mathbb{I}^{-1} \bar{\gamma}' (\Psi_{\rm hy}(t)) \Big)_X \\
= &\, 
- \Big( \frac{\d \BF_{c_*+ \bar{\gamma}(t)}}{\d \Psi_{\rm hy}}(\Psi_{\rm hy}(t)), 
    J \mathbb{I}^{-1} \bar{\gamma}' (\Psi_{\rm hy}(t)) \Big)_X \\
= &\, - \Big( \frac{\d \BF_{c_*+ \bar{\gamma}(t)}}{\d \Psi_{\rm hy}}(\Psi_{\rm hy}(t)) , 
    J \mathbb{I}^{-1} \bar{\gamma}' (\Psi_{\rm hy}(t)) \Big)_H \\
&\, 
- \Big( \p_x \frac{\d \BF_{c_*+ \bar{\gamma}(t)}}{\d \eta} (\Psi_{\rm hy}(t)), 
    \p_x J \mathbb{I}_{H^1}^{-1} \frac{\p \bar{\gamma}}{\p \eta} (\Psi_{\rm hy}(t)) \Big)_{L^2} .
\end{align*}
The first term is simply 
$ - \langle \BF_{c_*+ \bar{\gamma}(\BU)}'(\Psi_{\rm hy}(t)), J\mathbb{I}^{-1} \bar{\gamma}'(\Psi_{\rm hy}(t)) \rangle_{X^*,X} $. 
We then define $ \BN_{\bar{\gamma}} : \mathscr{O}_\e \to X $ by 
$ \BN_{\bar{\gamma}} (\BU) \equiv J \mathbb{I}^{-1} \bar{\gamma}' (\BU) 
- ( \p_x^2 J \mathbb{I}_{H^1}^{-1} \frac{\d \bar{\gamma}}{\d \eta} (\BU) , 0 ) \in X = H^1 \times L^2 $ 
(see the regularity shown for $ \bar{\gamma}'$ in Lemma \ref{Maeda3}), so that integration 
by parts yields
$$
\langle \bar{\gamma}' (\Psi_{\rm hy}(t)) , \p_t \Psi_{\rm hy}(t) \rangle_{X^*,X} 
= - \langle \BF_{c_*+ \bar{\gamma}(\BU)}'(\Psi_{\rm hy}(t)), \BN_{\bar{\gamma}}(\Psi_{\rm hy}(t)) \rangle_{X^*,X} .
$$
Inserting these relations into \eqref{OOdot} then gives
\begin{align*}
 \frac{d}{dt} \OO_\kappa ( \Psi_{\rm hy} (t) ) 
= - \Big\langle & \, \BF_{c_*+ \bar{\gamma}(t)}' (\Psi_{\rm hy} (t)) , \\
& \Big\{ J \Upsilon_\kappa ( \bar{\gamma}(t)) (\cdot + \bar{y}(t)) 
+ ( \Psi_{\rm hy}(t), \p_\gamma \Upsilon_\kappa (\bar{\gamma}(t)) (\cdot + \bar{y}(t) ) )_H 
\BN_{\bar{\gamma}} (\Psi_{\rm hy}(t))
 \Big\} \Big\rangle_{X^*,X} ,
\end{align*}
which is the desired equality.

If $ \Psi_{\rm hy} \in \BC^0 ( [ 0, T) , X ) $ is just a continuous in time solution 
to \eqref{hydrohamilto} that remains in $ \mathscr{O}_\e $, then the 
integrated relation
$$
\OO_\kappa ( \Psi_{\rm hy} (t) ) = \OO_\kappa ( \Psi_{\rm hy}^\ini ) 
+ \int_0^t \Xi_\kappa( \Psi_{\rm hy} (\tau) ) \ d \tau
$$
holds, as can be seen by using the continuity of the flow and the 
approximation of such a solution by smoother ones (see \cite{GaZhi}). \carre \\

We now compute the asymptotics of $ \Xi_\kappa ( \BW( \gamma ) ) $ for 
$ \gamma \to 0 $ and small $ \kappa $. 

%%%%%%%%%%%%%%
\begin{lemmas} 
\label{Maeda9} 
We have
$$
 \Xi_\kappa ( \BW( \gamma ) ) = - \frac{\gamma^2 \ddot{P}_*}{2 } 
  + o_{ ( \gamma , \kappa ) \to ( 0,0) }( \gamma^2 ) . 
$$
\end{lemmas}
%%%%%%%%%%%%

\noindent {\it Proof.} The proof follows the one of Lemma 5 in \cite{Maeda}. 
As a first step, notice that $ \bar{\gamma}(\BW(\gamma)) = \gamma $, 
$ \bar{y}(\BW(\gamma)) = 0 $, as it can be seen from the equality 
$ G( \BW(\gamma) , 0,0, 0 ) = 0 $ and the local uniqueness of the 
solution to $ G = 0 $. Therefore, since 
$ \BF'_{c_*+\gamma} (\BU_{c_*+\gamma}) = 0 $ and 
$ \s(\gamma) \sim - \gamma^2 \ddot{P}_* / (2 |\!| \BU_* |\!|^2_H) $,
\begin{align*}
 \BF'_{c_* + \bar{\gamma}(\BW(\gamma)) } ( \BW(\gamma) ) 
= & \, \BF'_{c_* + \gamma} ( \BU_{ c_* + \gamma} + \sigma(\gamma) \mathbb{B} \BU_{ c_* + \gamma} ) 
= \sigma(\gamma) \BF''_{c_* + \gamma} ( \BU_{ c_* + \gamma} ) [ \mathbb{B} \BU_{ c_* + \gamma} ] 
    + o_{ \gamma \to 0}( \gamma^2 ) \\ 
= & \, - \frac{\gamma^2 \ddot{P}_*}{ 2 |\!| \BU_* |\!|^2_H } 
\BF''_* ( \BU_* )[ \mathbb{B} \BU_* ] + o_{ \gamma \to 0}( \gamma^2 ) .
\end{align*}
In addition, since $\BU_c$ is even and $\Upsilon_\kappa (\gamma) $ is odd, 
we deduce
$$
 ( \BW(\gamma) ( \cdot + \bar{y}(\BW(\gamma))) , \p_\gamma \Upsilon_\kappa (\gamma) )_H
 = ( \BU_{c_* + \gamma } + \s ( \gamma) \mathbb{B} \BU_{c_* + \gamma } , \p_\gamma \Upsilon_\kappa (\gamma) )_H
= 0 .
$$
Consequently,
$$
 \Xi_\kappa (\BW(\gamma)) =
\frac{\gamma^2 \ddot{P}_*}{ 2 |\!| \BU_* |\!|^2_H } 
\langle \BF''_* ( \BU_* )[ \mathbb{B} \BU_* ] , J \Upsilon_\kappa ( \gamma ) \rangle_{X^*,X} 
+ o_{ \gamma \to 0}( \gamma^2 ) ,
$$
where ``$ o_{ \gamma \to 0}( \gamma^2 ) $'' does not depend on $\kappa$. Moreover, 
Lemma \ref{Maeda7} provides $ |\!| J \Upsilon_\kappa (\gamma ) - [\p_c \BU_c]_{|c = c_*+ \gamma} |\!|_X 
 \leq \kappa $ independently of $ \gamma \in ( - \gamma_1 , + \gamma_1 ) $, hence
\begin{align*}
 \Xi_\kappa (\BW(\gamma)) = & \, 
\frac{\gamma^2 \ddot{P}_*}{ 2 |\!| \BU_* |\!|^2_H } 
\langle \BF''_* ( \BU_* )[ \mathbb{B} \BU_* ] , [ \p_c \BU_c]_{|c = c_*+ \gamma} \rangle_{X^*,X} 
+ o_{ (\gamma , \kappa ) \to (0,0)}( \gamma^2 ) \\
= & \, 
\frac{\gamma^2 \ddot{P}_*}{ 2 |\!| \BU_* |\!|^2_H } 
\langle \BF''_* ( \BU_* )[ \mathbb{B} \BU_* ] , [ \p_c \BU_c]_{|c = c_*} \rangle_{X^*,X} 
+ o_{ (\gamma , \kappa ) \to (0,0)}( \gamma^2 ) .
\end{align*}
Finally, using once again the equality (for $ \phi \in X $) 
$ \langle \BF_*'' ( \BU_*) [ \p_c \BU_c]_{|c = c_*} , \phi \rangle_{X^*,X} 
= ( \mathbb{B} \BU_* , \phi )_H $ and that $ \BF_*'' $ is self-adjoint, we infer
$$
\langle \BF''_* ( \BU_* )[ \mathbb{B} \BU_* ] , [ \p_c \BU_c]_{|c = c_*} \rangle_{X^*,X} 
= \langle \BF''_* ( \BU_* )( [ \p_c \BU_c]_{|c = c_*} ) , \mathbb{B} \BU_* \rangle_{X^*,X} 
= \n \mathbb{B} \BU_* \n^2_H = \n \BU_* \n^2_H ,
$$
and reporting this into the previous expression gives the result. \carre \\ 

We now compute the asymptotics of $ \Xi_\kappa $ for more general functions.

%%%%%%%%%%%%%%
\begin{lemmas} 
\label{Maeda10} 
Let $\e > 0 $ be small enough. If $ \BU \in \mathscr{O}_\e $ satisfies 
$ P_{\rm hy}( \BU ) = P_{\rm hy}( \BU_* ) $ and $ \BF_*(\BU) - \BF_*(\BU_*) < 0 $, 
then, we have
$$
 \Xi_{ \kappa } ( \BU ) 
 = - \frac{\bar{\gamma}^2(\BU) \ddot{P}_*}{2 } + o( \bar{\gamma}^2 (\BU) ) 
$$
uniformly for $ 0 < \kappa \leq | \bar{\gamma} (\BU) |^3 $.
\end{lemmas}
%%%%%%%%%%%%

\noindent {\it Proof.} First, we may apply Lemma \ref{Maeda6} and infer that 
$ |\!| \vartheta( \BU) |\!|_X^2 + | \alpha ( \BU ) | = \BO( | \bar{\gamma}^3(\BU)| ) $. 
Then, we write $ \Xi_\kappa ( \BU ) = \Xi_\kappa ( \BU( \cdot - \bar{y}(\BU) ) ) = 
\Xi_\kappa ( \BW( \bar{\gamma}( \BU ) ) + \vartheta ( \BU ) 
+ \alpha ( \BU ) \mathbb{B} \BU_{c_* + \bar{\gamma}(\BU) } ) = 
\Xi_\kappa ( \BW( \bar{\gamma}( \BU ) ) + \vartheta ( \BU ) ) + \BO( | \bar{\gamma}(\BU) |^3 )$ 
and, recalling the expression
$$
 \Xi_\kappa (\BU) = 
- \Big\langle \BF_{c_*+ \bar{\gamma}(\BU)} ' (\BU) , 
\Big\{ J \Upsilon_\kappa ( \bar{\gamma}(\BU)) (\cdot + \bar{y}(\BU)) 
+ ( \BU, \p_\gamma \Upsilon_\kappa (\bar{\gamma}(\BU)) (\cdot + \bar{y}(\BU) ) )_H 
\BN_{\bar{\gamma}} (\BU) \Big\} \Big\rangle_{X^*,X} ,
$$
we wish to make a Taylor expansion. First, note that
\begin{align*}
\BF_{c_*+ \bar{\gamma}(\BU)} ' (\BU)  = & \, 
\BF_{c_*+ \bar{\gamma}(\BU)}' (\BW( \bar{\gamma}( \BU ) ) )
+ \BF_{c_*+ \bar{\gamma}(\BU)}'' (\BW( \bar{\gamma}( \BU ) ) ) [ \vartheta ( \BU ) ] 
+ \BO( | \bar{\gamma}(\BU) |^3 ) \\ 
= & \, \BF_{c_*+ \bar{\gamma}(\BU)}' (\BW( \bar{\gamma}( \BU ) ) ) 
+ \BF_{c_*+ \bar{\gamma}(\BU)}'' (\BU_{c_*+ \bar{\gamma} (\BU) } ) [ \vartheta ( \BU ) ] 
+ \BO( | \bar{\gamma}(\BU) |^3 ) 
\end{align*}
hence, since $ \BF_{c_*+ \bar{\gamma}}' ( \BU_{c_*+ \bar{\gamma}} ) = 0 $ 
and (Lemma \ref{Maeda1}) $ \sigma (\gamma ) = \BO(\gamma^2)$, we have 
$ \BW ( \bar{\gamma} ) = \BU_{c_*+ \bar{\gamma}} + \BO(\gamma^2)$, thus
\begin{align*}
 \Xi_\kappa (\BU) - & \, \Xi_\kappa (\BW( \bar{\gamma}( \BU ) )) = \BO( | \bar{\gamma}(\BU) |^3 )\\ & 
- \Big\langle  \BF_{c_*+ \bar{\gamma}(\BU)}'' (\BU_{c_*+ \bar{\gamma}(\BU) } ) [ \vartheta ( \BU ) ] , 
\Big\{ J \Upsilon_\kappa ( \bar{\gamma}(\BU)) (\cdot + \bar{y}(\BU)) 
+ ( \BU, \p_\gamma \Upsilon_\kappa (\bar{\gamma}(\BU)) (\cdot + \bar{y}(\BU) ) )_H 
\BN_{\bar{\gamma}} (\BU) \Big\} 
\Big\rangle_{X^*,X} .
\end{align*}
Now, in the bracket term, we may replace $ \BU $ by $ \BW( \bar{\gamma}( \BU ) ) 
+ \BO ( |\!| \vartheta ( \BU ) |\!|_{X} ) $ (since $ |\!| \vartheta ( \BU ) |\!|_{X}^2 = 
\BO( | \bar{\gamma}(\BU) |^3 ) $). By the computations of Lemma \ref{Maeda9} 
and the equalities $ \bar{\gamma}(\BW(\gamma)) = \gamma $, $ \bar{y}(\BW(\gamma)) = 0 $, 
this gives
\begin{align*}
 \Xi_\kappa (\BU) - \Xi_\kappa (\BW( \bar{\gamma}( \BU ) )) = & \, \BO( | \bar{\gamma}(\BU) |^3 ) 
 - \Big\langle  \BF_{c_*+ \bar{\gamma}(\BU)}'' (\BU_{c_*+ \bar{\gamma}(\BU) } ) [ \vartheta ( \BU ) ] ,  
 J \Upsilon_\kappa ( \bar{\gamma}(\BU)) \Big\rangle_{X^*,X} 
\\ = & \, \BO( | \bar{\gamma}(\BU) |^3 ) 
 - \Big\langle \BF_{c_*+ \bar{\gamma}(\BU)}'' (\BU_{c_*+ \bar{\gamma}(\BU) } ) [ \vartheta ( \BU ) ] ,  
\{ \p_c \BU_c\}_{| c= c_*+ \bar{\gamma}(\BU)} 
\Big\rangle_{X^*,X} + \BO(\kappa) 
 \\ 
 = & \, \BO( | \bar{\gamma}(\BU) |^3 ) 
 - \Big\langle \BF_{c_*+ \bar{\gamma}(\BU)}'' (\BU_{c_*+ \bar{\gamma}(\BU) } ) 
[ \{ \p_c \BU_c\}_{| c= c_*+ \bar{\gamma}(\BU)} ] , \vartheta ( \BU ) \Big\rangle_{X^*,X} + \BO(\kappa) 
\end{align*}
using Lemma \ref{Maeda7} and the self-adjointness of $ \BF_{c_*+ \bar{\gamma}(\BU)}'' $. 
Choosing $ 0 < \kappa \leq | \bar{\gamma}(\BU) |^3 $ and from the equality (for $ \phi \in X $) 
$ \langle \BF_c'' ( \BU_c) [ \p_c \BU_c] , \phi \rangle_{X^*,X} = ( \mathbb{B} \BU_c , \phi )_H $, 
we infer
$$
 \Xi_\kappa (\BU) - \Xi_\kappa (\BW( \bar{\gamma}( \BU ) )) = 
 \BO( | \bar{\gamma}(\BU) |^3 ) - ( \mathbb{B} \BU_{c_*+ \bar{\gamma}(\BU) } , 
\vartheta ( \BU ) )_H =  \BO( | \bar{\gamma}(\BU) |^3 ) ,
$$
by the orthogonality condition in Lemma \ref{Maeda4}. Inserting the expansion 
of $ \Xi_\kappa (\BW( \gamma )) $ given in Lemma \ref{Maeda9} yields the 
conclusion. \carre \\

\noindent {\bf Proof of Theorem \ref{stabcusp}.} We have to show that there 
exists $ \e > 0 $ such that, for any $ \d > 0 $, we can choose an initial 
datum at distance $ \leq \d $ from $ \BU_* $ but that escape from $ \mathscr{O}_\e $. 
Since $ \BW (\gamma ) \to \BU_* $ in $ X $, we shall take the initial datum to 
be $ \BW (\gamma ) $ for some small $\gamma $, and denote 
$ \Psi_{\rm hy} (t) $ the corresponding solution. In view of Lemma \ref{Maeda2}, 
we have $ \BF_* ( \BW(\gamma)) - \BF_*(\BU_*) \sim - \gamma^3 \ddot{P}_* / 6 $, 
hence we can choose $ \gamma $ with the sign of $ \ddot{P}_* \not = 0 $ so that
$$
 \BF_* ( \BW(\gamma)) - \BF_*(\BU_*) \sim - |\gamma|^3 |\ddot{P}_*| / 6 < 0 .
$$
We now assume that $ \Psi_{\rm hy} (t) $ is globally defined and remains in 
$ \mathscr{O}_\e $, where $ \e $ is as in Lemma \ref{Maeda6}. 
By conservation of energy and momentum and the construction of $\BW(\gamma) $, 
we deduce $ P_{\rm hy} ( \Psi_{\rm hy} (t) ) = P_{\rm hy} ( \BW(\gamma) ) = P_{\rm hy} ( \BU_*) $, 
and $  \BF_* ( \Psi_{\rm hy} (t) ) - \BF_*(\BU_*) = \BF_* ( \BW(\gamma)) - \BF_*(\BU_*) < 0 $. 
The first step is to have a control on $ \bar{\gamma}(t) \equiv \bar{\gamma}(\Psi_{\rm hy} (t)) $. 
We denote $ \alpha(t) = \alpha(\Psi_{\rm hy} (t)) $, $\bar{y}(t) = \bar{y}(\Psi_{\rm hy} (t)) $ 
and $\vartheta(t) =\vartheta(\Psi_{\rm hy} (t)) $. Applying Lemma \ref{Maeda6}, we obtain 
$ |\!| \vartheta( t ) |\!|_X^2 + | \alpha ( t ) | = \BO( | \bar{\gamma}^3 (t)| ) $. In addition, 
Lemma \ref{Maeda2} and Taylor expansion gives
\begin{align*}
 \BF_* ( \Psi_{\rm hy} (t) ) - & \, \BF_*(\BU_*) =  
 \BF_* ( \BW(\bar{\gamma} (t) + \vartheta (t ) + \alpha (t ) \mathbb{B} \BU_{c_* + \bar{\gamma}(t) } ) 
 - \BF_*(\BU_*) 
 \\
 = & \, \BF_* ( \BW(\bar{\gamma} (t) ) - \BF_*(\BU_*) 
 + \langle \BF_*' ( \BW(\bar{\gamma} (t) ) , \vartheta (t ) \rangle_{X^*,X}
  + \frac12 \langle \BF_*'' ( \BW(\bar{\gamma} (t) ) \vartheta (t ) , \vartheta (t ) \rangle_{X^*,X} 
  + o ( | \bar{\gamma}^3 (t)| )  \\
 = & \, - \frac{\bar{\gamma}^3 (t) \ddot{P}_*}{6} 
 + \langle \BF_*' ( \BW(\bar{\gamma} (t) ) , \vartheta (t ) \rangle_{X^*,X}
  + \frac12 \langle \BF_*'' ( \BW(\bar{\gamma} (t) ) \vartheta (t ) , \vartheta (t ) \rangle_{X^*,X} 
  + o ( | \bar{\gamma}^3 (t)| ) ,
\end{align*}
where we have used that $ \BF'_* ( \BW(\bar{\gamma} (t) ) = o(1) $ (for the terms 
involving $ \alpha (t ) $) and Lemma \ref{Maeda2}. Furthermore, by the 
orthogonality relations in Lemma \ref{Maeda3} and using that 
$ \sigma(\gamma) = \BO ( \gamma^2 ) $ and $ \BF_c'(\BU_c) = 0 $, it holds
\begin{align*}
  \langle \BF_*' ( \BW(\bar{\gamma} (t) ) , \vartheta (t ) \rangle_{X^*,X} 
 = & \, \langle \BF_{c_* + \bar{\gamma} (t)} ' ( \BW(\bar{\gamma} (t) ) , \vartheta (t ) \rangle_{X^*,X} 
  + \bar{\gamma} (t) ( \mathbb{B} \BW(\bar{\gamma} (t) ) , \vartheta (t ) )_H  \\
= & \, \langle \BF_{c_* + \bar{\gamma} (t)} ' ( \BU_{c_*+\bar{\gamma} (t)} + 
\sigma (\bar{\gamma} (t) ) ) \mathbb{B} \BU_{c_* + \bar{\gamma} (t)} , \vartheta (t ) \rangle_{X^*,X} 
= \BO (  | \bar{\gamma}^{7/2} (t)| ) .
\end{align*}
In addition, by Lemma \ref{Maeda4}, the before last term is $ \geq K_0 |\!| \vartheta (t ) |\!|_X^2 / 2 $. 
As a consequence, by conservation of $ \BF_* ( \Psi_{\rm hy} (t) ) $, we infer, for small 
$ \gamma $,
$$
0 > - |\gamma|^3 |\ddot{P}_*| / 3 >  \BF_* ( \BW(\gamma)) - \BF_*(\BU_*) 
 =  \BF_* ( \Psi_{\rm hy} (t) ) - \BF_*(\BU_*) \geq 
 - \frac{\bar{\gamma}^3 (t) \ddot{P}_*}{6} + o ( | \bar{\gamma}^3 (t)| ) ,
 $$
In particular, this forces $ \bar{\gamma} (t) $ to always be of the sign of 
$ \ddot{P}_* $ and to satisfy $ | \bar{\gamma} (t) | \geq | \gamma | / 2 $ 
(provided $ \e $ and $ \gamma $ are small enough). 

Since now, we have a good upper bound for $ | \bar{\gamma} (t) | $, we 
can choose $ \kappa = \kappa (\gamma ) \equiv \gamma^3 / 8 $, which is 
such that, for any $ t \geq 0 $, $  \kappa \leq | \bar{\gamma} (t) |^3 $. In 
particular, we can apply Lemma \ref{Maeda10} and get
$$
 \Xi_\kappa ( \Psi_{\rm hy} (t) ) =  - \frac{\bar{\gamma}(t)^2 \ddot{P}_*}{2 } 
 + o( \bar{\gamma}(t)^2 ) .
$$
With this choice $ \kappa = \kappa (\gamma ) $, we deduce from Lemma 
\ref{Maeda8} that
$$
 \frac{d}{dt} \OO_{\kappa(\gamma)} ( \Psi_{\rm hy} (t) ) 
 = \Xi_\kappa ( \Psi_{\rm hy} (t) ) =  - \frac{\bar{\gamma}(t)^2 \ddot{P}_*}{2 } 
 + o( \bar{\gamma}(t)^2 ) .
$$
Since $ | \bar{\gamma} (t) | \geq | \gamma | / 2 $, it follows that, when 
$ \ddot{P}_* < 0 $ (the case $ \ddot{P}_* > 0 $ is analogous), 
$$
 \frac{d}{dt} \OO_{\kappa(\gamma)} ( \Psi_{\rm hy} (t) ) 
 \geq  - \frac{ \gamma^2 \ddot{P}_*}{8 } > 0 ,
$$
hence $ \OO_{\kappa(\gamma)} ( \Psi_{\rm hy} (t) ) $ is unbounded as 
$ t $ goes to $ + \ii $. However, by definition of $ \OO_{\kappa} $, 
we have by Cauchy-Schwarz 
$ | \OO_{\kappa(\gamma)} (\BU) | \leq |\!| \BU |\!|_H |\!| \Upsilon_{\kappa(\gamma)} |\!|_H 
\leq C( \gamma ) $ for $ \BU \in \mathscr{O}_\e $. We have reached a contradiction. 
The proof of Theorem \ref{stabcusp} is complete. \carre

%%%%%%%%%%%%%%%%%%%%%%%%%%%%%%%%%%%%%%%%%%%%%%%%%%%%%%%%%
\section{The linear instability ($\bs{ 0 < c_* < \cs })$}

%%%%%%%%%%%%%%%%%%%%%%%%%%%%%%%%%%%%%%%%%%%%%%
\subsection{Proof of Theorem \ref{vpinstable}}
\label{sexins}

\noindent {\bf Existence of at least one unstable eigenvalue.} The proof 
of the existence of at least one unstable eigenvalue relies on Evans function 
technique, as in \cite{Zu}, \cite{BeGa}. We shall actually use Theorem 1 in 
\cite{BeGa} when observing (see {\it e.g.} \cite{Benzo}) that the Euler-Korteweg 
system
\be
\tag{EK}
\left\{\begin{array}{ll}
\ds{ \p_t \rho + 2 \p_{x} ( \rho u ) } = 0 \\ \ \\ 
\ds{ \p_{t} u + 2 u \p_{x} u - \p_x (f ( \rho )) 
- \p_x \Big( K( \rho) \p_{x}^2 \rho + \frac12 K'(\rho) (\p_x \rho )^2 \Big) } = 0 ,
\end{array}\right. 
\ee
where $ K : (0,+\ii) \to (0,+\ii) $ is the (smooth enough) capillarity, reduces 
to \eqref{MadTWk} (where, we recall, $ \Psi = A \ex^{ i \phi }$, $ \rho = A^2 $ %\eta = A^2 - r_0^2$ 
and $ u = \p_x \phi $), namely
$$
\left\{\begin{array}{ll}
\ds{ \p_t \rho + 2 \p_{x} ( \rho u ) } = 0 \\ \ \\ 
\ds{ \p_{t} u + 2 u \p_{x} u - \p_x (f ( \rho )) 
- \p_x \Big( \frac{ \p_{x}^2 (\sqrt{\rho}) }{\sqrt{\rho} } \Big) } = 0 ,
\end{array}\right. 
$$
for the capillarity $ K (\varrho) = 1 / (2 \varrho ) $, as can be shown by straightforward 
computations. The associated eigenvalue problem in the moving frame is
\be
\label{rouge}
\left\{\begin{array}{ll}
\ds{ \lambda \zeta - c_* \p_x \zeta 
+ 2 \p_{x} ( ( r_0^2 + \eta_*) \upsilon + \zeta u_* ) } = 0 \\ \ \\ 
\ds{ \lambda \upsilon - c_* \p_x \upsilon + 2 \p_{x} ( u_* \upsilon) 
- \p_x (f' ( r_0^2 + \eta_* ) \zeta ) 
- \p_x \Big\{ \frac{1}{2 \sqrt{r_0^2 + \eta_*} } \p_{x}^2 \Big( 
\frac{ \zeta }{ \sqrt{r_0^2 + \eta_*} } \Big) 
- \frac{ \zeta \p_{x}^2 (\sqrt{r_0^2 + \eta_*}) }{2 (r_0^2 + \eta_*)^{3/2} } \Big\} } = 0 .
\end{array}\right. 
\ee 
The link with the original eigenvalue problem \eqref{pbvp} is done 
through the formula
\be
\label{saumon}
 w = U_* \Big( \frac{\zeta}{2} + i \int_{-\ii}^x \upsilon \Big) ,
\ee
since this corresponds to $ \Psi = U_{c_*} + \psi = U_{c_*} + \ex^{\lambda t} w(x) 
= ( A_{c_*} + \ex^{\lambda t} \zeta(x) ) \exp( i \phi_{c_*} + i \ex^{\lambda t} \int_{-\ii}^x \upsilon  )$. 
Notice indeed that the second equation in \eqref{rouge} gives 
$ \int_{\R} \upsilon  \ dx = 0 $. It then follows from Theorem 1 in 
\cite{BeGa} that under the assumption $ \ds{ \frac{dP(U_c)}{dc}_{|c=c_*} > 0 }$, 
there exists at least one unstable eigenvalue $ \gamma_0 \in ( 0, +\ii ) $.

\bigskip

\noindent {\bf Existence of at most one unstable eigenvalue.} The fact 
that there exists at most one unstable eigenvalue follows from arguments 
as in \cite{BDDJ} (Appendix B) and is a direct consequence of Theorem 3.1 in 
\cite{PeWe}, that we recall now.

%%%%%%%%%%%%%%%
\begin{theorem} [\cite{PeWe}]
\label{atmost}
Let $\mathcal{J}$ and $\mathcal{L}$ be two two operators on a real Hilbert 
space $X$, with $\mathcal{L}$ self-adjoint and $\mathcal{J}$ skew-symmetric. 
Then, the number of eigenvalues, counting algebraic multiplicities, 
of $ [ \mathcal{J} \mathcal{L} ]_\C $ in the right-half plane 
$ \{ {\rm Re} > 0 \} $ is less than or equal to the number of negative 
eigenvalues of $ \mathcal{L} $, counting multiplicities.
\end{theorem}
%%%%%%%%%%%%

In order to apply this result to our problem, let us write the 
eigenvalue problem \eqref{rouge} under the form
$$
\lambda \left( \begin{array}{c} 
\zeta \\ \upsilon \end{array} \right) 
= - \p_x \left( \begin{array}{cc} 
0 & 1 \\ 1 & 0 \end{array} \right) 
\BL \left( \begin{array}{c} 
\zeta \\ \upsilon \end{array} \right) , 
$$
where $ \BM $ is the self-adjoint Sturm-Liouville operator
$$
 \BM \equiv - f' ( r_0^2 + \eta_* ) 
- \frac{1}{2 \sqrt{r_0^2 + \eta_*} } \p_{x}^2 \Big( 
\frac{ \cdot }{ \sqrt{r_0^2 + \eta_*} } \Big) 
+ \frac{ \p_{x}^2 (\sqrt{r_0^2 + \eta_*}) }{2 (r_0^2 + \eta_*)^{3/2} } 
$$
(which is bounded from below) on $ \BH \equiv L^2 \times L^2 $ and with
$$ 
\BL \equiv \left( \begin{array}{cc} 
\BM & 2 u_* - c_* 
\\ 
2 u_* - c_* 
& 
2 ( r_0^2 + \eta_*) 
 \end{array} \right) .
$$
We are in the setting of Theorem \ref{atmost} with $ \BJ = - \p_x \left( \begin{array}{cc} 
0 & 1 \\ 1 & 0 \end{array} \right) $ skew-symmetric and $ \BL $ self-adjoint. 
We thus show that $ \BL $ has at most one negative eigenvalue. Since $ r_0^2 + \eta_* $ 
remains bounded away from zero, it is clear that, for $ \s < 0 $ and 
$ ( \zeta , \upsilon ) $ given, $ \BL ( \zeta , \upsilon )^t = \s ( \zeta , \upsilon )^t $ 
if and only if
\be
\label{dupropre}
\BM^\dag \zeta - \frac{ ( c_* - 2 u_*)^2}{2(r_0^2 + \eta_*) } \cdot 
\frac{ \s }{2(r_0^2 + \eta_*) - \s } \zeta = \s \zeta , 
\quad \quad \quad {\rm with} \quad \quad \quad 
\BM^\dag \equiv \BM - \frac{ ( c_* - 2 u_*)^2 }{2(r_0^2 + \eta_*) } ,
\ee
since we may express $ \upsilon $ in terms of $\zeta$ with the second equation. 
We observe that the translation invariance shows that $ \p_x ( \eta_* , u_* )^t$ 
belongs to the kernel of $ \mathcal{L} $, that is, using once again the relation 
$ 2 u_c = 2 \p_x \phi_c = c \frac{\eta_c}{\eta_c + r_0^2} $, 
$ \BM^\dag \p_x \eta_* = 0 $. Furthermore, $ \BM^\dag $ has the same continuous 
spectrum as its constant coefficient limit as $x \to \pm \ii$, namely
$$
- \frac{1}{2r_0^2}\, \p_x^2 + \frac{\cs^2 - c_*^2}{2 r_0^2} ,
$$
that is $ \s_{\rm ess} (\BM^\dag) = [\cs^2 - c_*^2 , +\ii ) \subset ( 0 , +\ii )$, 
since $0 < c_* < \cs $. Since $ \p_x \eta_* $ has exactly one zero (at $x=0$), 
it follows from standard Sturm-Liouville theory that $ \BM^\dag $ has precisely 
one negative eigenvalue $ \mu < 0 $ and that the second eigenvalue is $0$. 
Taking the scalar product with \eqref{dupropre} yields
$$
 \langle \BM^\dag \zeta , \zeta \rangle_{L^2} 
- \int_\R \frac{ \s ( c_* - 2 u_*)^2 \zeta^2 }{2(r_0^2 + \eta_*) [ 2(r_0^2 + \eta_*) - \s ]} \ dx 
= \s |\!| \zeta |\!|_{L^2}^2 .
$$
%The integral is $ < 0 $ since $ \s < 0 $, hence we must have $ \mu < \s < 0 $. 
Now, for $ s \leq 0 $, we consider the self-adjoint operator
$$
\BM^\dag_s \equiv \BM^\dag - \frac{ ( c_* - 2 u_*)^2}{2(r_0^2 + \eta_*) } \cdot 
\frac{ s }{2(r_0^2 + \eta_*) - s } .
$$
Clearly, $ \BM^\dag_{s=0} = \BM^\dag $, $ \s_{\rm ess}( \BM^\dag_s ) \subset 
[\cs^2 - c_*^2 , +\ii ) \subset ( 0 , +\ii )$, and $ \R_- \ni s \mapsto \BM^\dag_s $ 
is decreasing. Let us assume now that the self-adjoint operator $ \BL $ has 
at least two negative eigenvalues. Then, we denote $ \s_1 < \s_2 < 0 $ the 
two smallest eigenvalues of $ \BL $ (necessarily simple), and $ \zeta_1 $, $ \zeta_2 $ 
two associated eigenvectors. Since $\BL$ is self-adjoint, 
$ \langle \zeta_1 , \zeta_2 \rangle_{L^2} = 0 $. Furthermore, 
$ \langle \BM^\dag_{s= \s_1} \zeta_1 , \zeta_1 \rangle_{L^2} 
= \s_1 |\!| \zeta_2 |\!|_{L^2}^2 < 0 $, hence, by monotonicity, 
$ \langle \BM^\dag_{s} \zeta_2 , \zeta_2 \rangle_{L^2} < 0 $ for any 
$ \s_1 \leq s \leq 0 $. Therefore, $ \BM^\dag_{s} $ has at least one negative 
eigenvalue for $ \s_1 \leq s \leq 0 $. We denote $ \lambda_{\rm min} (s) $ the smallest 
eigenvalue of $ \BM^\dag_s $. Then, $ \lambda_{\rm min} (s=0) = \mu < 0 $ and $ \lambda_{\rm min} $ 
decreases in $ [ \s_1 , 0 ] $. Moreover, we may choose a positive eigenvector 
$ \zeta_s $ for the eigenvalue $ \lambda_1 (s) $, with $\zeta_1 = \zeta_{\s_1} $. 
Since $ \s( \BM^\dag ) \cap \R_+ = \{ \mu, 0 \} $, it follows from the monotonicity 
that for any $ \s_1 \leq s < 0 $, we have $ \s( \BM^\dag_s ) \cap \R_- = \{ \lambda_{\rm min} (s) \} $. 
When $ s = \s_2 \in ( \s_1 ,0 ) $, we then have $ \s_2 \in \s( \BM^\dag_{s=\s_2} ) \cap \R_- $, 
and thus $ \s_2 = \lambda_{\rm min} ( \s_2 ) $, which implies that we may choose $ \zeta_2 > 0 $ 
without loss of generality. Similarly, if $ s = \s_2 $, we see that we may choose 
$ \zeta_2 > 0 $. We obtain a contradiction since then $ \langle \zeta_1 , \zeta_2 \rangle_{L^2} > 0 $ 
and thus $ \zeta_1 $ and $ \zeta_2 $ cannot be orthogonal in $L^2$. We have thus 
shown that $ \BL $ has at most one negative eigenvalue, and then Theorem \ref{atmost} 
shows that $ \mathcal{J} \mathcal{L} $ has at most one eigenvalue in $ \{ {\rm Re} > 0 \} $, 
as wished.

%%%%%%%%%%%%%%%%%%%%%%%%%%%%%%%%%%%%%%%%%%%%%%%%%%%%%%%%%%%%%%%%%%%%%%%%%%%%%%%%%%%%%%%%%%%
\subsection{Resolvent and semigroup estimates (proof of Corollary \ref{instabamelioree})}

In this section, we drop the ``$*$'' for the travelling wave 
we are considering. When linearizing the (NLS) equation in the moving 
frame with speed $c$, we obtain
\be
\label{eqlin}
i \frac{\p \psi}{\p t} - i c \p_x \psi + \p^2_x \psi + \psi f(|U|^2 ) 
+ 2 \langle \psi , U \rangle f'(|U|^2 ) U = 0 ,
\ee
or
\begin{align*}
\frac{\p }{\p t} \left( \begin{array}{c} \psi_1 \\ \psi_2 \end{array} \right) 
= & \, \left( \begin{array}{cc} 
 c \p_x - 2 f'(|U|^2 ) U_{1} U_{2} 
& 
- \p^2_x - f(|U|^2 ) - 2 f'(|U|^2 ) U_{2}^2  \\ 
\p^2_x + f(|U|^2 ) + 2 f'(|U|^2 ) U_{1}^2 
& 
 c \p_x + 2 f'(|U|^2 ) U_{1} U_{2}
\end{array} \right) 
\left( \begin{array}{c} \psi_1 \\ \psi_2 \end{array} \right) \\ 
= & \, \left( \begin{array}{cc} 0 & 1 \\ -1 & 0 \end{array} \right) 
\left( \begin{array}{cc} 
- \p^2_x - f(|U|^2 ) - 2 f'(|U|^2 ) U_{1}^2 
& 
- c \p_x - 2 f'(|U|^2 ) U_{1} U_{2} \\ 
 c \p_x - 2 f'(|U|^2 ) U_{1} U_{2} 
& 
- \p^2_x - f(|U|^2 ) - 2 f'(|U|^2 ) U_{2}^2 
\end{array} \right) 
\left( \begin{array}{c} \psi_1 \\ \psi_2 \end{array} \right) .
\end{align*}
We wish to show that this linear equation can be solved using a continuous 
semigroup. In order to handle later the nonlinear terms, we work in $H^1(\R,\C^2)$ 
instead of $ L^2(\R,\C^2) $. Therefore, we consider the unbounded operator 
$ \BA : D(\BA) = H^3(\R, \C^2 ) \subset H^1(\R, \C^2 ) \to H^1(\R, \C^2 ) $ 
on $ H^1(\R, \C^2 ) $ defined by
$$
\BA \equiv \left( \begin{array}{cc} 
 c \p_x - 2 f'(|U|^2 ) U_{1} U_{2} 
& 
- \p^2_x - f(|U|^2 ) - 2 f'(|U|^2 ) U_{2}^2  \\ 
\p^2_x + f(|U|^2 ) + 2 f'(|U|^2 ) U_{1}^2 
& 
 c \p_x + 2 f'(|U|^2 ) U_{1} U_{2}
\end{array} \right) .
$$
It follows easily that for $ \psi = \left( \begin{array}{c} \psi_1 \\ \psi_2 \end{array} \right) 
\in H^1(\R, \C^2 ) $,
\begin{align*}
 {\rm Re} \Big( \langle \BA \psi | \psi \rangle_{H^1(\R, \C^2 )} \Big) = & \, 
{\rm Re} \Big( 
\langle -2 f'(|U|^2 ) U_{1} U_{2} \psi_1 , \psi_1 \rangle_{H^1(\R, \C)} 
+ \langle -2 f'(|U|^2 ) U_{1} U_{2} \psi_2 , \psi_2 \rangle_{H^1(\R, \C)} 
\\ & \, 
+ \langle [ f(|U|^2 ) + 2 f'(|U|^2 ) U_{1}^2 ] \psi_1 , \psi_2 \rangle_{H^1(\R, \C)} 
- \langle [ f(|U|^2 ) + 2 f'(|U|^2 ) U_{1}^2 ] \psi_2 , \psi_1 \rangle_{H^1(\R, \C)} 
\Big) 
\\ \leq & \, 
K |\!| \psi |\!|^2_{H^1(\R, \C^2 )} .
\end{align*}
Moreover, the spectrum of $ \BA $ is included in the half-space $\{ {\rm Re} \leq \s_0 \}$, 
hence $ \BA $ generates a continuous semigroup $ \ex^{ t \BA } $ on $ H^1 (\R, \C^2 ) $.

In order to estimate the growth of the semigroup $ \ex^{ t \BA } $ on 
$ H^1 (\R, \C^2 ) $, we could try to use the same approach as \cite{dMG}, 
which relies on the proof of the spectral mapping theorem in \cite{GJLS}. 
However, our situation is slightly different since in these studies, the 
reference solutions is real-valued (it is a bound state in \cite{GJLS} 
and the kink in \cite{dMG}). Therefore, $ U_2 = 0 $ and $ \BA $ has no 
diagonal term, and the system is much more decoupled than in our situation. 
As a matter of fact, it is not very clear whether the arguments of \cite{GJLS} 
carry over to our problem. We thus have chosen to use the general approach 
given in Appendix B. We thus verify the assumptions of Theorem B.4 
(see also Corollary B.2) there, which are easy: $\BA$ generates a semigroup in $H^1(\R,\C^2)$ 
and the spectrum of $\BA$ is of the form $ i\R \cup \{ - \gamma_0 , + \gamma_0 \} $, 
where $ i\R $ is the essential spectrum and $ \pm \gamma_0 $ two simple eigenvalues. 
Moreover, the eigenvector associated with $ \gamma_0 $ belongs to $ H^3(\R,\C^2) = D(J) $. 
Therefore, Theorem B.4 in Appendix B applies and the growth estimate for the linearized 
problem follows. For the nonlinear instability result, we argue as for Corollary B.2 
in Appendix B, since the manifold $ \mathfrak{M} = \{ |U_*|( \cdot - y ), \ y \in \R \}$ 
is transverse to the curve $ \s \mapsto | U_* + \s w | $ in $ r_0 + H^1 (\R) $. 
Indeed, it follows from \eqref{saumon} that 
$ | U_* + \s w | = A_* + \s \zeta + \BO_{H^1} (\s^2) $. Assume that 
$ \zeta = \alpha \p_x |U_*| $, with $ \alpha \in \R $. Then, integration 
of the first equation of \eqref{rouge} provides
$$
 \lambda ( |U_*| - r_0 ) - c_* \p_x |U_*| + 2 \Big( (r_0^2 + \eta_*) \upsilon + u_* \zeta  \Big) = 0 ,
$$
hence, using that $ |U_*| = \sqrt{r_0^2 + \eta_*} $ and the equality 
$ 2 u_* = c \eta_* / (r_0^2 + \eta_*) $, we infer
$$
 \upsilon + \alpha \Big\{ \lambda \frac{ |U_*| - r_0 }{r_0^2 + \eta_* } 
+ \frac{c _* r_0^2}{ 4(r_0^2 + \eta_*)^{3/2}} \p_x \eta_* \Big\} = 0 .
$$
Since $ \int_\R \upsilon = 0 $ and $ |U_*| - r_0 $ has constant sign in $ \R $, 
integrating over $ \R $ then implies $ \alpha = 0 $, which in turn yields 
$ \zeta = \upsilon = 0 $ and $ w_*= 0 $, a contradiction. Consequently, 
$ \zeta \not \in \R \p_x |U_*| $ and the manifold 
$ \mathfrak{M} = \{ |U_*|( \cdot - y ), \ y \in \R \}$ is indeed transverse 
to the curve $ \s \mapsto | U_* + \s w | $ in $ r_0 + H^1 (\R) $.

%%%%%%%%%%%%%%%%%%%%%%%%%%%%%%%%%%%%%%%%%%%%%%%%%%%%%%%%%%%%%%%%%%%%
\section{Stability analysis for the kink $\bs{( c = 0)} $} 
\label{pfkink}

%%%%%%%%%%%%%%%%%%%%%%%%%%%%%%%%%%%%%%%%%%%%%%
\subsection{Proof of Lemma \ref{prolongement}}

Let us recall that the momentum $ P(U_c) $, for $ c > 0 $, has 
the expression
$$
 P(U_c) = c \int_{\xi_c}^0 \frac{\xi^2}{ r_0^2 + \xi} 
\frac{d\xi}{ \sqrt{ - \BV_c(\xi) } } ,
$$
since $  {\rm sgn} (\xi_c) = - 1 $. Therefore, we decompose 
$ P (U_c )$ with two integrals:
\be
\label{split}
 P(U_c) = c \int_{\xi_c}^0 \frac{ \xi^2 }{ r_0^2 + \xi} 
\frac{d\xi}{ \sqrt{ - \BV'_c(\xi_c)(\xi - \xi_c) } } 
+ c \int_{\xi_c}^0 \frac{ \xi^2 }{ r_0^2 + \xi} 
\Big( \frac{1}{ \sqrt{ - \BV_c(\xi) }} - \frac{1}{\sqrt{ -\BV'_c(\xi_c)(\xi -\xi_c)} } 
\Big) \ d \xi .
\ee
Using the change of variables $ \xi = t \xi_c $, the second integral in \eqref{split} 
is equal to
\begin{align*}
\xi_c^3 \int_{1}^0 \frac{ t^2 }{ r_0^2 + t \xi_c } &\, 
\Big( \frac{1}{ \sqrt{ - \BV_c(t \xi_c) }} - \frac{1}{\sqrt{ -\xi_c\BV'_c(\xi_c)(t - 1)} } 
\Big) \ d t \\ 
& = (-r_0^2)^3 \int_{1}^0 \frac{ t^2 }{ r_0^2 - t r_0^2 } 
\Big( \frac{1}{ \sqrt{ - \BV_0( -t r_0^2) }} - \frac{1}{\sqrt{ - 4 r_0^2 F(0) (t - 1)} } 
\Big) \ d t + o_{c \to 0} (1) \\ 
& = \int_{-r_0^2}^0 \frac{ \xi^2 }{ r_0^2 + \xi} 
\Big( \frac{1}{ \sqrt{ - \BV_0(\xi) }} - \frac{1}{\sqrt{ 4F(0)(\xi +r_0^2)} } 
\Big) \ d \xi + o_{c \to 0} (1) .
\end{align*}
The passage to the limit $ c \to 0 $ being justified by the dominated 
convergence theorem since the absolute value of the integrand is 
$ \leq K t $ for $ 0 \leq t \leq 1/2 $ for small $c$ and for 
$ 1/2 \leq t \leq 1 $, since $ \xi_c > -r_0^2 $, 
$ r_0^2 + t \xi_c \geq r_0^2 ( 1-t) $ and hence is equal to
\begin{align*}
& = \Big| \frac{ t^2 }{ r_0^2 + t \xi_c } \cdot 
\frac{\BV_c(t \xi_c) - \xi_c \BV'_c(\xi_c)(t - 1)}{ 
\sqrt{ - \BV_c(t \xi_c) } \sqrt{ -\xi_c\BV'_c(\xi_c)(t - 1)} 
\Big[ \sqrt{ - \BV_c(t \xi_c) } + \sqrt{ -\xi_c\BV'_c(\xi_c)(t - 1)} \Big] } \Big| \\ 
& \leq K \frac{( 1-t)^2 }{ ( 1-t )\sqrt{1-t} \sqrt{1-t} \sqrt{1-t} } 
= \frac{K}{\sqrt{1-t} } \in L^1( (1/2,1) ) .
\end{align*}
Furthermore, letting $\xi = \xi_c + (r_0^2 + \xi_c) t^2 $, $t \geq 0$, 
the first integral in \eqref{split} is equal to
\begin{align*}
\frac{1}{\sqrt{ r_0^2 + \xi_c}} & \, \int_0^{ \sqrt{- \xi_c /(r_0^2 + \xi_c)}} 
\frac{ (\xi_c + (r_0^2 + \xi_c) t^2)^2 }{1 + t^2} 
\cdot \frac{2 d t}{ \sqrt{ - \BV'_c(\xi_c) } } \\
& = \frac{2}{ \sqrt{ r_0^2 + \xi_c} \sqrt{ - \BV'_c(\xi_c) } } 
\Big\{ r_0^4 \Big[ \frac{\pi}{2} - \arctan \Big( \sqrt{ \frac{r_0^2 +\xi_c}{-\xi_c} } \Big) 
\Big] - 2 r_0^2 (r_0^2 + \xi_c ) \sqrt{ \frac{- \xi_c }{ r_0^2 + \xi_c} } 
\\ & \quad 
+ (r_0^2 + \xi_c )^2 \Big[ \sqrt{ \frac{- \xi_c }{ r_0^2 + \xi_c} }
+ \frac13 \sqrt{ \frac{- \xi_c }{ r_0^2 + \xi_c} }^3 \Big] \Big\} .
\end{align*}
by direct computation. Since $ \xi_c \simeq -r_0^2 $ is a simple 
zero of $ \BV_c (\xi) = c^2 \xi^2 - 4 ( r_0^2 + \xi ) F( r_0^2 + \xi ) $, 
we have
$$
 \xi_c = -r_0^2 + \frac{c^2r_0^4}{4F(0)} 
+ \frac{c^4r_0^6}{4 F(0) } \Big( \frac{r_0^2 f(0) }{ F(0) } - 2 \Big) 
+ o_{c\to 0} (c^4) 
= -r_0^2 + \frac{c^2r_0^4}{4F(0)} + \BO_{c\to 0}(c^4) ,
$$
thus
$$
 - \BV'_c(\xi_c) = 4 F(0) + \BO_{c\to 0} (c^2) 
$$
and
$$
\frac{2}{ \sqrt{ r_0^2 + \xi_c} \sqrt{ - \BV'_c(\xi_c) } } 
= \frac{2}{ r_0^2 c} + \BO_{c\to 0}(c) .
$$
As a consequence, the first integral  in \eqref{split} is equal to
\begin{align*}
 \frac{r_0^2 \pi}{c} + \Big\{ 
- \frac{ r_0^3}{ \sqrt{F(0)} } 
- \frac{2 r_0^3}{ \sqrt{F(0)} } 
+ \frac{ r_0^3}{ 3 \sqrt{F(0)} } 
\Big\} + \BO_{c\to 0} (c) 
=  \frac{r_0^2 \pi}{c} - \frac{8 r_0^3}{ 3 \sqrt{F(0)} } + \BO_{c\to 0} (c) .
\end{align*}
Gathering these two relations, we obtain
\begin{align*}
 P(U_c) = r_0^2 \pi + c \Big\{ - \frac{8 r_0^3}{ 3 \sqrt{F(0)} } 
+ \int_{-r_0^2}^0 \frac{ \xi^2 }{ r_0^2 + \xi} 
\Big( \frac{1}{ \sqrt{ - \BV_0(\xi) }} - \frac{1}{\sqrt{ 4F(0)(\xi +r_0^2)} } 
\Big) \ d \xi \Big\} + o_{c\to 0} (c ),
\end{align*}
as wished.

%%%%%%%%%%%%%%%%%%%%%%%%%%%%%%%%%%%%%%%%%%%%%%
\subsection{Proof of Theorem \ref{statiomini}}
\label{sexmini}

Since we have a kink solution $U_0$ for $c=0$, 
this implies that $ \BV_0 (\xi) = -4(r_0^2 + \xi) F (r_0^2 + \xi) $ is negative 
in $(-r_0^2 ,0)$ and that $ -r_0^2 $ is a simple zero of $ \BV_0 $, 
that is $ F (0) > 0 $. Then, $ F > 0 $ in $ [ 0 , r_0^2 ) $ and 
$ F( \varrho ) \simeq \frac{\cs^2}{4 r_0^2} ( \varrho - r_0^2 )^2 $ 
for $  \varrho \to r_0^2 $, it follows that there exists $ K_0 > 0 $ 
such that
$$
 F(\varrho) \geq \frac{1}{K_0} ( \varrho - r_0^2 )^2 
$$

We consider for $ \mu \geq 0 $ the quantity
$$
 \mathscr{K}_{\rm min}(\mu) \equiv \inf \Big\{ \mathscr{K} (u),\ 
u \in \mathcal{Z} , \ \inf_\R |u| = \mu \Big\} .
$$
The study of $ \mathscr{K}_{\rm min}(0) $ is easy.

%%%%%%%%%%%%%%%%%%%%
\begin{propositions} 
\label{minokink}
There holds
$$
\mathscr{K}_{\rm min}(0) = E(U_0) .
$$ 
More precisely, for any $ U \in \mathcal{Z} $,
$$ 
E(U) \geq 4 \int_{\inf_\R |U|}^{r_0} \sqrt{ F(s^2)} \ ds 
\quad \quad \quad {\it and} \quad \quad \quad 
E(U_0) = 4 \int_{0}^{r_0} \sqrt{ F(s^2)} \ ds .
$$
Finally, if $ U \in \mathcal{Z} $, $ \inf_\R |U| = 0 $ and 
$ \mathscr{K} (U ) = E(U_0) $, then there exists $ y \in \R $ 
and $ \theta \in \R $ such that $ U = \ex^{ i \theta } U_0( \cdot -y ) $.
\end{propositions}
%%%%%%%%%%%%%%%%%%%%

\noindent {\it Proof.} Taking $ U_0$ as a comparison map, we see that 
$ \mathscr{K}_{\rm min}(0) \leq E(U_0) $. Moreover, if $ U \in \mathcal{Z} $ 
and $ \inf_\R |U| = \mu \geq 0 $, we may assume, up to a translation, that 
$\mu = |U|(0) $. Then, denoting
$$
G ( r ) \equiv 2 \int_{r_0}^r \sqrt{ F(s^2)} \ ds ,
$$
we have the inequalities
\begin{align*}
\int_0^{+\ii} |\p_x U|^2 + F(|U|^2) \ dx 
\geq & \, 
\int_0^{+\ii} |\p_x |U| \, |^2 + F(|U|^2) \ dx 
\geq 2 \int_0^{+\ii} | \sqrt{ F(|U|^2)} \p_x |U| \, | \ dx 
\\ = & \, \int_0^{+\ii} | \p_x [G(|U|) ] | \ dx 
\geq \Big| \int_0^{+\ii} \p_x [ G(|U|) ] \ dx \Big| 
= \Big| G(|U|(+\ii)) - G(|U|(0)) \Big| 
\\ = & \,
\Big| G(r_0) - G(\mu) \Big| = 2 \int_\mu^{r_0} \sqrt{ F(s^2)} \ ds .
\end{align*}
Arguing similarly in $ ( -\ii, 0)$, we get
$$ 
E(U) \geq 4 \int_\mu^{r_0} \sqrt{ F(s^2)} \ ds .
$$
For the kink $U_0$, which is real-valued, we have the first integral 
$ |\p_x U_0|^2 = F(U_0^2) $, hence, using the change of variables $s = U_0(x)$,
$$
E(U_0) = 4 \int_0^{+\ii} F(U_0^2) \ dx = 4 \int_0^{r_0} \sqrt{ F(s^2)} \ ds .
$$
If $\mu = 0 $, we have then $ E(U) \geq E(U_0)$, hence 
$ \mathscr{K} (U) \geq E(U) \geq E(U_0)$ as wished. 

Assume finally that $ U \in \mathcal{Z} $ verifies $ \inf_\R |U| = 0 $ 
and $ \mathscr{K} ( U ) = E(U_0) $. Then $\mu=0$ and all the above 
inequalities are equalities. In particular, we must have 
$ |\p_x U| = |\p_x |U| \, | $ and equality in 
$ |\p_x |U| \, |^2 + F(|U|^2) \geq 2 | \sqrt{ F(|U|^2)} \p_x |U| \, | $, 
which means that $ |\p_x |U| | = \sqrt{ F(|U|^2)} $. Combining this 
ode with the condition $|U|(0) = 0$, we see that $ |U| = |U_0| $, since 
$ |U_0| $ solves $ \p_x U_0 = \sqrt{ F(U_0^2)} $. Finally, the fact 
that $ |\p_x U| = |\p_x |U| | $ implies that the phase is constant in 
$ (-\ii, 0) $ and in $(0,+\ii)$: there exist two constants 
$\theta_\pm \in \R $ satisfying $ U (x) = \ex^{i \theta_\pm } |U_0|(x) $ 
for $ \pm x \geq 0 $. Therefore, $ \mathfrak{P}(U) = r_0^2 ( \theta_+ - \theta_- ) $ 
mod $ 2 \pi r_0^2 $, and then
$$ 
 E(U_0) = \mathscr{K} (u ) = E(U_0) + 2 M r_0^4 \sin^2 \Big( \frac{ \theta_+ - \theta_- - \pi}{2} \Big) 
$$
implies $ \theta_+ - \theta_- = \pi $ mod $2\pi$, that is $ U = \ex^{i \theta_+ } U_0 $ 
in $\R$, which is the desired result. \carre \\

We recall the expansion $ P(U_s) = r_0^2 \pi + s \dot{P}_0 + o(s) $ 
as $ s \to 0 $, where $ \dot{P}_0 \equiv \ds{ \frac{dP(U_s)}{ds}_{|s=0} } $. 
From the Hamilton group relation $ \ds{ \frac{dE(U_s)}{ds} = s \frac{dP(U_s)}{ds} } $, 
we also infer by integration $ E(U_s) = E(U_0) + \ds{ \frac{s^2}{2} } \dot{P}_0 + o(s^2) $. 
As a first step, we define the small parameter $ \mu_* > 0 $. 
The key point is to prove the following result.

%%%%%%%%%%%%%%%%%%%%
\begin{propositions} 
\label{supercool} 
There exist some constant $K > 0 $ and a small $ \mu_* > 0 $ such that, 
for any $0 < \mu \leq \mu_* $,
$$
 \mathscr{K}_{\rm min}(\mu) = \inf \Big\{ \mathscr{K} (U),\ 
U \in \mathcal{Z} , \ \inf_\R |u| = \mu \Big\} \geq E(U_0) + \frac{\mu^2}{K} .
$$ 
\end{propositions}
%%%%%%%%%%%%%%%%%%%%

\noindent {\it Proof.} Notice first that for $ c > 0 $ small, 
there exists $ U_c $ travelling wave of speed $c$ and that $ \inf_\R |U_c| 
= \sqrt{r_0^2 + \xi_c } $ with $ \xi_c $ a smooth function in $c$ 
such that $ \xi_c = - r_0^2 + \ds{\frac{c^2 r_0^4}{4F(0)}} + \BO(c^4) $, 
hence $ \inf_\R |U_c| = \ds{ \frac{c r_0^2}{2 \sqrt{F(0)}} } + \BO(c^2)$ and is 
smooth. Therefore, there exists, for $0 \leq \mu \leq \mu_* $ small, 
a unique $\s_\mu$, with $\s_\mu= \ds{ \frac{2 \mu \sqrt{F(0)}}{r_0^2} + \BO(\mu^2)} $, 
such that $ \mu = \inf_\R |U_{\s_\mu} | $. In particular, taking 
$ U_{\s_\mu} $ as a comparison map in $ \mathscr{K}_{\rm min}(\mu) $, 
we have
\begin{align*}
 \mathscr{K}_{\rm min}(\mu) \leq & \, \mathscr{K}(U_{\s_\mu} ) 
= E ( U_{\s_\mu} ) + 2M r_0^4 \sin^2\Big( \frac{P(U_{\s_\mu}) - r_0^2 \pi}{2 r_0^2} \Big) \\
= & \, E(U_0) + \frac{\s^2_\mu }{2} \dot{P}_0 + o(\s^2_\mu) 
+ 2M r_0^4 \sin^2\Big( \frac{\s_\mu \dot{P}_0 + o(\s_\mu) }{2 r_0^2} \Big) \\
= & \, E(U_0) + \frac{\s^2_\mu }{2} \Big( \dot{P}_0 + M \dot{P}^2_0 \Big) + o(\s^2_\mu) .
\end{align*}
In particular, it follows that, for some positive constant $K$ and 
for $\mu_*$ small enough,
\be
\label{paraudessus}
 \mathscr{K}_{\rm min}(\mu) \leq E(U_0) + K \mu^2 
\leq \frac{11}{10} E(U_0) . 
\ee
Consider now $c$ small, a bounded open interval $ (x_- ,x_+) $ and 
$\eta$ a solution to the Newton equation
$$
2 \p^2_x \eta + \BV_c'(\eta) = 0 
$$
in $ (x_- ,x_+) $, with $ \p_x \eta(x_+) \leq 0 \leq \p_x \eta(x_-) $, 
$ \eta(x_+) \leq -r_0^2 + \mu_*^2 $ and $ \eta(x_-) \leq -r_0^2 + \mu_*^2 $. 
As $c \to 0 $, $ \BV_c $ converges to $ \BV_0 $ in $ \BC^1([ -r_0^2 ,0 ] ) $. 
Moreover, $\BV_0 $ is negative in $ ( -r_0^2 ,0 ) $ and has a simple zero at $-r_0^2 $. 
Therefore, if $c$ and $\mu_* > 0 $ are sufficiently small, we must have 
$ \int_{x_-}^{x_+ } F(r_0^2 + \eta) \ d x \geq \frac12 \int_\R F(U_0^2 ) \ d x $. 
Consequently, if $ v = A \ex^{i \vp} $ solves (TW$_c$) on a bounded interval $ (x_- ,x_+) $, 
satisfies $ 2 \p_x \vp = \ds{ \frac{ c \eta}{ r_0^2 + \eta} } $ ($\eta \equiv A^2 - r_0^2$) 
and if $ |v| $ is $ \leq \mu_* $ at $x_+$ and at $x_-$, 
with $ \p_x |v| (x_+) \leq 0 \leq \p_x |v| (x_-) $, then
\be
\label{chero}
\int_{x_-}^{x_+} |\p_x v|^2 + F(|v|^2) \ dx \geq \frac12 E(U_0) .
\ee
Here, we use that the Newton equation on the modulus $|V|$ actually holds true 
in $ (x_- ,x_+) $. Since $ F > 0 $ in $ [ 0 , r_0^2 ) $ and 
$ F( \varrho ) \simeq r_0^2 ( \varrho - r_0^2 )^2 $ 
when $ \varrho \to r_0^2 $, there exists $ K > 0 $ and $ \kappa > 0 $ such that 
$ F( \varrho ) \geq ( \varrho - r_0^2 )^2 /K $ for 
$ 0 \leq \varrho \leq r_0^2 (1 + \kappa)^2 $. Hence, if $ \inf_\R |v| \geq \mu > 0 $, 
then
\be
\label{garcia}
| P(v) | \leq \frac{K}{\mu} E (v) .
\ee
Moreover, arguing as in the proof of Proposition \ref{minokink}, we show 
that there exists $ \varkappa > 0 $ such that if $ U \in \mathcal{Z} $ and 
$ |U| $ takes values $ \leq \mu_* $ and $ \geq r_0 (1 + \kappa) $, then
$$
 E(U) \geq E(U_0) ( 1 + \varkappa ) .
$$
In particular, since $ \mathscr{K}_{\rm min}(\mu) \leq E(U_0) + \BO (\mu^2) $, 
we may choose $\mu_*$ sufficiently small so that if $ U \in \mathcal{Z} $ and 
$ \mathscr{K} (U) \leq \mathscr{K}_{\rm min} (\mu) + \mu_* $, then 
$ |U| \leq r_0 (1 + \kappa) $. This means that for the mappings we are 
considering, $ F( \varrho ) \geq ( \varrho - r_0^2 )^2 /K $.\\

\noindent {\bf Step 1: Construction of a suitable minimizing sequence.} 
There exists  a sequence $ ( V_n )_{n \geq 0} $ in $ \mathcal{Z} $ 
such that $\inf_\R |V_n | = \mu = |V_n|(0)$, 
$ V_n = A_n \ex^{ i \phi_n } $, $ P(V_n ) \in [0, \pi r_0^2 ] $,
$$ 
2 A_n^2 \p_x \phi_n = c_n ( A_n^2 - r_0^2 ) ,
\quad \quad \quad c_n \equiv M r_0^2 \sin \Big( \frac{r_0^2 \pi - P(V_n )}{2r_0^2} \Big) \geq 0
$$
and
$$
\lim_{n \to +\ii}  \mathscr{K} (V_n) = \mathscr{K}_{\rm min}(\mu) .
$$

Since $ \mu > 0 $, the maps $V$ we consider may be 
lifted $ V = A \ex^{i \phi } $. Therefore (with $u = \p_x \phi$), 
\begin{align}
\mathscr{K}_{\rm min}(\mu) = & \, \inf \Big\{ 
\int_\R ( \p_x A)^2 + F(A^2) \ dx  \nonumber \\
& \quad \quad \quad + \inf \Big\{ \int_\R A^2 u^2 \ dx 
+ 2M r_0^4 \sin^2 \Big( \frac{\int_\R (A^2 - r_0^2) u \ dx - r_0^2 \pi}{2r_0^2} \Big) , 
u \in L^2 ( \R , \R ) \Big\} , \label{bernardo} \\ 
& \quad \quad \quad A \in r_0 + H^1(\R, \R ) , 
\inf_\R A = \mu \Big\} .
 \nonumber
\end{align}
The infimum in $u$ may be written
$$
\inf_{ p \in \R} \inf \Big\{ \int_\R A^2 u^2 \ dx 
+ 2M r_0^4 \sin^2 \Big( \frac{p - r_0^2 \pi}{2r_0^2} \Big) , u \in L^2 ( \R , \R ) 
\ s.t. \ \int_\R (A^2 - r_0^2) u \ dx = p \Big\} .
$$
For each $ p \in \R $, we minimize in $u$ a quadratic 
functional on an affine hyperplane, with minimizer given 
by
$$
 u_p = p \Big( \int_\R \frac{(A^2-r_0^2)^2}{A^2} \ dx \Big)^{-1} 
 \frac{A^2-r_0^2}{A^2} .
$$
As a consequence, the infimum in $u$ in \eqref{bernardo} is
$$
\inf_{ p \in \R} \Big[ \int_\R A^2 u_p^2 \ dx 
+ 2M r_0^4v\sin^2 \Big( \frac{p - r_0^2 \pi}{2 r_0^2 } \Big) \Big] 
= 
\inf_{ p \in \R} \Big[ p^2 \Big( \int_\R \frac{(A^2-r_0^2)^2}{A^2} \ dx \Big)^{-1} 
+ 2M r_0^4 \sin^2 \Big( \frac{p - \pi r_0^2 }{2 r_0^2} \Big) \Big] .
$$
It is clear that this last infimum is achieved only for $ p $ inside 
$[ - \pi r_0^2 , + \pi r_0^2 ] $. Indeed, the second term is $2 \pi r_0^2 $-periodic and 
if $p > \pi r_0^2 $, then $ p - 2 \pi r_0^2 $ is a better competitor. Moreover, the 
function $ p \mapsto \sin^2 \Big( \frac{p - \pi r_0^2 }{2 r_0^2 } \Big) $ is 
continuous and even, hence we may consider some $ p \in [ 0 , \pi r_0^2 ] $ (depending on $A$), 
which is a minimizer for this last infimum. The corresponding $u_p$ 
is then a minimizer for the infimum in $u$ in \eqref{bernardo}. 
Writing that
$$
 \frac{d}{dp} \Big[ p^2 \Big( \int_\R \frac{(A^2-r_0^2)^2}{A^2} \ dx \Big)^{-1} 
+ 2M r_0^4 \sin^2 \Big( \frac{p - \pi r_0^2}{2 r_0^2 } \Big) \Big] = 
2 p \Big( \int_\R \frac{(A^2-r_0^2)^2}{A^2} \ dx \Big)^{-1} 
+ 2 M r_0^2 \sin \Big( \frac{ p - \pi r_0^2}{2r_0^2} \Big) = 0 , 
$$
we deduce the relations
$$
2 A^2 u = c ( A^2 - r_0^2 ) , \quad \quad \quad 
c \equiv M r_0^2 \sin \Big( \frac{ p - \pi r_0^2}{2r_0^2} \Big) .
$$
We conclude by considering a minimizing sequence $( A_n ) $ in 
\eqref{bernardo}, and translating in space so that 
$ \inf_\R A_n = \mu = |A_n|(0) $.\\

Since $F \geq 0$ in $ \R_+ $, we have
$$
 \int_\R |\p_x V_n |^2 \ dx \leq \mathscr{K} ( V_n ) \leq \frac{12}{10} E(U_0) 
$$
for $n$ large. Therefore, by compact Sobolev embedding 
$ H^1([-R,+R]) \hookrightarrow L^\ii([-R,+R])$, we may assume, up to 
a possible subsequence, that there exists $ V \in H^1_{\rm loc}(\R) $ 
such that for any $R > 0 $, $ V_n \rightharpoonup V $ in $ H^1([-R,+R]) $ 
and $ V_n \to V $ uniformly on $ [-R,+R] $. Moreover, by lower semicontinuity 
and Fatou's lemma, $ E(V) \leq \varliminf_{n \to +\ii} E(V_n)$. Since 
$ |V_n| \geq \mu > 0 $ in $\R$, $ |V| \geq \mu > 0 $ in $ \R $ and we 
may lift $ V = A \ex^{i \phi} $. Furthermore, $ \inf_\R A_n = \mu = |V_n|(0) $, 
hence $ \inf_\R A = \mu = |V|(0) $. We also know that $P(V_n) \in [ 0,  r_0^2 \pi] $ 
for all $n$, hence we may assume, up to another subsequence, that 
$P(V_n) $ converges to some $ P_\ii \in [ 0, r_0^2 \pi] $. We also set
$$
 c \equiv \lim_{n \to +\ii} c_n = M r_0^2 \sin \Big( \frac{ P_\ii - \pi r_0^2}{2r_0^2} \Big) . 
$$
In view of Step 1, and the convergence $ A_n \to A $ uniformly on 
any compact interval $ [-R,+R] $, it follows that
\be
\label{zorro}
2 A^2 \p_x \phi = c ( A^2 - r_0^2 ) 
\quad \quad \quad {\rm and} \quad \quad \quad 
\p_x \phi_n \to \p_x \phi \quad {\rm in} \quad L^\ii_{\rm loc} (\R).
\ee

Note that
$$
 \int_\R |\p_x V|^2 + \frac{1}{K} ( |V|^2 - r_0^2 )^2 \ dx 
 \leq E(V) < + \ii 
$$
hence $ |V| \to r_0 $ at $\pm \ii$. In particular, there exist 
$ - \ii < R_- \leq 0 \leq R_+ < +\ii $ such that 
$ |V| > \mu $ in $ ( - \ii , R_- ) $ and in $ ( R_+ , + \ii ) $ 
and $ |V|(R_\pm) = \mu $. \\

\noindent {\bf Step 2.} There exist $ -\ii < z_- \leq 0 \leq z_+ < + \ii $ 
such that
$$
 A(x) = A_c ( x - R_+ + z_+ ) \quad \quad {\rm for} \ x \geq R_+ 
\quad \quad \quad {\rm and} \quad \quad \quad 
 A(x) = A_c ( x - R_- + z_- ) \quad \quad {\rm for} \ x \leq R_- .
$$
We work for $ x \geq R_+ $, the other case being similar. We consider 
$ \chi \in \BC_c^1( (R_+, + \ii) , \C ) $, $ t \in \R $ small such 
that $ V^t_n \equiv v_n + t \chi $ verifies $ | V^t_n | > \mu $ in 
$ ( R_+ , + \ii ) $. This is possible since $ \inf_{{\rm Supp}(\chi)} | V_n | > \mu $. 
Then, $ | V^t_n | \geq \mu $ in $ \R $ and $ |V^t_n| (0) = \mu $, 
hence $ V^t_n $ is then a comparison map for $ \mathscr{K}_{\rm min}(\mu) $, 
and in view of the equality $ P(V^t_n ) = P(V_n ) + 2 t \int_{R_+}^{+\ii} 
\langle i \p_x V_n | \chi \rangle \ dx + \BO(t^2)$, it follows that
\begin{align*}
\mathscr{K}_{\rm min}(\mu) \leq \mathscr{K} ( V^t_n ) 
= & \,\mathscr{K}_{\rm min}(\mu) + o_{n \to +\ii} (1) 
+ 2 t \int_{R_+}^{+\ii} \langle \p_x V_n , \p_x \chi \rangle \ dx 
+ t^2 \int_{R_+}^{+\ii} | \p_x \chi |^2 \ dx 
\\ & \, 
- 2 t \int_{R_+}^{+\ii} f (|V_n|^2 ) \langle V_n , \chi \rangle \ dx 
+ M t \sin \Big( \frac{ P (V_n) - \pi r_0^2 }{r_0^2} \Big) 
\int_{R_+}^{+\ii} \langle i \p_x V_n , \chi \rangle \ dx + \BO_{t \to 0} (t^2) .
\end{align*}
Letting $ n \to +\ii $ and using the weak and strong convergences 
for $V_n$, we infer
$$
0 \leq 
2 t \int_{R_+}^{+\ii} \langle \p_x V, \p_x \chi \rangle \ dx 
- 2 t \int_{R_+}^{+\ii} f (|V|^2 ) \langle V , \chi \rangle \ dx 
- M t \sin \Big( \frac{r_0^2 \pi - P_\ii}{2 r_0^2} \Big) 
\int_{R_+}^{+\ii} \langle i \p_x V , \chi \rangle \ dx + \BO_{t \to 0} (t^2) .
$$
Dividing by $ t \not = 0 $ and letting $ t \to 0^+$ and then $ t \to 0^- $, 
we deduce that $ V $ solves (TW$_c$) in $ ( R_+ , + \ii ) $ and $V$ 
has finite energy. Moreover, $|V| (R_+) = \mu $ is small, thus 
$ V = \ex^{i\theta_+} U_c ( \cdot - R_+ +z_+ ) $ in $ ( R_+ , + \ii ) $ 
for some constants $ z_+ $ and $ \theta_+ $, and the speed $c$ is such 
that $ \inf_\R A_c = \sqrt{r_0^2 + \xi_c} \leq \mu $, hence 
$ c \leq \s(\mu) \leq K \mu $. 
Since $|V|$ has finite energy in $\R$ and solves (TW$_c$) in 
$ ( R_+ , + \ii ) $, $V$ is $\BC^1$ in $ [ R_+ , + \ii ) $. Moreover, 
$|V|$ reaches a minimum at $x = R_+$, thus we must have $ \p_x^+ |V| (R_+) \geq 0 $, 
which imposes $z_+ \geq 0 $. Note that $ A_c $ being even, it is possible 
to translate $V$ so that $ R \equiv R_+ = - R_- $ and 
$ z \equiv z_+ = - z_- $. Observe that $ \mu = A_c(z) \geq A_0 (z )$, 
hence $ z \leq K \mu $. This yields
\be
\label{cacoute}
 \int_{|x| \geq R } |\p_x V|^2 + F(|V|^2) \ dx 
= \int_{|x| \geq z } |\p_x U_c|^2 + F(|U_c|^2) \ dx 
\geq E(U_0) - K \mu .
\ee
In particular, we deduce from \eqref{paraudessus}
$$
 2 R F(\mu^2)  \leq \int_{|x| \leq R } |\p_x V|^2 + F(|V|^2) \ dx \leq K \mu ,
$$
hence $ R \leq K \mu $ for $\mu $ small ($F(0) > 0 $).

\bigskip

\noindent {\bf Step 3.} We prove that $ A = \mu $ in $ ( R_- , R_+ ) 
= ( -R , + R )$. \\

Indeed, if it is not the case, there exists a bounded interval $ ( x_- , x_+ )$ 
such that $ A = |V| > \mu $ in $ ( x_- , x_+ ) $ and $ |V| (x_\pm) = \mu $, 
with $ \p_x |V| (x_+ ) \leq 0 \leq \p_x |V| (x_- ) $. Therefore, we can 
make perturbations of the amplitude $A_n$ localized 
in $ (x_-, x_+ ) $. Hence, arguing as in Step 2, we see that then, 
$V$ solves (TW$_c$) in $ ( x_- , x_+ ) $, with $ 2 A^2 \p_x \phi = c ( A^2 - r_0^2 ) $ 
and $ |V| (x_\pm) = \mu $, $ \p_x |V| (x_+ ) \leq 0 \leq \p_x |V| (x_- ) $. 
We then are in position to apply \eqref{chero}, yielding
$$
\int_{x_-}^{x_+} |\p_x V|^2 + F(|V|^2) \ dx \geq \frac12 E(U_0) ,
$$
but the combination with \eqref{cacoute} provides
\begin{align*}
 \frac{11}{10} E(U_0) \geq \mathscr{K}_{\rm min}(\mu) 
\geq & \, 
\int_{x_-}^{x_+} |\p_x V|^2 + F(|V|^2) \ dx 
+ \int_{|x| \geq R } |\p_x V|^2 + F(|V|^2) \ dx \\ 
\geq & \, \frac12 E(U_0) + E(U_0) - K \mu_* 
= \frac32 E(U_0) - K \mu_*,
\end{align*}
which is not possible if $\mu_*$ is sufficiently small.\\

\noindent {\bf Step 4.} We have $ R = 0 $ or ($z= 0 $ and $c=\s_\mu$).\\

Indeed, assume $ R > 0 $, and consider $ \zeta \in \BC_c^1( (0 , +\ii) , \R ) $, 
$ \zeta \geq 0 $, $ t \geq 0 $ and $ V^t_n \equiv ( A_n + t \zeta ) \ex^{ i \phi_n } $, 
so that $ | V^t_n | = A_n + t \zeta \geq \mu $ in $ \R $. Since $ R > 0 $, 
we actually have $ \inf_\R | V^t_n | = \mu $ and $ V^t_n $ is a comparison 
map for $ \mathscr{K}_{\rm min}(\mu) $. Arguing as before, we thus have
\begin{align*}
\mathscr{K}_{\rm min}(\mu) \leq \mathscr{K} ( V^t_n ) 
= & \,\mathscr{K}_{\rm min}(\mu) + o_{n \to +\ii} (1) 
+ 2 t \int_0^{+\ii} \p_x A_n \p_x \zeta \ dx 
+ t^2 \int_0^{+\ii} ( \p_x \zeta )^2 \ dx 
\\ & \, 
+ 2 t \int_0^{+\ii} A_n \zeta ( \p_x \phi_n )^2 \ dx 
+ t^2 \int_0^{+\ii} \zeta^2 ( \p_x \phi_n )^2 \ dx 
- 2 t \int_0^{+\ii} f (A_n^2 ) A_n \zeta \ dx 
\\ & \, 
+ M r_0^2 t \sin \Big( \frac{P (V_n) - r_0^2 \pi}{r_0^2} \Big) 
\int_0^{R} 2 A_n \zeta \p_x \phi_n \ dx + \BO_{t \to 0} (t^2) .
\end{align*}
By \eqref{zorro}, we may pass to the limit as $ n \to +\ii $ in all 
the terms and deduce
\begin{align*}
0 \leq & \, 
2 t \int_0^{+\ii} \p_x A \p_x \zeta \ dx 
+ 2 t \int_0^{+\ii} A \zeta ( \p_x \phi )^2 \ dx 
- 2 t \int_0^{+\ii} f (A^2 ) A \zeta \ dx 
- 2 c t \int_0^{+\ii} A \zeta \p_x \phi \ dx 
+ \BO_{t \to 0} (t^2) .
\end{align*}
At this stage, we see the relevance of taking a minimizing sequence 
as chosen in Step 1, since it allows to pass to the limit in the 
nonlinear terms involving $ \p_x \phi_n $. 
As a consequence, using \eqref{zorro},
$$
 - \p^2_x A - A f(A^2) + \frac{c^2}{4} \frac{(A^2 - r_0^2)^2}{A^3} \geq 0 
$$
in the distributional sense in $ ( 0, +\ii )$. The term 
$ - A f(A^2) + \ds{ \frac{c^2}{4} \frac{(A^2 - r_0^2)^2}{A^3} } $ is continuous 
in $ \R $. However, since $ A(x) = \mu $ for $ 0 \leq x \leq R $ and 
$ A(x) = A_c ( x - R + z ) $ for $ x \geq R $, we infer 
$ - \p^2_x A = - \p_x A_c (z) \d_{x=R} $ plus a piecewise continuous 
function in the distributional sense in $ ( 0, +\ii )$. Since 
$ \p_x A_c (z) \geq 0 $ (recall that $z \geq 0 $), this forces 
to have $ \p_x A_c (z) = 0 $, that is $z=0$. Consequently, 
$ \mu = |V|(R) = A(R) = A_c(z) = A_c(0) $ and then $c=\s_\mu$.\\

In the next step, we take into account the loss in the weak 
convergence $ V_n \rightharpoonup V $.\\

\noindent {\bf Step 5.} There exists $ K > 0 $ such that
$$
E_\sharp \geq \frac{ P_\sharp}{K}, \quad \quad {\rm where} \quad \quad 
E_\sharp \equiv \varliminf_{n \to +\ii} E(V_n) - E(V) \geq 0, 
\quad \quad P_\sharp \equiv \lim_{n \to +\ii} P(V_n) - P(V) = P_\ii - P(V) .
$$

Let $ \epsilon > 0 $ be fixed but small, and pick some $ X > 0 $ large so that 
$$
\Big| E(V) - \int_{|x| \leq X} |\p_x V|^2 + F(|V|^2) \ dx \Big| 
\leq \epsilon \quad \quad \quad 
\Big| P(V) - \int_{|x| \leq X } (A^2 - r_0^2 ) u \ dx \Big| 
\leq \epsilon .
$$
We claim that there exists some small $ \bar{\mu} > 0 $, independent of 
$ \epsilon $ such that $ |V_n| \geq \bar{\mu} $ 
for $ |x| \geq X $ and $n$ large. Indeed, otherwise, we may argue 
as in Step 3 and show, as in the beginning of the proof there, that 
$ \int_{|x| \geq X} |\p_x V_n|^2 + F(|V_n|^2) \ dx \geq \frac12 E(U_0) $. 
This is not possible since
$$ 
\frac{12}{10} E(U_0) \geq \varliminf_{n \to +\ii} E(V_n) 
\geq \frac12 E(U_0) + \int_{|x| \leq X} |\p_x V|^2 + F(|V|^2) \ dx 
\geq \frac12 E(U_0) + E(V) - \epsilon ,
$$
and $ E(V) $ is close to $ E(U_0) $ as $ \mu \to 0 $. Therefore, 
as for \eqref{garcia},
$$
\Big| \int_{|x| \geq X} ( A_n^2 - r_0^2 ) u_n \ dx \Big|
 \leq \frac{K}{\bar{\mu}} \int_{|x| \geq X} |\p_x V_n|^2 + F(|V_n|^2) \ dx .
$$
Consequently,
\begin{align*}
E(V_n) - E(V) \geq & \, \int_{|x| \leq X} |\p_x V_n|^2 + F(|V_n|^2) \ dx 
- \int_{|x| \leq X} |\p_x V|^2 + F(|V|^2) \ dx 
+ \int_{|x| \geq X} |\p_x V_n|^2 + F(|V_n|^2) \ dx 
- \epsilon \\ 
\geq & \, \int_{|x| \leq X} |\p_x V_n|^2 + F(|V_n|^2) \ dx 
- \int_{|x| \leq X} |\p_x V|^2 + F(|V|^2) \ dx 
+ \frac{\bar{\mu}}{K} \Big| \int_{|x| \geq X} ( A_n^2 - r_0^2 ) u_n \ dx \Big| 
- \epsilon .
\end{align*}
Passing to the liminf and using the weak convergence in $ [-X ,+ X] $, we infer
\begin{align*}
\varliminf_{n \to + \ii } E(V_n) - E(V) \geq 
\frac{\bar{\mu}}{K} \varliminf_{n \to + \ii } 
\Big| P(V_n) - \int_{|x| \leq X} ( A_n^2 - r_0^2 ) u_n \ dx \Big| - \epsilon .
\end{align*}
However, \eqref{zorro} implies
$$
\int_{|x| \leq X} ( A_n^2 - r_0^2 ) u_n \ dx \to 
\int_{|x| \leq X} ( A^2 - r_0^2 ) u \ dx ,
$$
so that
$$
 E_\sharp \geq \frac{\bar{\mu}}{K} 
\Big| P_\ii - \int_{|x| \leq X} ( A^2 - r_0^2 ) u \ dx \Big| - \epsilon 
\geq  
\frac{\bar{\mu}}{K} \Big| P_\ii - P(V) \Big| - \Big( 1 + \frac{\bar{\mu}}{K} \Big) \epsilon 
= \frac{\bar{\mu}}{K} | P_\sharp | - \Big( 1 + \frac{\bar{\mu}}{K} \Big) \epsilon .
$$
Letting $ \epsilon \to 0 $, the conclusion follows.\\

\noindent {\bf Step 6.} There exists $K > 0$ such that, if $R > 0$, then
$$
\mathscr{K}_{\rm min}(\mu) \geq E(U_0) + \frac{\mu^2}{K} .
$$

We recall the expansion $ P(U_s) = r_0^2 \pi + s \dot{P}_0 + o(s) $ 
as $ s \to 0 $, where $ \dot{P}_0 \equiv \ds{ \frac{dP(U_s)}{ds}_{|s=0} } $. 
From the Hamilton group relation $ \ds{ \frac{dE(U_s)}{ds} = s \frac{dP(U_s)}{ds} } $, 
we also infer by integration $ E(U_s) = E(U_0) + \frac{s^2}{2} \dot{P}_0 + o(s^2) $. 
On the other hand, by definition of $c_n$,
\begin{align*}
2M r_0^4 \sin^2 \Big( \frac{P(V_n) - r_0^2 \pi}{2r_0^2} \Big) 
= & \, M r_0^4 \Big[ 1 - \cos \Big( \frac{P(V_n) - r_0^2 \pi}{r_0^2} \Big) \Big] 
\\ 
= & \, M r_0^4 \Big[ 1 - \sqrt{1 - \sin^2 \Big( \frac{P(V_n) - r_0^2 \pi}{r_0^2} \Big) } \Big] 
=  M r_0^4 \Big[ 1 - \sqrt{1 - \frac{c_n^2}{ M^2}} \Big] 
\end{align*}
for $n$ large. Here, we have used that 
$ M c_n = \sin ((r_0^2 \pi - P(V_n) ) / r_0^2 ) \to M c \in [ 0 , K \mu_*] $ 
({\it cf.} Step 2), thus $ \cos ( ( r_0^2 \pi - P(V_n)) / r_0^2) \geq 0 $, for 
otherwise, we would have, by Proposition \ref{minokink}
\begin{align*}
 \mathscr{K}(V_n) = E (V_n) + 2M r_0^4 \sin^2 \Big( \frac{P(V_n) - r_0^2 \pi}{2 r_0^2} \Big) 
 \geq & \, E(U_0 ) - K \mu + M r_0^4 \Big[ 1 + \sqrt{1 - \frac{c_n^2}{ M^2}} \Big] \\
 \geq & \, E(U_0 ) - K \mu_* + 2 M r_0^4 + \BO(\mu_*^2) , 
\end{align*}
but this contradicts \eqref{paraudessus} if $\mu_*$ is sufficiently small.

We assume $ R > 0 $, so that, by Step 4, $z=0 $ and $ c = \s_\mu $. 
We recall the expansion $\s _\mu = \ds{ \frac{2 \mu \sqrt{F(0)}}{r_0^2} + \BO(\mu^2) \sim 
\frac{2 \mu \sqrt{F(0)}}{r_0^2} } $. 
By definition of $E_\sharp $, one has
$$
E_\sharp + E(V) + M r_0^4 \Big[ 1 - \sqrt{1 - \frac{c^2}{ M^2}} \Big] 
\leq \varliminf_{n \to +\ii} E(V_n) + \lim_{n \to +\ii} 2M r_0^4 \sin^2 \Big( \frac{P(V_n) - r_0^2 \pi}{2 r_0^2} \Big) 
= \varliminf_{n \to +\ii} \mathscr{K} (V_n) = 
\mathscr{K}_{\rm min}(\mu) 
$$
since $ ( V_n ) $ is minimizing for $ \mathscr{K}_{\rm min}(\mu) $. 
Moreover, from the expression of $V$, it holds (for $R > 0$)
$$
 E(V) = E( U_{\s_\mu}) + 2 R \Big[ \frac{\s_\mu^2(r_0^2- \mu^2)^2}{4\mu^2} + F(\mu^2) \Big] 
\quad \quad \quad {\rm and} \quad \quad \quad 
 P(V) = P( U_{\s_\mu}) + R \s_\mu \frac{(r_0^2- \mu^2)^2}{ \mu^2 }.
$$
Furthermore, $ P_\sharp = P_\ii - P(V) $ and $ c = M r_0^2 \sin( ( r_0^2 \pi -P_\ii ) / r_0^2 ) $ 
with $ P_\ii \in [ 0, r_0^2 \pi ] $ and $ \cos ( ( r_0^2 \pi - P_\ii) / r_0^2) \geq 0 $, thus
$$ 
P_\sharp = P_\ii - P(V) = r_0^2 \pi - r_0^2 \arcsin \Big( \frac{c}{M r_0^2 } \Big) - 
P( U_{\s_\mu}) - R \s_\mu \frac{(r_0^2- \mu^2)^2}{ \mu^2 } . 
$$
Combining this with the expansion of $ E( U_\s) $ and $ P( U_\s) $ gives
\begin{align*}
 \mathscr{K}_{\rm min}(\mu) 
 \geq & \, 
 E(U_0) + E_\sharp + \frac{\s^2_\mu }{2} \dot{P}_0 + o(\s^2_\mu) 
 + M \Big[ 1 - \sqrt{1 - \frac{\s^2_\mu }{ M^2}} \Big] 
 + 2 R \Big[ \frac{\s_\mu^2(r_0^2- \mu^2)^2}{4\mu^2} + F(\mu^2) \Big] \\
 \geq & \, 
 E(U_0) + \frac{| P_\sharp | }{K} 
 + \frac{\s^2_\mu}{2} \Big[ \dot{P}_0 + \frac{1}{M} \Big] + o(\mu^2) 
 + 4 R F(0) \\
 \geq & \, 
 E(U_0) + \frac{1 }{K} \Big| r_0^2 \arcsin( \s_\mu / ( M r_0^2) ) + \s_\mu \dot{P}_0 
 + R \s \frac{(r_0^2- \mu^2)^2}{ \mu^2 } + o(\s_\mu) \Big|
 + \frac{\s^2_\mu}{2} \Big[ \dot{P}_0 + \frac{1}{M} \Big] + o(\mu^2) + 4 R F(0) \\
 \geq & \, 
 E(U_0) + \frac{1 }{K} \Big| \frac{ \s_\mu}{ M} + \s_\mu \dot{P}_0 + 
 R \s_\mu \frac{(r_0^2- \mu^2)^2}{ \mu^2 } + o(\s_\mu) \Big|
 + \frac{\s^2_\mu}{2} \Big[ \dot{P}_0 + \frac{1}{M} \Big] + o(\mu^2) + 4 R F(0) \\
 \geq & \, 
 E(U_0) + \frac{ \s_\mu}{K} \Big| \dot{P}_0 + \frac{1 }{ M} 
 + R \frac{(r_0^2- \mu^2)^2}{ \mu^2 } \Big|
 + \frac{\s^2_\mu}{2} \Big[ \dot{P}_0 + \frac{1}{M} \Big] + 4 R F(0) + o(\mu^2) .
\end{align*}
The right-hand side is a continuous piecewise affine function of 
$R$ (the "$o$" does not depend on $R$). Since 
$ \ds{ \frac{\s_\mu (r_0^2- \mu^2)^2}{ K \mu^2 } \simeq \frac{1}{\mu} \gg 4F(0) } $ 
and $ \ds{ \dot{P}_0 + \frac{1}{M} } < 0 $ (since $ M > - \dot{P}_0^{-1} $ by 
hypothesis), it follows that the right-hand side is a function of $R$ which is 
decreasing in $ [ 0, R_0 (\mu)] $ and increasing in $[ R_0(\mu), + \ii )$, with 
$$ R_0 (\mu) \equiv - \Big( \dot{P}_0 + \frac{1}{M} \Big) \frac{ \mu^2 }{ (r_0^2- \mu^2)^2} 
\sim - \Big( \dot{P}_0 + \frac{1}{M} \Big) \frac{\mu^2}{r_0^4} > 0 . $$
Therefore, using once again that $ \s^2_\mu \ds{ \sim \frac{4 \mu^2 F(0)}{r_0^4} } $,
\begin{align*}
 \mathscr{K}_{\rm min}(\mu) \geq 
 & \, 
E(U_0) + \frac{\s^2_\mu}{2} \Big[ \dot{P}_0 + \frac{1}{M} \Big] + 4 R_0(\mu) F(0) + o(\mu^2) 
\\ 
= & \, E(U_0) + \Big[ \dot{P}_0 + \frac{1}{M} \Big] \frac{2 \mu^2 F(0)}{r_0^4} 
- \Big[ \dot{P}_0 + \frac{1}{M} \Big] \frac{4 \mu^2 F(0)}{r_0^4} + o(\mu^2) 
\\ 
= & \, E(U_0) - \mu^2 \frac{2 F(0)}{r_0^4} \Big[ \dot{P}_0 + \frac{1}{M} \Big] 
+ o(\mu^2) .
\end{align*}
In view of our hypothesis $ \dot{P}_0 + \frac{1}{M} < 0 $, we infer that
$$
 \mathscr{K}_{\rm min}(\mu) \geq E(U_0) + \frac{\mu^2}{K} 
$$
for $ \mu_*$ sufficiently small and some positive constant $K$, as wished. 
If the assumption $ \dot{P}_0 + \frac{1}{M} < 0 $ is not verified, but if 
$ \dot{P}_0 + \frac{1}{M} > 0 $ for instance, then the function of $R$ above is 
increasing in $ [ 0 , +\ii ) $, with minimum value achieved at $R=0$ 
and equal to
$$ 
E(U_0) + \frac{\s^2_\mu}{2} \Big[ \dot{P}_0 + \frac{1}{M} \Big] + o(\mu^2) 
= E(U_0) + \frac{2 \mu^2 F(0)}{r_0^4} \Big[ \dot{P}_0 + \frac{1}{M} \Big] + o(\mu^2) 
\geq E(U_0) + \frac{\mu}{K} . 
$$
We then would have concluded a stronger estimate, which is actually in contradiction 
with \eqref{paraudessus}, hence we are necessarily in the case $ R > 0 $. The 
assumption $ \dot{P}_0 + \frac{1}{M} < 0 $ is however crucial for the last step. \\

\noindent {\bf Step 7.} We assume $ \dot{P}_0 + \frac{1}{M} < 0 $. Then, 
for $\mu_*$ sufficiently small, the case $ R = 0 $ does not occur.\\ 

We argue in a similar way, but since $ R = 0 $, the expressions for 
$E(V)$ and $P(V)$ are given by
$$
 E(V) = E( U_c) - 4 \int_0^z F(|U_c|^2) \ d x 
 \quad \quad \quad {\rm and} \quad \quad \quad 
 P(V) = P( U_c) - 2 \int_0^z \frac{c}{2} \frac{(r_0^2- A_c^2)^2}{ A_c^2 } \ d x .
$$
Here, we have used that $ |\p_x U_c|^2 = F(|U_c|^2) $ since 
$U_c$ solves (TW$_c$). Combining this here again with the expansion 
of $ E( U_c) $ and $ P( U_c) $ gives, using that $ 0 \leq c \leq K \mu $,
\begin{align*}
 \mathscr{K}_{\rm min}(\mu) 
 \geq & \, 
 E(U_0) + E_\sharp + \frac{c^2 }{2} \dot{P}_0 + o(c^2) 
 + M \Big[ 1 - \sqrt{1 - \frac{c^2 }{ M^2}} \Big] 
 - 4 \int_0^z F(|U_c|^2) \ d x
 \\
 \geq & \, 
 E(U_0) + \frac{| P_\sharp | }{K} 
 + \frac{c^2}{2} \Big[ \dot{P}_0 + \frac{1}{M} \Big] - 4 z F(0) 
 + o(\mu^2) 
\\
 \geq & \, 
 E(U_0) + \frac{1 }{K} \Big| \arcsin( c / M) + c \dot{P}_0 
 - c \int_0^z \frac{(r_0^2- A_c^2)^2}{ A_c^2 } \ d x + o(c) \Big| 
 + \frac{c^2}{2} \Big[ \dot{P}_0 + \frac{1}{M} \Big] + o(\mu^2) 
 - 4 z F(0) .
\end{align*}
Following the lines of the proof of Lemma \ref{prolongement}, 
we have
\be
\label{jdevlop}
c \int_0^z \frac{(r_0^2- A_c^2)^2}{ A_c^2 } \ d x 
= 2 \arctan \Big( \sqrt{\frac{\mu^2}{r_0^2 + \xi_c} } \Big) 
+ \BO(\mu^2) .
\ee
Indeed, noticing that $ A_c = \BO(\mu) $ in $[0, z ] $ with 
$ z \leq K \mu $, we write, expanding the square,
$$
\int_0^z \frac{(r_0^2- A_c^2)^2}{ A_c^2 } \ d x 
= \int_0^z \frac{r_0^4}{ A_c^2 } - 2 + A_c^2 \ d x 
= \int_0^z \frac{r_0^4}{ A_c^2 } \ dx + \BO(\mu) . 
$$
Then, using the change of variable $ \xi = \eta_c (x) $,
\begin{align}
\label{bienfoutue}
\int_0^z \frac{(r_0^2- A_c^2)^2}{ A_c^2 } \ d x 
= & \, \int_{\xi_c}^{\mu^2-r_0^2} \frac{r_0^4}{ (r_0^2 + \xi ) \sqrt{-\BV_c(\xi)} } \ d\xi 
+ \BO (\mu) \nonumber
\\ = &\, \int_{\xi_c}^{\mu^2-r_0^2} \frac{r_0^4}{ (r_0^2 + \xi ) \sqrt{-\BV_c'(\xi_c)(\xi - \xi_c)} } \ d\xi 
 \nonumber \\ & \, 
+ \int_{\xi_c}^{\mu^2-r_0^2} \frac{r_0^4}{ (r_0^2 + \xi ) } 
\Big( \frac{1}{\sqrt{-\BV_c(\xi)}} -  \frac{1}{ \sqrt{-\BV_c'(\xi_c)(\xi - \xi_c)}}\Big) \ d\xi 
+ \BO(\mu) \nonumber 
\\ = & \, \frac{2}{c} \arctan \Big( \sqrt{\frac{\mu^2}{r_0^2 + \xi_c} - 1} \Big)
+ \BO(\mu) ,
\end{align}
by computations similar to those for the proof of Lemma \ref{prolongement}. 
This proves \eqref{jdevlop}. Therefore,
\be
\label{pasmal}
 \mathscr{K}_{\rm min}(\mu) 
 \geq E(U_0) + \frac{1 }{K} \Big| c \Big[  \dot{P}_0 + \frac{1}{ M} \Big] 
 - 2 \arctan \Big( \sqrt{\frac{\mu^2}{r_0^2 + \xi_c} - 1 } \Big) + o(c) \Big| 
 + \frac{c^2}{2} \Big[ \dot{P}_0 + \frac{1}{M} \Big] + o(\mu^2) - 4 z F(0) . 
\ee
By \eqref{paraudessus}, the left-hand side is $ \leq E(U_0) + K \mu^2 $. 
Since $ \dot{P}_0 + \frac{1}{M} < 0 $, $c \leq K \mu $, $z \leq K \mu $ 
and $F(0) > 0 $, this implies
$$
 \Big| c \Big[  \dot{P}_0 + \frac{1}{ M} \Big] 
 - 2 \arctan \Big( \sqrt{\frac{\mu^2}{r_0^2 + \xi_c} } \Big) + o(c) \Big| 
 \leq K \mu ,
 $$
thus
$$
 \arctan \Big( \sqrt{ \frac{\mu^2}{r_0^2 + \xi_c} - 1} \Big) \leq K \mu
$$
and finally, for $\mu_*$ small enough,
$$
 0 \leq \frac{\mu^2}{r_0^2 + \xi_c} - 1\leq K \mu^2 . 
$$
Combining this with the equality $ r_0^2 + \xi_c = \ds{ \frac{c^2 r_0^4}{4 F(0)} } 
+ \BO(c^4) $ seen during the proof of Lemma \ref{prolongement}, we infer
$$
 c = \frac{2 \sqrt{F(0)}}{ r_0^2 } \mu + \BO(\mu^2) .
$$
In particular, going back to \eqref{bienfoutue} and since for 
$ 0 \leq x \leq z $,
$$
r_0^2 + \xi_c = A_c^2 (0) \leq A_c^2(x) \leq A_c^2 (z) = \mu^2 , 
$$
this implies
$$
\frac{z r_0^4}{ \mu^2 } \leq \int_0^z \frac{r_0^4}{ A_c^2 } \leq 
\frac{2}{c} \arctan \Big( \sqrt{\frac{\mu^2}{r_0^2 + \xi_c} - 1} \Big)
+ \BO(\mu) \leq \frac{K \mu}{c}+ K \mu  \leq K ,
$$
which provides (since $ c \approx \mu $)
$$
 z \leq K \mu^2 . 
$$
Inserting this into \eqref{pasmal} and keeping in mind that the left-hand 
side is $ \leq E(U_0) + K \mu^2 $, we deduce
$$
c \Big[ \dot{P}_0 + \frac{1}{ M} \Big] 
- 2 \arctan \Big( \sqrt{\frac{\mu^2}{r_0^2 + \xi_c} - 1 } \Big) = o ( \mu ) .
$$
However, since 
$ \arctan \Big( \ds{ \sqrt{\frac{\mu^2}{r_0^2 + \xi_c} - 1 } } \Big) \geq 0 $, 
this gives
$$
o ( \mu ) \leq c \Big[ \dot{P}_0 + \frac{1}{ M} \Big] 
\sim \frac{2 \mu \sqrt{F(0)}}{r_0^2} \Big[  \dot{P}_0 + \frac{1}{ M} \Big] ,
$$
yielding a contradiction for small $\mu$ since we have 
$ \dot{P}_0 + \frac{1}{ M} < 0 $ by assumption. Therefore, the case 
$ R = 0 $ does not occur for sufficiently small $\mu_*$. If we had 
$ \dot{P}_0 + \frac{1}{M} > 0 $, we would not have been able to show 
that $ \mathscr{K}_{\rm min}(\mu) $ gives a control on $ \mu$.\\ 

The proof of Proposition \ref{supercool} is complete. \carre \\

\noindent {\bf Proof of Theorem \ref{statiomini}.} Let $ U \in \mathscr{V}_{\mu_*} $. 
If $ \mu \equiv \inf_\R |U| > 0 $, then Proposition \ref{supercool} gives 
$ \mathscr{K} (U) \geq E(U_0) + \mu^2 /K > E(U_0) = \mathscr{K} (U_0) $. 
If $ \inf_\R |U| = 0 $, we deduce from Proposition \ref{minokink} that 
$ \mathscr{K} (U) \geq E(U_0) + 2 M r_0^4 \sin^2( (\mathfrak{P}(U) - \pi r_0^2 ) r_0^2 ) $. 
Hence $ \mathscr{K} (U) > E(U_0) $ except if $ \mathscr{K} (U) = E(U_0) $. 
From the study of the equality case in Proposition \ref{minokink}, it follows 
that $ U \in \{ \ex^{i\theta} U_0( \cdot -y), y \in \R, \theta \in \R \} $, as claimed.

%%%%%%%%%%%%%%%%%%%%%%%%%%%%%%%%%%%%%%%%%%%%
\subsection{Proof of Theorem \ref{stabmini}}

As a first step, we shall need a quantified version of Proposition \ref{minokink}.

%%%%%%%%%%%%%%%%%%%%
\begin{propositions} 
\label{minokinkameliore}
There exists $ \epsilon_0 > 0 $ and $ K > 0 $, depending only 
on $f$ such that, for any $ U \in \mathcal{Z} $ verifying
$$
\mathscr{K}(U) - E(U_0) \leq \epsilon_0 
\quad \quad \quad {\it and}\quad \quad \quad
\inf_\R |U| \leq \epsilon_0 ,
$$
there holds
$$ 
 \inf_{ \scriptsize{\begin{array}{c} y \in \R \\ \theta \in \R \end{array}}} 
d_\mathcal{Z} ( U , \ex^{i\theta} U_0 (\cdot - y ) ) 
\leq K \Big( \mathscr{K}(U) - E(U_0) + \inf_\R |U| \Big)^{1/4}.
$$
\end{propositions}
%%%%%%%%%%%%%%%%%%%%

\noindent {\it Proof.} First, we translate the problem in space 
so that $ \mu \equiv \inf_\R |U| = |U|(0) $ and shall choose the phase factor 
later. We follow the lines of the proof of 
Proposition \ref{minokink} and actually get (writing 
$ U = A \ex^{i\phi} $ locally in $ \{ |U| > 0 \} $)
\begin{align*}
\int_0^{+\ii} |\p_x U|^2 + F(|U|^2) \ dx 
= & \, 
\int_0^{+\ii} {\bf 1}_{|U| > 0 } A^2 (\p_x \phi )^2 \ dx 
+ \int_0^{+\ii} |\p_x |U| |^2 + F(|U|^2) \ dx 
\\ = & \, 
\int_0^{+\ii} {\bf 1}_{|U| > 0 } A^2 (\p_x \phi )^2 \ dx 
+ \int_0^{+\ii} \Big[ \sqrt{ F(|U|^2)} - | \p_x |U| | \Big]^2 \ dx 
\\ & \, 
+ 2 \int_0^{+\ii} | \sqrt{ F(|U|^2)} \p_x |U| | \ dx 
\\ \geq & \, 
\int_0^{+\ii} {\bf 1}_{|U| > 0 } A^2 (\p_x \phi )^2 \ dx 
+ \int_0^{+\ii} \Big[ \sqrt{ F(|U|^2)} - | \p_x |U| \, | \Big]^2 \ dx 
+ 2 \int_\mu^{r_0} \sqrt{ F(s^2)} \ ds .
\end{align*}
Arguing similarly in $ ( -\ii, 0)$, we get
\be
\label{merou}
E(U) \geq E(U_0) 
+ \int_\R {\bf 1}_{|U| > 0 } A^2 (\p_x \phi )^2 \ dx 
+ \int_\R \Big[ \sqrt{ F(|U|^2)} - | \p_x |U| | \Big]^2 \ dx 
- 4 \int_0^\mu \sqrt{ F(s^2)} \ ds .
\ee
The gradient of the phase is controlled using \eqref{merou}. We shall now 
estimate the modulus part. Let us denote $ A \equiv |U| $ and
$$
 h \equiv \p_x A - \sqrt{ F(A^2)} ,
$$
for which we have, by \eqref{merou},
\be
\label{bornh}
 \n h \n_{L^2(\R)}^2 \leq E(U) - E(U_0) + 4 \int_0^\mu \sqrt{ F(s^2)} \ ds 
\leq E(U) - E(U_0) + K \mu .
\ee
Recall that $U_0$ verifies $ (\p_x U_0)^2 = F(U_0^2) $ in $\R$, hence 
$ \p_x U_0 = \sqrt{ F(U_0^2)} $ in $\R_+$. Setting $ \Theta \equiv A - |U_0| $, we infer
$$
\p_x \Theta = \sqrt{ F(A^2)} - \sqrt{ F(U_0^2)} + h 
\quad \quad \quad {\rm in} \quad \R_+ . 
$$
We set, for $x \geq 0 $,
$$
G( x, \theta ) \equiv \sqrt{ F((U_0(x) + \theta)^2)} - \sqrt{ F(U_0^2(x))} 
+ \frac{ U_0 (x) f(U_0^2(x) )\theta }{\sqrt{ F(U_0^2(x))} }
$$
Since $U_0$ verifies $ \p^2_x U_0 + U_0 (x) f(U_0^2(x) ) = 0 $ and 
$ \p_x U_0 = \sqrt{ F(U_0^2(x))} $ in $ \R_+ $, it follows that
$$
G( x, \theta ) = \sqrt{ F((U_0(x) + \theta)^2)} - \sqrt{ F(U_0^2(x))} 
- \frac{ \p^2_x U_0(x)}{\p_x U_0(x)} \theta .
$$
Moreover, by Taylor expansion, we infer the existence of $K >0 $ and 
$\theta_0 > 0 $ such that, for $ |\theta | \leq \theta_0 $, $ x \in \R_+ $,
$$
| G( x, \theta ) | \leq K \theta^2 .
$$
The estimate is clearly uniform in view of the exponential decay of 
$ \p_x U_0 $ at infinity. Therefore,
\be
\label{jolieedo}
\p_x \Theta = \frac{ \p^2_x U_0(x)}{\p_x U_0(x)} \Theta + G( x, \Theta ) + h(x) . 
\ee
We view this ode as a linear ode with source term $ G( x, \Theta(x) ) + h (x) $. 
Since $ \p_x U_0 $ solves the homogeneous equation, we infer from Duhamel's 
formula and the fact that $ \Theta(0) = A(0) - U_0 (0) = |U(0)| = \mu $ that for $ x \geq 0 $,
\be
\label{Duhamel}
\Theta(x) = \mu + \p_x U_0 (x) \int_0^x \frac{ G( z, \Theta (z)) + h (z) }{ \p_x U_0(z) } \ dz .
\ee
We shall prove that this equation implies that if $\mu $ and $ |\!| h |\!|_{L^2(\R_+)} $ 
are sufficiently small, then
\be
\label{normeLinf}
 \n \Theta \n_{L^2(\R_+)} \leq K \Big( \n h \n_{L^2(\R_+)} + \mu \Big) .
\ee
We assume $ \mu < \theta_0 / 2 $. Note that since $U_0$ is a kink, we have 
the decays given in Proposition \ref{decroiss}. Hence, there exists two positive 
constants $K_1$ and $K_2$ such that
$$
\forall x \in \R_+, \quad \quad 
\frac{\ex^{-\cs x}}{K_1} \leq \p_x U_0(x) \leq K_2 \ex^{-\cs x} .
$$
In particular, if $ | \Theta (x) | \leq \theta_0 $ in the interval 
$ [ 0 , R] $, then \eqref{Duhamel} implies, for $ x \in [ 0 , R] $,
\begin{align*}
| \Theta(x) | \leq & \, \mu + K_1 K_2 \ex^{-\cs x} \int_0^x 
\ex^{\cs z} \Big[ K \n \Theta \n_{L^\ii([0,R])} |\Theta (z)| + |h| (z) \Big] \ dz \\
\leq & \, \mu + 
\frac{K K_1 K_2}{\cs} \n \Theta \n_{L^\ii([0,R])}^2  
+ \frac{K_1 K_2}{ \sqrt{ 2\cs} } \n h \n_{L^2(\R_+)}
\end{align*}
by Cauchy-Schwarz. We thus choose $|\!| h |\!|_{L^2(\R_+)} + \mu $ sufficiently 
small so that
$$
 4 \Big( \mu + \frac{K_1 K_2}{ \sqrt{ 2\cs} } \n h \n_{L^2(\R)} \Big) 
 \leq \tilde{\theta}_0 \equiv 
 \min \Big\{ \theta_0 , \frac{\cs}{2 K K_1 K_2} \Big\} .
$$
Then, we consider the set $\BR $ of all $ R > 0 $ such that 
$ | \Theta (x) | \leq \tilde{\theta}_0 $ in the 
interval $ [ 0 , R] $. Since $\Theta \in H^1( \R, \C ) $ is continuous 
by Sobolev embedding and $ | \Theta (0) | = \mu < \tilde{\theta}_0 $, 
$ \BR \not = \emptyset $ and is closed in $\R_+^*$. Moreover, the above 
estimate shows that for $R \in \BR $, 
$$
\n \Theta \n_{L^\ii([0,R])} \leq \mu + 
\frac{K K_1 K_2}{\cs} \n \Theta \n_{L^\ii([0,R])}^2  
+ \frac{K_1 K_2}{ \sqrt{ 2\cs} } \n h \n_{L^2(\R_+)} ,
$$
which gives
$$
\n \Theta \n_{L^\ii([0,R])} \Big( 1 - \frac{K K_1 K_2}{\cs} \n \Theta \n_{L^\ii([0,R])} \Big) 
\leq \mu + \frac{K_1 K_2}{ \sqrt{ 2\cs} } \n h \n_{L^2(\R_+)} ,
$$
and then
\be
\label{borninfi}
\n \Theta \n_{L^\ii([0,R])} \leq 2 \Big[ 
\mu + \frac{K_1 K_2}{ \sqrt{ 2\cs} } \n h \n_{L^2(\R)} \Big] 
\leq \frac{\tilde{\theta}_0}{2} < \tilde{\theta}_0 .
\ee
Consequently, $\BR $ is open in $\R_+^*$. By connexity, 
$ \BR = \R_+^* $, proving \eqref{normeLinf}. In what follows, 
we assume $ |\!| h |\!|_{L^2(\R_+)} + \mu $ sufficiently small 
so that $ |\!| \Theta |\!|_{L^\ii} \leq \tilde{\theta}_0 $, thus 
$ | G( x, \Theta ) | \leq K \Theta^2 $. In particular, 
$$
| \Theta(x) | \leq \mu + K_1 K_2 \int_0^x 
\ex^{ - \cs ( x- z)} \Big[ K \n \Theta \n_{L^\ii([0,R])} |\Theta (z)| + |h| (z) \Big] \ dz .
$$
For $ R > 0 $ to be determined later, we then deduce from classical 
convolution estimates that
$$
 \n \Theta \n_{L^2([0,R])} \leq \mu \sqrt{R} + K_3 \n \Theta \n_{L^\ii(\R_+)} 
 \n \Theta \n_{L^2([0,R])} + K_3 \n h \n_{L^2(\R_+)} .
$$
Imposing $ |\!| h |\!|_{L^2(\R_+)} + \mu $ smaller if necessary, 
we may assume that
$$
 K_3 \n \Theta \n_{L^\ii(\R_+)} \leq K_3 K ( \n h \n_{L^2(\R_+)} + \mu ) \leq \frac12 ,
$$
so that we get
$$
 \n \Theta \n_{L^2([0,R])} \leq K_4 \Big( \mu \sqrt{R} + \n h \n_{L^2(\R_+)} \Big) .
$$
Reporting this into \eqref{jolieedo} provides
$$
 \n \p_x \Theta \n^2_{L^2([0,R])} \leq K_5 \Big( \mu^2 R + \n h \n_{L^2(\R_+)}^2 \Big) . 
$$
Arguing similarly in $[ -R, 0] $ and using \eqref{bornh}, we obtain an $H^1$ estimate 
for $\Theta$ in $ [ -R,R] $:
\be
\label{near} 
\n \Theta \n_{H^1([-R,R])}^2 \leq K_6 \Big( E(U) - E(U_0) + \mu^2 R + \mu \Big)  . 
\ee
We now turn to the estimate in $ \{|x| \geq R \} $. For that purpose, we write
\begin{align}
\int_{|x| \geq R } ( \p_x |U| )^2 + \frac{1}{K} ( |U|^2 - r_0^2 )^2 \ dx 
\leq &\, E(U) - E(U_0) + \int_{ |x| \geq R } ( \p_x U_0 )^2 + F ( U_0^2 ) \ dx 
\nonumber \\ & \, 
\label{plouf}
- \int_{|x| \leq R } ( \p_x |U| |)^2 + F ( |U|^2 ) \ dx 
+ \int_{|x| \leq R } ( \p_x U_0 )^2 + F ( U_0^2 ) \ dx .
\end{align}
Since $U_0$ decays exponentially (see Proposition \ref{decroiss}), it follows that
$$
 \int_{|x| \geq R } | \p_x U_0 |^2 + F ( U_0^2 ) \ dx \leq K \ex^{- \cs R} . 
$$
Furthermore, by integration by parts,
\begin{align*}
& - \int_{|x| \leq R } ( \p_x |U| )^2 + F ( |U|^2 ) \ dx 
+ \int_{|x| \leq R } ( \p_x U_0 )^2 + F ( U_0^2 ) \ dx \\
= & \, - \int_{|x| \leq R } 2 \p_x U_0 \p_x \Theta - 2 U_0 f ( U_0^2 ) \Theta \ dx 
- \int_{|x| \leq R } ( \p_x \Theta )^2 + F ( [U_0 + \Theta]^2 ) - F ( U_0^2 ) -2U_0 F'(U_0^2) \Theta \ dx 
\\ 
\leq & \, 
\int_{|x| \leq R } 2 \Theta [ \p^2_x U_0 + U_0 f ( U_0^2 ) ] \ dx 
- 2 \Theta (+R) \p_x U_0 (+R) + 2 \Theta (-R) \p_x U_0 (-R) 
+ K \n \Theta \n_{H^1([-R,+R])}^2 \\ 
\leq & \, K \ex^{ - \cs R } + K \Big(  E(U) - E(U_0) + \mu^2 R + \mu \Big) .
\end{align*}
For the before last line, we have used that 
$ \theta \mapsto F ( [U_0 + \theta]^2 ) - F ( U_0^2 ) - 2 U_0 F' ( U_0^2 ) \theta $ is 
$ \BO (\theta^2 )$ as $ \theta \to 0 $ and for the last line, that $U_0$ solves 
$ \p^2_x U_0 + U_0 f ( U_0^2 ) = 0 $, the exponential decay of $\p_x U_0$ and the 
uniform bound on $ \Theta $. Reporting these estimates into \eqref{plouf} provides
\begin{align*}
\n \Theta \n_{H^1( \{|x| \geq R \} )}^2 
= & \, 
\int_{|x| \geq R } ( \p_x |U| - \p_x | U_0 | )^2 + ( |U| - | U_0 | )^2 \ dx 
\\ \leq & \, 
2 \int_{|x| \geq R } ( \p_x |U| )^2 + ( \p_x | U_0 | )^2 + ( |U| - r_0 )^2 
+ ( | U_0 | - r_0 )^2 \ dx 
\\ \leq & \, 
K \Big[ E(U) - E(U_0) + \ex^{- \cs R} + \mu^2 R + \mu \Big] .
\end{align*}
Combining this with \eqref{near}, we deduce that for any $R > 0$, there holds
$$
 \n \Theta \n_{H^1( \R )}^2 \leq K \Big[ E(U) - E(U_0) + \ex^{- \cs R} + \mu^2 R + \mu \Big] .
$$
We then choose $ R = \mu^{-1} $ if $ \mu > 0 $ or $ R \to +\ii $ if $\mu = 0 $, and get
$$
\n \Theta \n_{H^1( \R )} \leq K \sqrt{ E(U) - E(U_0) + \mu } .
$$
Notice that if $ f' < 0 $ everywhere, then we may give a quick proof of the 
above estimate, since, using here again integration by parts and that $\p^2_x U_0 + U_O f(U_O^2) = 0 $, 
we may deduce that
$$ 
E(U) - E(U_0) \geq - 4 \mu \p_x U_0 (0) + \int_0^{+\ii} (\p_x \Theta )^2 \ dx 
+ \int_\R F( ( U_0 + \Theta)^2 ) - F( U_0^2) - 2 U_0 \Theta F' ( U_0^2 ) \ dx ,
$$
and since $ f' < 0 $, $ F( ( U_0 + \theta)^2 ) - F( U_0^2) - 2 U_0 \theta F' ( U_0^2 ) \geq \theta^2 / K $ 
by Taylor expansion, providing the desired $H^1 $ bound on $ \Theta $.

Observe now that
$$
 \mathscr{K}(U) - E(U_0) \geq E(U) - E(U_0) \geq \int_\R {\bf 1}_{|U|>0} [ A \p_x \phi ]^2 \ dx ,
$$
hence
\begin{align}
\n \p_x U - \p_x U_0 \n_{L^2(\R)} = & \, \n \p_x ( |U_0| + \Theta ) \ex^{i\phi} {\bf 1}_{|U|>0}
+ i {\bf 1}_{|U|>0} A \p_x \phi \ex^{i\phi} - \p_x U_0 \n_{L^2(\R)} \nonumber \\ 
\leq  & \, \n \ex^{i\phi}  {\bf 1}_{|U|>0} \p_x |U_0| - \p_x U_0 \n_{L^2(\R)} 
+ \n {\bf 1}_{|U|>0} A \p_x \phi \n_{L^2(\R)} + \n \Theta \n_{L^2(\R)} \nonumber \\ 
\label{morue}
\leq  & \, \n \ex^{i\phi}  {\bf 1}_{|U|>0} \p_x |U_0| - \p_x U_0 \n_{L^2(\R)} 
%\n {\bf 1}_{|U|>0} \{ \ex^{i\phi} \p_x |U_0| - \p_x U_0 \} \n_{L^2(\R)} 
%+ \n {\bf 1}_{|U|=0} \p_x U_0 \n_{L^2(\R)} 
+ K \Big[ \mathscr{K}(U) - E(U_0) + \mu \Big]^{1/2} .
\end{align}
We distinguish now the cases $ \mu = 0 $ and $ \mu > 0 $, and begin with the assumption 
$ \mu > 0 $. Then, we have a global lifting $ U = A \ex^{i \phi } $ and
\begin{align*}
d_{\mathcal{Z}} ( U, U_0 ) = & \, 
\n \p_x U - \p_x U_0 \n_{L^2(\R)} + 
\n | U | - |U_0| \n_{L^2(\R)} + | U(0) - U_0 (0) | 
\\
= & \, 
\n \p_x U - \p_x U_0 \n_{L^2(\R)} + \n \Theta \n_{L^2(\R)} + \mu
\\ \leq & \, 
\n \ex^{i\phi} \p_x |U_0| - \p_x U_0 \n_{L^2(\R)} 
+ K \Big[ \mathscr{K}(U) - E(U_0) + \mu \Big]^{1/2} .
\end{align*}
Now, we notice that
\begin{align}
\n \ex^{i\phi} \p_x |U_0| - & \, \p_x U_0 \n_{L^2(\R)}^2
= 2 \int_\R \Big[ ( \p_x U_0 )^2 - \p_x U_0 \p_x |U_0| \cos(\phi) \ dx 
\nonumber \\ 
\label{espadon}
= & \, 2 \int_0^{+\ii} ( \p_x U_0 )^2 \Big( 1 - \cos(\phi) \Big) \ dx 
+ 2 \int_{-\ii}^0 ( \p_x U_0 )^2  \Big( 1 + \cos(\phi) \Big) \ dx 
\end{align}
and that
\be
\label{chapon}
 \mathscr{K}(U) - E(U_0) \geq 2 M r_0^4 \sin^2 \Big( \frac{\mathfrak{P}(U) - r_0^2 \pi }{2 r_0^2} \Big) 
\geq \frac{1}{K} ( P(U) - r_0^2 \pi \, {\rm mod} \, 2 \pi r_0^2 )^2 .
\ee
We define $ \d = ( \mathscr{K}(U) - E(U_0) + \mu )^{1/4} $. By Cauchy-Schwarz, we have
\begin{align*}
 \Big| \int_{ |x| \geq \d } ( A^2 - r_0^2 ) \p_x \phi \ dx \Big| 
\leq & \, 
\frac{K}{\inf_{|x| \geq \d} A } \Big( \int_{ |x| \geq \d } ( A^2 - r_0^2 )^2 \ dx \Big)^{1/2} 
\Big( \int_{|x| \geq \d } ( A \p_x \phi )^2 \ dx \Big)^{1/2} \\
\leq & \, 
\frac{K}{\inf_{|x| \geq \d} A } \Big( E(U) - E(U_0) + \mu \Big)^{1/2} .
\end{align*}
Inserting this into \eqref{chapon} gives
\begin{align*}
 \Big| \int_{|x| \leq \d } ( A^2 - r_0^2 ) \p_x \phi \ dx - r_0^2 \pi \, {\rm mod} \, 2 \pi r_0^2 \Big| 
\leq & \, K \Big[ ( \mathscr{K}(U) - E(U_0) )^{1/2}
+ \frac{1}{\inf_{|x| \geq \d} A } \Big( E(U) - E(U_0) + \mu \Big)^{1/2} \Big] \\ 
\leq & \, \frac{K}{\inf_{|x| \geq \d} A } \Big( \mathscr{K}(U) - E(U_0) + \mu \Big)^{1/2} .
\end{align*}
In addition, by Cauchy-Schwarz, 
$$
 \Big| \int_{ |x| \leq \d } A^2 \p_x \phi \ dx \Big| \leq \sqrt{2 \d} 
\Big( \sup_{|x| \leq \d} A \Big) \Big( \mathscr{K}(U) - E(U_0) + \mu \Big)^{1/2} .
$$
Consequently,
\begin{align}
\label{enphase}
 r_0^2 \Big| \phi (+\d) - \phi(-\d) - \pi \ {\rm mod} \, 2 \pi \Big| 
\leq & \, 
 \Big| \int_{|x| \leq \d } ( A^2 - r_0^2 ) \p_x \phi \ dx - r_0^2 \pi \, {\rm mod} \, 2 \pi r_0^2 \Big| 
\nonumber \\ 
& \, + \sqrt{2 \d} 
\Big( \sup_{|x| \leq \d} A \Big) \Big( \mathscr{K}(U) - E(U_0) + \mu \Big)^{1/2} \nonumber \\
\leq & \, 
\Big[ \frac{K}{\inf_{|x| \geq \d} A } + \sqrt{2 \d} \Big( \sup_{|x| \leq \d} A \Big) \Big] 
\Big( \mathscr{K}(U) - E(U_0) + \mu \Big)^{1/2} .
\end{align}
From our choice $ \d = ( \mathscr{K}(U) - E(U_0) + \mu )^{1/4} \ll 1 $ and since 
$ |\!| \Theta |\!|_{L^\ii(\R)} \leq K ( \mathscr{K}(U) - E(U_0) + \mu )^{1/2} = \BO(\d^2) $, 
we infer 
$ \inf_{|x| \geq \d} A \geq \inf_{|x| \geq \d} |U_0| - |\!| \Theta |\!|_{L^\ii(\R)} 
\geq \d / K $. Similarly, $ \sup_{|x| \leq \d} A \leq \sup_{|x| \leq \d} |U_0| 
+ |\!| \Theta |\!|_{L^\ii(\R)} \leq K \d $. Reporting this into \eqref{enphase} yields
$$
\Big| \phi (+\d) - \phi(-\d) - \pi \ {\rm mod} \, 2 \pi \Big| \leq K \d .
$$
We now freeze the gauge invariance by imposing $ \phi(+\d) = 0 $.  Note that then 
$ \phi(-\d) = \pi + \BO (\d) $. Furthermore, since $ \phi(+\d) = 0 $,
$$
 \int_{|x | \geq \d } ( \p_x \phi)^2 \ dx \leq \frac{K}{( \inf_{|x| \geq \d} A )^2} 
 \int_{|x | \geq \d } A^2 ( \p_x \phi)^2 \ dx \leq 
\frac{K}{ \d^2 } \d^4 = K \d^2 ,
$$
which implies, for $ x \geq \d $,
$$
 | 1 - \cos(\phi(x)) | \leq | 1 - \cos(\phi(0)) |  
+ \Big| \int_\d^x \p_x \phi \sin (\phi) \Big| \leq K \d \sqrt{x} 
$$
and similarly, since $ \cos( \phi(-\d) ) = \cos( \pi + \BO (\d) ) = -1 + \BO(\d^2) $, 
for $ x \leq - \d $,
$$
 | 1 + \cos(\phi(x)) | \leq K \d \sqrt{|x|} . 
$$
We turn back to \eqref{espadon} and infer
\begin{align*}
\n \ex^{i\phi} \p_x |U_0| - \p_x U_0 \n_{L^2(\R)}^2 
\leq & \, K \d + 2 \int_\d^{+\ii} ( \p_x U_0 )^2 \Big( 1 - \cos(\phi) \Big) \ dx 
+ 2 \int_{-\ii}^{-\d} ( \p_x U_0 )^2  \Big( 1 + \cos(\phi) \Big) \ dx 
 \\ 
\leq & \, K \d + K \d \int_\R ( \p_x U_0 )^2 \sqrt{|x|} \ dx = K \d .
\end{align*}
Inserting these estimates in \eqref{morue}, it follows that
\begin{align*}
d_{\mathcal{Z}} ( U, U_0 ) \leq K \d .
\end{align*}

We now turn to the case $ \mu = 0 $. Without loss of generality, we may assume 
that $ |U| > 0 $ in $ ( - \ii , 0 ) $ (since $|U| \to r_0 > 0 $ at $ \pm \ii $), and 
let $\ell \geq 0 $ be such that $ |U| (\ell) = 0 $ and $ |U| > 0 $ in $ ( \ell , + \ii ) $. 
We first estimate $ \ell $ by writing that
$$
 | U_0 | (\ell) = |U|(\ell) + \Theta (\ell) = \Theta (\ell) \leq 
\n \Theta \n_{L^\ii(\R)} \leq K \Big( \mathscr{K} (U) - E(U_0) + \mu \Big)^{1/2} 
= K \d^2 ,
$$
thus
$$
 \ell \leq K \d^2 .
$$
Moreover, we have two local liftings $ U = A \ex^{i \phi_+ } $ in $ [\ell , +\ii )$ 
and $ U = A \ex^{i \phi_- } $ in $ (-\ii , 0 )$. Going back to \eqref{morue}, we then 
deduce
\begin{align*} 
d_{\mathcal{Z}} ( U, U_0 ) \leq
 \n \ex^{i\phi_-} \p_x |U_0| - \p_x U_0 \n_{L^2(-\ii, 0 )} 
+ \n \ex^{i\phi_+} \p_x |U_0| - \p_x U_0 \n_{L^2(\ell , +\ii )} 
+ K \d 
+ K \Big[ \mathscr{K}(U) - E(U_0) \Big]^{1/2} .
\end{align*}
Arguing as for the case $ \mu > 0$, we obtain $|U| = A \geq \d / K $ in 
$ [\ell + \d , +\ii )$ and in $ (-\ii , -\d ) $. By definition of $ \mathfrak{P} $, 
we have
$$
 \mathfrak{P} (U) = \int_{-\d }^{\ell +\d } \langle i U | \p_x U \rangle 
+ \int_{\ell +\d }^{+\ii} ( A^2 - r_0^2 ) \p_x \phi_+ \ dx + r_0^2 \phi_+ (\ell + \d ) 
+ \int_{-\ii}^{ -\d } ( A^2 - r_0^2 ) \p_x \phi_- \ dx - r_0^2 \phi_+ ( - \d ) 
$$
in $ \R / ( 2 \pi r_0^2 \Z ) $, hence the same arguments as in the case $\mu > 0 $ 
provide
$$
 \Big| \phi_+ (\ell + \d ) - \phi_+ ( - \d ) - \pi \ {\rm mod} \, 2 \pi \Big| \leq K \d ,
$$
since the integral $ \int_{-\d }^{\ell +\d } \langle i U | \p_x U \rangle $ 
is bounded by $ K \sqrt{\d} $ by Cauchy-Schwarz. Imposing $ \phi_+ (\ell + \d ) $ 
for the gauge invariance, we infer $ 1 - \cos(\phi_+ (\ell + \d ) ) = 0 $ and 
$ \phi_+ ( - \d ) = \pi + \BO (\sqrt{\d}) \ {\rm mod} \, 2 \pi $, hence 
$ 1 + \cos(\phi_- ( - \d ) ) = \BO( \d) $. Therefore, we conclude as before that
$$
d_{\mathcal{Z}} ( U, U_0 ) \leq K \d ,
$$
which finishes the proof of the Proposition.\carre \\

In order to prove Theorem \ref{stabmini}, we use Proposition \ref{supercool}, 
which provides
$$
 \mathscr{K}(U) \geq E(U_0) + \frac{1}{K} \Big( \inf_\R |U| \Big)^2 ,
$$
thus
$$
\mu = \inf_\R |U| \leq K \sqrt{\mathscr{K}(U) - E(U_0)} .
$$
Inserting this bound in Proposition \ref{minokinkameliore} then gives
$$
d_{\mathcal{Z}} ( U, U_0 ) \leq K \Big[ \mathscr{K}(U) - E(U_0) + K \sqrt{\mathscr{K}(U) - E(U_0)} \Big]^{1/4}
\leq K \sqrt[8]{\mathscr{K}(U) - E(U_0)} ,
$$
and the proof is complete.

%%%%%%%%%%%%%%%%%%%%%%%%%%%%%%%%%%%%%%%%%%%%%%%%%%%%%%%%%%%%%%%%%%%%%%%%%%%%%%
\section{About the stability analysis for the sonic waves $\bs{( c = \cs )} $} 
\label{sexsonic}

We have left aside in our study the case of the sonic waves ($c = \cs $), but would 
like to say a few words on the difficulties associated with this critical case.

We note that if there exists a sonic nontrivial travelling wave, it does 
not vanish, hence we may use the hydrodynamical formulation \eqref{hydrohamilto} 
of (NLS) as in \cite{Lin}. The point is that the Sturm-Liouville 
operator (see section 4 in \cite{Lin})
$$
L \equiv - \frac{\p}{\p x} \Big( \frac{1}{4(r_0^2 - \eta)} \frac{\p}{\p x} \Big) 
+ q (x) ,
$$
with
$$
q(x) \equiv \frac{ (\p_x \eta)^2}{4(r_0^2 - \eta)^3} 
- \frac{\p}{\p x} \Big( \frac{ \p_x \eta }{4(r_0^2 - \eta)^2} \Big) 
- \frac12 f '(r_0^2 - \eta ) - \frac{c^2 r_0^4}{4(r_0^2 - \eta)^3 }
$$
has, by Weyl's theorem, essential spectrum $ \s_{\rm ess}(L) = [0 , +\ii ) $ 
when $ c = \cs $. Indeed, we know from Proposition \ref{decroiss} 
that $ \eta_{\cs} $ and its derivatives tend to zero at infinity, 
hence, as $ x \to \pm \ii $, 
$ q(x) \to - \frac12 f '(r_0^2 ) - \frac{c^2 }{4 r_0^2} = 0 $ 
since $ c^2 = \cs^2 = - 2 r_0^2  f '(r_0^2 ) $. Therefore, there 
does not exist $ \d > 0 $ such that 
$ \langle H p , p \rangle \geq \d |\!| p |\!|^2 $ for any $ p $ 
orthogonal to the subspace spanned by the negative and the zero eigenvalue, 
and thus the Grillakis-Shatah-Strauss theory does not apply.

In the case $ \ds{ \frac{dP}{dc}_{|c=\cs} < 0 } $, where it is natural 
to expect stability, a natural thing would be to try to work with the 
functional
$$
 \mathscr{L} (\psi ) \equiv E (\psi ) - \cs P (\psi ) 
+ \frac{M}{2} \Big( P(\psi ) - P (U_{\cs}) \Big)^2 
$$
and to follow the lines of the proof of Theorem \ref{statiomini}. Indeed, 
the spectral analysis shall not give positive definiteness of 
the Hessian due to presence of essential spectrum down to $0$. 
Therefore, we may study $ \mathscr{L} $ at fixed $ \mu = \inf_\R |\psi| $ 
close to $ \inf_\R |U_{\cs}| $. 
When $ 0 < c_* < \cs $ and $ \ds{ \frac{dP}{dc}_{|c=c_*} \not = 0 } $, the 
infimum of $|U_c|$ contains a neighborhood of $ \inf_\R |U_{c_*}| $ for 
$c$ close to $c_*$. For $ c_* = \cs $, this is no longer the case: we have 
only a one sided neighborhood of $ \inf_\R |U_{\cs}| $. It is plausible that 
the study for $\mu$ in this one sided neighborhood of $ \inf_\R |U_{\cs}| $ can 
be done as in the proof of Theorem \ref{statiomini}, but for the remaining 
values of $\mu$, we have to find a sharp ansatz, which is not very easy to find.

% However, our approach breaks down for establishing Step 7, which 
% works only when $ \inf_\R |U_{\cs}| $ is sufficiently small. 
% More generally, the way we prove Theorem \ref{statiomini} requires 
% some adaptation for Step 7 to be applied to the case of Theorem 
% \ref{Liapoupou}, if one wants to avoid spectral arguments.

Furthermore, for the linear instability which is expected if $ \ds{ \frac{dP}{dc}_{|c=\cs} > 0 } $, 
let us mention the following point. For the eigenvalue problem studied in \cite{BeGa}, 
the characteristic equation for the constant coefficient limit at infinity, namely
$$
 r^4 - (\cs^2 - c_*^2) r^2 - 2 c_* \lambda r + \lambda^2 = 0 .
$$
becomes, when $ c_* = \cs $,
\be
\label{karakcs}
 r^4 - 2 \cs \lambda r + \lambda^2 = 0 .
\ee
The behaviour of the roots for small $\lambda $ is then different from 
the case $ 0 < c_* < \cs $. Indeed, there exists a root $ \sim \lambda / ( 2 \cs ) $ 
for $ \lambda \to 0 $, and for the three other roots, we use the variable 
$r = \sqrt[3]{\lambda} z $, which transforms $ r^4 - 2 \cs \lambda r + \lambda^2 = 0 $ 
into $ z^4 - 2 \cs z + \lambda^{2/3} = 0 $. This last equation has, for $ \lambda \to 0 $, 
three roots $ \sim j^k \sqrt[3]{2\cs} $, where $ j = \ex^{ 2 i \pi /3 } $ and $k=0$, $1$, $2$. 
In particular, \eqref{karakcs} has three roots $ \sim j^k \sqrt[3]{2\cs \lambda } $, 
$k=0$, $1$, $2$. The value $ \lambda = 0 $ is then a branching point, and 
we shall have a smooth problem not in $ \lambda $ but in $ \sqrt[3]{ \lambda } $. 
Since analyticity is not necessary for our purpose,  we may define an Evans function 
$ \tilde{\mathbb{D}} $ in $\R_+ $, smooth, and such that, for $ \lambda> 0 $, 
$ \tilde{\mathbb{D}} (\sqrt[3]{ \lambda } ) = 0 $ if and only if $ \lambda $ is 
an unstable eigenvalue for \eqref{rouge}. another difficulty comes from the fact 
that it will be difficult to find an analytic extension of the Evans function $ \tilde{\mathbb{D}} $ 
near $0$ since, by Proposition \ref{decroiss}, for $ c_* = \cs $, $ u_* $ and $\eta_*$ decay 
only at an algebraic rate and not an exponential rate. Consequently, we can not use the 
Gap Lemma of \cite{GZ} and \cite{KS}. Finally, as a straightforward computation shows, 
the stable and unstable subspaces for the eigenvalue problem are transverse 
for $ \lambda > 0 $ but their continuous extensions at $ \lambda = 0 $ have a 
nontrivial intersection. Therefore, both stability and instability requires some 
further analysis, and the situation is then much more delicate than the one studied in 
subsection \ref{sexins}.

\appendix

\section*{Appendix A: construction of a Liapounov functional in the stable 
case in the Grillakis-Shatah-Strauss framework}
%%%%%%%%%%%%%%%%%%%%%%%%%%%%%%%%%%%%%%%%%%%%%%%%%%%%%%%%%%%%%%%%%%%%%%%%%%
\setcounter{equation}{0} \setcounter{subsection}{0}
\setcounter{lemma}{0}\setcounter{proposition}{0} \setcounter{theorem}{0}
\setcounter{corollary}{0} \setcounter{remark}{0} \setcounter{definition}{0}
%%%%%%%%%%%%%%%%%%%%%%%%%%%%%%%%%%%%%%%%%%%%%%%%%%%%%%%%%%%%%%%%%%%%%%%%%%

We work with the notations of \cite{GSS}, and recall them briefly. We 
consider a Hamiltonian equation in a real Hilbert space $\BX$, with scalar 
product $ ( \cdot , \cdot )_\BX $, under the form
\be
\tag{$\BH$}
\frac{\p u}{\p t} = J E'(u) ,
\ee
where $ J : \BX^* \to \BX $ is a closed linear operator with dense domain and 
skew-symmetric. Assume that $T$ is a $\BC_0$-group of unitary operators 
in $\BX$ generated by $ T'(0) $, which is skew-adjoint and with dense domain, 
and that $E$ is invariant by $T$, that is $E(T(s) u) = E(u) $ for any 
$ s \in \R $, $ u \in \BX $. Assume moreover that $ T(s) J = J T(-s)^* $ for 
any $ s \in \R $ and that there exists $ B : \BX \to \BX^* $ linear and bounded 
such that $ B^* = B $ and $JB$ is an extension of $ T'(0) $. We then set
$$
 Q(u) \equiv \frac12 \langle B u , u \rangle_{ \BX^* , \BX } .
$$
The basic assumptions of \cite{GSS} are the following ones.\\

\noindent {\it Assumption 1 (existence of solutions):} For any $ r > 0 $ 
there exists $ t_* > 0 $, depending only on $r$, such that for any $ u^\ini \in \BX $, 
there exists a $ u \in \BC( (- t_*, t_* ) , \BX ) $ with $ u (0) = u^\ini $ solution 
of ($\BH$) in the sense that for any $ \vp \in D ( J ) \subset \BX^* $,
$$
 \frac{d}{dt} \langle u(t) , \vp \rangle_{ \BX^* , \BX } = 
- \langle E'(u(t)) , J \vp \rangle_{ \BX^* , \BX } \quad \quad \quad 
{\rm in} \quad \BD'( (- t_*, t_* ) ) ,
$$
and verifying $ E(u(t))= E(u^\ini) $ and $ Q(u(t))= Q(u^\ini) $ for 
$ t \in (- t_*, t_* ) $.\\

\noindent {\it Assumption 2 (existence of ``bound states''):} There exists 
an interval $ \OO \subset \R $, not reduced to a singleton, and a mapping 
$ \OO \ni \oo \mapsto \phi_\oo \in \BX $ of class $ \BC^1 $ such that, for any 
$ \oo \in \OO $,
$$ 
E'(\phi_\oo) = \oo Q'(\phi_\oo) , 
\quad \quad 
\phi_\oo \in D(T'(0)^3) \cap D (JI T'(0)^2 ) ,
\quad \quad 
T'(0) \phi_\oo \not = 0 .
$$ 

\noindent {\it Assumption 3 (spectral decomposition):} For each $ \oo \in \OO $, 
the operator $ H_\oo \equiv E''(\phi_\oo) - \oo Q'' (\phi_\oo) : \BX \to \BX^* $ 
has its kernel spanned by $ T'(0) \phi_\oo $, has one negative simple 
eigenvalue and the rest of its spectrum is positive and bounded away from zero.\\

Under assumption 2, we consider some $ \oo_* \in \OO $ and the associated bound state 
$ \phi_{\oo_*} $, and then define, for $M > 0 $, the functional
$$
\mathscr{L}_{\oo_*} ( u ) \equiv E(u) - \oo_* Q(u) + \frac{M}{2} \Big( Q(u) - Q( \phi_{\oo_*} ) \Big)^2 .
$$
It is clear that $ \phi_{\oo_*} $ is a critical point of $ \mathscr{L}_{\oo_*} $: 
$ \mathscr{L}' ( \phi_{\oo_*} ) = E'(\phi_{\oo_*}) - \oo_* Q'(\phi_{\oo_*}) = 0 $. 
We denote by
$$
\Lambda \equiv \mathscr{L}_{\oo_*}'' ( \phi_{\oo_*} ) = 
H_{\oo_*} + M \langle Q'(\phi_{\oo_*}) , \cdot \rangle_{ \BX^* , \BX } Q'(\phi_{\oo_*}) 
$$
its second derivative, which is a self-adjoint operator. The main 
result of this appendix is the following.\\

%%%%%%%%%%%%%%%
\noindent {\bf Theorem A} {\it We make assumptions 2 and 3 and suppose 
moreover that the operator $ \langle Q'(\phi_{\oo_*}) , \cdot \rangle_{ \BX^* , \BX } Q'(\phi_{\oo_*}) $ 
is a compact perturbation of $ H_{\oo_*} $. If 
$ \ds{ \frac{dQ(\phi_{\oo})}{d \oo}_{|\oo = \oo_*} < 0 } $ and
$$ M > \frac{1}{ \ds{ - \frac{dQ(\phi_{\oo})}{d\oo}_{|\oo=\oo_*}}} , $$ 
there exists $ \d > 0 $ such that
$$
\forall v \in X, \quad {\it s.t.} \quad ( v , T'(0) \phi_{\oo_*} )_\BX = 0, 
\quad \quad \quad 
\langle \Lambda v , v \rangle \geq \d \n v \n^2 .
$$
In particular, for any $ u \in X $ with 
$ \inf_{ s \in \R } \n u - T(s) \phi_{\oo_*} \n^2 \leq \epsilon $, we have
$$
 \inf_{ s \in \R } \n u - T(s) \phi_{\oo_*} \n^2 \leq \frac{2}{\d} 
\Big( \mathscr{L} (u) - \mathscr{L} (\phi_{\oo_*}) \Big) .
$$
Therefore, when assumption 1 is moreover verified, the (global) solution 
$ u (t) $ to {\rm ($\BH$)} with initial datum $ u^\ini $ verifies
$$
 \sup_{t \in \R} \inf_{ s \in \R } \n u(t) - T(s) \phi_{\oo_*} \n^2 \leq \frac{2}{\d} 
\Big( \mathscr{L} (u^\ini) - \mathscr{L} (\phi_{\oo_*}) \Big) 
\leq K \n u^\ini - \phi_{\oo_*} \n^2 ,
$$
provided the right-hand side is sufficiently small. }
%%%%%%%%%%%%%

\bigskip

We point out that the condition that the operator 
$ \langle Q'(\phi_{\oo_*}) , \cdot \rangle_{ \BX^* , \BX } Q'(\phi_{\oo_*}) $ 
is a compact perturbation of $ H_{\oo_*} $ is not very restrictive, since in 
many cases coming from PDE's, it involves less derivatives than $ H_{\oo_*}$ 
and $ Q'(\phi_{\oo_*}) $ tends to zero at spatial infinity.\\

This type of Liapounov functional has been used by I. Barashenkov 
in \cite{Ba} to prove that the travelling waves of (NLS) in dimension 
one are stable when $ \ds{ \frac{dP}{dc} < 0 } $. The proof follows 
basically the one in \cite{Ba}, but some points have to be clarified. 
The interest of this type of Liapounov functional is that the saddle 
point $ \phi_{\oo_*} $ is now a non degenerate local minimum for 
$ \mathscr{L}_{\oo_*} $. This is a great advantage for numerical simulation 
of the ``bound states'', since a gradient flow method on $ \mathscr{L}_{\oo_*} $ 
can be used. This approach has been used, with a very similar functional, 
in \cite{PaSp} by N. Papanicolaou and P. Spathis for the numerical 
simulation of the travelling waves for a planar ferromagnets model. 
In the same spirit, in \cite{CS}, we also use a gradient flow method 
on this type of functional for the numerical simulation of the 
travelling waves for (NLS) in two dimensions.\\

\noindent {\it Proof of Theorem A.} Recall that the spectrum of $ H_{\oo_*} $ 
is, by assumption 3, such that $ - \lambda_*^2 \in \s (H_{\oo_*}) $, $ 0 \in \s (H_{\oo_*}) $ 
and $ \s (H_{\oo_*}) \setminus \{ - \lambda_*^2 , 0 \} \subset [ \d,+\ii ) $ for some $ \d > 0$. 
Since we assume that $ \langle Q'(\phi_{\oo_*}) , \cdot \rangle Q'(\phi_{\oo_*}) $ 
is a compact perturbation of $ H_{\oo_*} $, the essential spectrum of $ \Lambda $ 
is the same as the one of $ H_{\oo_*} $, hence is included in $ [ \d , + \ii ) $. 
Furthermore, $ 0 \in \s ( H_{\oo_*} ) $ and $ \ker (H_{\oo_*}) = \R T'(0) \phi_{\oo_*} $ 
by assumption 3. Since $ Q'(\phi_{\oo_*}) = B \phi_{\oo_*} $ and $ J B $ is an extension 
of $T'(0)$, we have $ \langle Q'(\phi_{\oo_*}) , T'(0) \phi_{\oo_*} \rangle_{ \BX^* , \BX } 
= \langle B \phi_{\oo_*} , JB \phi_{\oo_*} \rangle_{ \BX^* , \BX } = 0 $, hence 
$ \Lambda (T'(0) \phi_{\oo_*}) = 0 $. Noticing that 
$ \langle Q'(\phi_{\oo_*}) , \cdot \rangle_{ \BX^* , \BX } Q'(\phi_{\oo_*}) $ is a 
nonnegative operator, we infer that $ \ker ( \Lambda ) = \ker (H_{\oo_*}) = \R T'(0) \phi_{\oo_*} $ 
is one-dimensional. Therefore, it suffices to show that $ \Lambda $ has no 
eigenvalues in $ ( - \ii , 0 ) $. As we have seen that 
$ \langle Q'(\phi_{\oo_*}) , \cdot \rangle_{ X^* , X } Q'(\phi_{\oo_*}) $ is a nonnegative operator, 
we deduce that $ \s( \Lambda ) \subset [- \lambda_*^2, + \ii ) $. Let 
us first show that $ - \lambda_*^2 \not \in \s( \Lambda ) $ by contradiction. 
If $ - \lambda_*^2 $ is an eigenvalue of $ \Lambda $, then there exists 
$ v \in X $, $ v \not = 0 $, such that 
$ 0 = (\Lambda + \lambda_*^2 )v = ( H + \lambda_*^2 )v 
+ M \langle Q'(\phi_{\oo_*}) , v \rangle_{ \BX^* , \BX } Q'(\phi_{\oo_*}) $. 
Taking the duality product with $v$ yields 
$ 0 = \langle ( H_{\oo_*} + \lambda_*^2 )v , v \rangle_{ \BX^* ,\B X } 
+ M \langle Q'(\phi_{\oo_*}) , v \rangle_{ \BX^* , \BX }^2 $. 
Since the two terms in the sum are nonnegative, this implies 
$\langle Q'(\phi_{\oo_*}) , v \rangle_{ \BX^* ,\BX } = 0 $ and 
$ \langle ( H_{\oo_*} + \lambda_*^2 )v , v \rangle_{ \BX^* , \BX } = 0 $, 
which in turn implies $ v \in \ker ( H_{\oo_*} + \lambda_*^2 ) = \R \chi $ (here, 
$ \chi $ is a negative eigenvector of $H_{\oo_*}$ for the eigenvalue $ - \lambda_*^2 < 0 $. 
As a consequence, we must have $ \langle Q'(\phi_{\oo_*}) , \chi \rangle_{ \BX^* , \BX } = 0 $. On 
the other hand, differentiating the equality $E'(\phi_\oo) - \oo Q' (\phi_\oo) = 0 $ 
at $\oo = \oo_*$ yields $ Q'(\phi_{\oo_*}) = H_{\oo_*} \phi' $, where 
$ \phi' \equiv \frac{d \phi}{d\oo}_{|\oo = \oo_*} $. Thus we must have 
$ 0 = \langle H_{\oo_*} \phi' , \chi \rangle_{ \BX^* , \BX } 
= \langle H_{\oo_*} \chi , \phi' \rangle_{ \BX^* , \BX } 
= - \lambda_*^2 ( \chi , \phi' ) $. Therefore, $ \phi' $ is orthogonal to 
$ \chi $ and this gives $ \langle H_{\oo_*} \phi' , \phi' \rangle_{ \BX^* , \BX } \geq 0 $. 
However, this is not possible if $ \ds{ \frac{dQ(\phi_\oo)}{d \oo}_{|\oo = \oo_*} < 0 } $, 
since we have $ \ds{\frac{dQ(\phi_\oo)}{d \oo}_{|\oo = \oo_*} 
= - \langle H_{\oo_*} \phi' , \phi' \rangle_{ \BX^* , \BX } }$. 
As a consequence, if $ \lambda $ is a negative element of the spectrum of 
$ \Lambda $, then $ -\lambda_*^2 < \lambda < 0 $ and $\lambda$ is an 
eigenvalue: there exists $ v \in X$ such that $ v \not = 0 $ and
$$
\lambda v = \Lambda v = H_{\oo_*} v + M \langle Q'(\phi_{\oo_*}) , v \rangle_{ \BX^* , \BX } Q'(\phi_{\oo_*}) . 
$$
Since $ -\lambda_*^2 < \lambda < 0 $, we then infer
\be
\tag{A.1}
v = - M \langle Q'(\phi_{\oo_*}) , v \rangle_{ \BX^* , \BX } ( H_{\oo_*} - \lambda )^{-1} Q'(\phi_{\oo_*}) .
\ee
Since $ v \not= 0 $, we can not have $ \langle Q'(\phi_{\oo_*}) , v \rangle_{ \BX^* , \BX } = 0 $. 
Then, taking the scalar product of (A.1) with $ \mathbb{I}^{-1} Q'(\phi_{\oo_*}) $ 
(here, $ \mathbb{I} :\BX \to \BX^* $ is the usual Riesz isomorphism) gives
$$
 g (\lambda ) = 0, \quad \quad \quad {\rm where} \quad \quad \quad 
 g(t) \equiv 1 + M ( ( H_{\oo_*} - t )^{-1} Q'(\phi_{\oo_*}) , \mathbb{I}^{-1} Q'(\phi_{\oo_*}) )_\BX , 
\quad - \lambda_*^2 < t < 0 .
$$
It is clear that $g$ is smooth in $ ( - \lambda_*^2 , 0 ) $ and that 
$$
g'(t) = M ( ( H_{\oo_*} - t )^{-2} Q'(\phi_{\oo_*}) , \mathbb{I}^{-1} Q'(\phi_{\oo_*}) )_\BX 
= M \n ( H_{\oo_*} - t )^{-1} Q'(\phi_{\oo_*}) \n_\BX^2 > 0 .
$$
We now study the limit of $g$ at $ 0^- $. Let us recall that 
$ H_{\oo_*} \phi' = Q'(\phi_{\oo_*}) $ and that we have already seen that 
$ \langle Q'(\phi_{\oo_*}) , T'(0) \phi_{\oo_*} \rangle_{ \BX^* , \BX } = 0 $, {\it i.e.} 
$ \mathbb{I}^{-1} Q'(\phi_{\oo_*}) $ is orthogonal to $\ker (H_{\oo_*})$. Therefore, 
as $t \to 0^- $,
$$
( ( H_{\oo_*} - t )^{-1} Q'(\phi_{\oo_*}) , \mathbb{I}^{-1} Q'(\phi_{\oo_*}) ) \to 
( \phi' , \mathbb{I}^{-1} Q'(\phi_{\oo_*}) ) 
= \langle Q'(\phi_{\oo_*}) , \phi' \rangle_{ \BX^* , \BX } = \frac{d Q (\phi_\oo) }{d \oo}_{| \oo = \oo_*} 
$$ 
and it thus comes
$$
g(t) \to 1 + M \frac{d Q (\phi_\oo) }{d \oo}_{| \oo = \oo_*} 
\quad \quad \quad {\rm as} \quad t \to 0^- .
$$
Since $ \ds{\frac{d Q (\phi_\oo) }{d \oo}_{| \oo = \oo_*} < 0 } $ 
by hypothesis, it follows that if 
$ M > - 1 / ( \frac{d Q (\phi_\oo) }{d \oo}_{| \oo = \oo_*}) > 0 $, 
the function $g$ increases in $ (- \lambda_*^2 , 0 ) $ and tends to 
some negative limit at $ 0^- $. In particular, $g$ is negative, hence 
we can not have $ g( \lambda ) = 0 $ with $ \lambda \in (- \lambda_*^2 , 0 ) $. 
We have therefore shown that the spectrum of $ \Lambda $ consists in a 
simple eigenvalue $0$ with eigenspace spanned by $ T'(0) \phi_{\oo_*} $ and 
the rest of the spectrum is positive and bounded away from $0$. This 
concludes the proof. \carre \\

We would like to point out that in the proof of \cite{Ba}, the fact that 
$ - \lambda_*^2 \not \in \s( \Lambda ) $ was not shown, the kernel of 
$ \Lambda $ was not studied and the essential spectrum was not considered. 
Moreover, the functional spaces are not given, hence we do not know for which 
perturbations stability holds.

\appendix

\section*{Appendix B: from linear to nonlinear instability}
%%%%%%%%%%%%%%%%%%%%%%%%%%%%%%%%%%%%%%%%%%%%%%%%%%%%%%%%%%%%%%%%%%%%%%%%%%
\setcounter{equation}{0} \setcounter{subsection}{0}
\setcounter{lemma}{0}\setcounter{proposition}{0} \setcounter{theorem}{0}
\setcounter{corollary}{0} \setcounter{remark}{0} \setcounter{definition}{0}
%%%%%%%%%%%%%%%%%%%%%%%%%%%%%%%%%%%%%%%%%%%%%%%%%%%%%%%%%%%%%%%%%%%%%%%%%%

We still consider in this appendix an abstract Hamiltonian equation in the framework of \cite{GSS} 
\be
\tag{$\BH$}
\frac{\p u}{\p t} = J E'(u)   
\ee 
on the real Hilbert space $\BX$, with scalar product $ ( \cdot , \cdot )_\BX $. 
Here $ E : \BX \to \R $ is of class $ \BC^2 $ and $ J : \BX^* \to \BX $ is a closed 
linear operator with dense domain and skew-symmetric in the sense that 
$ ( u , Jw )_{\BX} = - \langle w , Ju \rangle_{\BX^*,\BX} $ for 
$ u \in \BX $, $ w \in \BX^* $. 

We assume that there exists a $\BC_0$-group $T$ of unitary operators in 
$\BX$ generated by $ T'(0) $, which is skew-adjoint and with dense domain, 
and that $E$ is invariant by $T$, that is $E(T(\oo) u) = E(u) $ for any 
$ \oo \in \R $, $ u \in \BX $. Assume moreover that $ T(\oo) J = J T(-\oo)^* $ 
for any $ \oo \in \R $ and that there exists $ B : \BX \to \BX^* $ linear 
and bounded such that $ B^* = B $ and $JB$ is an extension of $ T'(0) $. 
We then set
$$
 Q(u) \equiv \frac12 \langle B u , u \rangle_{ \BX^* , \BX } ,
$$
which is invariant by the flow ($\BH$) (see \cite{GSS}). By 
``bound state'', we mean a particular solution $ U $ of ($\BH$) of the 
form $ U(t) = T(\oo t ) \phi $ for some $ \oo \in \R $ and where $ \phi \in \BX $, 
$ \phi \not = 0 $. In other words, $ E'(\phi) = \oo Q'(\phi) $. 

There exists an open interval $ \OO \subset \R $, not reduced to a singleton, 
and a mapping $ \OO \ni \oo \mapsto \phi_\oo \in X $ of class $ \BC^1 $ such 
that, for any $ \oo \in \OO $,
$$ 
E'(\phi_\oo) = \oo Q'(\phi_\oo) , 
\quad \quad 
\phi_\oo \in D(T'(0)^3) \cap D (JI T'(0)^2 ) ,
\quad \quad 
T'(0) \phi_\oo \not = 0 .
$$

The solution $ U(t) = T(\oo t ) \phi $ is said to be stable in $\BX$ if 
for any $ \e > 0 $, there exists $ \d > 0 $ such that any solution to ($\BH$) 
with initial datum $ u^\ini \in B_\BX ( \phi , \d ) $ is global in time 
and remains in $ B_\BX (\phi , \e ) $ for $ t \geq 0 $. Otherwise, it is said unstable. 
This supposes some knowledge on the Cauchy problem for ($\BH$) (at least existence 
of solutions). If we are given some Banach 
space $ \BY \supset \BX $ with continuous imbedding $ \BX \hookrightarrow \BY $, 
we may also say that the solution $ U(t) = T(\oo t ) \phi $ is said to be stable 
from $\BX$ to $\BY $ if for any $ \e > 0 $, there exists $ \d > 0 $ such that any 
solution to ($\BH$) with initial datum $ u^\ini \in B_\BX ( \phi , \d ) $ remains in 
$ B_\BY (\phi , \e ) $ for $ t \geq 0 $. Clearly, a solution stable in $\BX$ 
is precisely a solution able from $\BX$ to $\BX$, and is also stable from 
$\BX$ to $\BY $, hence instability from $\BX$ to $\BY $ is a stronger statement 
that instability in $\BX$. 

In our framework, the notion of 
orbital stability is more relevant. Let us consider $ \mathbb{G} $ a group and 
$ \mathbb{T} :\R \times \mathbb{G} \to \mathscr{G}\mathscr{L}_c (\BX ) $ 
a unitary representation of $ \R \times \mathbb{G} $ on $ \BX $, extending 
$ T : \R \to \BX $ and leaving $E$ and $Q$ invariant. Then, 
$ U(t) = T(\oo t ) \phi $ is said to be orbitally stable in $\BX$ (for the 
group $ \mathbb{G} $) if for any $ \e > 0 $, there exists $ \d > 0 $ 
such that any solution to ($\BH$) with initial datum 
$ u^\ini \in B ( \phi , \d ) $ is global in time and remains in 
$ \cup_{ ( \oo, g ) \in \R \times \mathbb{G} } B (\mathbb{T} (\oo, g) \phi , \e)$ 
for $ t \geq 0 $. We may also define orbital instability from $\BX$ to 
$ \BY \supset \BX $ in a natural way.

In \cite{GSS, GSS2}, a general framework for the stability analysis for the 
``bound state'' has been given. In particular, the nonlinear orbital instability 
in proved in \cite{GSS} through the construction of a Liapounov type functional. 
However, this method does not give a clear understanding neither of how 
we get farther from the ``bound state'', nor on which timescale it occurs.

The need for allowing an additional group of invariances $ \mathbb{G} $ 
can be seen in the case of bound state solutions, that is 
$ U(t) = \ex^{i\oo t} \phi_\oo $ to the Nonlinear Schr\"odinger equation
\be
\tag{NLS}
 i \p_t \Psi + \Delta \Psi + \Psi f(|\Psi|^2) = 0 ,
\ee
or the Nonlinear Klein-Gordon equation in $ \R^d $
\be
\tag{NLKG}
 \p^2_t \Psi = \Delta \Psi + \Psi f(|\Psi|^2) ,
\ee
since then, the invariance by translation in space must be taken 
into account in the definition of orbital stability, and we are 
in a case where $ \mathbb{G} = \R^d $ acts naturally by translation. 
The translations are taken into account in \cite{CaLi}. In 
\cite{GSS} and \cite{GSS2}, the notion of orbital stability is 
for $ \mathbb{G} $ trivial. It is clear from the definition 
that orbital stability for $ \mathbb{G} = \{ 0 \} $ implies 
orbital stability for arbitrary $ \mathbb{G} $. For the 
instability in the Nonlinear Schr\"odinger equation or 
the Nonlinear Klein-Gordon equation, \cite{GSS} and \cite{ShaStrainsta} 
work with radial $H^1$ functions. The fact that this also 
implies the orbital instability with the action of 
$ \mathbb{G} = \R^d $ by translations follows immediately 
from the fact that for any $ \theta \in [ 0 , 2 \pi ] $ the manifold 
$ \mathfrak{M}_\theta \equiv \{ \ex^{i\theta} \phi ( \cdot -y) , 
\ y \in \R^d \} $ is orthogonal to $ H^1_{\rm rad}(\R^d) $.

For the stability analysis of a ``bound state'' 
$ U (t) = T(\oo_* t)( \phi_{\oo_*} ) $, it is natural to consider 
the linearization of ($\BH$) near $\phi$. More precisely, we 
linearize according to the ansatz $ u (t) = T(\oo_* t)( \phi_{\oo_*} + v(t) ) $, 
so that the ``bound state'' becomes stationary. The linearized problem then becomes
\be
\tag{$\BH_{\rm lin}$}
\frac{\p v}{\p t} = J (E''(\phi) - \oo Q''(\phi)) v = \BJ \BL v ,
\ee
where,  $ \mathbb{I} : \BX \to \BX^* $ denoting the Riesz isomorphism, 
$ \BJ \equiv J \mathbb{I}  : \BX \to \BX $ is skew-adjoint.

The purpose of this appendix is to give a general result, for Hamiltonian equations, 
showing that linear instability implies nonlinear (orbital) instability. By linear 
instability, we mean that the complexification of [$\BJ \BL $]$_\C$ has 
at least one eigenvalue in the right-half space $ \{ {\rm Re} > 0 \} $. The 
argument follows ideas from the works of F. Rousset and N. Tzvetkov \cite{RT2,RT1}.

Showing the existence of an unstable eigenvalue can be done through 
various techniques: \cite{GSS2} (in the framework of \cite{GSS} when 
$J$ is onto), \cite{Gri} (assuming a special structure of the Hamiltonian 
equation); for uses of the Vakhitov-Kolokolov function, see 
\cite{dBbubble}, \cite{dMG} or \cite{PeKe}. 
When $J$ is not onto, we quote \cite{Lopes}. 
For one dimensional partial differential equations, one may 
also use Evans' function (see the survey \cite{San}) as in 
\cite{PeWe}, \cite{GaZu}, \cite{KS}, \cite{Zu}. The paper 
\cite{Lin2} proposes another approach which allows to treat 
pseudo-differential equations, such as the BBM equation, the 
Benjamin-Ono equation, regularized Boussinesq equations, the 
Intermediate Long Wave equation...

In order to pass from linear to nonlinear instability, the following 
result is standard. We refer to the paper \cite{HPW} by D. Henry, 
J. Perez and W. Wreszinski. It can also be found in \cite{Gri} and \cite{ShaStra}.\\

%%%%%%%%%%%%%%%
\noindent {\bf Theorem B.1} (\cite{HPW}, \cite{Gri}, \cite{ShaStra})  
%\label{insta} 
{\it We assume that $ \BA $ generates a continuous semigroup on $X$ and 
that $ \s( \BA ) $ meets the right-half space $ \{ {\rm Re} > 0 \} $. We assume 
moreover that $ F : X \to X $ is locally Lipschitz continuous and verifies, for some 
$ \alpha > 0 $, $ |\!| F(v) |\!|_X = \BO( |\!| v |\!|_X^{1+\alpha} ) $ as $ v \to 0 $. 
Then, the solution $ \phi = 0 $ is unstable for the equation 
$ \p_t v = \BA v + F(v) $.}\\
%\end{theorem} 
%%%%%%%%%%%%%

In \cite{ShaStra}, it is claimed that an orbital instability result can also 
be established. Theorem B.1 shows nonlinear instability without 
assuming that the equation is Hamiltonian. However, if ($\BH_{\rm lin}$) 
can be solved using a semigroup, it does not give the growth of its norm. 
Moreover, it does not say that if the initial datum is in a most unstable 
direction, that is an eigendirection of $ \BA $ corresponding to an eigenvalue 
of maximal positive real part (plus the complex conjugate if necessary), 
then you can track the exponential growth of the solution. In particular, 
it does not explain the mechanism of instability and does not give any 
information on the timescale on which you see the instability. For 
instance, some strong instability results are shown by proving blow-up 
in finite time (see \cite{BeCa}), but the instability due to an 
exponentially growing mode holds on a much smaller timescale. 
We wish to provide here some results clarifying the 
instability mechanism by tracking the exponentially growing mode.

%%%%%%%%%%%%%%%%%%%%%%%%%%%%%%%%%%%%%%%%%%%%%%%%%%%%%%%%%%%%%%%%%%%%%%%%%%%%
\subsection*{B.1 A spectral mapping theorem for linearized Hamiltonian equations}

When we want to prove a nonlinear instability result from a linear instability one, 
one needs some information on the growth of the semigroup $ \BJ \BL $, 
when such a semigroup $ \ex^{ t \BJ \BL } $ exists, which we shall assume in 
this appendix. The growth estimate on $ \ex^{ t \BJ \BL } $ relies classically 
on the following spectral mapping result due to J. Pr\"uss \cite{Pruss}, which 
generalizes the work of L. Gearhart \cite{Gea}.\\

%%%%%%%%%%%%%%%
%\begin{theorem} 
\noindent {\bf Theorem B.2} (\cite{Pruss}) %\label{spectrum} 
{\it Let $X$ be a complex Hilbert space and $ \BA $ an unbounded operator on $X$ 
which generates a continuous semigroup $ \ex^{t \BA} $ on $X$. For 
$ t \in ( 0 ,+\ii )$, we have}
$$
\s(\ex^{t \BA}) \setminus \{ 0 \} = \Big\{ \ex^{\lambda t}, \ {\it either} \ 
\Big( \lambda + \frac{2i\pi }{t} \Z \Big) \cap \s( \BA ) \not = \emptyset , \ {\it or } \
\sup_{k \in \Z} \no \Big( \BA - \lambda - \frac{2i\pi k}{t} \Big)^{-1} \no_{\mathscr{L}_c(X)} = +\ii \Big\} .
$$
%\end{theorem} 
%%%%%%%%%%%%%%%

The following result is an immediate corollary.\\

%%%%%%%%%%%%%%%
%\begin{corollary} 
%\label{resol} 
\noindent {\bf Corollary B.1} 
{\it Let $X$ be a complex Hilbert space and $ \BA $ an unbounded operator on $X$ 
which generates a continuous semigroup $ \ex^{t \BA} $ on $X$. Assume that 
for any $ \gamma \in \R^* $, we have
$$
 \limsup_{|\tau | \to +\ii } 
\no \Big( \BA - \gamma - i \tau \Big)^{-1} \no_{\mathscr{L}_c(X)} < + \ii ,
$$
and that there exists $ \vartheta_0 \in [ 0, +\ii ) $ such that 
$ \s_{\rm ess} (\BA) = \{ i \vartheta , \ \vartheta \in \R , \ |\vartheta | \geq \vartheta_0 \} $. 
Then, for any $ t \in ( 0 ,+\ii )$, the spectral mapping holds: 
$ \s(\ex^{t \BA}) \setminus \{ 0 \} = \ex^{ t \s (\BA)} $. }\\
%\end{corollary}
%%%%%%%%%%%%%%%

\noindent {\it Proof.} Since 
$ \s_{\rm ess} (\BA) = \{ i \vartheta , \ \vartheta \in \R , \ |\vartheta | \geq \vartheta_0 \} $, 
we have $ \SS^1 \subset \ex^{t \s(\BA)} \subset \s(\ex^{t \BA}) $. If $ \lambda \in \C $ 
does not have modulus one, then note that when 
$( \lambda + \frac{2i\pi }{t} \Z ) \cap \s( \BA ) = \emptyset $, the supremum 
for $ k \in \Z $ in Theorem B.2 can be $ +\ii $ only when 
$ |k| \to +\ii $, and we conclude with our hypothesis. \carre \\
 
The fact that we exclude $0$ in the spectral mapping theorem just comes from 
the fact that we consider a semigroup and not a group. However, in most 
Hamiltonian PDE's, we have time reversibility and we have actually a continuous 
group and not only a semigroup. In most cases, we work with 
$ A : D(A) \subset Y \to Y $ where $Y$ is a real Hilbert space, thus for 
applying Theorem B.2 or Corollary B.1 we have to consider, as usual, the complexified 
operator $ A_\C : D(A_\C) \equiv D(A) \oplus i D(A) \subset Y_\C \equiv Y \oplus i Y \to Y_\C $ 
defined by $ A_\C( u + i v) = A u + i A v $.

It seems that the first time Theorem B.2 is used to prove a growth 
estimate on a semigroup was by T. Kapitula and B. Sandstede in \cite{KS}. 
Later, the work \cite{GJLS} by F. Gesztesy {\it et al.} also uses 
this result for bound states for (NLS). The bounds on the resolvent 
in \cite{KS} were proved using the particular structure of the 
linearized operator. In \cite{GJLS}, the computations are more involved 
and rely on suitable kernel estimates of some Hilbert-Schmidt operators. 
The same type of estimates have also been used in \cite{dMG}.

The main objective of this appendix is to provide a generalization of 
these results to a wide class of Hamiltonian equations. Indeed, the 
approaches in \cite{KS} and \cite{GJLS} seem specific to the problem. 
In addition, it is not clear whether the computations in \cite{GJLS} 
and \cite{dMG} can be extended to other types of equations. In 
particular, in \cite{C1d} and in the present paper, we have a situation 
similar to the one studied in \cite{dMG}, namely travelling wave 
solutions to a Nonlinear Schr\"odinger equation with nonzero 
condition at infinity, but for nonzero propagation speeds, the travelling 
wave is not real-valued (as it is the case in \cite{dMG} for stationary waves 
or for bound state solutions), and the bloc diagonal structure of the linearized 
Hamiltonian disappears. An additional difficulty is that in \cite{C1d} and the 
present work, the limits of the travelling waves at $ + \ii $ and $ - \ii $ differ.

The proof we give is based on ideas from \cite{RT2,RT1} and make very few 
spectral assumptions on $ \BL $.\\

\noindent {\it Assumption} {\rm (A)}: The spectrum of $ \BL $ consists in 
a finite number (possibly zero) of nonpositive eigenvalues 
$ - \mu_1 $, ... , $ - \mu_q $ in $ ( - \ii , 0 ] $, each one with finite 
multiplicity, and the rest of the spectrum is positive and bounded away from $0$. 
Furthermore that for any $ 1 \leq k \leq q $, we have 
$ \ker ( \BL + \mu_k ) \subset D( \BJ) $ and $ \BJ [ \ker ( \BL + \mu_k ) ] \subset D(\BL) $. 
Finally, there exists $ \vartheta_0 \in [ 0, +\ii ) $ such that 
$ \s_{\rm ess} (\BJ\BL) = \{ i \vartheta , \ \vartheta \in \R , \ |\vartheta | \geq \vartheta_0 \} $.\\

The first hypothesis on the location of the spectrum of $ \BL $ is quite weak, 
since it is verified when $ \BL $ is bounded from below and has essential 
spectrum positive and bounded away from zero. Indeed, if $ \d > 0 $ is such 
that $ \s_{\rm ess} (\BL) \subset [2 \d , +\ii ) $, then the eigenvalues 
of $ \BL $ in $ (-\ii , \d ] $ are isolated, of finite multiplicity, and 
are bounded from below by assumption. The second hypothesis 
$ \ker ( \BL + \mu_k ) \subset D( \BL \BJ ) $ is a regularity assumption 
on the eigenvectors.

Let us recall that Theorem \ref{atmost} ensures that the number of eigenvalues 
(with algebraic multiplicities) of $\mathcal{J} \mathcal{L} $ in the right-half 
space $ \{ {\rm Re} > 0 \} $ is less than or equal to the number of negative 
eigenvalues of $\BL$, hence is finite under assumption (A). Let us now state our 
main result, the proof of which is given in section B.3.1. \\

%%%%%%%%%%%%%%%
\noindent {\bf Theorem B.3} {\it 
%\label{spectrmap}
We make assumption {\rm (A)} and suppose that $ \BJ \BL $ generates 
a continuous semigroup. Then, for any $ t \in ( 0 ,+\ii )$, the spectral 
mapping holds: $ \s(\ex^{t [\BJ \BL]_\C }) \setminus \{ 0 \} = \ex^{ t \s ( [ \BJ \BL]_\C )} $. 
Furthermore, denoting
$$
 \gamma_0 \equiv \sup \Big\{ {\rm Re} (\lambda) , \ \lambda \in \s( [ \BJ \BL ]_\C ) \cap 
\{ {\rm Re} \geq 0 \} \Big\} \in [ 0 , + \ii) ,
$$
for any $ \beta > 0 $, there exists $ M (\beta ) > 0 $ such that, for any $ t \geq 0 $, we have
$$
 \n \ex^{ t \BJ \BL } \n_{\mathscr{L}_c(\BX)} \leq M ( \beta ) \ex^{ (\gamma_0 + \beta ) t } .
$$
Assume in addition $\gamma_0 > 0 $ and denote
$$
 m \equiv \max \Big\{ \text{algebraic multiplicity of } \lambda , 
\lambda \in \s( [ \BJ \BL ]_\C ) \ \text{s.t.} \ {\rm Re} \, \lambda = \gamma_0 \Big\} \in \N^* .
$$
Then, there exists $ M_0 > 0 $ such that, for any $ t \geq 0 $, we have}
$$
 \n \ex^{ t \BJ \BL } \n_{\mathscr{L}_c(\BX)} \leq M_0 (1+t)^{m-1} \ex^{ \gamma_0 t } .
$$
%%%%%%%%%%%%%%%

In particular, Theorem B.3 provides a very simple proof of the spectral 
mapping theorem used in \cite{GJLS} and \cite{dMG}. Indeed, the 
self-adjoint operator $ \BL $ involved in these papers is block diagonal:
$$
\BL = \left( \begin{array}{cc} 
\BL_1 & 0 \\ 0 & \BL_2 
\end{array}\right), 
$$
and both $ \BL_1 $, $ \BL_2 $ have at most two nonnegative eigenvalues. 
More generally, if $ \BL_1 $ and $ \BL_2 $ are closed self-adjoint operators 
on $ X $ verifying assumption {\rm (A)} and if $ \BN : X \to X $ is a linear bounded 
operator which is compact with respect to $\BL_1$ and $\BL_2$, then the self-adjoint 
operator
$$
\BL = \left( \begin{array}{cc} 
\BL_1 & \BN \\ \BN^* & \BL_2 
\end{array}\right) 
$$
also verifies assumption {\rm (A)}. Indeed, $ \BL $ is bounded from below (since 
$\BN$ is bounded) and its essential spectrum is 
$ \s_{\rm ess} (\BL_1 ) \cup \s_{\rm ess} (\BL_2 ) \subset [ \d ,+ \ii ) $ 
for some positive $ \d $, since $ \BN $ is compact with respect to $\BL_1$ and $\BL_2$). 
In \cite{KS} (section 7.1 there) and \cite{GeoOhta} (see Proposition 10 there), 
a spectral mapping theorem is used for such an operator. In \cite{KS}, the 
specific algebra of the problem was used, and for \cite{GeoOhta}, the proof 
relies on the arguments in \cite{GJLS}, but here again, in both cases, 
we may use Theorem B.3 to show the same result.

%%%%%%%%%%%%%%%%%%%%%%%%%%%%%%%%%%%%%%%%%%%%%%%%%%%%%
\section*{B.2 Passing from linear to nonlinear instability}

%%%%%%%%%%%%%%%%%%%%%%%%%%%%%%%%%%%%%%%%%%
\subsection*{B.2.1 Semilinear type models}

We start with a classical result for ``semilinear'' equations, proved in 
section B.3.3.\\

%%%%%%%%%%%%%%%
\noindent {\bf Theorem B.4} 
%\label{semili}
{\it Let $X$ be a real Hilbert space, and consider an evolution equation of the form
$$
 \frac{dv}{dt} = \BA v + \Phi(v) ,
$$
where $ \Phi : X \to X $ is a locally Lipschitz mapping satisfying 
$ \Phi (v) = \BO( |\!| v |\!|_X^2 ) $ as $ v \to 0 $ and $ \BA $ is a linear operator 
which generates a semigroup. We assume that $ \BA_\C : D(\BA_\C) \subset X_\C \to X_\C $ 
has an unstable eigenvalue in the right half plane $ \{ {\rm Re} > 0 \} $ and 
a finite number of eigenvalues in $ \{ {\rm Re} > 0 \} $. We denote
$$
 \gamma_0 \equiv \sup \Big\{ {\rm Re} (\mu) , \ \mu \in \s( [ \BJ \BL ]_\C ) \cap 
\{ {\rm Re} > 0 \} \Big\} \in ( 0 , + \ii) 
$$
and fix $ \lambda \in \s (  \BA_\C ) $ with $ {\rm Re}(\lambda) = \gamma_0 $  
and an associated eigenvector $ w_\C \in D(\BA_\C) $ such that 
$ |\!| {\rm Re}( w_\C ) |\!|_X = 1 $. Assume furthermore that there exists 
$ 0 \leq \beta < \gamma_0 $ and $ M_0 > 0 $ such that
$$
\n \ex^{t \BA } \n_{\mathscr{L}_c(X)} \leq M_0 \ex^{ ( \gamma_0 + \beta ) t } .
$$
Then, $0$ is an unstable solution. More precisely, there exist 
$ K > 0 $, $ \e_0 > 0 $ and $ \d_0 > 0 $ such that, for any $ 0 < \d < \d_0 $, 
the solution $v$ with initial datum $ v^\ini = \d {\rm Re}( w_\C ) \in D(\BA) $ 
exists at least on $ [0, \ln( 2 \e_0 / \d) / \gamma_0 ] $ and verifies, for 
$ 0 \leq t \leq \ln( 2 \e_0 / \d) / \gamma_0 $,
$$
 \n v(t) - \d {\rm Re}( \ex^{ t \lambda} w_\C ) \n_X 
\leq K \d^2 \ex^{ 2 t \gamma_0 } 
\quad \quad \quad  {\it and} \quad \quad \quad 
 \n v(t) \n_X \geq \d \ex^{ t \gamma _0} - K \d^2 \ex^{ 2 t \gamma_0 } ,
$$
In particular, for $ 0 < \e < \e_0 $, we see the instability for 
$ t = \ds{\frac{1}{\gamma_0} \ln \Big( \frac{2 \e}{\d} \Big)} $. If $Y$ 
is a Banach space containing $ X $ and with continuous imbedding $ X \hookrightarrow Y $, 
the trivial solution $ 0 $ is also unstable from $ X $ to $ Y $.}\\
%\end{theorem}
%%%%%%%%%%%%%

Let us observe that it is always possible to choose the (complex) eigenvector 
$w$ so that $ {\rm Re}( w_\C ) \not = 0 $ since for any $\theta \in \R $, 
$ \ex^{i\theta} w$ is also an eigenvector. The 
following corollary deals with the orbital instability. We recall that under 
assumption (A), $ [ \BJ \BL ]_\C $ has a finite number of eigenvalues in 
$ \{ {\rm Re} > 0 \} $.\\

%%%%%%%%%%%%%%%
\noindent {\bf Corollary B.2} 
%\label{semiliorbit} 
{\it We make assumption {\rm (A)} and suppose that $ \BJ \BL $ generates 
a continuous semigroup. Let $\BY$ be a Banach space containing $ \BX $ and 
with continuous imbedding $ \BX \hookrightarrow \BY $. Assume moreover that $ [ \BJ \BL ]_\C $ 
has at least one eigenvalue in $ \{ {\rm Re} > 0 \} $ and choose $ \lambda \in \C $ 
with
$$ {\rm Re}(\lambda) = \gamma_0 \equiv \max \Big\{ {\rm Re} (\mu) , \ \mu \in \s( [ \BJ \BL ]_\C ) 
\cap \{ {\rm Re} > 0 \} \Big\} \in (0 , + \ii) $$
and $ w_\C \in D(\BA_\C) $ an associated eigenvector such that 
$ |\!| {\rm Re}( w_\C ) |\!|_\BX = 1 $. We assume moreover 
that $ \mathfrak{M} \equiv \{ \mathbb{T}(\oo,g) \phi_{\oo_*} , \ \oo \in \R , \ 
g \in \mathbb{G} \} $ is a $ \BC^1 $ submanifold of $ \BX $. We finally 
suppose that the equation {\rm ($\BH$)} is semilinear in the sense that 
there exists $ \Phi : \BX \to \BX $ locally Lipschitz continuous such 
that $ \Phi (v) = \BO( |\!| v |\!|_\BX^2 ) $ as $ v \to 0 $ and
$$
 J (E' - \oo_* Q') ( \phi_{\oo_*} + v ) =  J (E'' - \oo_* Q'')(\phi_{\oo_*}) [v] + \Phi (v) . 
$$
Then, there exists $K > 0 $, $ \e_0 > 0 $ and $ \d_0 > 0 $, depending only 
on $ {\rm Re}( w_\C ) $ and $ \mathfrak{M} $, 
with the following properties. For any $ 0 < \d < \d_0 $, the 
solution $u$ to {\rm ($\BH$)} with initial datum 
$ u^\ini = \phi_{\oo_*} + \d {\rm Re}( w_\C ) \in D(\BA) $ exists at least on 
$ [0, \ln( 2 \e_0 / \d) / \gamma_0 ] $ and verifies, for 
$ 0 \leq t \leq \ln( 2 \e_0 / \d) / \gamma_0 $,
$$
 {\rm dist}_\BY ( u(t) , \mathfrak{M} ) \geq \frac{\d}{K} \ex^{ t \gamma _0} 
- K \d^2 \ex^{ 2 t \gamma_0 } .
$$
In particular, the ``bound state'' solution $ T(\oo_* t) \phi_{\oo_*} $ is 
nonlinearly orbitally unstable from $\BX$ to $\BY$ and, for $ 0 < \e < \e_0 / K $, 
we see the instability for 
$ t = \ds{\frac{1}{\gamma_0} \ln \Big( \frac{2 K \e}{ \d} \Big)} $.}\\
%\end{corollary}
%%%%%%%%%%%%%%%

In \cite{HPW}, a similar assertion is made for the orbital instability 
in the remark after Theorem 2 there, but with $\BY = \BX $. For applications 
to PDE's, the space $\BX$ may be a Sobolev space $H^s$, and $ \BY $ a space 
like $L^2$ or $L^\ii$ for instance. The framework of \cite{GSS} is the single 
energy space (for instance $ H^1 $), but an instability result established 
by tracking exponentially growing modes allows to prove instability from the 
regular space $\BX$ ($ H^1 $) to the nonregular space $\BY$ ($L^2$ or $L^\ii$). 
Here, we may obtain instability in $ L^2 $. \\

% This is an advantage of instability results 
% established by tracking exponentially growing modes: it allows to prove 
% instability from the regular space $\BX$ to the nonregular space $\BY$.\\
% 
% Moreover, let us 
% observe that \cite{GSS} shows the instability with $ \BX $ the energy space 
% . Here, we may obtain instability from $H^1$ to $ L^2 $. 

%%%%%%%%%%%%%%
\noindent {\bf Remark B.1} In the framework of \cite{GSS}, where a Lyapounov type 
functional is used, it follows that the instability is seen for a time at most equal 
to $ K \ds{\frac{\e}{\d^2}} $, where $K$ is some positive constant. 
This timescale is much larger than the natural one 
$ \ds{\frac{1}{\gamma_0} \ln \Big( \frac{2 K \e}{ \d} \Big)} $.\\
%\end{remark}
%%%%%%%%%%%%%

%%%%%%%%%%%%%%%%%%%%%%%%%%%%%%%%%%%%%%%%%%%%%%%%%%%%%%%%%%%%%%%%%%%%%%%%%%%%%%%%%%%%%%%%
\subsection*{B.2.2 Some applications}

We may apply our result to the Nonlinear Schr\"odinger Equation
\be
\tag{NLS}
 i \p_t \Psi + \Delta \Psi + \Psi f(|\Psi|^2) = 0 ,
\ee
or the Nonlinear Klein-Gordon equation
\be
\tag{NLKG}
 \p^2_t \Psi = \Delta \Psi + \Psi f(|\Psi|^2)  
\ee
in $ \R^d $. We shall consider a nonlinearity $f$ at least $\BC^1$, 
so that we are in the framework of \cite{GSS}. \\

\noindent $\bullet$ A bound state solutions for these two equations 
is a particular solution of the form $ U(t) = \ex^{i\oo t} \phi_\oo $. 
The instability is in general linked to the fact that
\be
\tag{B.1}
\frac{d}{d\oo} \int_{\R^d} |\phi_\oo|^2 \ dx < 0 
\quad \quad \quad {\rm for \ (NLS),} \quad \quad \quad 
{\rm resp.}\quad \frac{d}{d\oo} \Big( \oo \int_{\R^d} |\phi_\oo|^2 \ dx \Big) < 0
\quad \quad \quad {\rm for \ (NLKG).} 
\ee
The existence of at least one unstable eigenvalue has been shown 
under assumption (B.1) by \cite{Gri} for radial bound states 
with an arbitrary number of nodes and in \cite{GSS2} for radial ground 
states. Corollary B.2 may be applied with $ \BX = H^s (\R^d) $, 
where $s \in \N $, $ s > d/ 2 $ and assuming that the nonlinearity 
verifies $ f \in \BC^{s+2} $, and $\BY= L^2(\R^d)$ or $L^\ii(\R^d)$. 
The result in \cite{Mizu} shows the instability of linearly unstable bound 
states for (NLS) (in dimension $d=2$) with $f(\varrho) = \varrho^{\frac{p-1}{2}} $ 
by showing the exponential growth of an unstable eigenmode. Our 
result gives a simple proof of this result, but restricted to the 
sufficiently smooth cases, namely $ p $ an odd integer or $ p > 5 + 2s > 5 + d $. 
For non smooth nonlinearities, the situation is more delicate (see \cite{Mizu}). 
An alternative approach is to combine Strichartz estimates with the 
growth estimate on the semigroup $ \ex^{ t\BJ \BL } $ given in 
Theorem B.3, as in \cite{GeoOhta}.\\

\noindent $\bullet$ Corollary B.2 also applies 
to the Discrete Nonlinear Schr\"odinger Equation
\be
\tag{DNLS}
\forall n \in \Z, \quad \quad \quad 
i \p_t \Psi_n + \e \Big( \Psi_{n+1} - 2 \Psi_{n} + \Psi_{n-1} \Big) + \Psi_n f(|\Psi_n|^2) = 0 ,
\ee
as studied in \cite{MeChKeCu} with the saturated nonlinearity 
$ f ( \varrho ) = \ds{ \frac{\beta}{1+ \varrho} } $, $\beta > 0 $ 
(existence of travelling wave solution) and in \cite{FiKeSuFr} (defocusing cubic 
(DNLS), {\it i.e.} $ f ( \varrho ) = - \beta \varrho $ for some $\beta > 0 $). 
The numerical analysis in the paper \cite{FiKeSuFr} shows the existence of linearly 
unstable bound state solutions. The travelling wave solutions numerically obtained 
in \cite{MeChKeCu} are linearly stable, but it may happen that for other nonlinearities 
$f$, some are linearly unstable. 
%A specificity of equations on a lattice is that there is no 
%analog of the ``momentum'' or ``charge'', hence these models do not enter into 
%the Grillakis-Shatah-Strauss framework. 

%%%%%%%%%%%%%%%%%%%%%%%%%%%%%%%%%%%%%
\subsection*{B.2.3 Quasilinear PDE's}

For quasilinear problems, we shall not make restrictions on the smoothness 
of the nonlinearity. The result relies on the strategy of E. Grenier \cite{Grenier} 
and the works \cite{RT2,RT1}. We consider the evolution equation
\be
\tag{E}
 \frac{du}{dt} = J ( L_0 u + \nabla F (u) ) 
\ee
for $ u : \R^d \to \R^\nu $, where $ F \in \BC^\ii ( \R^\nu , \R ) $, with the 
following hypothesis. The operator $J$ is a Fourier multiplier, skew-symmetric on $L^2$, 
into and with domain containing $H^1$. There exists $ \s >0 $ such that the operator $ L_0 $ 
is a Fourier multiplier with domain containing $H^{2\s}$, symmetric and having a self-adjoint 
realization on $L^2( \R^d , \R^\nu ) $. Moreover, for some $ C > 0 $, $L_0$ verifies: 
$$ 
\frac{1}{C} |\!| u |\!|_{H^\s}^2 %- C |\!| u |\!|_{H^{\s_0}}^2 
\leq ( L_0 u , u )_{L^2} 
\leq C |\!| u |\!|_{H^\s}^2 .
%+ C |\!| u |\!|_{H^{\s_0}}^2 .
$$
The framework proposed in \cite{RT1} was for $ L_0 $ coercive in $ H^1 $, that is 
$ \s = 1 $. For the examples below, we shall have $ \s = 1/2 $ or $ \s = 2 $, which 
requires very few modifications to the proof of \cite{RT1}. We still assume that 
for some group $ \mathbb{G} $, there exists a unitary representation of 
$ \mathbb{G} $ on $ \BX $, $ \mathbb{T} : \mathbb{G} \to \mathscr{G}\mathscr{L}_c (\BX ) $, 
leaving the equation (E) invariant.

We consider a stationary solution of the evolution equation (E), 
that is some $ Q \in H^\ii (\R^d, \R^\nu ) $ such that $ L_0 Q + \nabla F (Q) = 0 $. We are 
interested in the stability of this solution. We assume that the commutator 
$ [ J, \nabla^2 F (Q) ] $ is bounded in $ L^2 $, which is the case when $ J $ 
is bounded in $ L^2 $ or when $d=1$ and $ J = \p_x $. We suppose that for the problem
$$
 \frac{ \p u }{ \p t } = J ( L_0 u + \nabla F( u^a + u ) - \nabla F ( u^a ) + G ) ,
$$
where $ u^a $ is smooth, bounded as well as its derivatives and 
$ G \in \BC( \R , H^s ) $ for every $ s $, we have local well-posedness for 
$s$ large enough: there exists a time $ T > 0 $ and a unique solution in 
$ \BC( [ 0 , T ] , H^s ) $. We moreover assume that for some continuous 
non-decreasing function $ \kappa : \R_+ \to \R_+ $ with $ \kappa (0) =0 $, the tame 
estimate
$$
 | ( \p_x^\alpha J \{ \nabla^2 F ( w +v )[ v] \}, \p_x^\alpha v )_{L^2} | \leq 
 \kappa( \n w \n_{W^{s+1,\ii}} + \n v \n_{ H^s} ) \n v \n_{H^s}^2 ,
$$
with $ | \alpha | \leq s $ holds true. In order to control high order derivatives, we 
finally require that for $s$ large enough, there exists a self-adjoint operator 
$ \mathbb{M}_s $ and $C_s$ such that
$$
 | ( \mathbb{M}_s u , v )_{L^2} | \leq C_s \n u \n_{H^s} \n v \n_{H^s} , 
\quad \quad \quad ( \mathbb{M}_s u , u )_{L^2} \geq \n u \n_{H^s}^2 - C_s \n u \n_{H^{s - \min(\s , 1)}}^2
$$
and
$$
 {\rm Re} ( J L u , \mathbb{M}_s u )_{L^2} \leq C_s \n u \n_{H^s} \n u \n_{H^{s - \min(\s , 1) }}
$$
(for a criterion which ensures the existence of such a multiplier, see Lemma 5.1 in 
\cite{RT1}).\\

%such that the we assume that for
%There exists $ C > 0 $ such that $ |\!|
%in the Sobolev spaces $ H^s ( \R^d ) $, $ s \in \N$. Here, $ \BA $ is a linear operator 
%in $ H^s $ which generates a semigroup and skew-adjoint in $ L^2 $. We assume 
%that $ \BA_\C : D(\BA_\C) \subset H^s_\C \to H^s_\C $ verifies assumption (A) in $ L^2 $.
%$ H^s $, $ s \in \N $.
%has an unstable eigenvalue in the right half-plane $ \{ {\rm Re} > 0 \} $ and 
%We assume that there exists $ \ell \in \N $ such that, for any $ s \in \N $,
%$$
 %\n \Phi (v) \n_{H^s} \leq C_s \n v \n_{H^{s+\ell}}^2 
 %\quad \quad \quad {\rm as} \quad v \to 0 \quad {\rm in} \quad H^{s+\ell}.
%$$

Adapting the strategy of \cite{RT2,RT1}, we may deduce the following result. 
Since the proof is very similar, we omit it.\\

%%%%%%%%%%%%%%%
\noindent {\bf Theorem B.5} {\it 
%\label{quasili}
We make the above assumptions and moreover that $ L_0 + \nabla^2 F (Q) $ 
verifies hypothesis {\rm (A)} in $ L^2 $. We assume furthermore that 
$ [ J (  L_0 + \nabla^2 F (Q) ) ]_\C $ has an unstable eigenvalue in the right 
half-plane $ \{ {\rm Re} > 0 \} $, denote
$$
 \gamma_0 \equiv \sup \Big\{ {\rm Re} (\lambda) , \ \lambda \in 
 \s( [ J (  L_0 + \nabla^2 F (Q) ) ]_\C ) \cap \{ {\rm Re} > 0 \} \Big\} \in ( 0 , + \ii) 
$$
and fix $ \lambda \in \s ( [ J ( L_0 + \nabla^2 F (Q) ) ]_\C ) $ with 
$ {\rm Re}(\lambda) = \gamma_0 $  and an associated eigenvector 
$ w_\C \in D([ J ( L_0 + \nabla^2 F (Q) ) ]_\C) $ such that 
$ |\!| {\rm Re}( w_\C ) |\!|_{H^s} = 1 $. There exists $ s_0 \in \N $ such that, if 
$ s \geq s_0 $, $Q$ is nonlinearly unstable from $ H^s $ to $ L^2 $ and to 
$L^\ii$: there exists $ K > 0 $, $ \e_0 > 0 $ and $ \d_0 > 0 $ such that, for any 
$ 0 < \d < \d_0 $, the $H^s$ solution $u$ to {\rm (E)} with initial datum 
$ u^\ini = Q + \d {\rm Re}( w_\C ) \in H^s $ exists at least on 
$ [0, \ln( 2 \e_0 / \d) / \gamma_0 ] $ and verifies, for 
$ 0 \leq t \leq \ln( 2 \e_0 / \d) / \gamma_0 $,
$$
 \n u(t) - Q - \d {\rm Re}( \ex^{ t \lambda} w_\C ) \n_{H^s} 
\leq K \d^2 \ex^{ 2 t \gamma_0 } ,
$$
hence
$$
\n u(t) - Q \n_{L^2} \geq \d \ex^{ t \gamma _0} - K \d^2 \ex^{ 2 t \gamma_0 } 
\quad \quad \quad  {\it and} \quad \quad \quad 
 \n u(t) - Q \n_{L^\ii} \geq \d \ex^{ t \gamma _0} - K \d^2 \ex^{ 2 t \gamma_0 } .
$$
If, in addition, $ \mathfrak{M} \equiv \{ \mathbb{T}( g) Q , \ g \in \mathbb{G} \} $ 
is a $ \BC^1 $ submanifold of $ H^s $, then we also have
$$
 {\rm dist}_{L^2} ( u(t) , \mathfrak{M} ) \geq K \d \ex^{ t \gamma _0} - K \d^2 \ex^{ 2 t \gamma_0 } 
\quad \quad \quad  {\it and} \quad \quad \quad 
 {\rm dist}_{L^\ii} ( u(t) , \mathfrak{M} ) \geq K \d \ex^{ t \gamma _0} - K \d^2 \ex^{ 2 t \gamma_0 } .
$$
In particular, for $ 0 < \e < \e_0 / K $, we see the nonlinear orbital instability for 
$ t = \ds{\frac{1}{\gamma_0} \ln \Big( \frac{2 K\e}{\d} \Big)} $. }
\subsection*{B.2.4 Some applications to nonlinear dispersive wave equations}

%$ \bullet $ {\bf Nonlinear dispersive wave equations.} 
Some model quasilinear equations are given by wave equations (in one 
space dimension) such as the generalized Korteweg-de Vries equation
\be
\tag{gKdV}
 \p_t u + \p_x(f(u)) + \p_x^3 u = 0 ,
\ee
the generalized regularized Korteweg-de Vries equation, also called 
Benjamin-Bona-Mahony equation or Peregrine equation when 
$ f(u) = u^2 /2 $,
\be
\tag{gBBM}
 \p_t u + \p_x u + \p_x(f(u)) - \p_t \p_x^2 u = 0 ,
\ee
the generalized regularized Boussinesq equation
\be
\tag{grBsq}
 \p_t^2 u - \p^2_x u - \p^2_x(f(u)) - \p_t^2 \p_x^2 u = 0 .
\ee
Each of these equations admit nontrivial solitary wave solutions 
$ u(t,x) = U_c( x - c t ) $ for $c$ in $ (0, +\ii ) $, 
$ (1, +\ii ) $ and $ ( - \ii , -1 ) \cup (1, +\ii ) $ respectively. 
For these solitary wave solutions, the momentum is, respectively,
$$
 P (U_c) = \int_\R U_c^2 \ dx = \n U_c \n_{L^2}^2 ,
 \quad \quad  
 P (U_c) = \int_\R U_c^2 + (\p_x U_c)^2 \ dx ,
 \quad \quad 
 P ( U_c) = c \int_\R U_c^2 + (\p_x U_c)^2 \ dx .
$$
The existence of exactly one unstable eigenvalue has been shown 
with the use of an Evans function by R. Pego and M. Weinstein in 
\cite{PeWe} for these three equations under the condition 
$ \ds{ \frac{dP(U_c)}{dc} < 0 } $. The paper \cite{Lopes} by O. Lopes 
also gives a linear instability result. Equations (gBBM) and (grBsq) 
turn out to be semilinear due to the regularization effect. Indeed, they may 
be written
$$
\p_t u + (1 - \p_x^2)^{-1} \p_x u + ( 1 - \p_x^2 )^{-1} \p_x (f(u)) = 0 , 
\quad \quad \quad 
 \p_t^2 u - ( 1 - \p_x^2 )^{-1} \p^2_x u - ( 1 - \p_x^2 )^{-1} \p^2_x(f(u)) = 0 .
$$
Therefore, Corollary B.2 applies to these two models and this shows the 
nonlinear instability when linear instability holds.

In \cite{Lin2}, some generalizations of the equations 
(gKdV), (gBBM) and (gBBM) have been proposed that take into account 
pseudodifferential operators. These are respectively
\be
\tag{I}
 \p_t u + \p_x(f(u)) - \p_x \BM u = 0 ,
\ee
\be
\tag{II}
 \p_t u + \p_x u + \p_x(f(u)) + \p_t \BM u = 0 
\ee
and
\be
\tag{III}
 \p_t^2 u - \p^2_x u - \p^2_x(f(u)) + \p_t^2 \BM u = 0 . 
\ee
Here, $\BM$ is a Fourier multiplier of symbol $ \hat{\BM} $: 
$\widehat{\BM w} = \hat{\BM} \hat{w} $ (here, $\hat{\cdot}$ denotes 
Fourier transform). We assume $ \hat{\BM} \geq 0 $ (otherwise, see \cite{Lin2}). 
When $ \BM = - \p_x^2 $, these equations reduce to (gKdV), (gBBM) and (gBBM) 
respectively. The Benjamin-Ono equation ($ \hat{\BM} = | \xi |$), 
the Smith equation ($\hat{\BM} = \sqrt{1+ \xi^2 } - 1 $) and the Intermediate Long Wave 
(or Whitham) equation ($ \hat{\BM} = \xi / {\rm tanh}(\xi H) - 1 / H $, for some 
constant $H> 0$) are common models of dispersive wave equations that are of type (I). 
We refer to \cite{Lin2} for references on these models and the existence of solitary waves. 
The associated momentum is
$$
 P_I (U_c) = \int_\R U_c^2 \ dx = \n U_c \n_{L^2}^2
 \quad \quad 
 P_{II} (U_c) = \n (1+ \BM)^{1/2} U_c \n_{L^2}^2 
 \quad \quad 
 P_{III} (U_c) = c \n (1+ \BM)^{1/2} U_c \n_{L^2}^2 .
$$
For these models, Evans' function type arguments do not work since 
we no longer have a differential equation (it is nonlocal). The paper 
\cite{Lin2} by Z. Lin proposes another approach than the Evans' function 
technique for establishing the existence of unstable eigenvalues. 
However, it is not completely clear whether this method extends easily to 
the case of systems such as the Euler-Korteweg system (EK) (given at the 
beginning of section \ref{sexins}). \\

%%%%%%%%%%%%%%%
\noindent {\bf Theorem B.6 (\cite{Lin2})}
%\label{instabLin2} 
{\it We consider one of the equations {\rm (I)}, {\rm (II)} or {\rm (III)} 
with $f $ of class $\BC^1$ satisfying $ f(0) = f'(0) = 0 $ and $| f(u) | \gg | u | $ 
for $ |u| \to + \ii $. We assume moreover that $ \hat{\BM} $ is even, nonnegative, 
and verifies, for some $ m \geq 1 $, 
$ \ds{ 0 < \varliminf_{+\ii} \frac{\hat{\BM}(\xi)}{\xi^m} \leq \varlimsup_{+\ii} 
\frac{\hat{\BM}(\xi)}{\xi^m} < \ii } $. Assume that $ c \mapsto \phi_c = U_c (x - ct ) $ 
is a $ \BC^1 $ branch of travelling wave solution to {\rm (I)}, {\rm (II)} or {\rm (III)} 
with $ U_c \in H^{m/2}(\R) $ defined near $ c_* $ and suppose that the linearized operator 
$ \BL $ has exactly one negative eigenvalue, that $ \ker{ \BL} $ is spanned by 
$ \p_x U_{c_*} $ and that $ \ds{\frac{dP(U_c)}{dc}_{|c = c_*} < 0 } $. Then, 
$ U_{c_*} $ is linearly unstable.}\\
%\end{theorem}
%%%%%%%%%%%%%

It is not easy to determine whether the hypothesis of Theorem 
B.6 hold true when $\BM$ is not a (differential) Sturm-Liouville 
operator. See however \cite{Al} on this question. It is clear that if the 
assumptions of Theorem B.6 are verified, then assumption (A) is 
also satisfied. As for the (gBBM) and the (grBsq) equations, the equations 
(II) and (II) turn out to be semilinear, thus we may prove nonlinear orbital 
instability by applying Corollary B.2.

\bigskip

The Kawahara equation (or fifth order (KdV) equation) 
\be
\tag{K}
 \p_t u + \p_x(f(u)) + \alpha \p_x^3 u + \beta \p_x^5 u = 0 ,
\ee
with $ \alpha $, $\beta \not = 0 $ two real constants, is another relevant 
dispersive model. For this equation, it may happen 
that the linearized equation around the solitary wave has more than 
one negative eigenvalues, in which case the works \cite{GSS,GSS2,Lopes,Lin2} 
do not give a clear necessary and sufficient condition for stability. 
The paper \cite{BriDe} by T. Bridges and G. Derks gives a sufficient 
condition for linear instability for solitary wave solutions, but also 
for other types of travelling solutions. This condition is probably not necessary 
since it may happen that there exist at least two unstable eigenvalues, 
or two complex conjugate eigenvalues.

Instead of stating a general result for nonlinear orbital instability, we shall 
consider several model cases on which we will verify the hypothesis of Theorem B.5, 
in particular the question of the existence of the multiplier $ \mathbb{M}_s $. \\

%%%%%%%%%%%%%%%
\noindent {\bf Proposition B.1} 
{\it We consider the equation {\rm (I)}, namely
$$
 \p_t u + \p_x(f(u)) - \p_x \BM u = 0 
 $$
with $f $ of class $\BC^1$ satisfying $ f(0) = f'(0) = 0 $ and $ | f(u) | \gg | u | $ 
for $ |u| \to + \ii $. We assume that $ \hat{\BM} $ is one of the following functions:
$$
-  \xi^2 \quad {\rm (KdV)}; \quad \quad \quad
 \xi^4 + \alpha \xi^2 \quad {\rm (Kawahara)}; \quad \quad \quad
 |\xi| \quad {\rm (Benjamin-Ono)};
 $$
 $$ 
 \frac{ \xi}{ {\rm tanh}(\xi H) } - \frac{1}{ H} 
 \quad {\rm (Intermediate \ Long\ Wave)} ; \quad \quad \quad
 \sqrt{1+ \xi^2 } - 1 \quad {\rm (Smith)} .
$$
There exists $ s_0 > 0 $ such that, if there exists $ c \in \R $ such that 
{\rm (I)} has a nontrivial solitary wave $ U_c \in L^2 $ which is linearly 
unstable, then for any $ s \geq s_0 $, it is also nonlinearly unstable from 
$ H^s $ to $ H^s $, to $ L^2 $ and to $ L^\ii $.
}\\
%%%%%%%%%%%%%%%

By application of Theorem B.5, we are thus able to show the nonlinear instability 
from $H^s$ to $L^2$ or $L^\ii$ by tracking the exponentially growing mode (this 
question was left open in \cite{Lin2} and also in \cite{Lopes}). In particular, we 
obtain the $L^2$ nonlinear instability of the linearly unstable solitary waves for 
these models.\\

\noindent {\it Proof.} All the assumptions for Theorem B.5 for these types of 
models are verified in section 8.1 in \cite{RT1}, except the existence of 
the multiplier $ \mathbb{M}_s $. 

For (KdV), where $ \s = 1 $, we shall take (for $ s \geq 2 $ integer)
$$
 \mathbb{M}_s \equiv (-1)^s \p_x^{2s} + 
\frac{1+2s}{3} (-1)^{s-1} \p_x^{s-1} \{ f'(Q) \p_x^{s-1} \cdot \} ,
$$
as the computations from \cite{RT1}, section 8.1 show. For the Kawahara equation, with $ \s = 2 $, 
we take (for $ s \geq 4 $ integer)
$$
 \mathbb{M}_s \equiv (-1)^s \p_x^{2s} + 
\frac{1+2s}{5} (-1)^{s-1} \p_x^{s-2} \{ f'(Q) \p_x^{s-2} \cdot \} 
$$
and since the computations are very similar, we omit them. 
For the Benjamin-Ono equation, we have $ \hat{\BM} (\xi )= | \xi | $ and $ \s = 1/2 $, and 
we will then have to deal with pseudo-differential operator which are Fourier multipliers 
with homogeneous symbol. For this type of operator, we shall need some commutator 
estimates. We denote by $ \mathscr{F} (w) $ or $ \hat{w} $ the Fourier transform of $w$, 
and $ \mathscr{H} $ the Fourier multiplier with symbol $ - i\, {\rm sgn}(\xi ) $ (this 
is the Hilbert transform).\\

%%%%%%%%%%%%%%%
\noindent {\bf Lemma B.1} 
{\it (i) Let $ h \in L^{\ii} (\R) $ with $ \mathscr{F} ( \BM^{\frac12} h ) \in L^1 (\R) $ 
(for instance, $ h \in H^\s (\R)$ for some $\s > 1$). Then, 
there exists $ C > 0 $ such that, for any $ v \in H^{\frac12} (\R)$,
$$
 \n \BM^{\frac12} ( h v ) - h \BM^{\frac12} v \n_{L^2(\R)} \leq C \n v \n_{L^2(\R)} .
$$
(ii) Let $ h \in L^{\ii}(\R) $ with $ \mathscr{F} ( \BM^{\frac32} h ) \in L^1 (\R)$ 
(for instance, $ h \in H^\s(\R) $ for some $\s > 2$). Then, 
there exists $ C > 0 $ such that, for any $ v \in H^{\frac32} (\R)$,
\begin{align*}
 \no \BM^{\frac32} \{ h v \} - h \BM^{\frac32} v 
 - \frac32 [ \p_x h] \BM^{\frac12} \mathscr{H} v \no_{L^2(\R)} 
  \leq & \, C \n v \n_{L^{2}(\R)} .
\end{align*}
(iii) Let $ h \in L^{\ii}(\R) $ with $ \mathscr{F} ( \p_x \BM^{\frac12} h ) \in L^1(\R) $ 
(for instance, $ h \in H^\s(\R) $ for some $\s > 2$). Then, 
there exists $ C > 0 $ such that, for any $ v \in H^{\frac32}(\R) $,
\begin{align*}
 \no \p_x \BM^{\frac12} \{ h v \} - h \p_x \BM^{\frac12} v 
 - \frac32 [ \p_x h] \BM^{\frac12} v \no_{L^2(\R)} \leq & \, C \n v \n_{L^{2}(\R)} .
\end{align*}
}\\
%%%%%%%%%%%%%%%

\noindent {\it Proof.} We have
$$ 
 \mathscr{F} \Big( \BM^{\frac12} ( h v ) - h \BM^{\frac12} v \Big) (\xi) = 
 \int_\R | \xi |^{\frac12} \hat{h} ( \xi - \zeta) \hat{v}(\zeta)   \ d \zeta 
 - \int_\R | \zeta |^{\frac12} \hat{h} ( \xi - \zeta) \hat{v}(\zeta)   \ d \zeta .
$$
Using the inequality 
$ \Big| | \xi |^{\frac12} - | \zeta |^{\frac12} \Big| \leq C | \xi - \zeta |^{\frac12} $, 
we thus obtain
$$
 \Big| \mathscr{F} \Big( \BM^{\frac12} ( h v ) - h \BM^{\frac12} v \Big) (\xi) \Big| 
 \leq C \int_\R | \xi - \zeta |^{\frac12} | \hat{h} ( \xi - \zeta) | \cdot 
|\hat{v}(\zeta) | \ d \zeta = C \{ | \mathscr{F} ( \BM^{\frac12} h ) | * | \hat{v} | \} (\xi) 
$$
and we conclude with the classical convolution estimate $ L^1 * L^2 \subset L^2 $. 
This argument does not provide the sharpest bound in $h$, since it involves 
$ \n \mathscr{F} ( \BM^{\frac12} h ) \n_{ L^1} $, whereas the use of paradifferential 
calculus will use only $ \n h \n_{\BC^{\frac12}} $. However, we shall to use this 
refinement here. 

The starting point for the second inequality is
$$
 \Big| | \xi |^{\frac32} - | \zeta |^{\frac32} - 
\frac32 | \zeta |^{\frac12} {\rm sgn}(\zeta) ( \xi - \zeta ) \Big| 
 \leq C | \xi - \zeta |^{\frac32} .
$$
To prove this, note that by homogeneity $ \xi = \theta \zeta $, $ \theta \in \R $, it 
suffices to prove
$$
 \Big| | \theta |^{\frac32} - 1 - \frac32 ( \theta - 1 ) \Big| \leq C | \theta - 1 |^{\frac32} ,
$$
which is easy. Therefore,
\begin{align*}
 \Big| & \, 
 \mathscr{F} \Big( \BM^{\frac32} \{ h v \} - h \BM^{\frac32} v 
 - \frac32 [ \p_x h ] \BM^{\frac12} \mathscr{H} v \Big) (\xi) \Big|
 \\ & 
 =  \Big| \int_\R | \xi |^{\frac32} \hat{h} ( \xi - \zeta) \hat{v}(\zeta) \ d \zeta 
 - \int_\R | \zeta |^{\frac32} \hat{h} ( \xi - \zeta) \hat{v}(\zeta) \ d \zeta  
 - \int_\R \frac32 | \zeta |^{\frac12} {\rm sgn}(\zeta) ( \xi - \zeta ) \hat{h} ( \xi - \zeta) 
\hat{v}(\zeta) \ d \zeta \Big| 
\\ & \leq C  \int_\R | \xi - \zeta |^{\frac32} | \hat{h} ( \xi - \zeta) | \cdot | \hat{v}(\zeta) | \ d \zeta 
\\ & = C  | \mathscr{F} ( \BM^{\frac32} h ) | * | \hat{v} | ,
\end{align*}
and we conclude as before. For the third inequality, we argue in a similar way 
with the estimate
$$
 \Big| i \xi | \xi |^{\frac12} - i \zeta | \zeta |^{\frac12} 
- i \frac32 | \zeta |^{\frac12} ( \xi - \zeta ) \Big| 
 \leq C | \xi - \zeta |^{\frac32} .
$$
The proof is complete. \carre \\

For the Benjamin-Ono equation, $ \hat{\BM}(\xi = |\xi| $, $\s = 1/2$ and 
the index $s$ will be half an integer: $ s \in \N / 2 $. Therefore, we set 
$ s = [s] + \{ s \} $, with $ [s] $ integer and $ \{ s \} \in \{0 ; 1/2 \} $. 
Let us define, for $ s \in \N / 2 $, $ s \geq 1 $,
$$
 \mathbb{M}_s \equiv 
 \left\{\begin{array}{ll}
\ds{ (-1)^s \p_x^{2s} 
 + \gamma_s \BM^{\frac12} \p_x^{s-1} \{ f'(Q) \p_x^{s-1} \BM^{\frac12} \cdot \} }  
& \quad {\rm if \ } \{ s \} = 0 
\\ \ \\ 
\ds{ (-1)^{[s]} \p_x^{2[s]} \BM  
 + \gamma_s \p_x^{[s]} \{ f'(Q) \p_x^{[s]} \cdot \} } 
& \quad {\rm if \ } \{ s \} = \frac12 ,
\end{array}\right. 
$$
for some real constant $ \gamma_s $ to be determined later.  It is clear that
$ \mathbb{M}_s $ is self-adjoint on $L^2$ and that there exists  $C_s > 0 $ such that
$$
 | ( \mathbb{M}_s u , v )_{L^2} | \leq C_s \n u \n_{H^s} \n v \n_{H^s} 
\quad \quad \quad {\rm and} \quad \quad \quad 
( \mathbb{M}_s u , u )_{L^2} \geq \n u \n_{H^s}^2 - C_s \n u \n_{H^{s-\frac12}}^2 .
$$
In order to verify the assumptions for the multiplier $  \mathbb{M}_s $, it remains to 
study $ {\rm Re} ( J ( L_0 + \nabla^2 F(Q) ) u , \mathbb{M}_s u )_{L^2}  $. 
When $ \{ s \} = 0 $, {\it i.e.} $ s \in \N $, this quantity is
\begin{align}
 {\rm Re} ( \p_x ( \BM + c + & \, f'(Q) ) u , \mathbb{M}_s u )_{L^2} 
 \nonumber \\
& =  {\rm Re} ( \p_x \BM u , (-1)^s \p_x^{2s} u )_{L^2} 
 + \gamma_s {\rm Re} ( \p_x \BM u , 
\BM^{\frac12} \p_x^{s-1} \{ f'(Q) \p_x^{s-1} \BM^{\frac12} u \} )_{L^2} 
 \nonumber \\ & \quad
 + {\rm Re} ( \p_x [ f'(Q) u ] , (-1)^s \p_x^{2s} u )_{L^2} 
 + \gamma_s {\rm Re} ( \p_x [ f'(Q) u ] , 
\BM^{\frac12} \p_x^{s-1} \{ f'(Q) \p_x^{s-1} \BM^{\frac12} u \} )_{L^2} 
 \nonumber \\ \tag{B.1}
& \quad + c  {\rm Re} ( \p_x u , \mathbb{M}_s u )_{L^2} .
\end{align}
By skew-adjointness, the first and last scalar products are zero. By integration 
by parts and Leibniz formula, we deduce, since $ Q \in H^\ii $, 
\begin{align*}
 {\rm Re} ( \p_x [ f'(Q) u ] , (-1)^s \p_x^{2s} u )_{L^2} 
 = & \, {\rm Re} ( \p_x^{s+1} [ f'(Q) u ] , \p_x^{s} u )_{L^2} 
 \\ 
 \leq & \, {\rm Re} ( f'(Q) \p_x^{s+1} u , \p_x^{s} u )_{L^2} 
 + (s + 1) {\rm Re} ( \p_x [f'(Q)] \p_x^{s} u , \p_x^{s} u )_{L^2} 
 + C_s \n u \n_{H^s} \n u \n_{H^{s-1}} 
  \\ 
 \leq & \, \Big( s + \frac12 \Big) {\rm Re} ( \p_x [f'(Q)] \p_x^{s} u , \p_x^{s} u )_{L^2} 
 + C_s \n u \n_{H^s} \n u \n_{H^{s-1}} .
\end{align*}
Similarly, using the easy estimates $ \n \BM^{\frac12} v \n_{L^2} \leq K \n v \n_{H^{\frac12}} $ 
and $ \n h v \n_{H^{\frac12}} \leq C(h) \n v \n_{H^{\frac12}} $ for $ h \in L^\ii $ with 
$ \mathscr{F} (\BM^{\frac12} h ) \in L^1 $ (this is an immediate consequence of Lemma B.1)
\begin{align*}
 \gamma_s {\rm Re} ( \p_x [ f'(Q) u ] , & \, \BM^{\frac12} \p_x^{s-1} 
\{ f'(Q) \p_x^{s-1} \BM^{\frac12} u \} )_{L^2} 
 = 
 \gamma_s (-1)^{s-1} {\rm Re} ( \BM^{\frac12} \p_x^{s} [ f'(Q) u ] , f'(Q) \p_x^{s-1} \BM^{\frac12} u )_{L^2} 
 \\ \leq & \, 
 \gamma_s (-1)^{s-1} {\rm Re} ( \BM^{\frac12} [ f'(Q) \p_x^s u ] , f'(Q) \p_x^{s-1} \BM^{\frac12} u )_{L^2} 
 + C \n u \n_{H^{s-\frac12} }^2 .
\end{align*}
Using Lemma B.1, we deduce 
$ \n \BM^{\frac12} [ f'(Q) \p_x^s u ] - f'(Q) \BM^{\frac12} \p_x^s u \n_{L^2} \leq 
C (Q) \n u \n_{ H^s } $, thus
\begin{align*}
 \gamma_s {\rm Re} ( \p_x [ f'(Q) u ] , & \, 
\BM^{\frac12} \p_x^{s-1} \{ f'(Q) \p_x^{s-1} \BM^{\frac12} u \} )_{L^2} 
 \\ \leq & \, 
 \gamma_s (-1)^{s-1} {\rm Re} ( f'(Q) \p_x^s \BM^{\frac12} u , f'(Q) \p_x^{s-1} \BM^{\frac12} u )_{L^2} 
 + C  \n u \n_{H^s} \n u \n_{H^{s-\frac12} } 
 \\ = & \, 
\frac{ \gamma_s}{2} (-1)^{s} 
{\rm Re} ( \p_x [ f'(Q) ] \p_x^{s-1} \BM^{\frac12} u , f'(Q) \p_x^{s-1} \BM^{\frac12} u )_{L^2} 
+ C  \n u \n_{H^s} \n u \n_{H^{s-\frac12} } 
\\ \leq & \, 
C \n u \n_{H^{s-\frac12} }^2 + C  \n u \n_{H^s} \n u \n_{H^{s-\frac12} } 
\leq C \n u \n_{H^s} \n u \n_{H^{s-\frac12} } .
\end{align*}
We now turn to the term
\begin{align*}
 \gamma_s {\rm Re} ( \p_x \BM u , \BM^{\frac12} \p_x^{s-1} \{ f'(Q) \p_x^{s-1} \BM^{\frac12} u \} )_{L^2} 
 = \gamma_s (-1)^{s-1} {\rm Re} ( \p_x^s u , \BM^{\frac32} \{ f'(Q) \p_x^{s-1} \BM^{\frac12} u \} )_{L^2} .
 \end{align*}
Using Lemma B.1, we write
\begin{align*}
 \no \BM^{\frac32} \{ f'(Q) \p_x^{s-1} \BM^{\frac12} u \} - f'(Q) \p_x^{s-1} \BM^{2} u 
 - \frac32 \p_x [f'(Q)] \BM^{\frac12} \mathscr{H} \{ \p_x^{s-1} \BM^{\frac12} u \} \no_{L^2} 
  \leq & \, C (Q) \n \p_x^{s-1} \BM^{\frac12} u  \n_{L^{2}} \\
  \leq & \, C (Q) \n u \n_{H^{s-\frac12}} ,
\end{align*}
which implies
\begin{align*}
\gamma_s {\rm Re} ( \p_x \BM u , & \, \BM^{\frac12} \p_x^{s-1} \{ f'(Q) \p_x^{s-1} \BM^{\frac12} u \} )_{L^2} 
 \\ & \leq 
 \gamma_s (-1)^{s-1} {\rm Re} ( \p_x^s u , f'(Q) \p_x^{s-1} \BM^{2} u )_{L^2} 
 + \frac32 \gamma_s (-1)^{s-1} {\rm Re} ( \p_x^s u , \p_x [f'(Q)] \BM^{\frac12} \mathscr{H} 
 \{ \p_x^{s-1} \BM^{\frac12} u \} )_{L^2} 
 \\ & \quad \quad + C \n u \n_{H^s} \n u \n_{H^{s-\frac12} } .
 \end{align*}
Noticing that $ \BM^2 = - \p_x^2 $ and $ \BM^{\frac12} \mathscr{H} \p_x^{s-1} \BM^{\frac12} 
= \p_x^{s-1} \BM \mathscr{H} = - \p_x^{s} $ (since $ \BM \mathscr{H} $ has symbol 
equal to $ - i \xi $), we infer
\begin{align*}
\gamma_s {\rm Re} ( \p_x \BM u , & \, \BM^{\frac12} \p_x^{s-1} \{ f'(Q) \p_x^{s-1} \BM^{\frac12} u \} )_{L^2} 
 \\ & \leq 
 \gamma_s (-1)^{s} {\rm Re} ( \p_x^s u , f'(Q) \p_x^{s+1} u )_{L^2} 
 + \frac32 \gamma_s (-1)^{s} {\rm Re} ( \p_x^s u , \p_x [f'(Q)] \p_x^{s} u )_{L^2} 
% \\ & \quad \quad 
+ C \n u \n_{H^s} \n u \n_{H^{s-\frac12} } 
  \\ & = \gamma_s (-1)^{s} {\rm Re} ( \p_x^s u , \p_x[ f'(Q) ] \p_x^{s} u )_{L^2} 
  + C \n u \n_{H^s} \n u \n_{H^{s-\frac12} } 
 \end{align*}
by integration by parts. 

Reporting these estimates into (B.1), we infer
\begin{align*}
 {\rm Re} ( \p_x ( \BM + c + f'(Q) ) u , \mathbb{M}_s u )_{L^2} 
 \leq & \Big( s + \frac12 \Big) {\rm Re} ( \p_x [f'(Q)] \p_x^{s} u , \p_x^{s} u )_{L^2} 
 + \gamma_s (-1)^{s} {\rm Re} ( \p_x^s u , \p_x[ f'(Q) ] \p_x^{s} u )_{L^2} \\
 & + C \n u \n_{H^s} \n u \n_{H^{s-\frac12} } .
\end{align*}
Therefore, the choice
$$
 \gamma_s \equiv (-1)^{s-1} \Big( s + \frac12 \Big)
$$
provides the desired control
$$
  {\rm Re} ( \p_x ( \BM + c + f'(Q) ) u , \mathbb{M}_s u )_{L^2} \leq 
  C \n u \n_{H^s} \n u \n_{H^{s-\frac12} } .
$$
When $ \{ s\} = 1/2 $, the computations are similar: (B.1) becomes now
\begin{align}
 {\rm Re} ( \p_x ( \BM + c + & \, f'(Q) ) u , \mathbb{M}_s u )_{L^2} 
 \nonumber \\
& =  {\rm Re} ( \p_x \BM u , (-1)^{[s]} \p_x^{2[s]} \BM u )_{L^2} 
 + \gamma_s {\rm Re} ( \p_x \BM u , \p_x^{[s]} \{ f'(Q) \p_x^{ [s]} u \} )_{L^2} 
 \nonumber \\ & \quad
 + {\rm Re} ( \p_x [ f'(Q) u ] , (-1)^{[s]} \p_x^{2[s]} \BM u )_{L^2} 
 + \gamma_s {\rm Re} ( \p_x [ f'(Q) u ] , \p_x^{[s]} \{ f'(Q) \p_x^{ [s]} u \} )_{L^2} 
 \nonumber \\ \tag{B.2}
& \quad + c  {\rm Re} ( \p_x u , \mathbb{M}_s u )_{L^2} ,
\end{align}
and the first and last scalar product still vanish. Moreover, by integration 
by parts and Leibniz formula, we deduce, since $ Q \in H^\ii $, 
\begin{align*}
 \gamma_s {\rm Re} ( \p_x [ f'(Q) u ] , & \, \p_x^{[s]} \{ f'(Q) \p_x^{[s]}u \} )_{L^2} 
 = 
 \gamma_s (-1)^{[s]} {\rm Re} ( \p_x^{ [s] + 1} [ f'(Q) u ] , f'(Q) \p_x^{[s]} u )_{L^2} 
 \\ \leq & \, 
 \gamma_s (-1)^{ [s] } {\rm Re} ( f'(Q) \p_x^{[s] + 1} u , f'(Q) \p_x^{[s]}u )_{L^2} 
 + C \n u \n_{H^{ [s]} }^2 
  \\ \leq & \, 
 \gamma_s (-1)^{ [s] -1} {\rm Re} (  \p_x [ f'(Q) ]  \p_x^{[s]} u  , f'(Q) \p_x^{[s]} u )_{L^2} 
 + C \n u \n_{H^{ [s]} }^2 
 \\ \leq & \, C \n u \n_{H^{ [s]} }^2 = C \n u \n_{H^{ s - \frac12 } }^2 .
\end{align*}
Furthermore,
\begin{align*}
 {\rm Re} ( \p_x [ f'(Q) u ] , (-1)^{[s]} \p_x^{2[s]} \BM u )_{L^2} 
 = & \, {\rm Re} ( \BM^{\frac12} \p_x^{[s]+1} [ f'(Q) u ] , \p_x^{[s]} \BM^{\frac12} u )_{L^2} 
 \\ 
 \leq & \, {\rm Re} ( \p_x \BM^{\frac12} \{ f'(Q) \p_x^{[s]} u \}, \p_x^{[s]} \BM^{\frac12} u )_{L^2} 
 \\  & \,
 + [s] {\rm Re} ( \p_x \BM^{\frac12} \{ \p_x [f'(Q)] \p_x^{[s]-1} u \} , \p_x^{[s]} \BM^{\frac12} u )_{L^2} 
 + C \n u \n_{H^{ [s] - \frac12}} \n u \n_{H^{ [s] + \frac12}} .
 \end{align*}
 For the second scalar product, we write, by Lemma B.1,
\begin{align*}
 {\rm Re} ( \p_x \BM^{\frac12} \{ \p_x [f'(Q)] \p_x^{[s]-1} u \} , & \, \p_x^{[s]} \BM^{\frac12} u )_{L^2} \\
 = & \, {\rm Re} ( \BM^{\frac12} \{ \p_x^2 [f'(Q)] \p_x^{[s]-1} u \} , \p_x^{[s]} \BM^{\frac12} u )_{L^2} 
 + {\rm Re} ( \BM^{\frac12} \{ \p_x [f'(Q)] \p_x^{[s]} u \} , \p_x^{[s]} \BM^{\frac12} u )_{L^2} 
 \\ \leq & \, C \n u \n_{H^{ [s] - \frac12}} \n u \n_{H^{ [s] + \frac12}} 
 + {\rm Re} ( \p_x [f'(Q)] \p_x^{[s]} \BM^{\frac12} u , \p_x^{[s]} \BM^{\frac12} u )_{L^2} 
 + C \n u \n_{H^{ [s] }} \n u \n_{H^{ [s] + \frac12}} 
  \\ \leq & \, {\rm Re} ( \p_x [f'(Q)] \p_x^{[s]} \BM^{\frac12} u , \p_x^{[s]} \BM^{\frac12} u )_{L^2} 
 + C \n u \n_{H^{ s- \frac12 }} \n u \n_{H^{ s }} .
\end{align*}
For the first scalar product, we use Lemma B.1 once again:
\begin{align*}
 {\rm Re} ( \p_x & \, \BM^{\frac12} \{ f'(Q) \p_x^{[s]} u \}, \p_x^{[s]} \BM^{\frac12} u )_{L^2} 
 \\ \leq & \, {\rm Re} ( f'(Q) \p_x \BM^{\frac12} \p_x^{[s] } u , \p_x^{[s]} \BM^{\frac12} u )_{L^2} 
 + \frac32 {\rm Re} ( \p_x [ f'(Q) ] \p_x^{[s]} \BM^{\frac12} u , \p_x^{[s]} \BM^{\frac12} u )_{L^2} 
 + C  \n u \n_{H^{ s- \frac12 }} \n u \n_{H^{ s }} 
 \\ \leq & \, {\rm Re} ( \p_x [ f'(Q) ] \p_x^{[s]} \BM^{\frac12} u , \p_x^{[s]} \BM^{\frac12} u )_{L^2} 
 + C  \n u \n_{H^{ s- \frac12 }} \n u \n_{H^{ s }} .
\end{align*}
As a consequence, since $ [s] = s - \frac12 $,
$$
 {\rm Re} ( \p_x [ f'(Q) u ] , (-1)^{[s]} \p_x^{2[s]} \BM u )_{L^2} \leq 
 \Big( s + \frac12 \Big) {\rm Re} ( \p_x [ f'(Q) ] \p_x^{[s]} \BM^{\frac12} u , 
\p_x^{[s]} \BM^{\frac12} u )_{L^2} 
 + C \n u \n_{H^{ s- \frac12 }} \n u \n_{H^{ s }} .
$$
We turn finally to the term
$$
 \gamma_s {\rm Re} ( \p_x \BM u , \p_x^{[s]} \{ f'(Q) \p_x^{ [s]} u \} )_{L^2} 
 = \gamma_s (-1)^{[s]} {\rm Re} ( \p_x^{[s]} \BM^{\frac12} u , 
\p_x \BM^{\frac12} \{ f'(Q) \p_x^{ [s]} u \} )_{L^2} ,
$$
and infer, by Lemma B.1,
\begin{align*}
\gamma_s {\rm Re} ( \p_x \BM u , \p_x^{[s]} \{ f'(Q) \p_x^{ [s]} u \} )_{L^2} 
\leq & \, 
\gamma_s (-1)^{[s]} {\rm Re} ( \p_x^{[s]} \BM^{\frac12} u , f'(Q) \p_x \BM^{\frac12} \p_x^{ [s]} u )_{L^2} 
\\ & \quad 
+ \frac32 \gamma_s (-1)^{[s]} {\rm Re} ( \p_x^{[s]} \BM^{\frac12} u , \p_x [  f'(Q) ] 
\BM^{\frac12} \p_x^{ [s]} u )_{L^2} + C \n u \n_{H^{ s- \frac12 }} \n u \n_{H^{ s }} \\
= & \, \gamma_s (-1)^{[s]} {\rm Re} ( \p_x^{[s]} \BM^{\frac12} u , \p_x [  f'(Q) ] 
\BM^{\frac12} \p_x^{ [s]} u )_{L^2} + C \n u \n_{H^{ s- \frac12 }} \n u \n_{H^{ s }} .
\end{align*}
Therefore,
$$
 {\rm Re} ( \p_x ( \BM + c + f'(Q) ) u , \mathbb{M}_s u )_{L^2}  \leq 
 \Big( s + \frac12 + \gamma_s (-1)^{[s]} \Big) {\rm Re} ( \p_x^{[s]} \BM^{\frac12} u , \p_x [  f'(Q) ] 
\BM^{\frac12} \p_x^{ [s]} u )_{L^2} + C \n u \n_{H^{ s- \frac12 }} \n u \n_{H^{ s }} ,
$$
hence choosing $ \gamma_s \equiv (-1)^{[s]- 1} ( s + \frac12 ) $ gives the result. 

\bigskip

It remains to study the cases of the Intermediate Long Wave equation and the Smith 
equation, for which $ \hat{ \BM} $ is, respectively,
$$
 \frac{ \xi}{ {\rm tanh}(\xi H) } - \frac{1}{ H} ; \quad \quad \quad
 \sqrt{1+ \xi^2 } - 1 .
$$
We denote $ \BM_0 $ the operator with symbol $ | \xi | $ (the one of the Benjamin-Ono 
equation), and define $ \mathbb{M}_s $ as for the Benjamin-Ono case (hence with 
"$\BM$"$= \BM_0$). We observe that in both cases, $ \tilde{\BM} \equiv \BM - \BM_0 $ 
is bounded on $ L^2 $. Indeed, its symbol is continuous in $ \R $ and, for $ \xi \to \pm \ii $,  
$$ 
\hat{\BM}(\xi ) =  \frac{ \xi}{ {\rm tanh}(\xi H) } - \frac{1}{ H} 
= \frac{ \xi}{ {\rm sgn}(\xi) + \BO( \ex^{ -2 | \xi | H }) } - \frac{1}{ H} 
= | \xi |  - \frac{1}{ H} + \BO( | \xi | \ex^{ -2 | \xi | H } )  
$$
and
$$ 
 \hat{\BM}(\xi ) = \sqrt{1+ \xi^2 } - 1 = | \xi | \sqrt{1+ \xi^{-2} } - 1 
 = | \xi | - 1 + \BO( | \xi|^{-1} ) 
$$
respectively. In the quantity $ {\rm Re} ( \p_x ( \BM + c + f'(Q) ) u , \mathbb{M}_s u )_{L^2} $, 
we then have to bound from above the extra term 
$ {\rm Re} ( \p_x ( \tilde{\BM} u ) , \mathbb{M}_s u )_{L^2} $, that is (using the 
skew-adjointness for the higher order derivatives in $ \mathbb{M}_s $),
$$
{\rm Re} ( \p_x ( \tilde{\BM} u ) , \gamma_s \BM_0^{\frac12} \p_x^{s-1} 
\{ f'(Q) \p_x^{s-1} \BM_0^{\frac12} u \} )_{L^2} 
= \gamma_s (-1)^{s-1} {\rm Re} ( \p_x^s \BM_0^{\frac12} ( \tilde{\BM} u ) , 
f'(Q) \p_x^{s-1} \BM_0^{\frac12}  u )_{L^2} \quad {\rm if \ } \{ s \} = 0 ;$$
$$
{\rm Re} ( \p_x ( \tilde{\BM} u ) ,  \gamma_s \p_x^{[s]} \{ f'(Q) \p_x^{[s]} u \} )_{L^2}  
= \gamma_s (-1)^{[s]} {\rm Re} ( \p_x^{[s]+1} ( \tilde{\BM} u ) , f'(Q) \p_x^{[s]} u )_{L^2}  
\quad \quad {\rm if \ } \{ s \} = \frac12 .
$$
We then note that in both cases, one may actually split 
$ \tilde{\BM} = \BM - \BM_0 = \tilde{\BM}_{c} + \tilde{\BM}_{h} $, where 
$ \tilde{\BM}_{c} $ is the multiplication by $ - 1/ H $ (resp. $-1$) and 
$ \tilde{\BM}_{h} $ has a symbol which is continuous in $ \R $ and 
$\BO( |\xi|^{-1}) $ at infinity, so that $ \tilde{\BM}_{h} $ is bounded from $H^\s$ to 
$H^{\s+1}$ if $ \s \geq 0 $. Therefore, when $ \{ s \} = 0 $, we easily get
\begin{align*}
{\rm Re} ( \p_x ( \tilde{\BM} u ) & \, , 
\gamma_s \BM_0^{\frac12} \p_x^{s-1} \{ f'(Q) \p_x^{s-1} \BM_0^{\frac12} u \} )_{L^2} \\ 
= & \, \gamma_s (-1)^{s-1} {\rm Re} ( \p_x^s \BM_0^{\frac12} ( \tilde{\BM}_c u ) , 
f'(Q) \p_x^{s-1} \BM_0^{\frac12}  u )_{L^2}  
+\gamma_s (-1)^{s-1} {\rm Re} ( \p_x^s \BM_0^{\frac12} ( \tilde{\BM}_h u ) , 
f'(Q) \p_x^{s-1} \BM_0^{\frac12}  u )_{L^2} 
\\ 
\leq & \, \frac12 \gamma_s (-1)^{s} {\rm Re} ( \p_x^{s-1} \BM_0^{\frac12} u , 
\tilde{\BM}_c \p_x [f'(Q) ] \p_x^{s-1} \BM_0^{\frac12}  u )_{L^2}  
+ C \n u \n_{H^{s - \frac12 }}^2 \leq C \n u \n_{H^{s - \frac12 }}^2 
\end{align*}
and similarly when $ \{ s \} = 1/2 $. Therefore, the estimate
$$
 {\rm Re} ( \p_x ( \BM + c + f'(Q) ) u , \mathbb{M}_s u )_{L^2}  \leq 
 C \n u \n_{H^{ s- \frac12 }} \n u \n_{H^{ s }} 
$$
remains true for the Intermediate Long Wave equation and the Smith equation. The 
proof of Proposition B.1 is thus complete by applying Theorem B.6. \carre

%%%%%%%%%%%%%%%%%%%%%%%%%%%%%%%%%%%%%%%%%%%%%%%%%%%%%%%%%%%%%%%%%%%%%%%%%%%
\section*{B.3 Proofs}

\subsection*{B.3.1 Proof of Theorem B.3}
%\noindent {\bf Proof of Theorem B.3.} 
We shall 
prove the resolvent estimate required in Corollary B.1. 
Let us consider $ \lambda = \gamma + i \tau \in \C $ with $ \gamma \not = 0 $ 
and the resolvent equation $ ( \BJ \BL - \lambda ) v = \Sigma $, or
\be
\tag{B.3}
 ( \gamma + i \tau ) v = \BJ \BL (v) - \Sigma .
\ee
By hypothesis, the essential spectrum of $ \BJ \BL $ is of the form 
$ i [ \R \setminus (- \vartheta_0 ,+ \vartheta_0 ) ] $. Moreover, we have seen 
that $ \BJ \BL $ has a finite number of eigenvalues in the half-space 
$ \{ {\rm Re} > 0 \} $, hence, for $|\tau| \geq \tau_0 $ sufficiently large, we 
know that there exists a unique solution $v$ to (B.3). By taking the scalar 
product with $ \BL (v) $, we deduce the conservation law
\be
\tag{B.4}
 \gamma ( v , \BL (v) )_\BX = - {\rm Re} ( \Sigma , \BL (v) )_\BX .
\ee
By our assumption, there exists a finite (possibly empty) number of eigenvalues 
in $ ( -\ii , 0 ] $, $ ( - \mu_1 , ... , - \mu_q ) $, each one of finite multiplicity. 
For any $ 1 \leq k \leq q $, we fix an orthonormal basis 
$ (\chi_{k,\ell})_{1 \leq \ell \leq n_k} $ of the eigenspace $ {\rm Ker} ( \BL + \mu_k ) $. 
By assumption (A), any eigenvector $ \chi_{k,\ell} $ is smooth in the sense that 
$ \chi_{k,\ell} \in D(\BJ) $ and $ \BJ \chi_{k,\ell} \in D(\BL ) $.

We then make a spectral orthogonal decomposition
$$
 v = \sum_{k=1}^q \sum_{\ell=1}^{n_k} \alpha_{k,\ell} \chi_{k,\ell} + v_+,
$$
where $ \BL ( \chi_{k,\ell} ) = \mu_k \chi_{k,\ell} $ and 
$ ( v_+ , \BL (v_+) )_\BX \geq \d |\!| v_+ |\!|_\BX^2 $ for some 
positive $ \d $. In the double sum, we have a finite number (independent of $v$) 
of terms. Inserting this into (B.4) yields
$$
 |\gamma| \d |\!| v_+ |\!|_\BX^2 \leq |\gamma| \d ( v_+ , \BL (v_+) )_\BX 
\leq \d \Big[ | {\rm Re} ( \Sigma , \BL (v) )_\BX | 
+ \sum_{k,\ell} \mu_k |\alpha_{k,\ell} |^2 \Big] 
\leq K |\!| \Sigma |\!|_\BX |\!| v |\!|_\BX + K \sum_{k,\ell} |\alpha_{k,\ell} |^2 
$$
Using the inequality $ a b \leq \e a^2 + b^2 / ( 4\e ) $ with 
$ a =  |\!| v |\!|_\BX $, $ b = K |\!| \Sigma |\!|_\BX $ and $\e = | \gamma |\d / 2 $, 
the equality $ |\!| v |\!|_\BX^2 = |\!| v_+ |\!|_\BX^2 + \sum_{k,\ell} |\alpha_{k,\ell}|^2 $ 
and incorporating the term $ |\gamma| \d |\!| v_+ |\!|_\BX^2 / 2 $ in the left-hand side, 
we infer
\be
\tag{B.5}
\frac{ |\gamma| \d }{2} |\!| v_+ |\!|_\BX^2 \leq 
K' \sum_{k,\ell} |\alpha_{k,\ell} |^2 + K'' |\!| \Sigma |\!|_\BX^2 .
\ee
On the other hand, since $ \chi_{k,\ell} \in D(\BJ) $ and 
$ \BJ \chi_{k,\ell} \in D(\BL ) $ by assumption (A), taking the 
scalar product of (B.1) with $ \chi_{k,\ell} $ provides
$$
 ( \gamma + i \tau ) \alpha_{k,\ell} = - ( v, \BL \BJ \chi_{k,\ell} )_\BX
- ( \Sigma , \chi_{k,\ell} )_\BX .
$$
Consequently,
$$
 ( |\gamma| + | \tau | ) | \alpha_{k,\ell} | \leq K_{k,\ell} |\!| v |\!|_\BX 
+ K |\!| \Sigma |\!|_\BX ,
$$
thus
$$
 ( |\gamma| + | \tau | )^2 \sum_{k,\ell} | \alpha_{k,\ell} |^2 
\leq K_0 |\!| v |\!|_\BX^2 + K |\!| \Sigma |\!|_\BX^2 
= K_0 \sum_{k,\ell} | \alpha_{k,\ell} |^2 + K |\!| v_+ |\!|_\BX^2 + K |\!| \Sigma |\!|_\BX^2 ,
$$
which implies, if $ | \tau | \geq 1 + \sqrt{K_0} - |\gamma| $,
$$
\sum_{k,\ell} | \alpha_{k,\ell} |^2 \leq K \frac{ |\!| v_+ |\!|_\BX^2 + |\!| \Sigma |\!|_\BX^2 }
{ ( |\gamma| + | \tau | )^2 - K_0 } .
$$
Reporting this into (B.5) gives
$$
\frac{ |\gamma| \d }{2} |\!| v_+ |\!|_\BX^2 \leq 
K'K \frac{ |\!| v_+ |\!|_\BX^2 + |\!| \Sigma |\!|_\BX^2 }
{ ( |\gamma| + | \tau | )^2 - K_0 } + K'' |\!| \Sigma |\!|_\BX^2 .
$$
If $ | \tau | \geq 1 + \sqrt{K_0 + 4 K K' / |\gamma| \d } - |\gamma| $, we deduce
$$
\frac{ |\gamma| \d }{4} |\!| v_+ |\!|_\BX^2 \leq 
\Big( K'' + \frac{K'K}{ ( |\gamma| + | \tau | )^2 - K_0 } 
\Big) \n \Sigma \n_\BX^2 \leq K_1 \n \Sigma \n_\BX^2 ,
$$
and it follows that
$$
\n v \n_\BX^2 = \n v_+ \n_\BX^2 + \sum_{k,\ell} |\alpha_{k , \ell } |^2 
\leq K_2 \n \Sigma \n_\BX^2 ,
$$
where $K_2$ does not depend on $ |\tau | $ (large enough) as wished.

\bigskip

The proof of the first semigroup estimate follows then easily, see, for instance, 
Proposition 2 in \cite{Pruss}.

\bigskip

\noindent {\bf Proof of the semigroup estimate when $\bs{ \gamma_0 > 0 }$.} 
Here, we assume $ \gamma_0 > 0 $. As a consequence, the spectrum of $ [\BJ \BL]_\C $ 
is of the form $ \s_{\rm s} \cup \s_{\rm u} $, where 
$ \s_{\rm ess} ( [\BJ \BL]_\C ) \subset \s_{\rm s} \subset \{ {\rm Re} \leq 0 \} $ 
and $ \emptyset \not = \s_{\rm u} \subset \{ {\rm Re} > 0 \} $ consists in a 
finite number of eigenvalues of finite algebraic multiplicities. 
Therefore, we may define (see, {\it e.g.} \cite{Kato}, \cite{HisSig}) 
the spectral Riesz projection
$$
 \mathbb{P} \equiv \frac{1}{2i\pi} \int_{\Gamma} ( [\BJ \BL]_\C - z )^{-1} \ dz ,
$$
where $ \Gamma $ is any simple (positively oriented) closed curve enclosing 
$ \s_{\rm u} $. As a consequence, $ \mathbb{P} $ is bounded, commutes with 
$ [ \BJ \BL]_\C $ on $ D([\BJ \BL]_\C ) $ and verifies $ \s( [\BJ \BL]_\C \mathbb{P} ) 
= \s_{\rm u} $, $ \s( [ \BJ \BL ]_\C ( {\rm Id} - \mathbb{P} ) ) = \s_{\rm s} $. Moreover, 
$ [\BJ \BL]_\C \mathbb{P} $ is bounded, hence generates a continuous semigroup, 
$ \ex^{t [\BJ \BL]_\C \mathbb{P} } $, given by the exponential series 
$$ \ex^{t [\BJ \BL]_\C \mathbb{P} } = 
\sum_{n=0}^{+\ii} \frac{t^n ([\BJ \BL]_\C \mathbb{P})^n}{n!} . $$
In addition, $ [\BJ \BL]_\C ( {\rm Id} - \mathbb{P} ) = [\BJ \BL]_\C - [\BJ \BL]_\C \mathbb{P} $ 
also generates a continuous semigroup and that we have $ \ex^{t[\BJ \BL]_\C } = 
\ex^{t \BJ \BL \mathbb{P} } \ex^{t [\BJ \BL]_\C ( {\rm Id} - \mathbb{P} ) } $. 

The semigroup generated by the bounded operator $ [\BJ \BL]_\C \mathbb{P} $ is 
easily analyzed. We shall now apply the spectral mapping theorem of J. Pr\"uss 
(Theorem B.2) to $ [\BJ \BL]_\C ( {\rm Id} - \mathbb{P} ) $ in order to 
control the growth of its norm. By Corollary B.1, it suffices to estimate 
its resolvent $ [ [\BJ \BL]_\C ( {\rm Id} - \mathbb{P} ) - ( \gamma + i \tau ) ]^{-1} $ 
for large $ | \tau | $ (note that $ \s( [\BJ \BL]_\C ( {\rm Id} - \mathbb{P} ) ) = \s_{\rm s} 
\subset \{ {\rm Re} \leq 0 \} $. If $ \Sigma \in \BX_\C $ and $| \tau | $ is large, 
it is clear that the solution $ u \in \BX_\C $ to 
$ [ [\BJ \BL]_\C ( {\rm Id} - \mathbb{P} ) - ( \gamma + i \tau ) ]u = \Sigma $ 
is given by
$$ u = [ [\BJ \BL]_\C - ( \gamma + i \tau ) ]^{-1} ( {\rm Id} - \mathbb{P} ) \Sigma 
- \frac{1}{ \gamma + i \tau } \mathbb{P} \Sigma , $$
thus, for $| \tau | $ large,
$$
 \no [ [\BJ \BL]_\C ( {\rm Id} - \mathbb{P}) -(\gamma + i \tau )]^{-1} \no_{\mathscr{L}_c(\BX_\C)} 
\leq 
\no [ [\BJ \BL]_\C - ( \gamma + i \tau ) ]^{-1} \no_{\mathscr{L}_c(\BX_\C)} 
 \n {\rm Id} - \mathbb{P} \n_{\mathscr{L}_c(\BX)} 
+ \frac{1}{ |\gamma + i \tau |} \n \mathbb{P} \n_{\mathscr{L}_c(\BX_\C)}  
$$
is bounded. Consequently, by Theorem B.2 and since $ \s_{\rm s} \subset \{ {\rm Re} \leq 0 \} $, 
$ \s ( \ex^{t [\BJ \BL]_\C ( {\rm Id} - \mathbb{P} ) } ) = 
\ex^{t \s([\BJ \BL]_\C ( {\rm Id} - \mathbb{P} ) )} = \ex^{t \s_{\rm s}} \subset \bar{D}(0,1) $. 
It follows that for any $ \epsilon > 0 $, there exists $ K_\epsilon > 0$ such that
$$
\forall \ t \geq 0, \quad \quad \quad 
\no \ex^{t [\BJ \BL]_\C ( {\rm Id} - \mathbb{P} ) } \n_{\mathscr{L}_c(\BX_\C)} 
\leq K_\epsilon \ex^{ \epsilon t } .
$$
Since $ \ex^{t [\BJ \BL]_\C \mathbb{P} } $ is given by the exponential series, 
we also have the optimal estimate
$$
\forall \ t \geq 0, \quad \quad \quad 
\no \ex^{t [\BJ \BL]_\C \mathbb{P} } \n_{\mathscr{L}_c(\BX_\C)} 
\leq K_0 (1 + t)^{m-1} \ex^{ \gamma_0 t } 
$$
by definition of $m$. We conclude by taking $ \epsilon = \gamma_0 / 2 $ for instance.

%%%%%%%%%%%%%%%%%%%%%%%%%%%%%%%%%%%%%%%%%%%%%%%%%%%%%%%%%%%%%%%%%
\subsection*{B.3.3 Proof of Theorem B.4}

%\noindent {\bf Proof of Theorem B.4.} 
Since $\BA$ generates a continuous 
semigroup, $v$ is a solution to $ \p_t v = \BA v + \Phi(v) $ if and only if it is a 
mild solution:
$$
 v(t) = \ex^{t \BA } v^\ini + \int_0^t \ex^{ (t-\tau) \BA } \Phi( v(\tau) ) \ d \tau .
$$
There exists $ r_0 > 0 $ such that $ |\!| \Phi (v) |\!|_X \leq M |\!| v |\!|_X^2 $ if 
$ |\!| v |\!|_X \leq r_0 $. We choose $ v^\ini = \d {\rm Re}\, w $, where 
$ |\!| {\rm Re}\, w |\!|_X = 1 $ and $ w $ is an eigenvector for the eigenvalue $ \lambda $, 
and write the solution under the form 
$ v = \ex^{t \BA } v^\ini + \tilde{v} = {\rm Re} ( \ex^{ t \lambda } w ) + \tilde{v} $. 
If $ \lambda \in \R $, we can choose $ w \in D( \BA ) \subset D( \BA_\C )$. Then,
$$
 \tilde{v}(t) = \int_0^t \ex^{ (t-\tau) \BA } \Phi( \d {\rm Re} ( \ex^{ t \lambda } w) 
+ \tilde{v}(\tau) ) \ d \tau .
$$
Let us denote $ r_1 \equiv \min \Big( r_0, \ds{\frac{\gamma_0 - \beta}{2 M M_0 } } \Big) $ 
and let $ T > 0$ be the maximal time such that $ T < \ln( r / (2 \d ) ) / \gamma_0 $ and 
$ |\!| \tilde{u}(\tau) |\!|_X < r_1/2 $ in $ [0, T) $, where $ 0 < r < r_1 $ will 
be determined later. We shall work for $0 \leq t < T $, so that 
$ |\!| \d {\rm Re}( \ex^{ t \lambda } w ) + \tilde{v}(\tau)|\!|_X 
< \d \ex^{ t \gamma_0 } + r_1/2 \leq r_1 \leq r_0 $. Then,
\begin{align*}
 \n \tilde{v}(t) \n \leq & \, \int_0^t 
\no \ex^{ (t-\tau) \BA } \no_{\mathscr{L}_c(X)} 
M \n \d {\rm Re} ( \ex^{ \tau \lambda } w ) + \tilde{v}(\tau) \n^2 \ d \tau \\
\leq & \, 
2 M_0 M \int_0^t \ex^{ ( \gamma_0 + \beta ) (t-\tau) } 
\Big( \d^2 \ex^{ 2 \tau \gamma_0 } + \n \tilde{v}(\tau) \n^2 \Big) \ d \tau 
\\ \leq & \, 
\frac{ 2 M_0 M}{\gamma_0 - \beta } \d^2 \ex^{ 2 t \gamma_0 } 
+ r_1 M_0 M \int_0^t \ex^{ ( \gamma_0 + \beta ) (t-\tau)  } \n \tilde{v}(\tau) \n_X \ d \tau ,
\end{align*}
since $ \beta < \gamma_0 $. Applying now the Gronwall inequality to 
$  \ex^{ -( \gamma_0 + \beta ) t } |\!| \tilde{v}( t) |\!|_X $ then gives, since 
$ M_0 M r_1 < \gamma_0 - \beta $,
$$
\n \tilde{u}(t) \n_X \leq \Big[ \frac{2 M_0 M }{ \gamma_0 - \beta } 
+ \frac{ 2 r_1 M_0^2 M^2 }{ ( \gamma_0 - \beta ) (  \gamma_0 - \beta - r_1 M_0 M)} \Big] 
\d^2 \ex^{ 2 t \gamma_0 } = K \d^2 \ex^{ 2 t \gamma_0 } .
$$
We now choose $ r \equiv \sqrt{r_1 / K} $, so that the right-hand side is 
$ \leq K r^2 / 4 < r_1 / 2 $, and this implies that $u$ exists at least on 
$ [0, \ln( r/ (2 \d)) / \gamma_0 ] $. In addition, for $ 0 \leq t < T $,
$$
 \n u(t) \n_X \geq \d \ex^{ t \gamma _0} - \n \tilde{u}(t) \n_X \geq 
 \d \ex^{ t \gamma _0} - K \d^2 \ex^{ 2 t \gamma_0 } ,
$$
as desired. We conclude choosing $\e_0 >0 $ so small that 
$ 2 \e_0 - K \e_0^2 \geq \e_0$.

%\bigskip

%%%%%%%%%%
\subsection*{B.3.3 Proof of Corollary B.2}

%\noindent {\bf Proof of Corollary B.2.} 
We pick some $ 0 < \beta < \gamma_0 $ 
(for instance $ \beta = \gamma_0 / 2 $) in order to have the semigroup estimate 
required in Theorem B.4. The solution $u(t) = T(\oo_* t ) ( \phi_{\oo_*} + v (t)) $ 
verifies, for $ 0 \leq t \leq \gamma_0^{-1} \ln (2 \e_0 / \delta ) $,
$$
\n v(t) \n_\BX = \n T( - \oo_* t) u(t) - (\phi_{\oo_*} + \d {\rm Re} ( \ex^{t\lambda} w) ) \n_\BX 
\leq K \d^2 \ex^{ 2 t \gamma_0 } .
$$
Hence, $ T( - \oo_* t) u(t) $ remains at distance $ \leq K \e_0$ from 
$ \phi_{\oo_*} \in \mathfrak{M} $ and therefore 
$$
 {\rm dist}_\BX ( u(t) , \mathfrak{M} ) \geq 
 {\rm dist}_\BX ( \d {\rm Re} ( \ex^{ t \lambda } w ) , \mathfrak{M} - \phi ) 
- K \d^2 \ex^{ 2 t \gamma_0 } .
$$
Assume $ \lambda \in \R $. Then, we observe that the straight line 
$ \R \ni \theta \mapsto \theta w $ is transverse to the tangent space T$_\phi \mathfrak{M}$ 
of the manifold $\mathfrak{M}$, since $ w $ is an eigenvector of $ \BJ \BL $ for 
$ \lambda \not = 0 $, hence does not belong to the kernel of $ \BL $. Therefore, 
$ {\rm dist}_\BX ( \theta w , \mathfrak{M} - \phi ) \geq | \theta | / K_1 $ 
for small $ | \theta | $. Thus,
$$
 {\rm dist}_\BX ( u(t) , \mathfrak{M} ) \geq \frac{1}{K_1} \d \ex^{ t \lambda } 
- K \d^2 \ex^{ 2 t \gamma_0 } .
$$
If $ \lambda \in \C \setminus \R $, the equation $ [\BJ \BL]_\C (w) = \lambda w $ 
splits as $ \BJ \BL ( {\rm Re}\, w ) = 
{\rm Re}(\lambda) {\rm Re} \, w -  {\rm Im} (\lambda) {\rm Im}\, w $ 
and $ \BJ \BL ( {\rm Im} \, w ) = 
{\rm Im}(\lambda) {\rm Re}\, w + {\rm Re} (\lambda) {\rm Im} \, w $. 
Therefore, $ {\rm Re} \, w $ and $ {\rm Im} \, w $ do not belong to $\ker (\BL) $. 
Consequently, the surface $ \C \ni \theta \mapsto {\rm Re} ( \theta w ) $ is transverse 
to the tangent space T$_\phi \mathfrak{M}$ of the manifold $\mathfrak{M}$, and 
we conclude as before that
$$
 {\rm dist}_\BX ( u(t) , \mathfrak{M} ) \geq \frac{1}{K_1} \d \ex^{ t \gamma_0 } 
 - K \d^2 \ex^{ 2 t \gamma_0 } .
$$

\bigskip

\noindent
{\bf Acknowledgement. } We acknowledge the support of the French ANR (Agence 
Nationale de la Recherche) under Grant ANR JC ArDyPitEq.

%%%%%%%%%%%%%%%%%%%%%%%%%%%%%%%%%%%%%%%%%%%%%%%%%%%%%%%%%%%%%%%%%%%%%%%%%%%%%%%%%%%%%%%%%

%%%%%%%%%%%%%%%%%%%%%%%%%%%%%%%%%%%%%%%%%%%%%

\end{document}